\documentclass[noinfoline]{imsart}

\usepackage[dvips,letterpaper,textwidth=16.5cm,vcentering,hcentering]{geometry}
\usepackage{amsmath,amssymb,amsthm,amsfonts,amsbsy}
\usepackage{graphicx}
\usepackage{wrapfig}
\usepackage[normalem]{ulem}
\usepackage[active]{srcltx}
\usepackage{enumerate}
\usepackage{appendix}
\usepackage{multirow}
\usepackage{afterpage}
\usepackage{cite}
\usepackage{fancyhdr}
\usepackage[crop=pdfcrop]{pstool}
\usepackage{subfig}
\usepackage{hyperref}

\renewcommand%RM 03-10-13\newcommand
{\psfragfig}{\includegraphics} %%%%%%%%%%%%% 

\providecommand{\texorpdfstring}[2]{#1}
\providecommand{\url}[1]{#1}

\newcommand{\ol}[1]{\overline{#1}}
\renewcommand{\bar}{\ol}
\renewcommand{\le}{\leqslant}
\renewcommand{\ge}{\geqslant}

\newcommand{\tsum}{\mathop{{\textstyle\sum}}}
\newcommand{\tsumnl}{\mathop{{\textstyle\sum}}\nolimits}
\newcommand{\lims}{\limits}
\newcommand{\nolim}{\nolimits}

\newcommand{\ffrown}{\text{\raisebox{3pt}[0pt][0pt]{$\frown$}}}
\renewcommand{\O}{\underset{\ffrown}{<}}

\renewcommand{\P}{\operatorname{\mathbb{P}}} % probability
\DeclareMathOperator{\E}{\mathbb{E}}         % expectation
\DeclareMathOperator{\var}{\mathbb{V}ar}     % variance
\DeclareMathOperator{\cov}{\mathbb{C}ov}     % covariance
\DeclareMathOperator{\I}{I}                  % indicator
\DeclareMathOperator{\tr}{tr}                % trace
\newcommand{\T}{\ensuremath{^\mathsf{T}}}    % transpose
\newcommand{\Id}{{I}}                        % identity matrix
\newcommand{\R}{\mathbb{R}}                  % real #'s
\newcommand{\Z}{\mathbb{Z}}                  % integers
\newcommand{\N}{\mathbb{N}}                  % natural #'s
\newcommand{\C}{\mathbb{C}}
\DeclareMathOperator{\Lam}{{\sf W}}          % Lambert's product-log
\DeclareMathOperator{\PU}{{\sf PU_{tail}}}   % Pinelis-Utev tail probability bound
\DeclareMathOperator{\PUalt}{{\sf PU_{alt}}}     % alternative Pinelis-Utev tail probability bound
\DeclareMathOperator{\PUthr}{{\sf PU_3}}
\DeclareMathOperator{\PUtwo}{{\sf PU_2}}
\DeclareMathOperator{\PUexp}{{\sf PU_{exp}}} % Pinelis-Utev exponential moment bound
\DeclareMathOperator{\BH}{{\sf BH_{tail}}}   % Bennett-Hoeffding tail probability bound
\DeclareMathOperator{\BHalt}{{\sf BH_{alt}}}     % alternative Bennett-Hoeffding tail probability bound
\DeclareMathOperator{\BHexp}{{\sf BH_{exp}}} % Bennett-Hoeffding exponential moment bound

\DeclareMathOperator{\intr}{int}             % interior
\DeclareMathOperator{\sgn}{{\sf Sgn}}        % sign (with Sgn(0)=1)
\newcommand{\nc}{\mathsf{nc}} %RM 02-21-13

\newcommand{\nov}{\operatorname{\mathsf{nov}}}
\newcommand{\as}{\operatorname{\mathsf{asymp}}}
\newcommand{\NZ}{\mathsf{N\!Z}}

\newcommand{\al}{\alpha}
\newcommand{\be}{\beta}
\newcommand{\ga}{\gamma}
\newcommand{\Ga}{\Gamma}
\newcommand{\de}{\delta}
\newcommand{\De}{\Delta}
\newcommand{\ep}{\epsilon}
\newcommand{\vp}{\varepsilon}
\renewcommand{\th}{\theta}
\newcommand{\vt}{\vartheta}
\newcommand{\Th}{\Theta}
\newcommand{\ka}{\kappa}
\newcommand{\vka}{\varkappa}
\newcommand{\la}{\lambda}
\newcommand{\La}{\Lambda}
\newcommand{\si}{\sigma}
\newcommand{\Si}{\Sigma}
\newcommand{\tSi}{\tilde\Sigma}
\newcommand{\vsi}{\varsigma}
\newcommand{\vpi}{\varphi}
\newcommand{\om}{\omega}
\newcommand{\Om}{\Omega}

\newcommand{\bs}[1]{\boldsymbol{#1}}
\newcommand{\bc}{{\bs{c}}}
\newcommand{\bu}{{\bs{u}}}
\newcommand{\x}{{\bs{x}}}

\newcommand{\tvp}{{\tilde{\vp}}}
\newcommand{\tth}{{\tilde{\theta}}}
\newcommand{\tsi}{{\tilde\sigma}}
\newcommand{\ta}{{\tilde{a}}}
\newcommand{\tc}{{\tilde{c}}}

\newcommand{\tx}{{\tilde{x}}}
\newcommand{\ty}{{\tilde{y}}} %RM 12.21.12 
\newcommand{\tw}{{\tilde{w}}}
\newcommand{\tA}{{\tilde{A}}}
\newcommand{\tS}{{\tilde{S}}}
\newcommand{\tX}{{\tilde{X}}}
\newcommand{\tY}{\tilde Y}

\newcommand{\uu}{{\mathfrak u}}
\newcommand{\vv}{{\mathfrak v}}
\renewcommand{\AA}{\mathfrak{A}}
\newcommand{\BB}{\mathfrak{B}}
\newcommand{\CC}{\mathfrak{C}}
\newcommand{\KK}{\mathfrak{K}}
\newcommand{\KKu}[1]{\KK_{\sf{u}#1}}

\newcommand{\XX}{\mathfrak{X}} 
\newcommand{\YY}{\mathfrak{Y}}
\newcommand{\ZZ}{\mathfrak{Z}}
\newcommand{\pp}{\mathfrak{p}}
\newcommand{\bb}{\mathfrak{b}}

\newcommand{\uhat}{{\hat{u}}}

\newcommand{\Ko}{\mathsf{Ko}}
\newcommand{\bW}{\mathsf{bW}}
\newcommand{\W}{\mathsf{W}}

\newcommand{\cc}{c_{\ast}}

\newcommand{\Mf}{M_\ep}
\newcommand{\ccc}{{\mathfrak K}}
\newcommand{\m}{\mathfrak{m}}
\newcommand{\res}{\mathcal{R}}
\newcommand{\tKKu}[1]{\tilde{\KK}_{{\sf u}#1}}
\newcommand{\Cucon}{\tKKu{0}}%{{\ccc}_{{\sf u}0}} %C_3
\newcommand{\Cuvsi}{\tKKu{1}}%{{\ccc}_{{\sf u}1}} %C_4
\newcommand{\Cuvtwo}{\tKKu{2}}%{{\ccc}_{{\sf u}2}} %C_5
\newcommand{\Cuvthr}{\tKKu{3}}%{{\ccc}_{{\sf u}3}} %C_6
\newcommand{\Cnvsi}{{\ccc}_{{\sf n}1}} %C_7
\newcommand{\Cnvtwoa}{{\ccc}_{{\sf n}21}} %C_{8,1}
\newcommand{\Cnvtwob}{{\ccc}_{{\sf n}22}} %C_{8,2}
\newcommand{\Cnvthra}{{\ccc}_{{\sf n}31}} %C_{8,3}
\newcommand{\Cnvthrb}{{\ccc}_{{\sf n}32}} %C_9
\newcommand{\Cnecon}{{\ccc}_{{\sf e0}}} %C_{10}
\newcommand{\Cnevsi}{{\ccc}_{{\sf e1}}} %C_{11}
\newcommand{\Cnevtwo}{{\ccc}_{{\sf e2}}} %C_{12}
\newcommand{\Cnevthr}{{\ccc}_{{\sf e3}}} %C_{13}
\newcommand{\LL}{{\mathfrak C}}

\newcommand{\Lnvsi}{{\LL}_{{\sf n}1}} %C_7
\newcommand{\Lnvtwoa}{{\LL}_{{\sf n}21}} %C_{8,1}
\newcommand{\Lnvtwob}{{\LL}_{{\sf n}22}} %C_{8,2}
\newcommand{\Lnvthra}{{\LL}_{{\sf n}31}} %C_{8,3}
\newcommand{\Lnvthrb}{{\LL}_{{\sf n}32}} %C_9
\newcommand{\Lnecon}{{\LL}_{{\sf e0}}} %C_{10}
\newcommand{\Lnevsi}{{\LL}_{{\sf e1}}} %C_{11}
\newcommand{\Lnevtwo}{{\LL}_{{\sf e2}}} %C_{12}
\newcommand{\Lnevthr}{{\LL}_{{\sf e3}}} %C_{13}
\newcommand{\A}{\mathcal{A}}
\newcommand{\B}{\mathcal{B}}
\newcommand{\EE}{\mathcal{E}}

\newcommand{\XXX}{\mathcal{X}}
\newcommand{\xx}{\x}%RM12.20\newcommand{\xx}{\mathbf{x}}
\newcommand{\X}{\bs{Y}}%RM12.20\newcommand{\X}{\mathbf{X}}

\newcommand{\Ybar}{\bar Y}
\newcommand{\Ysqbar}{\ol{Y^2}}
\newcommand{\YYTbar}{\ol{YY\T}}
\newcommand{\YYtensbar}{\ol{Y\otimes Y}}

\newcommand{\bigO}{O}

\newcommand{\tB}{B}

\newcommand{\ip}[2]{\langle#1,#2\rangle}

\providecommand{\bigprob}[1]{\P\bigl(#1\bigr)}
\providecommand{\Bigprob}[1]{\P\Bigl(#1\Bigr)}
\providecommand{\biggprob}[1]{\P\biggl(#1\biggr)}

\providecommand{\ind}[1]{\I\{#1\}}
\providecommand{\bigind}[1]{\I\bigl\{#1\bigr\}}

\providecommand{\abs}[1]{\lvert#1\rvert}
\providecommand{\bigabs}[1]{\bigl\lvert#1\bigr\rvert}
\providecommand{\Bigabs}[1]{\Bigl\lvert#1\Bigr\rvert}
\providecommand{\biggabs}[1]{\biggl\lvert#1\biggr\rvert}

\providecommand{\norm}[1]{\lVert#1\rVert}
\providecommand{\bignorm}[1]{\bigl\lVert#1\bigr\rVert}
\providecommand{\Bignorm}[1]{\Bigl\lVert#1\Bigr\rVert}

\providecommand{\regp}[1]{(#1)}
\providecommand{\bigp}[1]{\bigl(#1\bigr)}
\providecommand{\Bigp}[1]{\Bigl(#1\Bigr)}
\providecommand{\biggp}[1]{\biggl(#1\biggr)}

\providecommand{\bigb}[1]{\bigl\{#1\bigr\}}
\providecommand{\Bigb}[1]{\Bigl\{#1\Bigr\}}
\providecommand{\biggb}[1]{\biggl\{#1\biggr\}}

\newtheorem{thm}{Theorem}
\newtheorem{cor}[thm]{Corollary}
\newtheorem{lem}[thm]{Lemma}
\newtheorem{prop}[thm]{Proposition}
\theoremstyle{remark}
\newtheorem{remark}[thm]{Remark}
\newtheorem{ex}[thm]{Example}

\numberwithin{equation}{section}
\numberwithin{thm}{section}

\begin{document}

\begin{frontmatter}

%IP
%\title{Berry-Esseen bounds for general nonlinear statistics, with applications to Pearson's and non-central Student's and Hotelling's}
%\runtitle{Berry-Esseen bounds for nonlinear statistics}
%
%\begin{aug}
% \author{\fnms{Iosif}  \snm{Pinelis}\thanksref{a,e1}\ead[label=e1,mark]{ipinelis@mtu.edu}}
% \and
% \author{\fnms{Raymond} \snm{Molzon}\thanksref{a}}
% \runauthor{Iosif Pinelis and Raymond Molzon}
% \affiliation{Michigan Technological University}
% \address[a]{Department of Mathematical Sciences\\Michigan Technological University\\Hough\-ton, Michigan 49931\\\printead{e1}}
% \end{aug}
 
%IP
%IP16 \title{Berry-Esseen bounds for general nonlinear statistics, with applications to Pearson's and non-central Student's and Hotelling's\protect\thanksref{T1}}
%\runtitle{Berry-Esseen bounds for nonlinear statistics}
%\thankstext{T1}{Supported in part by NSF grant DMS-0805946 and NSA grant H98230-12-1-0237.}

\title{Optimal-order bounds on the rate of convergence to normality in the multivariate delta method}
\runtitle{Convergence rate in delta method}

%\begin{aug}
%\author{\fnms{Iosif} \snm{Pinelis}\thanksref{t2}\ead[label=e1]{ipinelis@mtu.edu}}
%  \thankstext{t2}{Supported in part by NSF grant DMS-0805946 and NSA grant H98230-12-1-0237}
%\runauthor{Iosif Pinelis}
%
%%IP
%\begin{aug}
% \author{\fnms{Iosif}  \snm{Pinelis}\thanksref{a,e2}\ead[label=e2,mark]{ipinelis@mtu.edu}}
% \and
% \author{\fnms{Raymond} \snm{Molzon}\thanksref{a}}
% \runauthor{Iosif Pinelis and Raymond Molzon}
% \affiliation{Michigan Technological University}
% \address[a]{Department of Mathematical Sciences\\Michigan Technological University\\Hough\-ton, Michigan 49931\\\printead{e2}}
% \end{aug}

%IP 
\begin{aug}
  \author{\fnms{Iosif}  \snm{Pinelis}\corref{}%\thanksref{t2}
  \ead[label=e1]{ipinelis@mtu.edu}}
  \and
  \author{\fnms{Raymond} \snm{Molzon}\ead[label=e2]{remolzon@mtu.edu}}
%  \and
%  \author{\fnms{Third}  \snm{Author}%
%  \ead[label=e3]{third@somewhere.com}%
%  \ead[label=u1,url]{http://www.foo.com}}

 % \thankstext{t2}{Footnote to the first author with the `thankstext' command.}

  \runauthor{Iosif Pinelis and Raymond Molzon}

  \affiliation{Michigan Technological University}

  \address{Department of Mathematical Sciences\\Michigan Technological University\\Hough\-ton, Michigan 49931\\ 
          \printead{e1,e2}}

%  \address{Address of the Third author,\\
%          \printead{e3,u1}}

\end{aug} 

\begin{abstract} 
Uniform and nonuniform Berry--Esseen (BE) bounds of optimal orders on the closeness to normality for general abstract nonlinear statistics are given, which are then used to obtain optimal bounds on the rate of convergence in the delta method for vector statistics. 
Specific applications to Pearson's, non-central Student's and Hotelling's statistics, sphericity test statistics, a regularized canonical correlation, and maximum likelihood estimators (MLEs)  
are given; all these uniform and nonuniform BE bounds appear to be the first known results of these kinds, except for  uniform BE bounds for MLEs. 
When applied to the well-studied case of the central Student statistic, 
our general results compare well with known ones in that case, obtained previously by specialized methods. 
The proofs use a Stein-type method developed by Chen and Shao, a Cram\'er-type of tilt transform,  exponential and Rosenthal-type inequalities for sums of random vectors established by Pinelis, Sakhanenko, and Utev, as well as a number of other, quite recent results motivated by this study.  
The method allows one to obtain bounds with explicit and rather moderate-size constants, at least as far as the uniform bounds are concerned. 
For instance, one has the uniform BE bound $
3.61 %RM 02-21-13 4.06
\E(Y_1^6+Z_1^6)\,(1+\sigma^{-3})/\sqrt n$ for the Pearson sample correlation coefficient based on independent identically distributed random pairs $(Y_1,Z_1),\dots,(Y_n,Z_n)$ with $\E Y_1=\E Z_1=\E Y_1Z_1=0$ and $\E Y_1^2=\E Z_1^2=1$, where $\sigma:=\sqrt{\E Y_1^2Z_1^2}$. 
\end{abstract}

%IP16 \begin{keyword}[class=AMS]
%\kwd[Primary ]{60F05}
%\kwd{60E15}
%\kwd{62F12}
%\kwd[; secondary ]{60E10}
%\kwd{62F03}
%\kwd{62F05}
%\kwd{62G10}
%\kwd{62G20}
%\end{keyword}
%
%\begin{keyword}
%\kwd{Berry-Esseen bound}
%\kwd{delta method}
%\kwd{rates of convergence}
%\kwd{Cram\'er's tilt}
%\kwd{exponential inequalities}
%\kwd{non-central Hotelling's statistic}
%\kwd{non-central Student's statistic}
%\kwd{nonlinear statistics}
%\kwd{Pearson's correlation coefficient}
%\kwd{Rosenthal's inequality}
%\kwd{Stein method}
%\end{keyword}



\end{frontmatter}

\tableofcontents

\section{Introduction}
\label{sec:intro}

Initially, we were interested in studying certain properties of the Pitman asymptotic relative efficiency (ARE) between Pearson's, Kendall's, and Spearman's correlation coefficients. 
As is well known (see e.g.\ \cite{noe55}), the standard expression for the Pitman ARE is applicable when the distributions of the corresponding test statistics are close to normality uniformly over a neighborhood of the null set of distributions. 
Such uniform closeness can usually be provided by Berry-Esseen (BE) type of bounds. 

%IP2 %RM12.28
BE bounds, especially in the special case of linear statistics, constitute a well-established area of research, which originated mainly in work by Scandinavian authors, who were to a large degree concerned with applications in insurance industry and published many of their results on the accuracy of the normal approximation in actuarial journals. For a small sample of recent uses of BE bounds in various areas of sciences and engineering (again for linear statistics), 
see e.g.\ \cite{li14,kuechler-tappe,zeifman-etal,horgan-murphy,kotevski-mitrevski}. 

Kendall's and Spearman's correlation coefficients are instances of $U$-statistics, for which BE bounds are well known; see e.g.\ \cite{kor94}. 
As for the Pearson statistic (say $R$), we have not been able to find a BE bound in the literature. 

This may not be very surprising, considering that an optimal BE bound for the somewhat similar (and, perhaps, somewhat simpler) Student's statistic was obtained only in 1996, by Bentkus and G\"otze \cite{bent96} for independent identically distributed (i.i.d.) random variables (r.v.'s) and by Bentkus, Bloznelis and G\"otze \cite{bbg96} in the general, non-i.i.d.\ case. 
\big(A necessary and sufficient condition, in the i.i.d.\ case, for the Student statistic to be asymptotically standard normal was established only in 1997 by Gin\'e, G\"otze and Mason \cite{ggm}, 
and Hall and Wang \cite{hall04} derive the leading term in the convergence rate in this general setting.\big)  %RM12.28
For more recent developments concerning the Student statistic, see e.g.\ %RM 03.22.13 the 2005 paper by Shao \cite{shao05} and the 2011 preprint by Pinelis \cite{pin11}. 
Shao \cite{shao05} and Pinelis \cite{pin11}. 

Employing such simple and standard tools as a delta-method type linearization together with the Chebyshev and Rosenthal inequalities, we quickly obtained (in the i.i.d.\ case) a uniform bound of the form $\bigO(n^{-1/3})$ for the Pearson statistic. 
Indeed, Pearson's $R$ can be expressed as $f(\bar V)$, a smooth nonlinear function of the sample mean $\bar V=\frac1n\sum_{i=1}^nV_i$, where the $V_i$'s are independent zero-mean random vectors constructed based on the observations of a random sample; cf.\ \eqref{eq:Rdef}. 
A natural approximation to $f(\bar V)-f(0)$, obtained by the delta method, is the linear statistic $L(\bar V)=\sum_{i=1}^n L(\frac1n V_i)$, where $L$ is the linear functional that is the first derivative of $f$ at the origin. 
Since BE bounds for linear statistics is a well-studied subject, we are left with estimating the closeness between $f(\bar V)$ and $L(\bar V)$. 
Assuming $f$ is smooth enough, one will have $|f(\bar V)-L(\bar V)|$ on the order of $\norm{\bar V}^2$, and so, demonstrating the smallness of this remainder term becomes the main problem. 

Using (instead of the mentioned Rosenthal inequality) exponential inequalities for sums of random vectors due to Pinelis and Sakhanenko \cite{pinsak85} or Pinelis \cite{pin94,pin95}, for each $p\in(2,3)$, under the assumption of the finiteness of the $p$th moment of the norm of the $V_i$'s, one can obtain a uniform bound of the form $\bigO(1/n^{p/2-1})$, which is similar to the BE bound for a linear statistic with a comparable moment restriction. 
However, the corresponding constant factor in the $\bigO(1/n^{p/2-1})$ will then explode to infinity as $p\uparrow3$. 
As for $p\ge3$, this method produces bounds of order $\bigO((\ln n)^{3/2}/\sqrt n)$ (for $p=3$) and $\bigO((\ln n)/\sqrt n)$ (for $p>3$), with the extra logarithmic factors. 

While any of these bounds would have sufficed as far as the ARE is concerned, we became interested in obtaining an optimal-rate BE bound for the Pearson statistic. 
Soon after that, we came across the remarkable paper by Chen and Shao \cite{chen07}. 
Suppose that $T$ is any nonlinear statistic and $W$ is any linear one, and let $\De:=T-W$; then make the simple observation that
\begin{equation*}\label{eq:concentr}
-\P(z-|\De|\le W\le z)\le\P(T\le z)-\P(W\le z)\le\P(z\le W\le z+|\De|)
\end{equation*}
for all $z\in\R$. 
Chen and Shao \cite{chen07} offer a Stein-type method to provide relatively simple bounds on the two concentration probabilities in the above inequality, hence bounding the distance between $T$ and $W$; the reader is referred e.g.\ to \cite{barbour05} for illustrations of the elegance and power of Stein's method to a wide array of problems. 
Chen and Shao provided a number of applications of their general results. 

However, in the applications that we desired, such as to Pearson's $R$, it was difficult to deal with $\De=T-W$, as defined above. 
The simple cure applied here was to allow for any $\De\ge|T-W|$, so that, for $T=f(\bar V)$, $W=L(\bar V)$, and smooth enough $f$, the random variable $\De$ could be taken as $\norm{\bar V}^2$ (up to some multiplicative constant). 
This allowed for a BE bound of order $\bigO(1/\sqrt n)$, though under the excessive moment restriction that $\E\norm{V_i}^4<\infty$. 

To obtain a BE bound of the ``optimal'' order $\bigO(1/\sqrt n)$ using only the assumption $\E\norm{V_i}^3<\infty$, we combine the Chen-Shao technique with a Cram\'er-type tilt transform. 
Yet another modification was made by introducing a second level of truncation, to obtain a bound of order $\bigO(1/n^{p/2-1})$ in the case when $\E\norm{V_i}^p<\infty$ for $p\in(2,3)$. 
Thus we obtain our first group of main results (presented in Section~\ref{sec:chen.mod}), on the closeness in distribution of general abstract nonlinear statistics to linear ones. 
These results may be represented by Theorem~\ref{thm:nub}, which provides a ``nonuniform'' upper bound on $|\P(T>z)-\P(W>z)|$ (that is, an upper bound which decreases to $0$ in $|z|$), for a general abstract nonlinear statistic $T$ and a general linear statistic $W$; a ``uniform'' bound on $|\P(T>z)-\P(W>z)|$ is given by Theorem~\ref{thm:ub}. 

The other kind of main results, based on Theorems~\ref{thm:ub} and \ref{thm:nub}, is presented in Section~\ref{sec:f(S)}. 
For instance, Theorem~\ref{thm:f(S).nub} provides a nonuniform upper bound on $|\P(f(S)>z)-\P(L(S)>z)|$  
and thus may be considered as a bound on the rate of convergence in the delta method for vector statistics; it is the latter bound that took more of our time and effort. 
%RM12.28
The reader is referred to \cite{small10} for a rather detailed description of the delta method and its applications; see \cite{kosorok08,romisch05} for a more modern treatment of the delta method applied to infinite-dimensional random vectors. 

Finally, as applications of the delta-method bounds given in Section~\ref{sec:f(S)}, we present (in Section~\ref{sec:apps}) %IP16 
uniform and nonuniform BE-type bounds for the Pearson statistic, %RM12.28 as well as for 
the noncentral Student and Hotelling 
statistics, various statistics commonly used in testing hypotheses about a population covariance matrix, the largest eigenvalue of a certain linear operator on an infinite-dimensional Hilbert space, and %IP16 the maximum likelihood estimator. %RM12.28 ones. 
maximum likelihood estimators. 
No such BE bounds appear to be previously known, %IP16 
except for uniform BE bounds for MLEs. 
As for the known BE bounds for the central Student statistic (obtained by specialized methods, targeting this specific statistic), it turns out that our bounds (even though based on the mentioned results for general nonlinear statistics) compare well with the former ones.  

%RM12.28
Our general BE bounds in the multivariate delta method can of course be used in applications other than the ones considered here; we mention %IP2 several 
a number of other potential applications in Subsections~\ref{subsec:cov} and \ref{subsec:pca}. 
In fact, a result from an earlier arXiv version of this paper, similar to %IP16 Corollary
Theorem~\ref{cor:iid}, was already used in \cite{gamboa13}. Of course, our results cannot perfectly cover the entire variety of uses of the delta method; they may require modification or use of different ideas; see e.g.\ \cite[pages~1198 and 1211]{wass13}. 

To obtain the delta-method bounds stated in Sections~\ref{sec:f(S)} and their applications presented in Section~\ref{sec:apps}, we use a number of previously known results, including precise exponential and Rosenthal-type bounds developed by Pinelis, Sakhanenko, and Utev \cite{pin80,pinsak85,pin86,pinutev86,pin94} and also a number of other known results due to Bennett \cite{bennett62}, Hoeffding~\cite{hoeff63}, de Acosta and Samur \cite{acosta79}, Michel \cite{michel81}, 
%RM 03.20.13 Ibragimov and Sharakhmetov \cite{ibrag97}, 
and Shevtsova \cite{shev11}.  
There we also use the recent results developed in  \cite{%IP 03.11.13
student-mono-v2,pinBH,pin09.winsor,pin11tilt,pin11,pin11.BE,re-center,pin12.tilt.symm
,pin12-2smooth,pin13.de,pin13.rosen%RM 03.20.13
}. 

As for the requirement that the observations be identically distributed, it may (and will) be dispensed in general; that is, $\bar V$ will in general be replaced by a sum $S$ of independent but not necessarily identically distributed random vectors.  

The paper is organized as follows. 
\begin{enumerate}[{-}]
\item
In Section~\ref{sec:chen.mod}, we state and discuss the mentioned upper bounds on $|\P(T>z)-\P(W>z)|$ for general $T$ and $W$. 
\item
In Section~\ref{sec:f(S)}, the mentioned Theorem~\ref{thm:f(S).nub} and other results are stated, providing general bounds on the rate of convergence in the vector delta method, that is, bounds on $|\P(f(S)>z)-\P(L(S)>z)|$. 
\item
Applications to several commonly used statistics, namely the non-central Student $T$, the Pearson $R$, %RM12.28 and 
the non-central Hotelling $T^2$, 
various test statistics constructed from a sample covariance matrix, the largest eigenvalue of a certain linear operator, and %IP16 the maximum likelihood estimator, %RM12.28 
maximum likelihood estimators 
are stated in Section~\ref{sec:apps}. 
\item Proofs of results from Sections~\ref{sec:chen.mod} and \ref{sec:f(S)}, as well as selected results from Section~\ref{sec:apps}, are deferred to Section~\ref{sec:proofs}.
\end{enumerate}
Certain results and proofs are relegated to appendices. 
\begin{enumerate}[{-}]
\item The statement and proof of an explicit (and quite complicated in appearance) nonuniform bound on the distance to normality of $f(\bar{V})$ in an i.i.d.\ setting is provided in Appendix~\ref{app:nonunif}. 
\item The nonuniform bounds developed in this paper are valid under the restriction that $z=\bigO(\sqrt{n})$ (in the i.i.d.\ case); in Appendix~\ref{app:z<=sqrt.n} we prove that this restriction cannot generally be discarded or even relaxed.
\item Appendix~\ref{app:num.proofs} contains the proofs of bounds from Section~\ref{sec:apps} which, for practical purposes, make the use of a computer algebra system (CAS) preferable. 
\item In Appendix~\ref{app:z}, we discuss the potential application of the bounds presented in Section~\ref{sec:f(S)} to the Fisher $z$-transform of the Pearson statistic. 
%RM12.28
\item In Appendix~\ref{sec:compact} we provide a short, self-contained proof of the compactness of the covariance operator for a random vector taking values in a separable Hilbert space and possessing a finite second moment; this is used in one of our applications on the principal component of a certain linear operator.
\item In Appendix~\ref{sec:inf-dim} we outline the proof of the existence of the spectral decomposition for the covariance operator of a random vector taking values in an infinite-dimensional separable Hilbert space. 
\end{enumerate}

\section{Approximation of the distributions of general abstract nonlinear statistics by the distributions of linear ones}
\label{sec:chen.mod}

Let $X_1,\dotsc,X_n$ be independent r.v.'s with values in some measurable space $\XX$, and let $T\colon\XX^n\to\R$ be a Borel-measurable function. 
For brevity, let $T$ also stand for $T(X_1,\dotsc,X_n)$, the statistic of the random sample $(X_i)_{i=1}^n$. 
Further let
\begin{equation}\label{eq:xi,eta}
 \xi_i:=g_i(X_i)
 \quad\text{and}\quad
 \eta_i:=h_i(X_i)
\end{equation}
for $i=1,\dotsc,n$, where $g_i\colon\XX\to\R$ and $h_i\colon\XX\to\R$ are Borel-measurable functions.
Assume that 
\begin{equation}\label{eq:=1}
 \E\xi_i=0\text{ for all $i=1,\dotsc,n$, and }\tsum_{i=1}^n\E\xi_i^2=1.
\end{equation}
Consider the linear statistic
\begin{equation}\label{eq:W}
 W:=\tsum\lims_{i=1}^n\xi_i.
\end{equation}
Further, take an arbitrary $\cc\in(0,1)$ and let $\de$ be any real number such that
\begin{equation}\label{eq:de}
 \tsum\lims_{i=1}^n\E|\xi_i|\bigp{\de\wedge|\xi_i|}\ge\cc;
\end{equation}
note that such a number $\de$ always exists (because the limit of the left-hand side of \eqref{eq:de} as $\de\uparrow\infty$ is 1). Necessarily, $\de>0$. 

\begin{thm}\label{thm:ub}
Let $\De$ be any r.v.\ such that $\abs{\De}\ge\abs{T-W}$ almost surely (a.s.), and for each $i=1,\dotsc,n$, let $\De_i$ be any r.v.\ such that $X_i$ and $(\De_i,W-\xi_i)$ are independent. Take any real number $w>0$, and let $\bar\De$ be any r.v.\ such that
\begin{equation}\label{eq:barDe}
 \bar\De=\De\text{ a.s.\ on the event }\Bigb{\max_{1\le i\le n}\eta_i\le w}.
\end{equation}
Then for all $z\in\R$
\begin{equation}\label{eq:ub}
 \bigabs{\P(T>z)-\P(W>z)}
 \le\frac1{2\cc}
  \Bigp{ 4\de + \E\bigabs{W\bar\De} + \tsum\lims_{i=1}^n\E\bigabs{\xi_i(\bar\De-\De_i)} }
  +\bigprob{\max\nolim_i\eta_i>w},
\end{equation}
where $\de$ %RM 03.22.13 satisfies \eqref{eq:de}.
is any number satisfying \eqref{eq:de}. 
\end{thm}

\begin{remark}\label{re:chen.improvement} 
Sacrificing some simplicity in appearance, one can improve the bound in \eqref{eq:ub} by replacing the term $4\de$ there with 
\begin{equation}\label{eq:ub.imp} 
 2\de+\frac{\de^2}{\cc}+2\de\sqrt{\frac1{2\cc}
  \Bigp{2\de+\frac{\de^2}{2\cc}+\E\bigabs{W\bar\De}+\tsum\nolim_i\E\bigabs{\xi_i(\bar\De-\De_i)}}};
\end{equation}
the validity of \eqref{eq:ub} after such a replacement will be shown in the proof of Theorem~\ref{thm:ub}. 
Evidently, when the upper bound in \eqref{eq:ub} is small, the expression \eqref{eq:ub.imp} will behave like $2\de$, in place of $4\de$ in \eqref{eq:ub}. 
\end{remark}

\begin{remark}\label{re:ub} 
Inequality \eqref{eq:ub} above is a rather straightforward generalization of the result (2.3) in Theorem~2.1 by Chen and Shao \cite{chen07}. 
The modifications we have made are as follows. 
First, $\De$ was defined in \cite{chen07} as simply equal to $T-W$. Then, in the applications given in our present paper, it becomes problematic to bound the term $\E|\xi_i(T-W-\De_i)|$ \big(which would arise in place of the term $\E|\xi_i(\bar\De-\De_i)|$ in \eqref{eq:ub}\big). 
Using the more general condition $\abs{\De}\ge\abs{T-W}$ instead of $\De=T-W$ allows one to choose a possibly larger $\De$ so that $\E|\xi_i(\De-\De_i)|$ be more amenable to analysis. 
However, if that $\De$ should happen to be ``too large,'' our second generalization allows one to truncate $\De$ to within acceptable constraints by using the additional truncation level $w$, as well as $\bar\De$ and $\P(\max_i\eta_i>w)$. 
The third difference is that in \cite{chen07} $c_*$ was chosen to be $\frac12$; the more general condition $c_*\in(0,1)$ results in improved explicit constants in the applications. 
\end{remark}

Before stating the ``nonuniform'' counterpart of Theorem~\ref{thm:ub}, let us introduce some notation. 
For any real $a$ and $b$, let $a\wedge b$ and $a\vee b$ denote the minimum and maximum, respectively, of $a$ and $b$; use also the notation $a_+:=a\vee0$. 
For any real-valued r.v.\ $\xi$ and any $p\in[1,\infty)$, let $\norm{\xi}_p:=\E^{1/p}|\xi|^p$. 
For the $\xi_i$'s as in \eqref{eq:xi,eta}, also let
\begin{equation}\label{eq:si_p}
 \si_p:=\biggp{\tsum_{i=1}^n\norm{\xi_i}_p^p}^{1/p} = \biggp{\tsum_{i=1}^n\E|\xi_i|^p}^{1/p}.
\end{equation}

In proving, and even stating, the forthcoming results of the current paper, we will need several tools for estimating moments and tail probabilities. 
Let here $\zeta:=(\zeta_1,\dotsc,\zeta_n)$, where $\zeta_1,\dotsc,\zeta_n$ are independent real-valued r.v.'s, $S:=\sum_i\zeta_i$, and 
\begin{equation}\label{eq:G}
 G_\zeta(z):=\tsum\lims_{i=1}^n\P(\zeta_i>z)
 \quad\text{for all}\quad z\in\R.  
\end{equation}

If the $\zeta_i$'s are zero-mean, then for each real $\al\ge2$ there exist positive constants $\AA_\R(\al)$ and $\BB_\R(\al)$, depending only on $\al$, such that 
\begin{equation}\label{eq:rosen1}   
%IP 03.10.13 here and elsewhere, replaced ^\al(\al) by (\al)^\al
%  \norm{S}_\al^\al\le\AA_\R^\al(\al)\tsumnl_i\norm{\zeta_i}_\al^\al+\BB_\R^\al(\al)\Bigp{\tsumnl_i\norm{\zeta_i}_2^2}^{\al/2}.  
\norm{S}_\al^\al\le\AA_\R(\al)^\al\tsumnl_i\norm{\zeta_i}_\al^\al+\BB_\R(\al)^\al\Bigp{\tsumnl_i\norm{\zeta_i}_2^2}^{\al/2}.  
\end{equation}
Such a result will be referred to in this paper as a Rosenthal-type inequality, since it was first obtained by Rosenthal in \cite[Theorem~3]{rosenthal70}; however, the constants there were too large, 
%IP 03.10.13
as e.g.\ compared with ones in \cite{pin12-2smooth,pin13.rosen}; cf.\ also \eqref{eq:A1,B1}. 
%RM 02-21-13 %% added remainder of paragraph concerning non-central Rosenthal-type inequality %%
If the $\zeta_i$'s are not centered, a similar inequality can be obtained. Namely,  
\begin{equation}\label{eq:rosen.noncentral}
 \norm{S-\E S}_\al^\al\le
 \AA_{\R,\nc}(\al)^\al%RM 03.20.13 \AA_{\R,\nc}^\al
 \tsumnl_i\norm{\zeta_i}_\al^\al+
 \BB_{\R,\nc}(\al)^\al%RM 03.20.13 \BB_{\R,\nc}^\al(\al)
 \Bigp{\tsumnl_i\norm{\zeta_i}_2^2}^{\al/2}
\end{equation}
for any $\al\ge2$ and some positive constants $\AA_{\R,\nc}(\al)$ and $\BB_{\R,\nc}(\al)$; see e.g.\ \cite[Corollary~4]{pin12-2smooth}. 

Next, we shall need upper bounds on the tail probabilities. 
Suppose now that $G_\zeta(y)=0$ for some $y>0$, i.e.\ each of the $\zeta_i$'s is bounded %RM 03.22.13 
from above a.s.\ by $y$. 
Then \cite[Theorem~2]{pinutev86} implies that for any $\la\ge0$ 
\begin{equation}\label{eq:PUexp}
 \E\exp\bigb{\la(S-m)}\le
  \PUexp(\la,y,B,\vp):=\exp\Bigb{\frac{\la^2}{2}\,B^2(1-\vp)+\frac{e^{\la y}-1-\la y}{y^2}\,B^2\vp},
\end{equation}
where $B=(\sum_i\E\zeta_i^2)^{1/2}<\infty$, $m=\E S$, $\vp=\sum_i\E(\zeta_i)_+^p/(B^2y^{p-2})\in(0,1)$, and $p\in[2,3]$. 
Further, an application of the Markov inequality and \eqref{eq:PUexp} yield
\begin{equation}\label{eq:PU}
 \P(S\ge x)\le\PU(x,y,B,m,\vp):=\inf_{\la\ge0}e^{-\la(x-m)}\PUexp(\la,y,B,\vp)\quad\text{for any $x\in\R$.}
\end{equation}
As functions of the real numbers $\la\ge0$, $y>0$, $B>0$, $\vp\in(0,1)$, $x$, and $m$, the bounds $\PUexp$ and $\PU$ possess certain monotonicity properties: $\PU$ is clearly nondecreasing in $m\in\R$, and from the inequality $e^t-1-t-t^2/2\ge0$ for all $t\ge0$ it follows that 
\begin{equation}\label{eq:incr}
 \text{$\PUexp$, and hence $\PU$, are nondecreasing in $B$ and in $\vp$.}
\end{equation}
Thus, we see the inequalities in \eqref{eq:PUexp} and \eqref{eq:PU} hold under the relaxed (and more convenient) conditions
\begin{equation}\label{eq:PUcond}
 \tsumnl_i\P(\zeta_i>y)=0,\quad
 \bigp{\tsumnl_i\E\zeta_i^2}^{1/2}\le B,\quad
 \E S\le m,\quad\text{and}\quad
 \frac{\sum_i\E(\zeta_i)_+^p}{B^2y^{p-2}}\le\vp\in(0,1];
\end{equation}
that \eqref{eq:PUexp} is true when $\vp=1$ is a result by Bennett \cite{bennett62} and Hoeffding \cite{hoeff63}, and we let $\BHexp(\la,y,B):=\PUexp(\la,y,B,1)$ and $\BH(x,y,B,m):=\inf_{\la>0}e^{-\la(x-m)}\BHexp(\la,y,B)$.  
The bounds $\PUexp$ and $\PU$ can be much less than $\BHexp$ and $\BH$, respectively, when $\vp$ is significantly less than 1. 
Expressions for $\PU$ are given in \cite[Corollary~1]{pinutev86} and \cite[Proposition~3.1]{pinBH}, and Lemma~\ref{lem:PU} will present these in a manner useful for the applications considered in the present paper. 
We remark here that an exponential bound on $\E e^{\la(S-m)}$ (and hence also $\P(S\ge x)$) which incorporates the moments $\E(\zeta_i)_+^p$ with  $p>3$ is stated in \cite[Theorem~6]{pinutev86}, though the resulting expression is considerably more complicated in appearance than the bound in \eqref{eq:PUexp}. 

In the proof of Theorem~\ref{thm:nub} stated below, we shall also have cause to find a lower bound for the exponential moment of a Winsorized r.v. 
Particularly, suppose that $\xi$ is a zero-mean r.v.\ with $\sqrt{\E\xi^2}\le B$ for some $B\in(0,\infty)$. 
Then for any $c>0$, \cite[Theorem~2.1]{pin09.winsor} states that
\begin{equation}\label{eq:L_W}
 \E\exp\bigb{c\bigp{1\wedge\xi}}
  \ge L_{W;\,c,B}:=\frac{a_{c,B}^2e^c+B^2e^{-ca_{c,B}}}{a_{c,B}^2+B^2}, 
\end{equation}
where $a_{c,B}$ is the unique positive root of the function $a\mapsto\frac ac\,\bigp{2(e^{c+ac}-1)-ac}-B^2$. 
In fact, as shown in \cite{pin09.winsor}, $L_{W;\,c,B}$ is the exact lower bound on $\E\exp\bigb{c\bigp{1\wedge\xi}}$ over all zero-mean r.v.'s $\xi$ with $\sqrt{\E\xi^2}\le B$, and hence $L_{W;\,c,B}$ is nonincreasing in $B\in(0,\infty)$. 

\begin{thm}\label{thm:nub}
Let $\De$ be any r.v.\ such that $\abs{\De}\ge\abs{T-W}$ a.s.
For each $i=1,\dotsc,n$, let $\De_i$ be any r.v.\ such that $X_i$ and $(\De_i,(X_j\colon j\neq i))$ are independent, and assume that the mentioned Borel-measurable functions $g_i$ and $h_i$ are such that $g_i\le h_i$, so that $\xi_i\le\eta_i$.
Take any real $p\in[2,3]$ and let $q:=\frac p{p-1}$, so that $\frac1p+\frac1q=1$; also take any real numbers 
\begin{equation}\label{eq:pi's}
 \cc\in(0,1),\ \th>0,\ w>0,\ \de_0\in(0,w],\ %RM 03.22.13
 \text{and }\pi_1>0,\ \pi_2>0,\ \text{and}\ \pi_3>0\text{ such that }\pi_1+\pi_2+\pi_3=1.
\end{equation}
Then for all $z\ge0$
\begin{equation}\label{eq:nub}
 \bigabs{\hat\P(T>z)-\hat\P(W>z)}\le\ga_z+\tau e^{-(1-\pi_1)z/\th},
\end{equation}
where
\begin{gather}
\label{eq:Phat}
 \hat\P(E):=\P\bigp{E\cap\bigb{\abs{\De}\le\pi_1z}}\text{ for any event $E$,}
\\
\label{eq:ga_z}
 \ga_z:=G_\xi\bigp{\pi_2z}+\tsum_{i=1}^n\bigprob{W-\xi_i\ge\pi_3z}\bigprob{\eta_i>w},
\\
\label{eq:tau}
 \tau:=c_1\tsum_{i=1}^n\bignorm{\xi_i}_p\bignorm{\bar\De-\De_i}_q+c_2\bignorm{\bar\De}_q+c_3\de,
\end{gather}
$G_\xi$ is defined by \eqref{eq:G}, $\bar\De$ is any r.v.\ satisfying \eqref{eq:barDe}, $\de$ is any number such that \eqref{eq:de} holds,
\begin{gather}
 \label{eq:c1} c_1:=\tfrac1{\cc}\PUexp\bigp{\tfrac{p}{\th},w,\tfrac1{\sqrt{p}},\vp_1}e^{\de_0/\th},\\
 \label{eq:c2} c_2:=c_1\Bigp{
 \AA_{\R,\nc}(p) %RM 02-21-13 \AA_\R(p)
 \bigp{
 %RM 02-21-13 1.32
  a_1e^{pw/\th}}^{1/p}\si_p+
  \BB_{\R,\nc}(p) %RM 02-21-13 \BB_\R(p)
  \bigp{a_1e^{pw/\th}}^{1/2}+\bigp{e^{pw/\th}-1}/w},\\
 \label{eq:c3} c_3:=\Bigp{2c_2+\tfrac1{\cc}\sqrt2\PUexp\bigp{\tfrac{2}{\th},w,\tfrac1{\sqrt2},\vp_1}}\vee\Bigp{\tfrac1{\de_0}\,\PUexp\bigp{\tfrac1\th,w,1,\vp_1}},\\
\label{eq:vp1}
 \vp_1:=\frac{\si_{p}^{p}}{w^{p-2}}\wedge1,\\
\label{eq:a1} 
 a_1:=1/L_{W;\,pw/\th,\,\max_i\norm{\xi_i}_2/w}.
\end{gather}
\end{thm}

\begin{remark}\label{re:Phat}
We shall use \eqref{eq:nub} in conjunction with the obvious inequality
\begin{equation*}
 \bigabs{\P(T>z)-\P(W>z)}\le\P\bigp{\abs{\De}>\pi_1\abs{z}}+\bigabs{\hat\P(T>z)-\hat\P(W>z)}.
\end{equation*}
Thus, the use of the measure $\hat\P$ in \eqref{eq:nub} will allow us to avoid a ``double counting'' of the probability $\P(|\De|>\pi_1|z|)$ when Theorem \ref{thm:nub} is used to obtain Theorem~\ref{thm:f(S).nub}. 
\end{remark}

\begin{remark}\label{re:z<0} 
The bound \eqref{eq:nub} (as well as other nonuniform bounds presented later in this paper) is stated only for $z\ge0$, which allows for one-tail expressions $G_\xi\bigp{\pi_2z}$ and $\bigprob{W-\xi_i\ge\pi_3z}$ to be used in \eqref{eq:ga_z}. 
In order to obtain the corresponding bound for $z<0$, all that is needed is to replace $T$ and $g_i$ with $-T$ and $-g_i$, respectively, where the $g_i$'s are as in \eqref{eq:xi,eta}. 
\end{remark}

%RM 02-21-13
\begin{remark}\label{re:de}
A simple modification of \cite[Remark~2.1]{chen07} shows us that \eqref{eq:de} is satisfied when
\begin{equation}\label{eq:c4}
 \de=\Bigp{\frac{(p-2)^{p-2}}{(p-1)^{p-1}(1-\cc)}}^{1/(p-2)}\si_p^{p/(p-2)} %RM 02-21-13
%  \de=c_4\si_p^{\tq},\quad\text{where}\quad
%  c_4:=\Bigp{\frac{(p-2)^{p-2}}{(p-1)^{p-1}(1-\cc)}}^{1/(p-2)},\quad\text{and}\quad\tq:=\frac{p}{p-2}
\end{equation}
for any $p\in(2,3]$. 
%IP 03.10.13
% An even better (i.e.\ smaller) choice for $\de$ exists when we assume that $\si_3<\infty$: according to \cite[Theorem~1]{pin13.de}, we have
% \eqref{eq:de} satisfied when
A smaller choice of $\de$ exists for $p=3$: according to \cite[Theorem~1]{pin13.de}, \eqref{eq:de} holds if $\si_3<\infty$ and 
\begin{equation}\label{eq:de-imp}
 \de=\begin{cases}\cc\si_3^3&\text{if }0<\cc\le\frac12,\\[5pt]\displaystyle\frac{\si_3^3-(2\cc-1)^2/\si_1}{4(1-\cc)}&\text{if }\frac12\le\cc<1.\end{cases}
\end{equation}
\end{remark}

\begin{remark}\label{re:tau}
It is easy to see that the expressions $c_1$, $c_2$, and $c_3$ in \eqref{eq:c1}--\eqref{eq:c3} can be bounded by finite positive constants depending only on the values of the parameters $p$, $\cc$, $\th$, $w$, and $\de_0$ (and not on the distributions of the $X_i$'s).
This follows because $\PUexp$ is nondecreasing in $\vp$ (recall \eqref{eq:incr}) and $a_1\le1/L_{W;pw/\th,1/w}$ (since $\max_i\norm{\xi_i}_2\le\norm{W}_2=1$ and $L_{W;\,c,B}$ is nonincreasing in $B$). 
Thus, one may refer to $c_1$, $c_2$, and $c_3$ as \emph{pre-constants}. 
\end{remark}

\begin{remark}\label{re:symm}
If we add the assumption that the $\xi_i$'s are all symmetric(ally distributed) to the %RM 03.22.13 hypotheses 
assumptions of Theorem~\ref{thm:nub}, then, according to the main result of \cite{pin12.tilt.symm}, $(e^{pw/\th}-1)/w$ in \eqref{eq:c2} may be replaced by the smaller quantity $\sinh(pw/\th)/w$. 
This sharpening of the inequality \eqref{eq:nub} allows for smaller absolute constants to be obtained in applications of Theorem~\ref{thm:nub}; cf.\ the nonuniform bound for the self-normalized sum in Corollary~\ref{cor:centralT.nub} and Remark~\ref{re:symm.T1}.
\end{remark}

For $p=2$, the result of Theorem~\ref{thm:nub} is similar to that by Chen and Shao \cite[Theorem 2.2]{chen07}. 
The bound given by \eqref{eq:nub} turns out to be more precise in the applications given in this paper. 
In particular, it allows one to weaken conditions on moments. 
Indeed, in Theorem~\ref{thm:f(S).nub} one will have $|\bar\De|$ on the order of $\norm{S}^2$ and $|\bar\De-\De_i|$ on the order of $\norm{X_i}^2+\norm{X_i}\,\norm{S-X_i}$, where $S:=\sum_{i=1}^nX_i$ and the $X_i$'s are independent random vectors. 
So, using Theorem~\ref{thm:nub} with $p=3$ (and hence $q=\frac32$) in order to obtain a bound of the classical form $\bigO(\frac1{\sqrt n(|z|+1)^3})$, one will need only the third moments of $\norm{X_i}$ to be finite. 
On the other hand, using \eqref{eq:tau} with $p=2$ to get the same kind of bound would require the finiteness of the fourth moments of $\norm{X_i}$. 

Expressions in Theorem~\ref{thm:nub} are complicated, especially the ones for $c_1$, $c_2$, and $c_3$. However, this may be considered as just another instance of the usual trade-off between accuracy and complexity of the bounds. 

Bounds \eqref{eq:ub} and \eqref{eq:nub} on the closeness of the distribution of the linear approximation $W$ to that of the original statistic $T$ are to be complemented by any number of well-known BE-type bounds on the closeness of the distribution of the linear statistic $W$ to the standard normal distribution; the reader may be referred to Petrov's monograph \cite[Chapter~V]{pet75} or the paper \cite{pin11.BE}. 
For the linear statistic $W$ as in \eqref{eq:W} with i.i.d.\ $\xi_1,\dots,\xi_n$ as in \eqref{eq:=1}, results due to Shevtsova \cite{shev11} and Michel \cite{michel81} imply 
\begin{equation}\label{eq:tyurin,michel}
 \bigabs{\P(W\le z)-\Phi(z)}\le %RM12.21.12 n\norm{\xi_1}_3^3\Bigp{0.4748\wedge\frac{30.2211}{\abs{z}^3+1}}.
  n\Bigp{0.33554\bigp{\norm{\xi_1}_3^3+0.415\norm{\xi_1}_2^3}\wedge\frac{30.2211\norm{\xi_1}_3^3}{\abs{z}^3+1}}.
\end{equation}

\section{Berry-Esseen bounds for smooth nonlinear functions of sums of independent random vectors}
\label{sec:f(S)}

In this section, we shall state applications of results of Section~\ref{sec:chen.mod}.
Assume from hereon that $(\XX,\|\cdot\|)$ is a separable Banach space of type 2; for a definition and properties of such spaces, see e.g.\ \cite{hoff76,pin86}. 
Let $X_1,\dotsc,X_n$ be independent random vectors in $\XX$ with $\E X_i=0$ for $i=1,\dotsc,n$, and also let
\begin{gather}
\notag
 S:=\tsum_{i=1}^nX_i,
\\
\notag
 \|X\|_p:=\E^{1/p}\|X\|^p,
\\
\label{eq:s_al}
 s_p:=\biggp{\tsum_{i=1}^n\norm{X_i}_p^p}^{1/p}=\biggp{\tsum_{i=1}^n\E\norm{X_i}^p}^{1/p},
\\
\label{eq:G_X}
 G_X(z):=\tsum_{i=1}^n\bigprob{\norm{X_i}>z},
\end{gather}
for any $p\ge1$ and $z\ge0$; compare \eqref{eq:s_al} and \eqref{eq:G_X} to \eqref{eq:si_p} and \eqref{eq:G}, respectively. 

Note that the results of \cite[Theorem~1]{pinsak85} (see also the remark in \cite[p.~343]{pinutev86}) may be used to derive bounds analogous to those given in \eqref{eq:PUexp} and \eqref{eq:PU} when the $\zeta_i$'s take values in a separable Banach space. Particularly, 
\begin{equation}\label{eq:PU,vector}
 \text{\eqref{eq:PU} and \eqref{eq:PUexp} hold under \eqref{eq:PUcond} when $S$ and $\zeta_i$ are replaced by $\norm{S}$ and $\norm{X_i}$, respectively.} 
\end{equation}

Since $\XX$ is of type 2 and the $X_i$'s are zero-mean, there exists a constant $D:=D(\XX)\in(0,\infty)$ such that
\begin{equation}\label{eq:D}
 \norm{S}_2\le D s_2.
\end{equation}
We shall assume that $D$ is chosen to be minimal with respect to this property; so, $D=1$ with the equality in \eqref{eq:D} whenever $\XX$ is a Hilbert space. 
By \cite[Theorem~2]{pin86} %IP 03.23.13
or \cite{pin12-2smooth}, one also has the Rosenthal-type inequality
\begin{equation}\label{eq:rosen}
%IP 03.11.13 \norm{S}_\al^\al\le \AA_\XX^\al(\al)s_\al^\al+\BB_\XX^\al(\al)s_2^\al
\norm{S}_\al^\al\le \AA_\XX(\al)^\al s_\al^\al+\BB_\XX(\al)^\al s_2^\al
\end{equation}
for any $\al\ge2$ and some pair of constants $(\AA_\XX(\al),\BB_\XX(\al))$; note that \eqref{eq:rosen} generalizes \eqref{eq:rosen1}. 
%RM 02-21-13 %% moved Rosenthal constants for p=3 to discussion introducing Section 4 %%
% In the particular case when $\XX$ is a Hilbert space and $\al=3$, %IP12.13.12\cite[Theorem~1]{pin80}
% inequality \cite[(12)]{pin12-2smooth} 
% allows us to use the values
% \begin{equation}\label{eq:A1,B1}
%  \AA_\XX(3)=1\quad\text{and}\quad \BB_\XX(3)=%IP12.13.12 3
%  2^{1/3};
% \end{equation}
% cf.\ also \cite{pin80}. 

\begin{remark}\label{rem:2-smooth}
The results of this section hold for vector martingales taking values in a 2-smooth separable Banach space; in such a case, one can apply results of \cite{pin94} instead of the ones of \cite{pinsak85} used in the present paper. 
By \cite{hoff76,pin94}, every 2-smooth Banach space is of type 2. 
It is known that $L^p$ spaces are $2$-smooth, and hence of type 2, for all $p\ge2$ \cite[Proposition~2.1]{pin94}. 
\end{remark}

Let next $f\colon\XX\to\R$ be a Borel-measurable functional with $f(0)=0$, satisfying the following smoothness condition: there exist %RM 03.22.13 $\ep>0$, $\Mf>0$, 
$\ep\in(0,\infty)$, $\Mf\in(0,\infty)$, and a nonzero continuous linear functional $L\colon\XX\to\R$ such that
\begin{align}\label{eq:smooth}
 \bigabs{f(x)-L(x)}\le\frac \Mf2\,\|x\|^2\text{ for all $x\in\XX$ with }\norm{x}\le\ep;
\end{align}
thus, $L$ necessarily coincides with the first Fr\'echet derivative, $f'(0)$, of the function $f$ at $0$. 
Moreover, for the smoothness condition \eqref{eq:smooth} to hold, it is enough that the second derivative $f''(x)$ exist and be bounded (in the operator norm) by $\Mf$ over all $x\in\XX$ with $\|x\|\le\ep$. 

%RM12.26 
\begin{remark}\label{re:compos}
A fact useful in applications is that the smoothness condition \eqref{eq:smooth} continues to hold over compositions of functions. 
Specifically, suppose that $\XX$, $\YY$, and $\ZZ$ are separable Banach spaces with respective norms $\norm{\cdot}_\XX$, $\norm{\cdot}_\YY$, and $\norm{\cdot}_\ZZ$, and let $h\colon\XX\to\YY$ and $g\colon\YY\to\ZZ$ be functions such that 
\begin{equation}\label{eq:h 2smooth}
 \norm{h(x)-L_h(x)}_\YY\le\tfrac{M_h}{2}\,\norm{x}_\XX^2\text{ for all $x\in\XX$ with $\norm{x}_\XX\le\ep_h$}
\end{equation}
and
\begin{equation}\label{eq:g 2smooth}
 \norm{g(y)-L_g(y)}_\ZZ\le\tfrac{M_g}{2}\,\norm{y}_\YY^2\text{ for all $y\in\YY$ with $\norm{y}_\YY\le\ep_g$}
\end{equation}
for some continuous linear operators $L_h\colon\XX\to\YY$, $L_g\colon\YY\to\ZZ$ and positive real numbers $M_h$, $\ep_h$, $M_g$, $\ep_g$. 
Then the composition $f:=g\circ h\colon\XX\to\ZZ$ satisfies \eqref{eq:smooth} with $\ZZ$ in place of $\R$, $L=L_g\circ L_h$, $\Mf=M_h\norm{L_g}+M_gm_h^2$, $m_h:=\norm{L_h}+M_h\ep_h/2$, and $\ep=\ep_h$, provided that $\ep_h$ is chosen small enough to ensure $m_h\ep_h\le\ep_g$. Such a statement can of course be generalized to the composition of any finite number of functions. 
We shall prove this assertion in Section~\ref{sec:proofs}. 
\end{remark}

%RM 02-21-13
Given a function $f$ which satisfies the smoothness condition \eqref{eq:smooth}, let us define
\begin{equation}\label{eq:si}
 \si:=\norm{L(S)}_2=\Bigp{\tsumnl_i\norm{L(X_i)}_2^2}^{1/2};
\end{equation}
further assume that $\si\in(0,\infty)$. %IP 03.10.13 Then let $g_i$ be defined by $x\mapsto L(x)/\si$ for each $i=1,\dotsc,n$, so that
In \eqref{eq:xi,eta}, take $g_i(x)\equiv L(x)/\si$ for each $i=1,\dotsc,n$, so that 
\begin{equation}\label{eq:xi=L(X)/si}
 \xi_i=\frac{L(X_i)}{\si}; 
\end{equation}
%IP 03.10.13 follows from \eqref{eq:xi,eta}; 
it is clear then that \eqref{eq:=1} is satisfied, and $W=L(S)/\si$ according to \eqref{eq:W}. 

The following bound for the distribution of $f(S)$ may still look rather abstract and complicated. However, especially in such applications to specific statistics as the ones presented in Corollaries~\ref{cor:centralT2} and \ref{cor:Runif}, it leads to comparatively simple BE type bounds of a ``correct'' order of magnitude and with explicit numerical constants of rather moderate sizes. 

\begin{thm}\label{thm:f(S).ub}
Let $f\colon\XX\to\R$ satisfy \eqref{eq:smooth}, %IP 03.23.13
and let $X_1,\dotsc,X_n$ be independent zero-mean random vectors in $\XX$. %IP 03.10.13 was already assumed after {eq:si}
% , and assume that $\si\in(0,\infty)$. %RM 02-21-13
%RM 02-21-13 % \begin{equation}\label{eq:si}
%  \si:=\norm{L(S)}_2=\Bigp{\tsum\nolim_i\norm{L(X_i)}_2^2}^{1/2}\in(0,\infty).
% \end{equation}
Further, take any $p\in(2,3]$, $\cc\in(0,1)$, $w>0$, and let $q:=\frac{p}{p-1}$, so that $\frac1p+\frac1q=1$. 
Then for all $z\in\R$ 
\begin{equation}
\label{eq:f(S).ub}
\begin{split}
 \biggabs{\Bigprob{\frac{f(S)}{\si}>z}-\Bigprob{\frac{L(S)}{\si}>z}}
  \le\bigprob{\norm{S}>\ep}
   +\frac{4 
   \de %RM 02-21-13 c_4\si_p^\tq
   +\bigp{\AA_\R(p)\si_p+\BB_\R(p)}\uu+\si_p\vv}{2\cc}
   +G_\eta(w),
\end{split}
\end{equation}
where 
$\de$ is any number satisfying \eqref{eq:de}, %RM 02-21-13
%RM 02-21-13 $c_4$ and $\tq$ are as defined in \eqref{eq:c4}, %RM 12.21.12
% \begin{equation}\label{eq:c4}
%  c_4:=\Bigp{\frac{(p-2)^{p-2}}{(p-1)^{p-1}(1-\cc)}}^{1/(p-2)},\quad\tq:=\frac{p}{p-2},
% \end{equation}
$\si_p$ and $G_\eta$ are as in \eqref{eq:si_p} and \eqref{eq:G} with
\begin{equation}\label{eq:eta.ub}
 %RM 02-21-13 \xi_i=\frac{L(X_i)}{\si},\quad
 \eta_i=\frac{\norm L\norm{X_i}}{\si}\,\ind{2<p<3},
\end{equation}
\begin{equation}\label{eq:u}
 \uu:=\frac{\Mf\si}{2\norm{L}^2}\times
 \begin{cases}
  \bigp{
  \AA_\XX(3)^2%RM 03.20.13 \AA_\XX^2(3)
  \la_{3}^2+
  \BB_\XX(3)^2%RM 03.20.13 \BB_\XX^2(3)
  \la_2^2}&\text{if }p=3,\\
  5w^2\bigp{
  \AA_\XX(2q)^2%RM 03.20.13 \AA_\XX^2(2q)
  \la_p^{p-1}+
  \BB_\XX(2q)^2%RM 03.20.13 \BB_\XX^2(2q)
  \la_2^2+\la_p^{2p}}&\text{if }p\in(2,3),
 \end{cases}
\end{equation}
\begin{equation}\label{eq:v}
 \vv:=\frac{\Mf\si}{2\norm{L}^2}\times
 \begin{cases}
  \bigp{\la_{3}^2+2D\la_2\la_{3/2}}&\text{if }p=3,\\
  w^2\bigp{\la_p^{p-1}+4D\la_2\la_q+2\la_q\la_p^p}&\text{if }p\in(2,3),
 \end{cases}
\end{equation}
\begin{equation}\label{eq:la_p}
 \la_\al:=\norm{L}\,\frac{s_\al}{\si}\times
 \begin{cases}
  1&\text{if }p=3,\\
  w^{-1}&\text{if }p\in(2,3).
 \end{cases}
\end{equation}
\end{thm}

\begin{remark}\label{re:|S|}
The term $\P(\norm{S}>\ep)$ in \eqref{eq:f(S).ub} can be bounded in a variety of ways. 
For instance, using Chebyshev's inequality and \eqref{eq:D}, one can write
\begin{equation}\label{eq:|S|>ep}
 \P(\norm{S}>\ep)
  \le\frac{\norm{S}_2^2}{\ep^2}
  \le\frac{D^2s_2^2}{\ep^2}.
\end{equation}
Alternatively, one can write 
\begin{equation*}
 \P(\norm{S}>\ep)
  \le\frac{\norm{S}_p^p}{\ep^p}
    %IP 03.10.13
  \le\frac{\AA_\XX(p)^{p}s_p^p+\BB_\XX(p)^{p}s_2^p}{\ep^p},
%   \le\frac{\AA_\XX^p(p)s_p^p+\BB_\XX^p(p)s_2^p}{\ep^p},
\end{equation*}
using a Rosenthal-type inequality \eqref{eq:rosen}. 
An exponential inequality as described in \eqref{eq:PU,vector} can also be used. 
\end{remark}

\begin{remark}\label{re:s<infty}
The expressions $\uu$ and $\vv$ in \eqref{eq:u} and \eqref{eq:v} are finite for any given $p\in(2,3]$ whenever $s_p<\infty$, whereas $\la_{2q}$ may be infinite for $p\in(2,3)$ even when the condition $s_p<\infty$ holds. 
It is the additional truncation, with $\bar\De$ instead of $\De$, in the bounds of Section~\ref{sec:chen.mod} that allows one to use $\la_p$ instead of $\la_{2q}$ in the terms $\uu$ and $\vv$ when $p<3$; cf.\ Remark~\ref{re:ub}. 
\end{remark}

The hardest to obtain result of this section is the nonuniform bound in Theorem~\ref{thm:f(S).nub} below. 

\begin{thm}\label{thm:f(S).nub}
Assume that the conditions of Theorem~\ref{thm:f(S).ub} are satisfied, and take any real numbers $\th$, 
%IP 03.10.13
$w$, $\de_0$, $\pi_1$, $\pi_2$, $\pi_3$, and $\om$ such that the conditions \eqref{eq:pi's} hold and 
\begin{equation}\label{eq:om}
 \om\in\Bigl(0,\frac{\Mf\ep^2}{2\pi_1}\Bigr].
\end{equation}
Let
\begin{equation}\label{eq:eta.nub}
 %RM 02-21-13 \xi_i:=\frac{L(X_i)}{\si},\quad\text{and}\quad
 \eta_i:=\frac{\norm L\norm{X_i}}{\si}\,\ind{2<p<3}+\frac{L(X_i)}{\si}\,\ind{p=3}.
\end{equation} 
Then for all 
\begin{equation}\label{eq:z}
 z\in(0,\om/\si] 
\end{equation}
one has
\begin{equation}\label{eq:f(S).nub}
 \Bigabs{\Bigprob{\frac{f(S)}{\si}>z}-\Bigprob{\frac{L(S)}{\si}>z}}\le \tilde\ga_z+\tilde\tau e^{-(1-\pi_1)z/\th},
\end{equation}
where
\begin{equation}\label{eq:tga_z}
 \tilde\ga_z:=\biggprob{\norm{S}>\sqrt{\frac{2\pi_1\si z}{\Mf}}}+\ga_z,
\end{equation}
\begin{equation}\label{eq:tau2}
 \tilde\tau:=c_1\si_p\vv+c_2\uu+c_3
  \de, %RM 02-21-13 c_4\si_p^{\tq}, 
\end{equation} 
and $\ga_z$, $c_1$, $c_2$, $c_3$ are as in Theorem~\ref{thm:nub}. 
\end{thm}

\begin{remark}\label{re:nuisance}
The restriction \eqref{eq:z} is of essence. 
Indeed, if $z>>\frac1\si$ (that is, if $z$ is much greater than $\frac1\si$) and the event $\{\frac{L(S)}{\si}>z\}$ in \eqref{eq:f(S).nub} occurs, then $L(S)>>1$ and hence $\|S\|>>1$, and in this latter zone, of large deviations of $S$ from its zero mean, the linear approximation of $f(S)$ by $L(S)$ will usually break down; cf.\ e.g.\ \eqref{eq:De}, in which $\si\De$, measuring the difference between $\si T=f(S)$ and $\si W=L(S)$, is on the order of magnitude of $\|S\|^2$ and thus much greater than $L(S)$ when $\|S\|>>1$. 
This heuristics will be implicitly used in Proposition~\ref{prop:z} in Appendix~\ref{app:z<=sqrt.n}, which shows that the upper bound $\frac{\om}{\si}$ on $z$ in \eqref{eq:z} is indeed the best possible up to a constant factor, even when the Banach space $\XX$ is one-dimensional. 
Note also that \eqref{eq:om} can be satisfied for any given $\om\in(0,\infty)$ by (say) taking $\pi_1$ to be small enough.  
\end{remark}

While the expressions for the upper bounds given in Theorems~\ref{thm:f(S).ub} and \ref{thm:f(S).nub} are quite explicit, they may seem complicated (as compared with the classical uniform and nonuniform BE bounds). However, one should realize that here there are a whole host of players: those associated with the function $f$ and the space $\XX$ (like $\|L\|$, $\Mf$, $\ep$, and $D$), the parameters we are free to choose (namely, $\cc$, $\th$, $w$, $\de_0$, $\pi_1$, $\pi_2$, $\pi_3$, and $\om$), and more traditional terms (as $s_p$, $\si$, and $G_\xi$) -- each with a significant and rather circumscribed role to play.

One should note that the bounds in Theorems~\ref{thm:f(S).ub} and \ref{thm:f(S).nub} do not depend on the dimension of the space $\XX$ but only on the choice of the norm $\|\cdot\|$ on $\XX$. One can exercise this choice to an advantage, as e.g.\ will be done in the application considered in Section~\ref{subsec:quadratic}. 
The only restriction on the norm is that the space $\XX$ (possibly even infinite-dimensional) be of type 2; in particular, the bounds will depend on the ``smoothness'' constant $D$ for the norm and on the corresponding Rosenthal-type inequality constants $(\AA_\XX(\cdot),\BB_\XX(\cdot))$. 

Another advantage of the bounds in \eqref{eq:f(S).ub} and \eqref{eq:f(S).nub} is that they do not explicitly depend on $n$. 
Indeed, $n$ is irrelevant when the $X_i$'s are not identically distributed (because one could e.g.\ introduce any number of %RM 03.22.13 extra 
additional zero summands $X_i$). 
In fact, \eqref{eq:f(S).ub} and \eqref{eq:f(S).nub} remain valid when $S$ is the sum of an infinite series of independent  zero-mean r.v.'s, i.e.\ $S=\sum_{i=1}^\infty X_i$, provided that the series converges in an appropriate sense; see e.g.\ Jain and Marcus \cite{jain75}.

On the other hand, for i.i.d.\ r.v.'s $X_i$ our bounds have the correct order of magnitude in $n$. 
Indeed, let
\begin{equation*}
 V,V_1,\dotsc,V_n\text{ be i.i.d.\ random vectors}
\end{equation*}
in $\XX$, with $\E V=0$. Here we shall use 
\begin{equation*}
 \bar V:=\frac{1}{n}\sum_{i=1}^nV_i
\end{equation*}
in place of $S$ (and hence $\frac1nV_i$ in place of $X_i$). 

\begin{%IP16 cor
thm}\label{cor:iid}
Take any $p\in(2,3]$. Suppose that \eqref{eq:smooth} holds,
\begin{align*}
 \tsi :=\|L(V)\|_2>0,
\end{align*}
and $\norm{V}_p<\infty$. 
Then for all $z\in\R$
\begin{equation}\label{eq:f(S).iid}
 \Bigabs{\Bigprob{\frac{f(\bar V)}{\tsi /\sqrt n}\le z}-\Phi(z)}\le\frac{\CC}{n^{p/2-1}};
\end{equation}
moreover, for any $\om\in(0,\infty)$, $\tth\in(0,\infty)$, and for all 
\begin{equation}\label{eq:z,iid}
z\in\Big(0,\,\frac{\om}{\tsi}\,\sqrt n\Big]
\end{equation}
one has
\begin{equation}\label{eq:f(S).iid.nu}
 \Bigabs{\Bigprob{\frac{f(\bar V)}{\tsi /\sqrt n}\le z}-\Phi(z)} \le\CC\Bigp{n\bigprob{\norm{V}>\CC
 z\sqrt{n}%RM 03.20.13 z
 }+\frac{n\P(\norm{V}>\CC\sqrt{n})}{z^p}+\frac{1}{(
 z\sqrt{n}%RM 03.20.13 \sqrt{n}z
 )^p}+\frac{1}{e^{z/\tth}n^{p/2-1}
 %RM 03.20.13 e^{z/\tth}
 }}.
\end{equation}
Each instance of $\CC$ above is a finite %IP 03.10.13
positive expression that depends only upon $p$, the space $\XX$ (through the constants $D$ in \eqref{eq:D} and $(\AA_\XX(\cdot),\BB_\XX(\cdot))$ in \eqref{eq:rosen}), the function $f$ (through \eqref{eq:smooth}), the moments $\tsi$, $\norm{L(V)}_p$, $\norm{V}_q$, $\norm{V}_2$, and $\norm{V}_p$, with $\CC$ in \eqref{eq:f(S).iid} also depending on $\om$ and $\tth$. 
Also, \eqref{eq:f(S).iid} and \eqref{eq:f(S).iid.nu} both hold when $\P(\sqrt{n}L(\bar V)/\tsi\le z)$ replaces $\Phi(z)$. 
\end{%IP16 cor
thm}

%RM12.27
\begin{thm}\label{thm:iid,uni-dim}
Suppose that a function $f\colon\R\to\R$ is twice continuously differentiable in a neighborhood of $0$, with $f(0)=0$ and $f'(0)\ne0$. 
Let $Y,Y_1,Y_2,\dotsc$ be a sequence of 
i.i.d.\ zero-mean unit-variance real-valued r.v.'s with $\norm{Y}_3<\infty$, and let $\bar Y_n:=\frac1n\sum_{i=1}^nY_i$. 
Then there exists a real number $\CC>0$ such that for all $n\in\N$ and all $z\in\R$ 
\begin{equation}\label{eq:iid,uni-dim}
  \Bigabs{\Bigprob{\frac{f(\bar Y_n)}{\abs{f'(0)}/\sqrt n}\le z}-\Phi(z)}\le\frac{\CC}{\sqrt{n}}.  
\end{equation}
Moreover, for any $\om\in(0,\infty)$ there exists a real number $\CC>0$ such that for all $n\in\N$ and all $z$ as in \eqref{eq:z,iid} 
\begin{equation}
  \Bigabs{\Bigprob{\frac{f(\bar Y_n)}{\abs{f'(0)}/\sqrt n}\le z}-\Phi(z)}\le\frac{\CC}{z^3\,\sqrt{n}}.  
\end{equation}
\end{thm}

%RM12.27 %% paragraph copied from ``ejs2.tex''
Theorem~\ref{thm:iid,uni-dim}, a straightforward consequence of 
%IP16 Corollary
Theorem~\ref{cor:iid}, %RM12.27 Theorem~\ref{thm:iid}, 
is stated here to provide an example of uniform and nonuniform BE bounds for the ``classical'', ``univariate'' delta method; even this very simple case appears to be new to the literature. 
Just as with the BE bound for linear statistics, we see that the moment restriction $\norm{Y}_3<\infty$ is sufficient to obtain a bound on the order of $\bigO(1/\sqrt{n})$. 
That bounds such as \eqref{eq:iid,uni-dim} are useful in applications was suggested to us by E.\ MolavianJazi \cite{MolavianJazi}, who needed such a result in his research in electrical engineering. 

In applications to problems of the asymptotic relative efficiency of statistical tests, usually it is the closeness of the distribution of the test statistic to a normal distribution (in $\R$) that is needed or most convenient; in fact, as mentioned before, obtaining uniform bounds on such closeness was our original motivation for this work.

On the other hand, there have been a number of deep results on the closeness of the distribution of $f(S)$, not to the standard normal distribution, but to that of $f(N)$, where $N$ is a normal random vector with the mean and covariance matching those of $S$. 
In particular, G\"otze \cite{gotze91} provided an upper bound of the order $\bigO(1/\sqrt n)$ on the uniform distance between the d.f.'s of the r.v.'s $f(S)$ and $f(N)$ under comparatively mild restrictions on the smoothness of $f$; however, the bound increases to $\infty$ with the dimension $k$ of the space $\XX$ (which is $\R^k$ therein). 
Bhattacharya and Holmes \cite{bhatt} obtained a constant which is $\bigO(k^{5/2})$, and Chen and Fang \cite[Theorem~3.5]{chenfang} %RM12.26 chen-fang} 
recently improved this to $\bigO(k^{1/2})$.

One should also note here such results as the ones obtained by G\"otze \cite{gotze86} (uniform bounds) and Zalesski{\u\i} \cite{zal88,zal89} (nonuniform bounds), also on the closeness of the distribution of $f(S)$ to that of $f(N)$. 
There (in an i.i.d.\ case), $\XX$ can be any type 2 Banach space, but $f$ is required to be at least thrice differentiable, with certain conditions on the derivatives. 
Moreover, Bentkus and G\"otze \cite{bent93} provide several examples showing that, in an infinite-dimensional space $\XX$, the existence of the first three derivatives (and the associated smoothness conditions on such derivatives) cannot be relaxed in general. 

\section{Applications}
\label{sec:apps}

Here we shall apply the results of Section \ref{sec:f(S)} to present several novel bounds on the rate of convergence to normality for some commonly used statistics. 
For the sake of simplicity and brevity, assume throughout this section that 
\begin{equation*}
	p=3 
\end{equation*} 
and 
$V,V_1,\dotsc,V_n$ are i.i.d.\ $\XX$-valued r.v.'s, where $\XX$ is a Hilbert space; also adopt the notation
\begin{equation}\label{eq:tsi,ga3,v_al}
 \tsi:=\norm{L(V)}_2,\quad\vsi_\al:=\frac{\norm{L(V)}_\al}{\tsi},\quad\text{and}\quad v_\al:=\norm{V}_\al\quad\text{for $\al\ge1$},
\end{equation}
where $L$ is as in \eqref{eq:smooth}. 
%RM 02-21-13 moved, and updated, from Section 3
Under these assumptions we then can choose the smallest (to our knowledge) constants for the Rosenthal-type inequalities in \eqref{eq:rosen1}, \eqref{eq:rosen.noncentral}, and \eqref{eq:rosen}. Namely, 
\begin{equation}\label{eq:A1,B1}
\begin{split}
 \bigp{\AA_\R(3),\BB_\R(3)}&=\bigp{1,(8/\pi)^{1/6}},\\
 \bigp{\AA_{\R,\nc}(3),\BB_{\R,\nc}(3)}&=\bigp{1.316^{1/3},2^{1/3}},\\
 \bigp{\AA_{\XX}(3),\BB_{\XX}(3)}&=\bigp{1,2^{1/3}},
\end{split}
\end{equation}
according to \cite[(5)]{pin13.rosen} (set $x=0$ there), \cite[Corollary~4]{pin12-2smooth} and \cite[(12)]{pin12-2smooth}, respectively; cf.\ also \cite{pin80}. 

Essentially two types of results will be presented in this section. 
Theorems~\ref{thm:t}, \ref{thm:R}, %RM12.30 and 
\ref{thm:T^2}, 
\ref{thm:cov}, \ref{thm:rspcc}, \ref{thm:MLE}, %RM12.30
containing uniform and nonuniform BE-type bounds for specific statistics (namely, Student's, Pearson's, %RM12.30 and 
noncentral Hotelling's, 
certain statistics used to test hypotheses about a covariance operator, a type of canonical correlation, and %IP16 the maximum likelihood estimator) %RM12.30 
maximum likelihood estimators 
are straightforward applications of %IP16 Corollary
Theorem~\ref{cor:iid}, in each specific instance with its own space $\XX$, function $f$, and random vector $V$. 
Of course, these results inherit from %IP16 Corollary
Theorem~\ref{cor:iid} the not quite explicit constants $\CC$, which, recall, were finite expressions depending only upon $p$, the function $f$, and the distribution of $V$, with $\CC$ in the nonuniform bounds also depending on $\om$; however, in contrast with %IP16 Corollary
Theorem~\ref{cor:iid}, the $\CC$'s in Theorems~\ref{thm:t}, \ref{thm:R}, and \ref{thm:T^2} will no longer depend on the space $\XX$, since one can use the same constants $D$ in \eqref{eq:D} and $(\AA_\XX(\cdot),\BB_\XX(\cdot))$ in \eqref{eq:rosen} for all Hilbert spaces $\XX$. 

On the other hand, Theorem~\ref{thm:iid,p=3,unif} will provide a uniform BE-type bound for a normalized statistic $\sqrt{n}f(\bar V)/\tsi$, with explicit coefficients on each of the terms in the bound. These coefficients, denoted by $\ccc$ with two or three subscripts, will in specific applications be variously bounded from above by finite explicit constants which do not depend on $n$ or $z$; so, such coefficients may be referred to as \emph{pre-constants}. 
The corresponding nonuniform bound is much more complicated and therefore will be relegated to Appendix~\ref{app:nonunif}, where it is stated (and proved) as Theorem~\ref{thm:iid,p=3,nonunif}. 
To help the reader follow our indexing of the pre-constants, let us say that the subscript of a pre-constant $\ccc$ will be $\sf u$ or $\sf e$ or $\sf n$, depending on whether the pre-constant appears in a \underline{{\sf u}}niform BE-type bound or in an \underline{{\sf e}}xponentially (in $z$) decreasing term of a nonuniform BE-type bound or in a power-like decreasing term of a \underline{{\sf n}}onuniform BE-type bound, respectively; the remaining subscripts refer to the moments of which the pre-constant is a coefficient.    

We then apply the inequalities of Theorems~\ref{thm:iid,p=3,unif} and \ref{thm:iid,p=3,nonunif} to obtain BE-type bounds for the self-normalized sum and Pearson's correlation coefficient containing only absolute constants and moments of relevant r.v.'s, with a simple (and optimal) dependence on $n$ and $z$; these latter bounds are given in Corollaries~\ref{cor:centralT}, \ref{cor:centralT.nub}, and \ref{cor:Runif}. 
The proofs of these three corollaries are somewhat lengthy and technical, and so are placed in Appendix~\ref{app:num.proofs}. 

\begin{thm}\label{thm:iid,p=3,unif}
Let $\XX$ be a Hilbert space, let $f$ satisfy \eqref{eq:smooth} for some real $\ep>0$, and assume that $\E V=0$, $\tsi>0$, and $v_3<\infty$. 
Take any real numbers
\begin{equation}\label{eq:ka's}
 \cc\in
 \bigl[\tfrac12,1\bigr),\ %RM 02-21-13 (0,1),\ 
 \ka_{2,0}>0,\ \ka_{3,0}>0,\ \ka_{2,1}>0,\text{ and }\ka_{3,1}>0.
\end{equation}
Then 
\begin{align}
\label{eq:iid,p=3}
 \Bigabs{\Bigprob{\frac{f(\bar V)}{\tsi /\sqrt{n}}\le z}-\Phi(z)}
 &\le\frac{\KKu{0}+%RM 12.21.12
 \KKu{1}\vsi_3^3+(\KKu{20}+\KKu{21}\vsi_3)v_2^2+(\KKu{30}+\KKu{31}\vsi_3)v_3^2+\KKu{\ep}}{\sqrt{n}}
 \\
\label{eq:iid,p=3,young's}
 &\le\frac{\Cucon+\Cuvsi\vsi_3^3+\Cuvtwo v_2^3+\Cuvthr v_3^3}{\sqrt n}
\end{align}
for all $z\in\R$ and $n\in\N$, where
\begin{equation}
\begin{split}
\label{eq:Ku}
 \KKu{0}:=0.13925-\frac{(2\cc-1)^2}{2\cc(1-\cc)},\quad %RM 12.21.12
 \KKu{1}:=0.33554 %RM 12.21.12 0.4748
  +\frac1{2\cc(1-\cc)},\\
   \bigp{\KKu{20},\KKu{21},\KKu{30},\KKu{31}}
  :=\frac{\Mf}{4\cc\tsi}\biggp{2
  \Bigp{\frac2\pi}^{1/6}%RM 03.22.13
  ,\, %RM 02-21-13 
  2+\frac{2^{2/3}}{n^{1/6}}
  ,\,%RM 03.22.13 
  \frac{(8/\pi)^{1/6}} %RM 02-21-13 2^{1/3}}
  {n^{1/3}} %RM 03.22.13
  ,\,
  \frac{2}{n^{1/2}}}, 
  %RM12.21.12 3,2+\frac{3^{2/3}}{n^{1/6}},\frac{3^{1/3}}{n^{1/3}},\frac{2}{n^{1/2}}},
\end{split}
\end{equation}
\begin{equation}\label{eq:Kuep}
 \KKu{\ep}:=\frac{v_2^2}{\ep^2n^{1/2}}\bigwedge\frac{
  2 %RM 02-21-13 3
  v_2^3+v_3^3/n^{1/2}}{\ep^3n},
\end{equation}
\begin{equation}\label{eq:tKu}
\begin{aligned}
 \Cucon&:=\KKu{0}+ %RM 12.21.12 
  \frac{1}{3\ka_{2,0}^3}\Bigp{\KKu{20}+\frac{1}{\ep^2n^{1/2}}}+\frac{1}{3\ka_{3,0}^3}\,\KKu{30}, &
 \Cuvsi&:=\KKu{1}+\frac{1}{3\ka_{2,1}^3}\,\KKu{21}+\frac{1}{3\ka_{3,1}^3}\,\KKu{31}, \\
 \Cuvtwo&:=\frac{2\ka_{2,0}^{3/2}}{3}\Bigp{\KKu{20}+\frac{1}{\ep^2n^{1/2}}}+\frac{2\ka_{2,1}^{3/2}}{3}\,\KKu{21}, &
 \Cuvthr&:=\frac{2\ka_{3,0}^{3/2}}{3}\,\KKu{30}+\frac{2\ka_{3,1}^{3/2}}{3}\,\KKu{31}.
\end{aligned}
\end{equation}
\end{thm}

%RM 02-21-13 
% \begin{remark}\label{re:symm.unif}
% Under the additional assumption that $L(V)$ is symmetric, Theorem~\ref{thm:iid,p=3,unif} remains true if the quadruple $(\KKu{20},\KKu{21},\KKu{30},\KKu{31})$ is replaced by the following one, with smaller values: 
% \begin{equation}\label{eq:KKu.symm}
% \frac{\Mf}{4\cc\tsi}\Bigp{3^{2/3}C_3^*,2,\frac{C_3^*}{n^{1/3}},\frac1{n^{1/2}}},\quad\text{where}\quad C_3^*:=\Bigp{1+\sqrt{\tfrac8\pi}}^{1/3}.
% \end{equation}
% Indeed, \cite[Theorem~1]{ibrag97} implies $\norm{W}_3\le C_3^*(1\vee\si_3)$ (recall here \eqref{eq:W}, \eqref{eq:=1}, and \eqref{eq:si_p}). 
% Since $\KKu{1}>1$ for any $\cc\in(0,1)$, we may assume w.l.o.g.\ that $\si_3^3=\vsi_3^3/\sqrt{n}<1$ and hence $\norm{W}_3\le C_3^*$. 
% So, following the lines of the proof of Theorem~\ref{thm:f(S).ub}, we see that the upper bound $\AA_\R(p)\si_p+\BB_\R(p)$ on $\norm{W}_p$, which appears in \eqref{eq:f(S).ub}, can be replaced there by $C_3^*$. 
% \end{remark}

\begin{remark}\label{re:nonunif} 
One can have a ``nonuniform'' counterpart to Theorem~\ref{thm:iid,p=3,unif}. 
Indeed, assume that the conditions of Theorem~\ref{thm:iid,p=3,unif} take place; in particular, let $\ep$ and $\Mf$ be any positive real numbers such that \eqref{eq:smooth} holds. 
Take any positive real numbers $z_0$, $\tth$, $K_1$, $K_2$, and $K_3$. 
Then%RM 03.22.13
, by Theorem~\ref{thm:iid,p=3,nonunif}, there exist some finite positive constants $\om$, $\Lnvsi$, $\Lnvtwoa$, $\Lnvtwob$, $\Lnvthra$, $\Lnvthrb$, $\Lnecon$, $\Lnevsi$, $\Lnevtwo$, and $\Lnevthr$, each depending only on $\ep$, $\Mf$, $z_0$, $\tth$, $K_1$, $K_2$, and $K_3$, such that 
\begin{equation}\label{eq:re:nonunif}
 \Bigabs{\P\Big(\frac{f(\bar V)}{\tsi/\sqrt{n}}\le z\Big)
 -\Phi(z)}
 \le 
  \frac{\Lnvsi\vsi_3^3+\bigp{(\Lnvtwoa\vee \Lnvtwob)v_2^4}\vee\bigp{\Lnvthra v_3^3}+\Lnvthrb v_3^3}{z^3\sqrt n}
  +\frac{\Lnecon+\Lnevsi\vsi_3^3+\Lnevtwo v_2^3+\Lnevthr v_3^3}{e^{z/\tth}\sqrt n}
\end{equation}
for all $z\in\R$ and $n\in\N$ such that
\begin{equation}\label{eq:z.n.conds}
 z_0\le z\le\frac{\om}{\tsi}\,\sqrt{n},\quad
 \frac{K_1\vsi_3^3}{\sqrt{n}}\le1,\quad
 \frac{K_2v_2^4}{\tsi^3z^3\sqrt{n}}\le1,\quad\text{and}\quad
 \frac{K_3v_3^3}{\tsi^3z^3\sqrt{n}}\le1.
\end{equation}
The constants $\Lnvsi,\dots,\Lnevthr$ in \eqref{eq:re:nonunif} are upper bounds on certain corresponding pre-constants $\Cnvsi,\dots,\Cnevthr$, explicit expressions for which are given in Theorem~\ref{thm:iid,p=3,nonunif}. 
Concerning the conditions in \eqref{eq:z.n.conds}, note the following: 
\begin{enumerate}
	\item The condition $z\ge z_0$ does not diminish generality, in view of uniform bounds \eqref{eq:iid,p=3} and \eqref{eq:iid,p=3,young's}. 
	\item The condition $z\le\frac{\om}{\tsi}\,\sqrt{n}$ is essential and even optimal, up to a constant factor, as shown in Appendix~\ref{app:z<=sqrt.n}. 
	\item The other three conditions in \eqref{eq:z.n.conds}, involving the constants $K_1$, $K_2$, and $K_3$, will be satisfied when $n$ and $z$ are large enough. As mentioned above, the case when $z$ is not large can be covered using a uniform bound. Finally, the remaining case with ``large'' $z$ and ``small'' $n$ can be dealt with based on an appropriate upper bound on large deviation probabilities. 
In fact, the proof (given in Appendix~\ref{app:num.proofs}) of the nonuniform bound in  Corollary~\ref{cor:centralT.nub} is conducted right along such lines.  
\end{enumerate}
\end{remark}

The mentioned pre-constants in Theorems~\ref{thm:iid,p=3,unif} and \ref{thm:iid,p=3,nonunif} are complicated in appearance. 
However, in particular applications -- presented in Corollaries~\ref{cor:centralT}, \ref{cor:centralT.nub}, and \ref{cor:Runif} -- these statements will result in bounds of much simpler structure, with explicit numerical constants, which are also rather moderate in size, especially in the uniform bounds.  
The following corollary shows that the asymptotic behavior of the uniform and nonuniform BE-type bounds given in Theorems~\ref{thm:iid,p=3,unif} and \ref{thm:iid,p=3,nonunif} is quite simple as well, and the corresponding constants are again moderate in size. 

\begin{cor}\label{cor:asymp}
Assume that the conditions of Theorem~\ref{thm:iid,p=3,unif} hold, and also that $f''$ is twice continuously differentiable in a neighborhood of the origin. 
Then 
\begin{equation}\label{eq:asymp.unif}
 \limsup_{n\to\infty}\,\sup_{z\in\R}\sqrt{n}\Bigabs{\Bigprob{\frac{f(\bar V)}{\tsi/\sqrt{n}}\le z}-\Phi(z)}
  \le0.63925+0.83554\vsi_3^3+\frac{y_*}{2}+\frac12\,\sqrt{(\vsi_3^3-1)(\vsi_3^3-1+2y_*)}, %RM 12.21.12
% \le 0.13925+0.83554\vsi_3^3 %RM 12.21.12 0.9748\vsi_3^3
%  +\frac{y_*}4
%  +\frac12\sqrt{\vsi_3^3\bigp{\vsi_3^3+y_*}}
%  %RM 12.21.12 \le1.4748\vsi_3^3+\frac{y_*}2,
\end{equation}
where 
\begin{equation}\label{eq:y_*}
y_*:=\tfrac{\norm{f''(0)}}{\tsi}
\Bigp{\bigp{\tfrac2\pi}^{1/6}+\vsi_3}v_2^2. %RM 02-21-13 (1+\vsi_3)v_2^2. 
%RM 12.21.12 (3+2\vsi_3)v_2^2. 
\end{equation}

Also, for any positive increasing unbounded function $g$ on $\N$ 
\begin{equation}\label{eq:asymp.nonunif}
 \limsup_{n\to\infty}\sup_{g(n)\le z\le\sqrt{n}/g(n)}z^3\sqrt{n}\bigabs{\bigprob{\tfrac{f(\bar V)}{\tsi /\sqrt{n}}\le z}-\Phi(z)}\le30.2211\vsi_3^3; 
\end{equation}
in fact, here it will be possible to replace the factor $30.2211$ by any improved constant factor that one will be able to obtain in place of $30.2211$ in the nonuniform BE inequality \eqref{eq:tyurin,michel} for linear statistics. 
\end{cor}

As one can see, in the expressions of the asymptotic uniform bounds in \eqref{eq:asymp.unif} the higher moment $v_3$ disappears, and in the asymptotic nonuniform bound in \eqref{eq:asymp.nonunif} the moment $v_2$ disappears as well; 
however, Corollary~\ref{cor:asymp} inherits the condition $v_3<\infty$ from Theorems~\ref{thm:iid,p=3,unif} and \ref{thm:iid,p=3,nonunif} -- where, as seen from Remark%IP 03.11.13
s~\ref{re:moments,T1} and \ref{re:R.6th.moments}, this condition is essential; cf.\ also Remark~\ref{re:shao-compar}.  

For the remainder of the results in this section, $\XX$ will be the Euclidean space $\R^k$ for some natural number $k$, and the nonlinear functional $f\colon\XX\to\R$ will be continuously twice differentiable in some neighborhood about the origin. 
Thus, for a given (small enough) $\ep$, the smoothness condition \eqref{eq:smooth} will hold when
\begin{equation}\label{eq:L,M,iid}
 L=f'(0)\quad\text{and}\quad \Mf=\sup_{\norm{\x}\le\ep}\norm{f''(\x)},
\end{equation}
where $f'(\x)$ and $f''(\x)$ are identified with the gradient vector and the Hessian matrix, respectively, of $f$ at some point $\x\in\XX$, and then $\norm{f''(\x)}$ denotes the spectral norm of the matrix $f''(\x)$. 
Upon specifying the function $f$ and the relevant r.v. $V$, the results of Theorems~\ref{thm:t}, \ref{thm:R}, and \ref{thm:T^2} (uniform and nonuniform bounds without explicit coefficients) will be proved by invoking %RM 03.22.13 the results of
%IP16 Corollary
Theorem~\ref{cor:iid}. 

\subsection{``Quadratic'' statistic}\label{subsec:quadratic} %% supporting calculations found in notebook '02-26-13.novak.quadratic.compare.nb' %%

The first application we consider involves a particularly simple nonlinear statistic investigated by Novak in \cite[Section~3]{novak05}. 
Let $V=(Y,Z),V_1=(Y_1,Z_1),\dotsc,V_n=(Y_n,Z_n)$ be i.i.d.\ r.v.'s with $\E V=0$, $\E Y^2=\E Z^2=1$. 
Take any real $\th>0$ and let $\XX$ be $\R^2$ with the norm defined by the formula $\|
\x%RM 03.20.13 x
\|:=\sqrt{x_1^2+x_2^2/\th^2}$ for $
\x%RM 03.20.13 x
=(x_1,x_2)\in\XX$. 
Next, take any real $c_0\ge0$ and let $f\colon\R^2\to\R$ be defined by $f(x_1,x_2)=x_1+c_0x_2^2$. 
Then $f$ satisfies the smoothness condition \eqref{eq:smooth} with $L(x_1,x_2)=x_1$ and $\Mf=2c_0\th^2=\|f''(0)\|$, for any $\ep>0$. 
Consider the statistic
\begin{equation}\label{eq:Q}
 Q:=\bar{Y}+c_0\bar{Z}^2=f(\bar{V})\quad\text{with}\quad\bigp{\bar{Y},\bar{Z}}=\bar{V}=\frac1n\,\tsum_{i=1}^nV_i,  
\end{equation}
so that the statistic $\sqrt{n}Q=\sum_i(Y_i/\sqrt{n})+c(\sum_iZ_i/\sqrt{n})^2$ with $c:=c_0/\sqrt{n}$ coincides with the quadratic statistic studied in \cite{novak05}; the $X_i$'s and $Y_i$'s in \cite{novak05} are replaced here by $Y_i/\sqrt{n}$ and $Z_i/\sqrt{n}$, respectively. 
One may also note that in \cite{novak05} the condition $\E Z^2=1$ was not assumed; however, it can be assumed (as we do) without loss of generality, by adjusting the choice of the factor $c_0$. 

Now one can use the inequalities $\norm{Y}_1\le\norm{Y}_2=1$, $\norm{Z}_1\le\norm{Z}_2=1$, $\norm{YZ}_1\le\norm{Y}_2\norm{Z}_2=1$, $\norm{\sqrt{n}\bar{Z}}_1\le\norm{\sqrt{n}\bar{Z}}_{3/2}\le\norm{\sqrt{n}\bar{Z}}_2=1$, and $\norm{\sqrt{n}\bar{Y}}_3\le\norm{Y}_3/n^{1/6}+
(8/\pi)^{1/6} %RM 02-21-13 3^{1/3}
$ (cf.\ \eqref{eq:rosen1} and \eqref{eq:A1,B1}) in conjunction with \cite[Theorem~2]{novak05} to obtain
\begin{equation}\label{eq:quad.Novak}
\begin{split}
 \limsup_{n\to\infty}\sqrt{n}\bigabs{\P(\sqrt{n}Q\le z)-\Phi(z)} &\le2+\Bigp{\tfrac9{\sqrt{2\pi}}+\sqrt{\tfrac\pi8}+1}\norm{Y}_3^3+\Bigp{\sqrt{\tfrac\pi2}+4}c_0
\\
 &<2+5.218\norm{Y}_3^3+5.254c_0.
\end{split}
\end{equation}
On the other hand, Corollary~\ref{cor:asymp} implies
\begin{equation}\label{eq:quad.nonlinear}
\begin{aligned} 
 \limsup_{n\to\infty}\sqrt{n}\bigabs{\P(\sqrt{n}Q\le z)-\Phi(z)}
 &\le0.63925+0.83554\norm{Y}_3^3+\frac{\ty_*}2+\frac12\sqrt{(\norm{Y}_3^3-1)(\norm{Y}_3^3-1+2\ty_*)}, %RM 12.21.12
%  &\le0.13925+0.83554\norm{Y}_3^3+\frac{\ty_*}4+\frac12\sqrt{\norm{Y}_3^3(\norm{Y}_3^3+\ty_*)},
  %RM 12.21.12 &\le\inf_{\th>0}\big[1.4748\norm{Y}_3^3+c_0\th^2\bigp{3+2\norm{Y}_3}\bigp{1+1/\th^2}\big] \\ 
  % &=1.4748\norm{Y}_3^3+\bigp{3+2\norm{Y}_3}c_0. 
\end{aligned}   
\end{equation} 
where
\begin{equation*}
 \ty_*=\inf_{\th>0}y_*=2c_0
 \bigp{\bigp{\tfrac2\pi}^{1/6}+\norm{Y}_3}<2c_0\bigp{0.928+\norm{Y}_3}. %RM 02-21-13 (1+\norm{Y}_3). 
 %RM 12.21.12 2c_0\bigp{3+2\norm{Y}_3}.
\end{equation*}
%IP12.13.12
Note that, in contrast with \cite[Theorem~2]{novak05}, which only required that $\E|Y|^3+\E Z^2< \infty$, one needs the condition $\E|Y|^3+\E|Z|^3< \infty$ to deduce \eqref{eq:quad.nonlinear} immediately from Corollary~\ref{cor:asymp}. 

\begin{wrapfigure}{L}{.3\textwidth}
%  \vspace*{-8pt} 
  \psfragfig[width=0.28\textwidth]{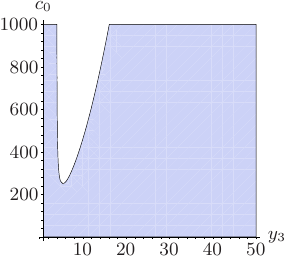}
  %RM 12.21.12 \includegraphics[width=.28\textwidth]{asymp-scr1.pdf}
  \caption{\eqref{eq:quad.Novak} vs.\ \eqref{eq:quad.nonlinear}}
% \vspace*{-12pt} 
 \label{fig:asymp-scr1}
\end{wrapfigure}
%RM 12.21.12 \noindent 
Figure~\ref{fig:asymp-scr1} shows the set (shaded) of all points $(y_3,c_0)\in[1,50]\times[0,1000]$ %RM 12.21.12 [1,10]\times[0,50]$ 
with $y_3:=\norm{Y}_3$ for which the asymptotic bound in \eqref{eq:quad.nonlinear} is less than  
that in \eqref{eq:quad.Novak}.   
It is seen that \eqref{eq:quad.nonlinear} works better than \eqref{eq:quad.Novak} un\-less the ``nonlinearity coefficient'' $c_0$ in \eqref{eq:Q} is very large. 
In particular, for \eqref{eq:quad.Novak} to be better than \eqref{eq:quad.nonlinear} it is necessary that 
$\norm{Y}_3>4.07$ and $c_0>249+3.06(\norm{Y}_3-5.44)^2$. %RM 02-21-13 $\norm{Y}_3>3.97$ and $c_0>203+2.67(\norm{Y}_3-5.27)^2$. 
%RM 12.21.12 $\norm{Y}_3>3.47$ and $c_0>16.37 + 1.871 (\norm{Y}_3 - 1)^2$. 
This and discussion in subsequent 
Subsubsection~\ref{subsec:centralT} 
suggest that bounds developed in this paper for general nonlinear statistics are competitive with bounds obtained earlier by specialized methods, tailored to a specific statistic or a specific class of statistics. 
\\
\\
\\
\\

\subsection{Student's \texorpdfstring{$T$}{T}}

Let $Y,Y_1,\dotsc,Y_n$ be i.i.d.\ real-valued r.v.'s, with
\begin{equation*}
 \mu:=\E Y\quad\text{and}\quad\var Y\in(0,\infty).
\end{equation*}
Consider the statistic commonly referred to as Student's $T$ (or simply $T$):
\begin{align*}
 T:=\frac{\Ybar}{S_Y/\sqrt n}=\frac{\sqrt{n}\ \Ybar}{\bigp{\Ysqbar-\Ybar^2}^{1/2}},
\end{align*}
where
\begin{equation*}
 \Ybar:=\tfrac1n\tsum\nolim_iY_i,\quad
 \Ysqbar:=\tfrac1n\tsum\nolim_iY_i^2,\quad\text{and}\quad
 S_Y:=\Bigp{\tfrac1n\tsum\nolim_i(Y_i-\Ybar)^2}^{1/2}=\Bigp{\Ysqbar-\Ybar^2}^{1/2};
\end{equation*}
let %IP12.13.12 $T$ take an arbitrary value $t_0$ 
$T:=0$ when $\Ysqbar=\Ybar^2$. 
Note that $S_Y$ is defined here as the empirical standard deviation of the sample $(Y_i)_{i=1}^n$, rather than the sample standard deviation $(\frac n{n-1}(\Ysqbar-\Ybar^2))^{1/2}$. 

Let us call $T$ ``central'' when $\mu=0$ and ``non-central'' when $\mu\ne0$. 

As $T$ is invariant under the transformation $Y_i\mapsto aY_i$ for arbitrary $a>0$, let us assume without loss of generality (w.l.o.g.) that   
\begin{equation*}
 \var Y=1.
\end{equation*}

Now let $\XX=\R^2$, and for $\x=(x_1,x_2)\in\XX$ such that $1+x_2-x_1^2>0$, let $f\colon\XX\to\R$ be defined by
\begin{equation*}
 f(\x)=f(x_1,x_2)=\frac{x_1+\mu}{\sqrt{1+x_2-x_1^2}}-\mu;
\end{equation*}
let $f(\x):=%IP12.13.12 t_0/\sqrt{n}
-\mu$ for all other $\x\in\XX$%IP12.13.12, where $t_0$ is the ``exceptional'' value chosen above for $T$
. 
Since
\begin{equation}\label{eq:d(ep)}
 \min_{x_1^2+x_2^2\le\ep^2}(1+x_2-x_1^2)=
 \begin{cases}
  1-\ep&\text{if }0<\ep\le\tfrac12,\\
  \tfrac34-\ep^2&\text{if }\ep\ge\tfrac12,
 \end{cases}
\end{equation}
it is easy to see that $f''$ is continuous (and hence uniformly bounded) on the closed ball $\{\x\in\XX\colon\norm{\x}\le\ep\}$ for any fixed $\ep\in(0,\sqrt{3}/2)$. 
Then the smoothness condition \eqref{eq:smooth} is satisfied, with $L(\x)=f'(0)(x_1,x_2)=x_1-\mu x_2/2$ for $\x=(x_1,x_2)\in\XX$, and upon letting
\begin{equation}\label{eq:V for t}
 V=\bigp{Y-\mu,(Y-\mu)^2-1}
\end{equation}
we see that $\sqrt{n}f(\bar V)=T-\sqrt{n}\mu$. 
Then %IP16 Corollary
Theorem~\ref{cor:iid} %IP12.13.12
and
Markov's inequality immediately yield%s 

\begin{thm}\label{thm:t}
Take any $\om>0$ and assume that $\tsi>0$ and $v_3<\infty$, for $\tsi$ and %RM 03.22.13 $v_p$ 
$v_\al$ defined in \eqref{eq:tsi,ga3,v_al}. 
Then for all $z\in\R$ and $n\in\N$ 
\begin{equation}\label{eq:t.bound.uni}
 \Bigabs{\Bigprob{\frac{T-\sqrt{n}\mu}{\tsi}\le z}-\Phi(z)}\le\frac{\CC}{\sqrt n},
\end{equation} 
where $\CC$ is a finite expression depending only on the distribution of $Y$; 
also, for all real $z>0$ and $n\in\N$ satisfying \eqref{eq:z,iid} 
\begin{equation}\label{eq:t.bound}
 \Bigabs{\Bigprob{\frac{T-\sqrt{n}\mu}{\tsi}\le z}-\Phi(z)}\le\frac{\CC}{z^3\,\sqrt n}, 
\end{equation}
where $\CC$ is a finite expression depending only on $\om$ and the distribution of $Y$. 
\end{thm}

\begin{remark}\label{re:t.si=0}
If $\mu=0$ then $\tsi\neq0$, and otherwise $\tsi=0$ only if $Y$ has a 2-point distribution, which depends only on $\mu$. 
Indeed, if $\mu\ne0$ then $\tsi=0\Leftrightarrow L(V)=0$ a.s.\ $\Leftrightarrow Y-\mu=(1\pm\sqrt{1+\mu^2})/\mu$ a.s. 
That is, $\tsi=0$ if and only if $Y=2\sqrt{p(1-p)}/(1-2p)+B_p$ a.s., where $B_p$ is a standardized Bernoulli($p$) r.v. with $p\in(0,1)\setminus\{\frac12\}$.
\end{remark}

%RM01.12
\begin{remark}\label{re:t.asymp}
The upper bound in \eqref{eq:t.bound.uni} is optimal in its dependence on $n$ for the noncentral $T$. 
Indeed, suppose that a function $f\colon\R^k\to\R$ is twice continuously differentiable in a neighborhood of the origin (so that $f$ satisfies the smoothness condition \eqref{eq:smooth}), and let $L$ and $H$ denote here the gradient vector and Hessian matrix of $f$ at 0. 
Further assume, in addition to the assumptions $\tsi>0$ and $v_3<\infty$, that $V$ satisfies the Cram\'er-type condition $\limsup_{\norm{t}\to\infty}\abs{\E e^{it\T V}}<1$. 
Then a calculation of the asymptotic distribution of $\sqrt{n}f(\bar V)/\tsi$ using \cite[Theorem~2]{bhatt78} implies
\begin{equation}\label{eq:BG}
 \sup_{z\in\R}\biggabs{\Bigprob{\frac{f(\bar V)}{\tsi/\sqrt{n}}\le z}-\Phi(z)-\frac{\De(z)}{\sqrt{n}}}=o\Big(\frac1{\sqrt{n}}\Big),
\end{equation}
where
\begin{equation}\label{eq:BG.De}
 \De(z):=-\Bigb{\Bigp{\frac{\E[(L\T V)^3]}{6\tsi^3}+a_3}(z^2-1)+a_1}\varphi(z),
\end{equation}
\begin{equation*}
 a_1:=\frac1{2\tsi}\,\tr(H\Si),\quad a_3:=\frac1{4\tsi^3}(L\T\Si L-\tsi^2)\tr(H\Si)+\frac1{2\tsi^3}\,L\T\Si H\Si L,
\end{equation*}
$\Si$ denotes the covariance matrix of $V$, and $\varphi$ is the standard normal density. 

In the conditions of Theorem~\ref{thm:t}, take the simple case where $Y$ is symmetric about its non-zero mean $\mu$, unit-variance, and has an absolutely continuous distribution; let $\nu_k:=\E(Y-\mu)^k$ denote the $k^\text{th}$ central moment of $Y$, so that $\nu_k=0$ for odd natural $k$. 
Then, for $\De(z)$ as in \eqref{eq:BG.De}, 
\begin{equation*}
 \De(1)=-\frac{\mu(1+3\nu_4)}{8\tsi}\,\varphi(1)\quad\text{and}\quad\tsi=1+\frac{\mu^2}4(\nu_4-1).
\end{equation*}
That is, $\tsi>0$ and $\De(1)\ne0$, and we see that the dependence of the upper bound in \eqref{eq:t.bound.uni} (when $\mu\ne0$) on $n$ is optimal.
\end{remark}

Much work has been done rather recently concerning the distribution of the central $T$; see some references in this regard in Subsubsection~\ref{subsec:centralT} below. 

On the other hand, the bounds in \eqref{eq:t.bound.uni} and \eqref{eq:t.bound} appear to be new for the non-central $T$. 
Bentkus, Jing, Shao, and Zhou \cite{bent07} recently showed that if $\norm{Y}_4<\infty$, then (after some standardization) $T$ has a limit distribution which is either the standard normal distribution or the $\chi^2$ distribution with one degree of freedom; the latter will be the case if and only if $Y$ has the two-point distribution described above in Remark~\ref{re:t.si=0} concerning the degeneracy condition $\tsi =0$. 

\begin{remark}\label{re:moments,T1}
The condition $\norm{Y}_4<\infty$ in \cite{bent07} is equivalent to $\norm{V}_2<\infty$, where $V$ is as in \eqref{eq:V for t}. 
Therefore, it appears natural to require that $\norm{V}_3<\infty$%RM 03.22.13 , or equivalently $\norm{Y}_6<\infty$, 
\ or, equivalently, $\norm{Y}_6<\infty$ 
in order to obtain a bound of order $\bigO(1/\sqrt{n})$; cf.\ the classical BE bound for linear statistics, where the finiteness of the third moment of the summand r.v.'s is usually imposed to achieve a bound of order $\bigO(1/\sqrt{n})$. 
In fact, the asymptotic expansion for the distribution of $T$ up to the order of $O(1/\sqrt n)$ (which follows from the general results for nonlinear statistics obtained by Bhattacharya and Ghosh~\cite{bhatt78}) indeed contains $\norm{Y}_6$ whenever the mean $\mu$ is nonzero. 

The ``central'', or ``null'', case when $\mu=0$ is in this sense exceptional, as discussed in Remark~\ref{re:T1-excep}. In this case, it is well known that the finiteness of the $\E|Y|^3$ is enough for a uniform BE bound for $T$. 
On the other hand, %IP12.13.12
it follows from the remark by 
Novak %IP12.13.12 shows 
at the end of \cite{novak05} that no nonuniform bound of the form $\E|Y|^3g(z)/\sqrt n$ for the self-normalized sum %IP12.13.12, a statistic closely related to 
or, equivalently, for the central $T$ 
can hold for any positive function $g$ such that $g(z)\downarrow0$ as $z\uparrow\infty$. 
%IP 03.10.13
Also, even for $\mu=0$, the presence of the higher order moments can be overcome by an appropriate truncation of the underlying distribution, as suggested by \cite[Corollary~1.5]{pin11} and the discussion therein following it; further details on this can be found in Remark~\ref{re:shao-compar2} below. 
\end{remark}

\subsubsection{Central \texorpdfstring{$T$}{T} and the self-normalized sum}
\label{subsec:centralT}

The central $T$ is very close to the self-normalized sum
\begin{equation}\label{eq:T1}
 T_1:=\frac{Y_1+\dotso+Y_n}{\sqrt{Y_1^2+\dotso+Y_n^2}}=\frac{\sqrt n\,\Ybar}{\sqrt{\Ysqbar}}=\frac{T}{\sqrt{1+T^2/n}}. 
\end{equation} 
In particular, letting $z_n:=z/\sqrt{1+z^2/n}$, one has $\P(T\le z)=\P(T_1\le z_n)$ for all $z\in\R$ and hence  
\begin{multline*}
	\Big|\sup_{z\in\R}|\P(T\le z)-\Phi(z)|-\sup_{z\in\R}|\P(T_1\le z)-\Phi(z)|\Big| \\ 
	\le\sup_{z\in\R}|\Phi(z_n)-\Phi(z)|
	\le\sup_{u\in\R}|u^3\Phi'(u)|/(2n)=(3/(2e))^{3/2}/(n\sqrt{\pi})<0.24/n, 
\end{multline*}
which is much less than $1/\sqrt n$; 
cf.\ \cite[Proposition~1.4]{pin11} and its proof, where Student's $T$ was defined using the sample standard deviation (as opposed to the empirical standard deviation) of the random sample $(Y_i)_{i=1}^n$. 

Slavova \cite{slav85} appears to have first produced a uniform BE-type bound for $T$ of the optimal order in $n$, namely of the form $C/\sqrt n$, where $C$ depends only on $\E|Y|^3$. 
It was only in 1996 that Bentkus and G\"otze \cite[Theorem~1.2]{bent96} obtained a uniform BE-type bound of the optimal order in $n$ \emph{and} with the ``correct'' dependence on the moments; namely, they showed that there exists an absolute constant $A$ such that
\begin{equation}\label{eq:bg}
 \bigabs{\P(T\le z)-\Phi(z)}
  \le An\E\Big[\Big(\frac Y{\sqrt{n}}\Big)^2\wedge\Big|\frac Y{\sqrt{n}}\Big|^3\Big]
\end{equation}
for all $z\in\R$; note that the above bound is no greater than $A\,\E|Y|^p/n^{p/2-1}$ for any  $p\in[2,3]$.  
Bentkus, Bloznelis, and G{\"o}tze \cite{bbg96} provided a similar bound when the $Y_i$'s are not necessarily identically distributed (i.d.).  
Shao \cite[Theorem~1.1]{shao05} obtained a version of \eqref{eq:bg} with explicit absolute constants (and also without the i.d.\ assumption), which in particular implies that in the i.i.d.\ case for all $z\in\R$ 
\begin{align}%RM12.21.12  %% added extra line to equation 
\label{eq:shao.exact} \bigabs{\P(T_1\le z)-\Phi(z)}&\le
  10.2n\E Y^2\ind{\abs{Y}>\sqrt{n}/2}+\frac{25\E\abs{Y}^3\ind{\abs{Y}\le\sqrt{n}/2}}{\sqrt{n}}\\
\label{eq:shao.selfnorm}
&\le\frac{25\norm{Y}_3^3}{\sqrt n}.
\end{align}

Novak \cite{novak00,novak05} obtained BE-type bounds for $T_1$; however, the structure of those bounds is rather complicated. 

Nagaev \cite[Theorem~1 and (1.18)]{nag02}, stated that for all $z\in\R$
\begin{equation}\label{eq:nag.selfnorm}
 \bigabs{\P(T_1\le z)-\Phi(z)}\le
  \frac{36\norm{Y}_3^3+9}{\sqrt n}\wedge
  \frac{4.4\norm{Y}_3^3+\norm{Y}_4^4/\norm{Y}_3^3+\norm{Y^2-1}_3^3}{\sqrt n}
\end{equation}
when the $Y_i$'s are i.i.d. However, there are a number of mistakes of various kinds in the proof in \cite{nag02}; see \cite{pin11} for details. 

\begin{remark}\label{re:pin-shao-compar}
Pinelis \cite[Theorem~1.2]{pin11} obtained a bound of the form
\begin{equation}\label{eq:pin,T}
 \bigabs{\P(T_1\le z)-\Phi(z)}\le\frac1{\sqrt{n}}\Bigp{A_3\norm{Y}_3^3+A_4\norm{Y^2-1}_2+A_6\,\frac{\norm{Y^2-1}_3^3}{\norm{Y}_3^%RM 12.21.12 3
 9
 }}
\end{equation}
for all $z\in\R$, where the triple $(A_3,A_4,A_6)$ depends on several parameters whose values may be freely chosen within certain ranges. 
For instance, a specific choice of the parameters yields $(A_3,A_4,A_6)=(1.53,1.52,1.28)$. Thus, all the constant factors $A_3,A_4,A_6$ in \eqref{eq:pin,T} can be made rather small. 
A bound for the general, non-i.d.\ case, similar to \eqref{eq:pin,T} but with slightly greater constants, was also obtained in \cite{pin11}; as shown there, that bound in \cite{pin11} compares well with %RM 03.22.13 \eqref{eq:shao.selfnorm}, 
\eqref{eq:shao.exact}, especially after truncation.  
\end{remark} 

A number of important advances concerning limit theorems for the central $T$ and/or $T_1$ have been made rather recently. 
For instance, Hall \cite{hall87} obtained an Edgeworth expansion of the distribution of $T$. 
It was only in 1997 that Gin{\'e}, G{\"o}tze, and Mason \cite{ggm} found a necessary and sufficient condition for the {S}tudent statistic to be asymptotically standard normal. 
Shao \cite{shao97,shao99}, Nagaev \cite{nag05}, Jing, Shao, and Wang \cite{jing03}, and Wang and Hall \cite{wang09} studied the probabilities of large deviations. 
Chistyakov and G\"otze \cite{chist03,chist04} and Jing, Shao, and Zhou \cite{jing08} considered the  probabilities of moderate deviations.  
See Gin\'e and Mason \cite{gine98} and Pang, Zhang, and Wang \cite{pang08} concerning the law of the iterated logarithm, and Wang and Jing \cite{wang99} and Robinson and Wang \cite{robinson05} for exponential nonuniform BE bounds.
This is of course but a sampling of the recent work done concerning asymptotic properties of the central $T$ and the related self-normalized sums; for work done somewhat earlier, the reader may be referred to the bibliography in \cite{bent96}. 

\begin{remark}\label{re:T1-excep}
The central $T$ (as compared with the noncentral one) is special for two reasons: (i) when $\mu=0$, then $L(V)=Y$ and, to be finite, $\tsi$ needs only the second moment of $Y$ (rather than the fourth) to exist; and (ii) while in general $\De$ is rather naturally of the order $\sqrt{n}\norm{\bar{V}}^2$, $\De$ is significantly smaller for the central $T$. 
Moreover, the first term, $\sqrt n L(\bar V)/\tsi$, in a formal stochastic expansion of the central $T$ is precisely $\sqrt n\bar Y$ and thus linear in the $Y_i$'s, whereas for the noncentral $T$ this term contains $\bar{Y^2}$. 
This heuristics is reflected in Corollary~\ref{cor:centralT2} below, which is derived using Theorem~\ref{thm:ub}, with a better choice of $\De$ for this specific case than that for the general results of Section~\ref{sec:f(S)}. 
\end{remark}

\begin{cor}[to Theorem~\ref{thm:ub}]\label{cor:centralT2} %RM %% supporting calculations found in notebook '03-08-13.unifbound.4thmoments.nb' %% 
Let $Y,Y_1,\dotsc,Y_n$ be i.i.d.\ r.v.'s, with $\E Y=0$ and $\norm{Y}_2=1$. Then
\begin{equation}\label{eq:centralT.4th}
 \bigabs{\P(T_1\le z)-\Phi(z)}\le\frac{1}{\sqrt n}\Bigp{A_3\norm{Y}_3^3+A_4\norm{Y}_4^4-A_0}
\end{equation}
for all $z\in\R$ and any triple
\begin{equation}\label{eq:Ajlist.4th}
 (A_3,A_4,A_0)\in
 \bigb{(3.00,4.66,4.33),%RM 03.22.13 
 \,(3.17,2.04,1.07),%RM 03.22.13 
 \,(3.48,1.27,-1.43)}. %RM 03-05-13 \bigb{(3.00,5.03,4.70),(3.18,2.09,1.01),(3.60,1.30,-1.44)}. %RM12.21.12 \bigb{(3.14,5.16,4.89),(3.34,2.20,1.28),(3.79,1.31,-1.32)}.
\end{equation}
\end{cor} 

It appears that the bound in \eqref{eq:centralT.4th} may in certain cases be competitive with the bound in \eqref{eq:pin,T} \big(say with $(A_3,A_4,A_6)=(1.53,1.52,1.28)$, as before\big), even though the bound in \eqref{eq:pin,T} was obtained by methods specifically designed for $T_1$. 
Therefore, by Remark~\ref{re:pin-shao-compar}, the bound in \eqref{eq:centralT.4th} may also in certain cases compare well with that in %RM 03.22.13 \eqref{eq:shao.selfnorm}; 
\eqref{eq:shao.exact}; see Remark%RM 03.22.13 
s~\ref{re:shao-compar} and \ref{re:shao-compar2} for some details. 

The uniform and nonuniform bounds presented in Corollaries~\ref{cor:centralT} and \ref{cor:centralT.nub}, respectively, involve the sixth moments of $Y$, as they are based on the general results of Theorems~\ref{thm:iid,p=3,unif} and \ref{thm:iid,p=3,nonunif}, with $\De$ being on the order of magnitude of $\norm{S}^2/\si=\sqrt{n}\norm{\bar V}^2$. 

\begin{cor}[to Theorem~\ref{thm:iid,p=3,unif}]\label{cor:centralT} %RM %% supporting calculations found in notebook '02-24-13.unifbound.selfnormsum.nb' %%
Let $Y,Y_1,\dotsc,Y_n$ be i.i.d. r.v.'s, with $\E Y=0$ and $\norm{Y}_2=1$. 
Then
\begin{equation}\label{eq:centralT}
 \bigabs{\P(T_1\le z)-\Phi(z)}\le\frac{1}{\sqrt n}\Bigp{A_3\norm{Y}_3^3+A_4\norm{Y}_4^6+A_6\norm{Y^2-1}_3^3}
\end{equation}
for all $z\in\R$ and either triple
\begin{equation}\label{eq:Ajlist.Tunif}
 (A_3,A_4,A_6)\in\bigb{(2.99,2.99,0.15),%RM 03.22.13 
 \,(4.46,1.12,0.22)}. %RM 02-21-13
 %RM 12.21.12 \bigb{(3.33,3.33,0.17),(5.79,1.45,0.26)}. 
\end{equation} 
\end{cor}

The two triples $(A_3,A_4,A_6)$ in \eqref{eq:Ajlist.Tunif} are the result of trying to approximately minimize $A_3\vee(A_4/w_4)\vee(A_6/w_6)$, with weights $(w_4,w_6)\in\{(1,0.05),%RM 03.22.13 
\,(0.25,0.05)\}$. 

One can see that the constants in \eqref{eq:centralT}--\eqref{eq:Ajlist.Tunif} are not much worse than those in \eqref{eq:centralT.4th}--\eqref{eq:Ajlist.4th}. 

\begin{cor}[to Theorem~\ref{thm:iid,p=3,nonunif}]\label{cor:centralT.nub} %RM %% supporting calculations found in notebook '03-08-13.nonunifbound.selfnormsum.nb' %%
Let $\om\in\{0.1,0.5\}$, $w_g\in\{0,1\}$, and 
\begin{equation}\label{eq:g}
 g(z):=\frac1{z^3}+\frac{w_g}{e^{z/2}}.
\end{equation}
Then under the assumptions of Corollary~\ref{cor:centralT}, for all
\begin{equation}\label{eq:z,T}
 z\in(0,\om\sqrt n\,]
\end{equation}
one has
\begin{equation}\label{eq:centralT2}
 \bigabs{P(T_1\le z)-\Phi(z)}\le\frac{g(z)}{\sqrt{n}}\,\bigp{\hat A_3\norm{Y}_3^3+\hat A_4\norm{Y}_4^8+\hat A_6\norm{Y^2-1}_3^3},
\end{equation}
where, for any given pair $(\om,w_g)\in\{0.1,0.5\}\times\{0,1\}$, the triple $(\hat A_3,\hat A_4,\hat A_6)$ is either one of the two triples given in the corresponding block of Table~\ref{tab:Tbound} below.

\begin{table}[ht]\small
\begin{center}
\begin{tabular}{c||ccc|ccc}
&\multicolumn{3}{c|}{$\om=0.1$}&\multicolumn{3}{c}{$\om=0.5$}\\
&$\hat A_3$&$\hat A_4$&$\hat A_6$&$\hat A_3$&$\hat A_4$&$\hat A_6$\\\hline\hline
\multirow{2}{*}{$w_g=1$}
&38 &36 &36 &48 &48 &42\\
&39 &20 &7  &66 &33 &13\\\hline
\multirow{2}{*}{$w_g=0$}
&151&148&147&166&166&165\\
&169%RM 03.20.13 170
&85 &29 &229&115&45
\end{tabular}
\end{center}
\caption{Constants associated with nonuniform bound in \eqref{eq:centralT2}\label{tab:Tbound}}
\end{table}
\end{cor}

One can see that, especially in the case when $w_g=1$ and $\om=0.1$, the sum of the constants $\hat A_3$, $\hat A_4$, and $\hat A_6$ is comparable with the constant factor $30.2211$ in the nonuniform BE inequality \eqref{eq:tyurin,michel} for linear statistics; recall here also the asymptotic bound in \eqref{eq:asymp.nonunif}, with the same constant $30.2211$. 
One may also note that the constants $\hat A_3$, $\hat A_4$, and $\hat A_6$ in the case when $w_g=0$ are significantly greater than those for $w_g=1$. This reflects the fact that, whereas $\frac1{e^{z/2}}$ is much smaller than $\frac1{z^3}$ for very large $z>0$, the maximum of the ratio $\frac1{e^{z/2}}/\frac1{z^3}$ over %RM 03.22.13 
all $z>0$ is (attained at $z=6$ and) quite large, about $10.75$. 
Whereas at least some of the constants $\hat A_3$, $\hat A_4$, and $\hat A_6$ are rather large when $w_g=0$, one can put this into a perspective by recalling that, even in the much simpler case of sums of independent identically distributed r.v.'s, the first explicit constant in the nonuniform BE bound (obtained in \cite{pad78}) was greater than $1955$. 

Similarly to their counterparts in \cite{pin11}, the proofs of Corollaries~\ref{cor:centralT} and \ref{cor:centralT.nub} demonstrate a method by which one may obtain a variety of specific numerical constants for the bounds of the form \eqref{eq:centralT} and \eqref{eq:centralT2}. 
In particular, the introduction of the numerous parameters in Theorems~\ref{thm:iid,p=3,unif} and \ref{thm:iid,p=3,nonunif} allows one to account more accurately for the relations between the possible sizes of the various moments (cf.\ e.g.\ the ideas represented by \cite[Theorems~5.2, 6.1, 6.2]{pin94}). 
On the other hand, such an approach rather understandably results in significantly more complicated expressions. 

%RM 02-21-13 %% moved from end of subsection %%
\begin{remark}\label{re:symm.T1} %RM %% supporting calculations found in notebook '03-08-13.nonunifbound.selfnormsum.nb' %%
Suppose here that, in addition to the other condition of Corollary~\ref{cor:centralT}, the r.v.\ $Y$ is symmetric. 
Then, 
%RM 02-21-13 by Remark~\ref{re:symm.unif}, the triples $(A_3,A_4,A_6)$ of the constants in Corollary~\ref{cor:centralT} can be replaced by either one of the triples (of smaller constants) in the set 
% \begin{equation}\label{eq:Ajlist.symm}
%  \bigb{(3.09,3.09,0.16),(5.33,1.34,0.27)}. 
% \end{equation}
% Similarly, 
by Remark~\ref{re:symm.nonunif}, Table~\ref{tab:Tbound} can 
%RM 02-21-13 then 
be replaced by Table~\ref{tab:Tbound.symm}%RM 03.22.13 . 
, with somewhat better constants. 
The proof of this remark is contained in the proof of Corollary~\ref{cor:centralT.nub}. %RM 02-21-13 
\begin{table}[ht]\small
\begin{center}
\begin{tabular}{c||ccc|ccc}
&\multicolumn{3}{c|}{$\om=0.1$}&\multicolumn{3}{c}{$\om=0.5$}\\
&$\hat A_3$&$\hat A_4$&$\hat A_6$&$\hat A_3$&$\hat A_4$&$\hat A_6$\\\hline\hline
\multirow{2}{*}{$w_g=1$}
&35 &32 &31 &48 &48 &41\\
&37 &19 &5  &57 &29 &12\\\hline
\multirow{2}{*}{$w_g=0$}
&124&123%RM 03.20.13 124
&121&141&138&138\\
&145&73 &22 &205%RM 03.20.13 206
&103&42
\end{tabular}
\end{center}
\caption{Constants associated with nonuniform bound in \eqref{eq:centralT2}\label{tab:Tbound.symm} when $Y$ is symmetric.}
\end{table}
%RM 02-21-13 Proofs of these statements are provided in the proofs of the corresponding corollaries.
\end{remark}

\begin{remark}\label{re:shao-compar} %RM %% supporting calculations found in notebook `03-08-13.shao.compare.selfnormsum.nb' %%
The uniform bounds in \eqref{eq:centralT.4th} and \eqref{eq:centralT} (as well as the nonuniform one in \eqref{eq:centralT2}) involve moments of orders higher than $3$, in contrast with the uniform bound in \eqref{eq:shao.selfnorm}, say. 
However, it appears that the effect of the smaller constants in \eqref{eq:centralT.4th}--\eqref{eq:Ajlist.4th} and \eqref{eq:centralT}--\eqref{eq:Ajlist.Tunif} will oftentimes more than counterbalance the ``defect'' of the higher-order moments. 
For instance, suppose that $Y\sim\tilde t_d$, where $\tilde t_d$ denotes the standardized $t$ distribution with $d$ degrees of freedom, where $d$ is any positive real number. 
This distribution is symmetric. Its tails vary from very heavy ones for small $d$ to the very light tails of the standard normal distribution, corresponding to the limit case $d=\infty$. 
The absolute moments, say $m_s(d)$, of order $s$ of the distribution $\tilde t_d$ will be infinite for all $s\in[d,\infty)$.  
Then, in particular, the bound in \eqref{eq:shao.selfnorm} will be infinite if $d\le3$. 
On the other hand, 
one can show that for $Y\sim\tilde t_d$ the bound in \eqref{eq:centralT.4th} \big(say with the choice of the triple $(A_3,A_4,A_0)=(3.48,1.27,-1.43)$  %RM 03-05-13 (3.60,1.30,-1.44)%RM 12.21.12 (3.79,1.31,-1.32)$ 
in \eqref{eq:Ajlist.4th}\big) will be smaller than that in \eqref{eq:shao.selfnorm} for all real $d\ge 4.15$; %RM 03-05-13 4.16$; 
this can be checked using monotonicity properties of $m_3(d)$ and $m_4(d)$. 
Namely, $m_4(d)=\frac{3 (d-2)}{d-4}$ clearly decreases in $d>4$. As for $m_3(d)$, one can write $\sqrt{\frac{\pi }8} m_3(d)=r(\frac{d-3}2)$ for $d>3$, where $r(x):=\frac{\sqrt{x+\frac{1}{2}} \Gamma (x)}{\Gamma \left(x+\frac{1}{2}\right)}$. So, reasoning as in the proof of \cite[Lemma~2.1]{student-mono-v2}, one has $(\ln m_3(d)\big)'_d\,(d-3)=-\int_0^1\frac{t^{d-3}}{(t+1)^2}\, dt-\frac{1}{2 (d-2)}<0$ for all $d>3$, whence $m_3(d)$ decreases in $d>3$. 
Note also here that the bound in \eqref{eq:shao.selfnorm} will be nontrivial (that is, less than $1$) for some $d\in(0,4.15)$ %RM 03-05-13 4.16)$
only if $n>25^2\|Y\|_3^6=25^2m_3(d)>25^2m_3
(4.15)>4439$. %RM 03-05-13 (4.16)>4408$, 
%RM 03.20.13 where $m_3(d)$ stands for the third absolute moment of $\tilde t_d$. 

Similarly, the bound in \eqref{eq:centralT.4th} \big(again with $(A_3,A_4,A_0)=(3.48,1.27,-1.43)$ %RM 03-05-13 (3.60,1.30,-1.44)%RM 12.21.12 (3.79,1.31,-1.32)$
\big) will be smaller than that in \eqref{eq:shao.selfnorm} 
when $Y$ has any standardized two-point distribution which is not too skewed -- it is enough that $\P\big(Y=\sqrt{q/p}\,\big)=p$, $\P\big(Y=-\sqrt{p/q}\,\big)=q$, $0<p<1$, $q:=1-p$, and $p\wedge q\ge
0.0035$; %RM 03-05-13  0.00383$; 
moreover, if $p\wedge q<
0.0035$ %RM 03-05-13 0.00383$ 
then the bound in \eqref{eq:shao.selfnorm} will be nontrivial only if $n>25^2\|Y\|_3^6=25^2\big(\frac1{\sqrt{pq}}-2\sqrt{pq}\,\big)^2>
176707$. %RM 03-05-13 161322$. 
Note that any zero-mean distribution is a mixture of zero-mean two-point distributions \cite{disintegr}, so that 
such distributions appear to be of particular interest. 

%RM 12.21.12
% Discussion in Pinelis \cite{pin11} shows that one could utilize appropriate truncation to further alleviate the presence of the higher-order moments and thus make the comparison even more favorable to the bounds with the smaller constants. 
\end{remark}

%RM 12.21.12 
\begin{remark}\label{re:shao-compar2} %RM %% supporting calculations found in notebooks '2-point.compare2.nb', 'shaocompare.plots4.nb', 'student-dist.shaocompare3.nb', and 'pareto-dist.shaocompare3.nb'
\begin{figure}[h!t]
\centering
\subfloat[$n=10$]{\psfragfig[width=0.4\textwidth]{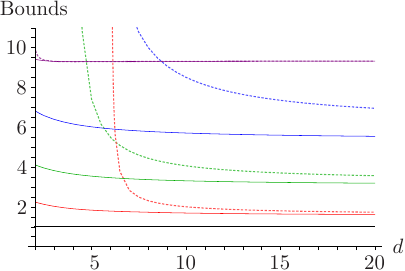}}\qquad\qquad
\subfloat[$n=100$]{\psfragfig[width=0.4\textwidth]{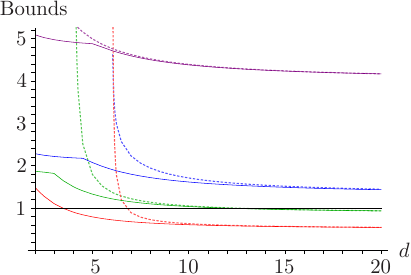}}\\
\subfloat[$n=1000$]{\psfragfig[width=0.4\textwidth]{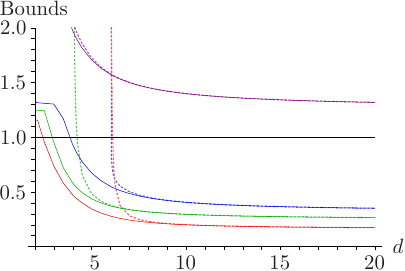}}\qquad\qquad
\subfloat[$n=10,000$]{\psfragfig[width=0.4\textwidth]{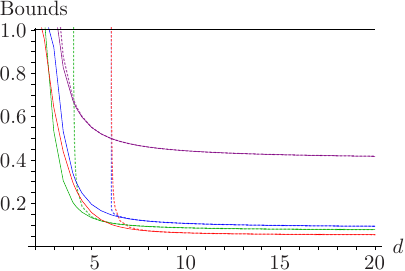}}
\caption{Comparison of bounds 
%IP 03.11.13 \eqref{eq:4thmoms} (green), \eqref{eq:iid,p=3} (blue), \eqref{eq:pin,T} (red), and \eqref{eq:shao.exact} (purple) 
\eqref{eq:pin,T} (red), \eqref{eq:4thmoms} (green), \eqref{eq:iid,p=3} (blue), and \eqref{eq:shao.exact} (purple)
when $Y$ has %IP 03.11.13 
the Student distribution with $d$ degrees of freedom 
%IP 03.11.13 (solid) or without (dotted) truncation}
-- using the postfactum truncation (solid) or not using it (dotted). }  
\label{fig:student}
\end{figure}

\begin{figure}[ht]
\centering
\subfloat[$n=10$]{\psfragfig[width=0.4\textwidth]{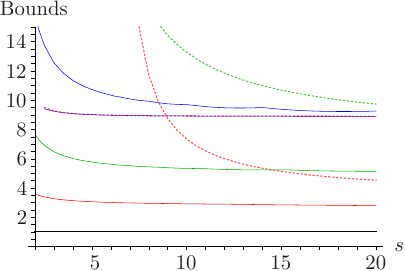}}\qquad\qquad
\subfloat[$n=100$]{\psfragfig[width=0.4\textwidth]{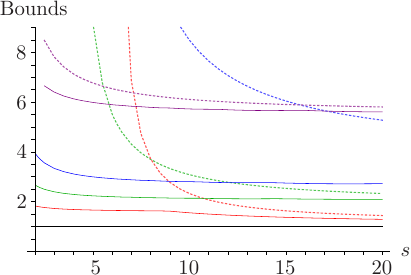}}\\
\subfloat[$n=1000$]{\psfragfig[width=0.4\textwidth]{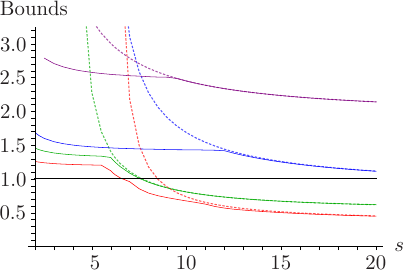}}\qquad\qquad
\subfloat[$n=10,000$]{\psfragfig[width=0.4\textwidth]{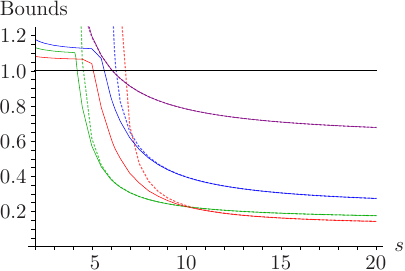}}
\caption{Comparison of bounds 
%IP 03.11.13 \eqref{eq:4thmoms} (green), \eqref{eq:iid,p=3} (blue), \eqref{eq:pin,T} (red), and \eqref{eq:shao.exact} (purple) 
\eqref{eq:pin,T} (red), \eqref{eq:4thmoms} (green), \eqref{eq:iid,p=3} (blue), and \eqref{eq:shao.exact} (purple)
when $Y$ has %IP 03.11.13 
the centered Pareto distribution with shape parameter $s$ %IP 03.11.13 
-- using the postfactum truncation (solid) or not using it (dotted). }
\label{fig:pareto}
\end{figure}

%IP 03.11.13 the main point of this remark: 
%
As was mentioned in Remark~\ref{re:moments,T1}, the potential ``defect'' caused by using higher order moments in our bounds for the central Student statistic \big(such as \eqref{eq:centralT.4th} and \eqref{eq:centralT}\big) can be eliminated or, at least, further reduced by an appropriate truncation, as suggested by \cite[Corollary~1.5]{pin11} and the discussion following it. Such a truncation may be referred to as ``postfactum truncation'' -- since it is done to the moments in the resulting bounds,  rather than in the proof of the bounds (which latter is the kind of technique usually employed to get rid of higher order moments). Looking at the comparisons made in \cite{pin11} and below in the present remark, it appears that the postfactum truncation may generally be more effective than the in-the-proof truncation; one possible reason for this advantage of the postfactum truncation is that it is sensitive to the underlying distribution of the observations, which seems to make sense, as knowledge of this distribution is needed anyway in order to compute the bounds. 

%IP 03.11.13
Yet another way to further improve the bounds in Corollaries~\ref{cor:centralT2} and \ref{cor:centralT} is to note the following. The last step in the proof of these bounds is the use of Young's inequality,  in order to eliminate products of different moments, and so, this step entails some loss in the accuracy. 
%
%The bounds of Corollaries~\ref{cor:centralT2} and \ref{cor:centralT} 
%are nearly direct applications of the inequalities \eqref{eq:ub} and \eqref{eq:iid,p=3,young's} to the self-normalized sum, respectively, though the use of Young's inequality (so as to obtain bounds without terms containing the product of different moments) entails some loss in the accuracy of our bounds. 
More accurate (and complicated in appearance) bounds on $\abs{\P(T_1\le z)-\Phi(z)}$ are given by \eqref{eq:4thmoms} (from which \eqref{eq:centralT2} is derived) and \eqref{eq:iid,p=3} (with $\tsi=1$, $\vsi_3=\norm{Y}_3$, $v_2=\norm{Y}_4^2$, $v_3=\norm{Y^4-Y^2+1}_{3/2}$, and $\Mf$ as defined in \eqref{eq:C1.selfnorm}). 
%IP 03.11.13
%Given any particular distribution for $Y$ along with a value of $n$, it is a straightforward task to numerically minimize either of the above two bounds in the two parameters $\cc$ and $\ep$%IP 03.11.13; we
%. 

Let us make a few graphical comparisons of the bounds \eqref{eq:4thmoms} and \eqref{eq:iid,p=3} to either of the bounds in \eqref{eq:shao.exact} or \eqref{eq:pin,T} (using the triple $(A_3,A_4,A_6)=(2,1,0.19)$, as found in the table at the end of the proof of Corollary~\ref{cor:centralT.nub}). 
Here let us consider the case  
%IP 03.11.13 such a distr. was consider in previous remark
%
%Consider first a two-point distribution, so that $Y$ takes the values $b$ and $-1/b$, with probabilities $1/(1+b^2)$ and $b^2/(1+b^2)$, respectively, for any $b\ge 1$. 
%Assume w.l.o.g.\ that $\sqrt{n}\ge2b$, since otherwise the bound in \eqref{eq:shao.exact} exceeds the trivial bound of $1$. Then the bound in \eqref{eq:shao.exact} is $25\E\abs{Y}^3/\sqrt{n}=25(b^4+1)/(b(1+b^2)\sqrt{n})$. Numerical minimization (in $\cc$ and $\vp$) suggests that the bound in \eqref{eq:4thmoms} or \eqref{eq:iid,p=3} is smaller than that in \eqref{eq:shao.exact} when $b\in[1,40.79]$ or $b\in[1,8.25]$, respectively (after using the previously assumed lower bound of $4b^2$ on $n$). That is, when the distribution is not too skewed, either of the two uniform bounds of Theorems~\ref{thm:ub} or \ref{thm:f(S).ub} applied to the self-normalized sum is smaller than the bound of Shao in \eqref{eq:shao.exact} for large enough $n$. 
%
%Next consider the bounds \eqref{eq:4thmoms}, \eqref{eq:iid,p=3}, \eqref{eq:pin,T}, and \eqref{eq:shao.exact} 
when the distribution of $Y$ has the Student distribution with $d$ degrees of freedom or the centered Pareto distribution with shape parameter $s$%IP 03.11.13. 
; the latter distribution has the density function $x\mapsto s (x + \frac s{s - 1})^{-(s + 1)}\I\{x>- \frac1{s - 1}\}$ for $s>1$. 
Plots of any of these four bounds (represented by the dotted curves) are found in Figures~\ref{fig:student} and \ref{fig:pareto} for $n\in\{10,100,1000,10000\}$ and $d\in[2,20]$ (or $s\in[2,20]$). 
%IP 03.11.13 
%While the bound of \eqref{eq:shao.exact} needs only $d>2$ (or $s>2$) in order to be finite, the finiteness of the bound in \eqref{eq:4thmoms} requires $d>4$ (or $s>4$), and the even greater restriction that $d>6$ (or $s>6$) is needed to assert the finiteness of the bounds in \eqref{eq:iid,p=3} or \eqref{eq:pin,T}. 
%This potential moment ``defect'' can be overcome by an appropriate truncation of the distribution in question, as suggested by \cite[Corollary~1.5]{pin11} and the discussion following it. 
The solid lines in Figures~\ref{fig:student} and \ref{fig:pareto} represent these four bounds after 
%IP 03.11.13 an optimal truncation 
a (numerically optimized) postfactum truncation is performed. 
The %IP 03.11.13 curves representing the 
bounds in \eqref{eq:4thmoms} and \eqref{eq:iid,p=3} have also been numerically minimized in $\cc$ and $\ep$. 
The remarks made in \cite{pin11} are also applicable here. Particularly, the effect of truncation in decreasing any of the bounds is most significant when the tails of the distribution are heavy (i.e.\ $d$ or $s$ is small). A general pattern to be found is that, when $n$ is large enough 
%IP 03.11.13 
for the (truncated or non-truncated) bounds to be 
%so that the (truncated or non-truncated) bounds are 
smaller than the trivial bound $1$, the smallest bound is that of \eqref{eq:pin,T}, followed by \eqref{eq:4thmoms}, then \eqref{eq:iid,p=3}, with %Shao's bound 
the bound in \eqref{eq:shao.exact} typically being the largest of the bounds under consideration. 
Again, we have the somewhat surprising result that the bounds presented in this paper, developed for a very general class of nonlinear statistics of which the self-normalized sum is but a single example, compare quite competitively with other bounds in the literature that were proven using methods %IP 03.11.13 which were 
tailored for the self-normalized sum.
\end{remark}

\begin{remark}\label{re:asymp-unif-compar}  %RM %% supporting calculations found in notebook 'novakcompare.asymp.selfnormsum.03-08-13.nb' %%
One may also want to compare, in the case of the statistic $T_1$, the asymptotic behavior of our bounds described in Corollary~\ref{cor:asymp} with the corresponding known asymptotic results. 
In particular, it follows from \cite[($\ast$)]{novak05} that 
\begin{equation}\label{eq:T1.as.unif.novak}
 %IP12.11.12
 d_{T_1,\as}:=\limsup_{n\to\infty}\,\sqrt n\,\sup_{z\in\R}|\P(T_1\le z)-\Phi(z)|  \le6.4\norm{Y}_3^3+2\norm{Y}_1  
\end{equation} 
whenever $\norm{Y}_3<\infty$. 
On the other hand, taking any real $\th>0$, $f(x_1,x_2):=\frac{x_1/\th}{\sqrt{1+x_2}}$ for $(x_1,x_2)\in\R\times(-1,\infty)$, and $V:=(\th\,Y,Y^2-1)$, one has $\sqrt{n}f(\bar V)=T_1$. 
Choose now 
$\th=\sqrt{\norm{Y}_4^4-1}$ %RM 02-21-13 $\th=\|Y\|_4^4-1$ 
(assuming that $\|Y\|_4\ne1$ and hence $\|Y\|_4>1$; the case $\|Y\|_4=1$ can then be treated by continuity, say). 
Then, by \eqref{eq:y_*}, $y_*=((2/\pi)^{1/6}+\norm{Y}_3)%RM 02-21-13 (3+2\|Y\|_3)
\sqrt{\|Y\|_4^4-1}$. 
Using this expression for $y_*$, 
one can show that 
the bound %RM 02-21-13 both bounds 
in \eqref{eq:asymp.unif} will be smaller than that in \eqref{eq:T1.as.unif.novak} \big(and even smaller than $6.4\norm{Y}_3^3$\big) 
whenever $\norm{Y}_4\le2.189\norm{Y}_3$. 
%IP 03.11.13 
In view of \cite[Corollary 1.3 (ii)]{student-mono-v2}, this 
%This 
will be the case %RM 03-05-13 -- say, in the case 
when $Y$ has the standardized $t$ distribution with $d$ degrees of freedom, for any real $d>6$.  
The same %IP 03.11.13
conclusion about the bounds 
in \eqref{eq:asymp.unif} and \eqref{eq:T1.as.unif.novak} 
will be true when $Y$ has any standardized two-point distribution which is not too skewed -- it is enough that $\P\big(Y=\sqrt{q/p}\,\big)=p$, $\P\big(Y=-\sqrt{p/q}\,\big)=q$, $0<p<1$, $q:=1-p$, and $p\wedge q\ge8.3\times 10^{-5}$. %RM 03-05-13 0.00272$. %RM 02-21-13 .00906$. 
Note also that in the case of the  statistic $T_1$ one can get an asymptotic bound better than the one just obtained based on Corollary~\ref{cor:asymp} (which latter is derived from Theorem~\ref{thm:f(S).ub}, which in turn is a corollary to Theorem~\ref{thm:ub}) -- 
if instead one uses Theorem~\ref{thm:ub} directly;  
cf.\ 
Corollary~\ref{cor:centralT2} (to Theorem~\ref{thm:ub}) vs.\ 
Corollary~\ref{cor:centralT} (to Theorem~\ref{thm:iid,p=3,unif}). 

%IP12.11.12

% \parbox{.3\textwidth}{
% \rule{0pt}{0pt}\hspace*{-10pt}\psfragfig[width=0.3\textwidth]{novak00}%RM 02-21-13 \includegraphics[width=.3\textwidth]{novak00compar.pdf}
% \vspace*{4pt}
% }
% \parbox{.67\textwidth}{
% In a paper preceding \cite{novak05}, Novak \cite[(5.8)]{novak00} obtained a bound which, \break 
% taken together with \eqref{eq:tyurin,michel}, implies that 
% \begin{equation}\label{eq:novak}
%  d_{T_1,\as}\le
%    C_{\nov,3}\norm{Y}_3^3+C_{\nov,4}\sqrt{\norm{Y}_4^4-1}, 
% \end{equation} 
% where $C_{\nov,3}:=0.4748+\frac8{e\sqrt{2\pi}}=1.648\dots$ and $C_{\nov,4}:=(\frac2{e\pi})^{1/4}=0.695\dots$. 
% %RM 12.21.12 One can show that the bound in \eqref{eq:asymp.unif} in the case of the statistic $T_1$ is always greater than the bound in \eqref{eq:novak}. 
% In the picture on  
% }

%RM 03-05-13 
\begin{minipage}{0.34\textwidth}
\vspace*{5pt}
% \begin{center}
 \includegraphics[width=\textwidth]{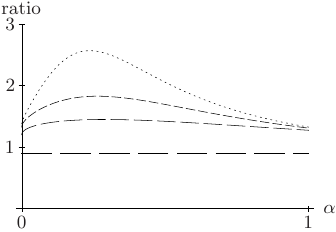}
%  \psfragfig*[width=0.33\textwidth]{novak00}
% \end{center}
 {\small Ratio of bound in \eqref{eq:asymp.unif} to \eqref{eq:novak}, where $\al:=\frac{\norm{Y}_3-1}{\norm{Y}_4-1}$ and $\norm{Y}_4=1+$\break (long-dash), 2 (medium-dash), 4 (short-dash), and 10 (dotted)}
\vspace{8pt}
\end{minipage}
\hfill
\begin{minipage}{0.62\textwidth}
\ \\
In a paper preceding \cite{novak05}, Novak \cite[page 424]{novak00} %RM 03-05-13 \cite[(5.8)]{novak00} 
obtained a bound which, %\break 
taken together with \eqref{eq:tyurin,michel}, implies that 
\begin{equation}\label{eq:novak}
 d_{T_1,\as}\le
   C_{\nov,3}\norm{Y}_3^3+C_{\nov,4}\sqrt{\norm{Y}_4^4-1}, 
\end{equation} 
where 
$C_{\nov,3}:=0.4748+\frac8{e\sqrt{2\pi}}=1.648\dots$ and $C_{\nov,4}:=(\frac2{e\pi})^{1/4}=0.695\dots$. 
%RM 12.21.12 One can show that the bound in \eqref{eq:asymp.unif} in the case of the statistic $T_1$ is always greater than the bound in \eqref{eq:novak}. 
In the picture on the left, one can see the graphs of the ratio%IP 03-11-13 s 
\ of these two bounds %in \eqref{eq:asymp.unif} for $T_1$ to that in \eqref{eq:novak} 
for $\E Y^4=1+$, $2^4$, $4^4$, and $10^4$ as functions of $\al:=(\|Y\|_3-1)/(\|Y\|_4-1)\in(0,1)$. %RM 02-21-13 ; these four graphs are shown, respectively, as long-dash, medium-dash, short-dash, and dotted lines
The limit of this ratio is approximately $0.8944$ when $\E Y^4$ approaches 1 from above. %RM 12.21.12
Thus, the bound in \eqref{eq:asymp.unif} for $T_1$ is usually moderately greater than the bound in \eqref{eq:novak}. 
%IP12.11.12 ???
On the other hand, in contrast with the general methods used in the present paper, 
the methods used in \cite{novak00} were specialized to target certain ratio-type statistics. Also, the non-asymptotic bounds in \cite[(5.6) and 
(5.7) %RM 03-05-13 (5.6*)
]{novak00} on which the asymptotic bound in \eqref{eq:novak} is based, were of a very complicated structure, 
%IP 03.11.13   
with further work needed to bound the various terms and choose explicit values of the parameters. 
\end{minipage}
\end{remark}

%IP12.11.12
% \vspace*{13pt}

\begin{remark}\label{re:asymp-nonunif-compar} 
Consider now the asymptotic behavior of the nonuniform bound for $T_1$. 
Novak \cite[Theorem~10]{novak00} provides an explicit, though complicated in appearance, nonuniform BE-type bound for this statistic. 
Using \cite[(5.10)]{novak00} and \eqref{eq:tyurin,michel} (and still assuming that $\E Y=0$ and $\E Y^2=1$, as well as $\E Y^4<\infty$) one can show that
\begin{equation}\label{eq:novak_nu,asymp}
 \limsup_{n\to\infty}\,\sup_{g(n)\le z\le n^{1/6}}z^3\sqrt{n}\bigabs{\P(T_1\le z)-\Phi(z)}\le30.2211\norm{Y}_3^3;  
\end{equation}
here $g$ stands for any positive increasing unbounded function on $\N$. 
Thus, for the specific statistic $T_1$, the asymptotic bound in \eqref{eq:novak_nu,asymp} coincides with that in \eqref{eq:asymp.nonunif}, obtained for general nonlinear statistics of the form $f(\bar V)$. 
Note also that the bound in \eqref{eq:asymp.nonunif} 
holds for $z$ in the zone $[g(n),n^{1/2}/g(n)]$, which is much wider than the zone $[g(n),n^{1/6}]$ in \eqref{eq:novak_nu,asymp} if $g$ is taken to grow slowly enough. 
On the other hand, Theorem~\ref{thm:iid,p=3,unif} and then Corollary~\ref{cor:asymp} contain the moment condition $v_3<\infty$, which is equivalent, in the specific case of $T_1$, to $\E Y^6<\infty$, which is more stringent than the corresponding condition $\E Y^4<\infty$ used here to derive \eqref{eq:novak_nu,asymp}.  
\end{remark}  

\subsection{Pearson's \texorpdfstring{$R$}{R}}

Let $(Y,Z),(Y_1,Z_1),\ldots,(Y_n,Z_n)$ be a sequence of i.i.d.\ random points in $\R^2$, with
\begin{equation*}
 \var Y\in(0,\infty)\quad\text{and}\quad\var Z\in(0,\infty).
\end{equation*}
Recall the definition of Pearson's product-moment correlation coefficient: 
\begin{equation}\label{eq:Rdef}
 R:=\frac{\sum_{i=1}^n(Y_i-\bar Y)(Z_i-\bar Z)}
  {\sqrt{\sum_{i=1}^n(Y_i-\bar Y)^2}\,
   \sqrt{\sum_{i=1}^n(Z_i-\bar Z)^2}}
  =\frac{\ol{YZ}-\bar Y\ \bar Z}{\sqrt{\ol{Y^2}-\bar Y^2}\sqrt{\ol{Z^2}-\bar Z^2}},
\end{equation}
where 
\begin{equation*}
 \bar Y:=\tfrac1n\tsum\nolim_iY_i,\quad
 \bar Z:=\tfrac1n\tsum\nolim_iZ_i,\quad
 \bar{Y^2}:=\tfrac1n\tsum\nolim_iY_i^2,\quad
 \bar{Z^2}:=\tfrac1n\tsum\nolim_iZ_i^2,\quad\text{and}\quad
 \bar{YZ}:=\tfrac1n\tsum\nolim_iY_iZ_i;
\end{equation*}
let %IP12.13.12 us allow $R$ to take an arbitrary value $r_0$ 
$R:=0$ 
if the denominator in \eqref{eq:Rdef} is $0$. 
Note that $R$ is invariant under all affine transformations of the form $Y_i\mapsto a+bY_i$ and $Z_i\mapsto c+dZ_i$ with positive $b$ and $d$; so, in what follows we may (and shall) assume that the r.v.'s $Y$ and $Z$ are standardized: 
\begin{equation*}
 \E Y=\E Z=0\quad\text{and}\quad\E Y^2=\E Z^2=1,\quad\text{and we let}\quad\rho:=\E YZ=\operatorname{\mathbb{C}orr}(Y,Z).
\end{equation*}

Let $\XX=\R^5$, and for $\x=(x_1,x_2,x_3,x_4,x_5)\in\XX$ such that $(1+x_3-x_1^2)(1+x_4-x_2^2)>0$, let 
\begin{equation}\label{eq:f.R}
 f(\x)=f(x_1,x_2,x_3,x_4,x_5)=\frac{x_5+\rho-x_1x_2}{\sqrt{1+x_3-x_1^2}\sqrt{1+x_4-x_2^2}}-\rho;
\end{equation}
let $f(\x):=%IP12.13.12 r_0
-\rho$ for all other $\x\in\XX$. 
Recall \eqref{eq:d(ep)} to see that $f''(\x)$ exists and is continuous on the closed $\ep$-ball about the origin for any fixed $\ep\in(0,\sqrt{3}/2)$; 
then the smoothness condition \eqref{eq:smooth} holds, with $L(\x)=f'(0)(x_1,x_2,x_3,x_4,x_5)=-\rho x_3/2-\rho x_4/2+x_5$. 
Letting $V=\bigp{Y,Z,Y^2-1,Z^2-1,YZ-\rho}$, so that $L(V)=YZ-\tfrac\rho2\bigp{Y^2+Z^2}$, we see that $f(\bar V)=R-\rho$. 
Then %IP16 Corollary
Theorem~\ref{cor:iid} immediately yields

\begin{thm}\label{thm:R}
Take any $\om>0$ and assume that $\tsi>0$ and $v_3<\infty$. 
Then for all $z\in\R$ and $n\in\N$ 
\begin{equation}\label{eq:R.bound.uni}
 \Bigabs{\P\Bigp{\frac{R-\rho}{\tsi/\sqrt{n}}\le z}-\Phi(z)}\le\frac{\CC}{\sqrt n}, 
\end{equation}
where $\CC$ is a finite expression depending only on the distribution of the random point $(Y,Z)$; also, for all real $z>0$ and $n\in\N$ satisfying \eqref{eq:z,iid} 
\begin{equation}\label{eq:R.bound}
 \Bigabs{\P\Bigl(\frac{R-\rho}{\tsi/\sqrt{n}}\le z\Bigr)-\Phi(z)}\le\frac{\CC}{z^3\,\sqrt n}, 
\end{equation}
where $\CC$ is a finite expression depending only on $\om$ and the distribution of $(Y,Z)$. 
\end{thm}

\begin{remark}\label{re:R.si=0}
Note that the \emph{degeneracy} condition $\tsi =0$ is equivalent to the following: there exists some $\ka\in\R$ such that the random point $(Y,Z)$ lies a.s.\ on the union of the two straight lines through the origin with slopes $\ka$ and $1/\ka$ (for $\ka=0$, these two lines should be understood as the two coordinate axes in the plane $\R^2$). 
Indeed, if $\tsi =0$, then $YZ-\tfrac\rho2(Y^2+Z^2)=0$ a.s.; solving this equation for the slope $Z/Y$, one obtains two roots, whose product is $1$. 
Vice versa, if $(Y,Z)$ lies a.s.\ on the union of the two lines through the origin with slopes $\ka$ and $1/\ka$, then $YZ=\frac r2(Y^2+Z^2)$ a.s.\ for $r:=2\ka/(\ka^2+1)$ and, moreover, $r=\E\frac r2(Y^2+Z^2)=\E YZ=\rho$. 

For example, let the random point $(Y,Z)$ equal $(cx,\ka cx)$, $(-cx,-\ka cx)$, $(\ka cy,cy)$, $(-\ka cy,-cy)$ with probabilities $\frac p2$, $\frac p2$, $\frac q2$, $\frac q2$, respectively, where $x\ne0$, $y\ne0$, $\ka\in\R$, $c:=\sqrt{\dfrac{x^{-2}+y^{-2}}{\ka^2+1}}$, $p:=\dfrac{y^2}{x^2+y^2}$, and $q:=1-p$; then $\tsi=0$ (and the r.v.'s $Y$ and $Z$ are standardized). 
In particular, one can take here $x=y=1$, so that $p=q=\frac12$. 
\end{remark}

\begin{remark}\label{re:R.6th.moments}
In order to get a uniform bound of order $\bigO(1/\sqrt{n})$ in Theorem~\ref{thm:R}, it is necessary to assume that $v_3<\infty$, which is equivalent to $\norm{Y}_6+\norm{Z}_6<\infty$. 
This moment condition might seem overly restrictive, since only third absolute moments are required to obtain a BE-type bound of the same order for linear statistics (or even for the central Student statistic). 
However, the moments $\norm{Y}_6$ and $\norm{Z}_6$ do appear in an asymptotic expansion (up to an order $n^{-1/2}$) of the distribution of $R$ when $\rho\ne0$; cf.\ Remark~\ref{re:moments,T1}; for details, one can see \cite{corr.asymp}. 
When $\rho=0$, the most restrictive moment assumption for the existence of the asymptotic expansion is that $\norm{YZ}_3<\infty$. 
\end{remark}

\begin{remark}\label{re:R.asymp}
%RM01.12
% The upper bound in \eqref{eq:R.bound.uni} is optimal in its dependence on $n$. 
% Indeed, suppose that a function $f\colon\R^k\to\R$ is twice continuously differentiable in a neighborhood of the origin (so that $f$ satisfies the smoothness condition \eqref{eq:smooth}), and let $L$ and $H$ denote here the gradient vector and Hessian matrix of $f$ at 0. 
% Further assume, in addition to the assumptions $\tsi>0$ and $v_3<\infty$, that $V$ satisfies the Cram\'er-type condition $\limsup_{\norm{t}\to\infty}\abs{\E e^{it\T V}}<1$. 
% Then a calculation of the asymptotic distribution of $\sqrt{n}f(\bar V)/\tsi$ using \cite[Theorem~2]{bhatt78} implies
% \begin{equation}\label{eq:BG}
%  \sup_{z\in\R}\biggabs{\Bigprob{\frac{f(\bar V)}{\tsi/\sqrt{n}}\le z}-\Phi(z)-\frac{\De(z)}{\sqrt{n}}}=o\Big(\frac1{\sqrt{n}}\Big),
% \end{equation}
% where
% \begin{equation}\label{eq:BG.De}
%  \De(z):=-\Bigb{\Bigp{\frac{\E[(L\T V)^3]}{6\tsi^3}+a_3}(z^2-1)+a_1}\varphi(z),
% \end{equation}
% \begin{equation*}
%  a_1:=\frac1{2\tsi}\,\tr(H\Si),\quad a_3:=\frac1{4\tsi^3}(L\T\Si L-\tsi^2)\tr(H\Si)+\frac1{2\tsi^3}\,L\T\Si H\Si L,
% \end{equation*}
% $\Si$ denotes the covariance matrix of $V$, and $\varphi$ is the standard normal density. 
Recall the asymptotic distribution results of Bhattacharya and Ghosh \cite{bhatt78} as outlined in Remark~\ref{re:t.asymp}. %RM01.12
In the conditions of Theorem~\ref{thm:R}, take now the very simple case when $Y$ and $Z$ are zero-mean, unit-variance, absolutely continuous r.v.'s independent of each other. 
Then a straightforward calculation shows that $a_1=0$, $a_3=0$, and hence $\De(z)=-\frac16\E Y^3\E Z^3(z^2-1)\varphi(z)$. 
So, the bound in \eqref{eq:R.bound.uni} has an optimal dependence on $n$ whenever $\E Y^3\ne0$ and $\E Z^3\ne0$. 
Moreover, since $\De(z)$ is real-analytic in $z$, $L$, $H$, and moments of $V$, we see that generally $\De(z)\ne0$ and hence the bound in \eqref{eq:R.bound.uni} is generally of the optimal order in $n$. 
\end{remark}

The bounds in \eqref{eq:R.bound.uni} and \eqref{eq:R.bound} appear to be new. 
In fact, we have not been able to find in the literature any uniform (or nonuniform) bound on the closeness of the distribution of $R$ to normality. 
Note that such bounds are important in considerations of the asymptotic relative efficiency of statistical tests; see e.g.\ Noether \cite{noe55}. 
Shen~\cite{shen07} recently provided results concerning probabilities of large deviations for $R$ in the special case when $(Y,Z)$ is a bivariate normal r.v. 
Formal asymptotic expansions for the density of $R$ follow from the paper by Kollo and Ruul~\cite{kollo03}.

We next state one particular simplification of the uniform bound in \eqref{eq:iid,p=3} when applied to the Pearson statistic in the case when $\rho=0$. 

\begin{cor}[to Theorem~\ref{thm:iid,p=3,unif}]\label{cor:Runif}
Assume that $\E YZ=0$ and $\tsi=\norm{YZ}_2>0$. Then for all $z\in\R$ and $n\in\N$
\begin{equation}\label{eq:Runif}
 \Bigabs{\Bigprob{\frac{R}{\tsi/\sqrt n}\le z}-\Phi(z)}\le\frac{B_0+B_3/\tsi^3}{\sqrt{n}}\bigp{\norm{Y}_6^6+\norm{Z}_6^6},
\end{equation}
where $(B_0,B_3)$ is any ordered pair in the set
\begin{equation}\label{eq:Bjlist,Runif}
\bigb{(3.61,3.61),%RM 03.22.13 
\,(1.12,8.94),%RM 03.22.13 
\,(13.33,1.69),%RM 03.22.13 
\,(0.56,14.97),%RM 03.22.13 
\,(36.32,1.37)}. %RM 02-21-13
%RM 12.21.12 \bigb{(4.08,4.08),(1.38,11.02),(14.73,1.85),(0.71,19.16),(39.22,1.47)}.
\end{equation}
\end{cor}

Similarly to the proof of Corollary~\ref{cor:centralT}, that of Corollary~\ref{cor:Runif} gives a method by which one may obtain a variety of values for the pair $(B_0,B_3)$. 
The specific pairs listed in \eqref{eq:Bjlist,Runif} are obtained by trying to minimize $B_0\vee B_3/\tsi^3$ for $\tsi\in\{1,2,1/2,3,1/3\}$. 

%RM 02-21-13
% \begin{remark}\label{re:R.imp}
% By employing the improvement described in Remark~\ref{re:chen.improvement}, it is possible to slightly decrease the constants in \eqref{eq:Bjlist,Runif}; namely, \eqref{eq:Runif} holds for any ordered pair $(B_0,B_3)$ in the set
% \begin{equation*}
%  \bigb{(4.06,4.06),(1.38,10.99),(14.62,1.83),(0.71,19.14),(38.70,1.44)}. 
% \end{equation*}
% \end{remark}

%RM 02-21-13
% \begin{remark}\label{re:R.symm}
% By Remark~\ref{re:symm.unif}, in the case when the r.v.\ $YZ$ is symmetrically distributed, 
% the constants in \eqref{eq:Bjlist,Runif} can be slightly improved; namely, then \eqref{eq:Bjlist,Runif} can be replaced by  
% \begin{equation}\label{eq:R.symm}
%  \bigb{(3.74,3.74),(1.31,10.47),(13.19,1.65),(0.68,18.27),(36.40,1.35)}. %RM 12.21.12 \bigb{(3.72,3.72),(1.31,10.44),(13.49,1.69),(0.68,18.24),(37.24,1.38)}.
% \end{equation}
% \end{remark}

\begin{remark}\label{re:fisher-z} 
Bounds similar to the ones in Corollary~\ref{cor:Runif} can be obtained, e.g., for other statistics related to Pearson's $R$, including the Fisher $z$ transform. 
However, for reasons discussed in Appendix~\ref{app:z} and because the paper is already quite long, we chose not to present such results here. 
\end{remark}

\subsection{Non-central Hotelling's \texorpdfstring{$T^2$}{T2} statistic}
\label{subsec:Hotelling}

Let $k\ge2$ be an integer, and let $Y,Y_1,\dotsc,Y_n$ be i.i.d.\ r.v.'s in $\R^k$, with finite 
\begin{equation*}
 \mu:=\E Y\quad\text{and}\quad\cov Y=\E YY\T-\mu\mu\T\text{ strictly positive definite}.
\end{equation*}
Consider Hotelling's $T^2$ statistic
\begin{equation}\label{eq:T^2}
 T^2:=\Ybar\T(S_Y^2/n)^{-1}\Ybar
  =n\Ybar\T\Bigp{\YYTbar-\Ybar\,\Ybar\T}^{-1}\Ybar,
\end{equation}
where
\begin{equation*}
 \Ybar:=\tfrac1n\tsum\nolim_iY_i,\quad
 \YYTbar:=\tfrac1n\tsum\nolim_iY_iY_i\T,\quad\text{and}\quad
 S_Y^2:=\tfrac1n\tsum\nolim_i\bigp{Y_i-\Ybar}\bigp{Y_i-\Ybar}\T
  =\YYTbar-\Ybar\,\Ybar\T;
\end{equation*}
the generalized inverse is often used in place of the inverse in \eqref{eq:T^2}, though here we may %IP12.13.12 allow $T^2$ to take any value $t_0^2$ 
just let $T^2:=0$ whenever $S_Y^2$ is singular. 
Also note that $S_Y^2$ is defined as the empirical covariance matrix of the sample $(Y_i)_{i=1}^n$, rather than the sample covariance matrix $\frac{n}{n-1}S_Y^2$.
Call $T^2$ ``central'' when $\mu=0$ and ``non-central'' otherwise. 

For any nonsingular matrix $B$, $T^2$ is invariant under the invertible transformation $Y_i\mapsto BY_i$, so let us assume w.l.o.g.\ that
\begin{equation*}
 \cov Y=\Id,
\end{equation*}
the $k\times k$ identity matrix. 

%RM12.26 copied paragraph from ``ejs2.tex''
Now let $\XX=\bigb{(x_1,x_2)\colon x_1\in\R^k,x_2\in\R^{k\times k}}$ be equipped with the norm 
\begin{equation}\label{eq:k,k-by-k norm}
  \norm{(x_1,x_2)}:=\sqrt{\norm{x_1}^2+\norm{x_2}_F^2}, 
\end{equation}
where $\norm{x_2}_F:=\sqrt{\tr(x_2x_2\T)}$ is the Frobenius norm. 
For $\x=(x_1,x_2)\in\XX$ such that $\Id+x_2-x_1x_1\T$ is nonsingular, let
\begin{equation*}
 f(\x)=(x_1+\mu)\T\bigp{\Id+x_2-x_1x_1\T}^{-1}(x_1+\mu)-\mu\T\mu,
\end{equation*}
and let $f(\x):=-\mu\T\mu$ for all other $\x\in\XX$. 
The Fr\'echet derivative of $f$ at the origin is the linear functional defined by $L(\x)=f'(0)(x_1,x_2)=2x_1\T\mu-\mu\T x_2\mu$. 
Let us recall a couple of other useful facts (found in, say, the monograph \cite{horn85}): the spectral norm $\norm{B}$ of any $k\times k$ matrix $B$ does not exceed $\norm{B}_F$, and $\norm{B}<1$ implies $\Id-B$ is nonsingular and $\norm{(\Id-B)^{-1}}\le1/(1-\norm{B})$. 
In particular, 
\begin{equation*} \norm{x_1x_1\T-x_2}\le\norm{x_1x_1\T-x_2}_F\le\norm{x_1x_1\T}_F+\norm{x_2}_F=\norm{x_1}^2+\norm{x_2}_F<1
\end{equation*}
for any $\x$ in the closed $\ep$-ball about the origin and any fixed $\ep\in(0,\sqrt{3}/2)$ (which again follows from \eqref{eq:d(ep)}), so that the smoothness condition \eqref{eq:smooth} holds. 
Upon letting 
\begin{equation*}
 V=\bigp{Y-\mu,\,(Y-\mu)(Y-\mu)\T-\Id}, 
\end{equation*}
we see that $nf(\bar V)=T^2-n\mu\T\mu$. 
Then %IP16 Corollary
Theorem~\ref{cor:iid} immediately yields 

\begin{thm}\label{thm:T^2}
Take any $\om>0$ and assume that $\tsi>0$ and $v_3<\infty$. 
Then for all $z\in\R$ and $n\in\N$ 
\begin{equation}\label{eq:T^2.bound.uni}
 \Bigabs{\Bigprob{\frac{T^2-n\mu\T\mu}{\tsi\sqrt{n}}\le z}-\Phi(z)}\le\frac{\CC}{\sqrt n},
\end{equation}
where $\CC$ is a finite expression depending only on the distribution of $Y$; also, for all real $z>0$ and $n\in\N$ satisfying \eqref{eq:z,iid} 
\begin{equation}\label{eq:T^2.bound}
 \Bigabs{\Bigprob{\frac{T^2-n\mu\T\mu}{\tsi\sqrt{n}}\le z}-\Phi(z)}\le\frac{\CC}{z^3\,\sqrt n}, 
\end{equation}
where $\CC$ is a finite expression depending only on $\om$ and the distribution of $Y$. 
\end{thm}

\begin{remark}\label{re:T^2.si=0}
The non-degeneracy condition $\tsi >0$ immediately implies that $\mu\ne0$, so that Theorem~\ref{thm:T^2} is applicable only to the non-central $T^2$. 
If $\mu\ne0$, then $\tsi =0$ if and only if $(Y-\mu)\T\mu=1\pm\sqrt{1+\|\mu\|^2}$ a.s., that is, if and only if $\P(Y\T\mu=x_1)=1-\P(Y\T\mu=x_2)=p$, where
\begin{equation*}
 x_1=1+\|\mu\|^2+\sqrt{1+\|\mu\|^2},\quad x_2=1+\|\mu\|^2-\sqrt{1+\|\mu\|^2},\quad
 p=\frac12\Bigl(1-\frac{1}{\sqrt{1+\|\mu\|^2}}\Bigr);
\end{equation*}
in other words, $\tsi =0$ if and only if $Y$ lies a.s.\ in the two hyperplanes defined by $Y\T\mu=x_1$ or $Y\T\mu=x_2$. 
Note the similarity to the degeneracy condition of Student's $T$ statistic described in Remark~\ref{re:t.si=0}. 
Recalling the conditions $\E Y=\mu$ and $\cov Y=I$, we have $\tsi =0$ if and only if
\begin{equation*}
 Y=\xi\,\frac{\mu}{\norm{\mu}}+\tilde Y\ \ \text{a.s.},
\end{equation*}
where
\begin{equation*}
 \xi=\frac{2\sqrt{p(1-p)}}{1-2p}+B_p\ \text{for some }p\in(0,\tfrac12),
\end{equation*}
and $\tilde Y$ is a random vector in $\R^k$ such that $\E\tilde Y=0$, $\E\xi\tilde Y=0$, $\tilde Y\T\mu=0$ a.s., and $\cov\tilde Y$ is the orthoprojector onto the hyperplane $\{\mu\}^\perp:=\{\x\in\R^k\colon\x\T\mu=0\}$.
\end{remark}

%RM01.12
\begin{remark}\label{re:T^2.asymp}
Using again the asymptotic expansion results of Bhattacharya and Ghosh \cite[Theorem~2]{bhatt78} (cf.\ Remark~\ref{re:R.asymp}), we can show that generally 
the upper bound in \eqref{eq:T^2.bound.uni} has an optimal dependence on $n$ as well. 
For instance, 
consider the simple case when $Y=(Y_1,Y_2)$, where $Y_1$ and $Y_2$ have absolutely continuous distributions and are independent of one another; further suppose that $\E Y_2=0$, $\E Y_1=\mu_1\ne0$, and that $Y_1$ is symmetric, so that $\E(Y_1-\mu_1)^m=0$ for odd natural $m$. 
Then, for $\De(z)$ as in \eqref{eq:BG.De}, 
\begin{equation*}
 \De(1)=-\frac{\mu_1^2(\nu_4+1)+2}{\tsi}\,\varphi(1)\quad\text{and}\quad\tsi=|\mu_1|\,\sqrt{(\nu_4-1)\mu_1^2+4}, 
\end{equation*}
with $\nu_4:=\E(Y_1-\mu_1)^4\ge\E^2(Y_1-\mu_1)^2=1$. So, $\tsi>0$ and $\De(1)\ne0$. 
Thus, the dependence of the upper bound in \eqref{eq:T^2.bound.uni} on $n$ is optimal.
\end{remark}

Again, the bounds in \eqref{eq:T^2.bound.uni} and \eqref{eq:T^2.bound} appear to be new; we have found no mention of BE bounds for $T^2$ in the literature. 
Probabilities of moderate and large deviations for the central Hotelling $T^2$ statistic (when $\mu=0$) were considered by Dembo and Shao \cite{dem06}.
Asymptotic expansions for the generalized $T^2$ distribution for \emph{normal populations} were given by It\^o \cite{ito56} (for $\mu=0$), and by It\^o \cite{ito60}, Siotani \cite{sio71}, and Muirhead \cite{muir72} (for any $\mu$); Kano \cite{kano95} and Fujikoshi \cite{fuji97} give an asymptotic expansion for the distribution of the central $T^2$ for non-normal populations, and Kakizawa and Iwashita \cite{kaki08} do this for the noncentral $T^2$ statistic.

%RM12.26
\subsection{Covariance test statistics}
\label{subsec:cov}

For any natural $k\ge2$, let $Y,Y_1,\dotsc,Y_n$ be i.i.d.\ r.v.'s in $\R^k$ with
\begin{equation*}
 \E Y=0\quad\text{and}\quad\Sigma:=\cov Y=\E YY\T >0. 
\end{equation*} 
Further let
\begin{equation}\label{eq:det,tr}
 \al:=\tr(\Si)/k,\quad\be:=\det(\Si)^{1/k},\quad\de:=\sqrt{\tr[(\Si-\al I)^2]/k}
\end{equation}
be the arithmetic mean, geometric mean, and standard deviation, respectively, of the eigenvalues of $\Si$; the assumption that $\Si>0$ implies $\al>0$ and $\be>0$. 

We consider here a few statistics used to test either the null hypothesis of sphericity \break ($H_{0,1}\colon\Si=\si^2I$ for some unknown $\si^2>0$) or the null hypothesis of the identity covariance ($H_{0,2}\colon\Si=I$).
Each of these statistics is a smooth function of the sample covariance matrix 
\begin{equation*} S:=\YYTbar-\Ybar\,\Ybar\T,\quad\text{where}\quad\Ybar:=\tfrac1n\tsum\nolim_iY_i,\quad\text{and}\quad\YYTbar:=\tfrac1n\tsum\nolim_iY_iY_i\T. 
\end{equation*}
In turn, $S$ is a smooth function of the zero-mean r.v.
\begin{equation*}
 \bar V:=\tfrac1n\,\tsum\nolim\nolimits_i V_i=(\Ybar,\YYTbar-\Si), \quad\text{where}\quad V_i:=(Y_i,Y_iY_i\T-\Si). 
\end{equation*}

Let $\XX=\R^k\times\R^{k\times k}$, $\YY=\R^{k\times k}$, and $\ZZ=\R$, where $\XX$ has the norm defined by \eqref{eq:k,k-by-k norm} and $\YY$ is equipped with the spectral norm. 
Then the function $h\colon\XX\to\YY$ defined by the formula $h(x_1,x_2)=x_2-x_1x_1\T$ satisfies the smoothness condition \eqref{eq:h 2smooth} with $L_h(x_1,x_2)=x_2$, $M_h=2$, and any $\ep_h\in(0,\infty)$. 
Moreover, $h(\bar V)=S-\Si$. 

The likelihood-ratio tests of $H_{0,1}$ and $H_{0,2}$ against their negations, based on a normal population, reject for small values of the statistics
\begin{equation*}
 \La_1=\frac{\det(S)}{(\tr(S)/k)^k}\quad\text{and}\quad\La_2=\frac{\det(S)}{e^{\tr(S)}},
\end{equation*}
respectively; see e.g.\ Muirhead \cite[Theorems~8.3.2 and 8.4.2]{muirhead09}; one can also find in \cite{muirhead09} asymptotic properties of these tests, including expansions of their distributions under both null and nonnull distributions. 
Associate with $\La_1$ the functions $g\colon\YY\to\ZZ$ and $L_g\colon\YY\to\ZZ$ defined by
\begin{equation*}
 g(x)=
   \frac{\det(x+\Si)}{[\tr(x+\Si)/k]^k}\ind{\tr(x)>-k\al}-\Bigp{\frac\be\al}^k
 \quad\text{and}\quad
 L_g(x)=\Bigp{\frac{\be}{\al}}^k\,\tr[(\Si^{-1}-\al^{-1}I)x]. 
\end{equation*}
Similarly, with the statistic $\La_2$ associate the functions $g$ and $L_g$ defined by
\begin{equation*}
 g(x)=\frac{\det(x+\Si)}{e^{\tr(x+\Si)}}-\Bigp{\frac{\be}{e^{\al}}}^k\quad\text{and}\quad L_g(x)=\Bigp{\frac{\be}{e^{\al}}}^k\tr[(\Si^{-1}-I)x]. 
\end{equation*}
It is clear that, for either of the two functions $g$ defined above, $L_g=g'(0)$ and $g$ satisfies \eqref{eq:g 2smooth} for small enough $\ep_g$. 
Hence $f:=g\circ h$ satisfies \eqref{eq:smooth}, for both versions of the function $g$, and so, 
%IP16 Corollary
Theorem~\ref{cor:iid} %RM12.27 Theorem~\ref{thm:iid} 
may be applied to $f(\bar V)=\La_1-(\be/\al)^k$ and $f(\bar V)=\La_2-(\be/e^{\al})^k$. 

\newcommand{\tV}{\widetilde{V}}

For the case when the dimension $k$ is large, Nagao \cite{nagao73} proposes the test statistics 
\begin{equation*}
 U:=\frac1k\,\tr\Bigb{\Bigp{\frac{S}{\tr(S)/k}-I}^2}\quad\text{and}\quad \tV:=\frac1k\,\tr\bigb{(S-I)^2}
\end{equation*}
in place of the statistics $\La_1$ and $\La_2$, respectively. 
%RM22 %IP22
John \cite{john71} shows that the test of $H_{0,1}$ based on $U$ is locally most powerful (assuming a normal population). 
%RM22 See \cite{nagao73,nagao74} for asymptotic expansions of both null and nonnull distributions of $U$ and $\tV$.  
Associate with $U$ the functions 
\begin{equation*}
 g(x)=\frac1k\,\tr\Bigb{\Bigp{\frac{x+\Si}{\tr(x+\Si)/k}-I}^2}-\frac{\de^2}{\al^2}
 \quad\text{and}\quad
 L_g(x)=\frac2{k^2\al^3}\tr\Bigb{\bigp{\Si-\al I}\bigp{k\al I-\Si}x}
\end{equation*}
and with $\tV$ the functions
\begin{equation*}
 g(x)=\frac1k\,\tr\Bigb{\bigp{x+\Si-I}^2}-\de^2-(1-\al)^2
 \quad\text{and}\quad
 L_g(x)=\frac2k\,\tr\bigb{(\Si-I)x}.
\end{equation*}
It is straightforward to verify that either of the above functions $g$ satisfy the smoothness condition \eqref{eq:g 2smooth}, and hence that 
%IP16 Corollary
Theorem~\ref{cor:iid} %RM12.27 Theorem~\ref{thm:iid} 
may be applied to either of the functions $f(\bar V)=U-\de^2/\al^2$ or $f(\bar V)=\tV-\de^2-(1-\al)^2$. 

Yet one more variation on these tests we consider is the ``large-dimensional'' case. 
Ledoit and Wolf \cite{ledoit02} investigate the asymptotic behavior of both $U$ and $\tV$ when $k/n\to c\in(0,\infty)$ as $n\to\infty$, as opposed to the ``fixed-dimensional'' case (where $n\to\infty$ while $k$ is assumed a constant). 
They show that the test of $H_{0,1}$ based on $U$ remains consistent in the large-dimensional setting, whereas the test of $H_{0,2}$ based on $\tV$ is not necessarily consistent. 
By not dropping terms like $k/n$ in investigations of the asymptotics of $\tV$, the authors propose the statistic
\begin{equation}\label{eq:Wcov}%RM12.26 \label{eq:W}
  W:=\tV-\frac kn\Bigp{\Bigb{\frac{\tr(S)}k}^2-1}
  =\frac1k\,\tr\bigb{(S-I)^2}-\frac kn\Bigp{\Bigb{\frac{\tr(S)}k}^2-1}
\end{equation}
as an alternative to $\tV$ in the test of $H_{0,2}$. 
It is shown that $W$ has the same limiting distribution as $\tV$ in the fixed-dimensional setting while also being consistent in a large-dimensional framework. %RM22; asymptotic behavior of $U$, $\tV$, and $W$ when $k/n\to\infty$ as $n\to\infty$ is investigated in \cite{birke05}. 
We see that $f(\bar V)=W-\de^2-(1-\al)^2+\frac kn(\al^2-1)$ when $f=g\circ h$ and $g$ is defined by
\begin{equation*}
 g(x)=\frac1k\,\tr\Bigb{\bigp{x+\Si-I}^2}-\frac kn\biggp{\Bigp{\frac{\tr(x+\Si)}{k}}^2-1}-\de^2-(1-\al)^2+\frac kn(\al^2-1);
\end{equation*}
moreover, $g$ satisfies \eqref{eq:g 2smooth} with $L_g(x)=\frac2k\,\tr[(\Si-I-\frac kn\,\al I)x]$. 

\begin{thm}\label{thm:cov}
Take any $t\in\{\La_1,\La_2,U,\tV,W\}$, and let $f=g\circ h$ and $L=L_g\circ L_h$ for the functions $g$ and $L_g$ paired with the statistic $t$ as described above. 
Assume that $\tsi>0$ and $v_3<\infty$, for $\tsi$ and $v_p$ defined in 
\eqref{eq:tsi,ga3,v_al}. %RM12.27 \eqref{eq:tsi,v_p,vsi_p}. 
Then for all $n\in\N$ and $z\in\R$, 
\begin{equation}\label{eq:covbound.uni}
 \Bigabs{\Bigprob{\frac{f(\bar V)}{\tsi/\sqrt{n}}\le z}-\Phi(z)}\le\frac{\CC}{\sqrt n},
\end{equation}
where $\CC$ is a finite expression depending only on the distribution of $Y$; also, for any $\om>0$ and all real $z>0$ and $n\in\N$ satisfying \eqref{eq:z,iid},  
\begin{equation}\label{eq:covbound.nonuni}
 \Bigabs{\Bigprob{\frac{f(\bar V)}{\tsi/\sqrt{n}}\le z}-\Phi(z)}\le\frac{\CC}{z^3\,\sqrt n}, 
\end{equation}
where $\CC$ is a finite expression depending only on $\om$ and the distribution of $Y$. 
\end{thm}

\begin{remark}\label{re:cov.si=0}
The non-degeneracy condition $\tsi>0$ immediately implies that Theorem~\ref{thm:cov} -- and the delta method itself -- are applicable only to non-null distributions of the statistics $\La_1$, $\La_2$, $U$, and $\tV$, since $L_g=0$ for any of these statistics under the assumption of their respective null hypotheses. 
This should hardly be surprising, as it is known that these statistics (or some normalizing function of them) all have a limiting $\chi^2$ distribution under the null hypothesis. 
However, one can fix the null-degeneracy of the statistics $\La_1$, $\La_2$, $U$, or $\tV$ and thus 
make the delta method and our BE bounds applicable even to the null distributions by using essentially the same trick as in the definition of the statistic $W$ in \eqref{eq:Wcov}, %RM12.26 \eqref{eq:W}, 
that is, by adding a term of the form $\al\,\bigp{\bigb{\frac{\tr(S)}k}^2-1}$ for some nonzero real $\al$. 

By diagonalization of $\Si$, we can simply characterize the degeneracy condition $\tsi=0$ for any of the above statistics in this subsection. 
Indeed, by the spectral decomposition, $\Si=Q\T DQ$, where $D$ is the diagonal matrix with the eigenvalues $\la_1,\dotsc,\la_k$ of $\Si$ on its diagonal and $Q$ is an orthogonal matrix whose columns are corresponding orthonormal eigenvectors of $\Si$. 
Let $Z=(Z_1,\dotsc,Z_k)\T:=QY$. 
Then, for the statistic $\La_1$,
\begin{align*}
 \bigp{\tfrac{\al}{\be}}^kL(V)&=\tr\bigb{(\Si^{-1}-\al^{-1}I)(Y\,Y\T-\Si)}
  =\tr\bigb{Q\T(D^{-1}-\al^{-1}I)QQ\T(Z\,Z\T-D)Q}\\
 &=\tr\bigb{(D^{-1}-\al^{-1}I)(Z\,Z\T-D)}
  =\tr\bigb{(D^{-1}-\al^{-1}I)Z\,Z\T}-\tr\bigb{I-\al^{-1}D}\\
  &=\tsum_{j=1}^k\bigp{\tfrac1{\la_j}-\tfrac1{\al}}Z_j^2.
\end{align*}
Since $\tsi=0$ means precisely that $L(V)=0$ a.s., it follows that for any non-null alternative, $\tsi=0$ for the statistic $\La_1$ if and only if the support of 
the distribution of the random vector $Y$ degenerates so as to lie entirely on a certain quadric conical surface in $\R^k$.  
Similar work shows that for one of the statistics $\La_2$, $U$, $\tV$, and $W$ we have $\tsi=0$ if and only if the respective one of the random (homogeneous or not) quadratic forms  
\begin{multline*}
 \tsum_{j=1}^k\bigp{\tfrac1{\la_j}-1}Z_j^2,\quad
 \tsum_{j=1}^k(\la_j-\al)(k\al-\la_j)(Z_j^2-\la_j),\\ 
 \tsum_{j=1}^k(\la_j-1)(Z_j^2-\la_j),\quad
 \tsum_{j=1}^k\bigp{\la_j-1-\tfrac{k\al}{n}}(Z_j^2-\la_j)
\end{multline*}
equals $0$ a.s. 
In particular, whenever the random vector $Y$ is absolutely continuous, one has $\tsi>0$ for all these statistics in the non-null case, and then $\tsi>0$ for the statistic $W$ even in the null case provided that $(1-\frac kn)\Si\ne I$. 
\end{remark}

\begin{remark}\label{re:cov=Si}
Let $\Si_0$ be any given positive definite symmetric matrix. Then the hypotheses $\Si=\si^2\Si_0$ (with an unknown $\si^2>0$) and $\Si=\Si_0$ on the common covariance matrix $\Si$ of i.i.d.\ random vectors $Y_i$ are obviously equivalent to the respective hypotheses $\tSi=\si^2 I$ (with an unknown $\si^2>0$) and $\tSi=I$ on the common covariance matrix $\tSi$ of the i.i.d.\ random vectors $\tY_i:=\Si_0^{-1/2}Y_i$. So, the results in this subsection can be obviously extended to the more general case of the null hypotheses $\Si=\si^2\Si_0$ and $\Si=\Si_0$.  
\end{remark}

It appears certain that the bounds in Theorem~\ref{thm:cov} are all new to the literature; indeed, any of the results concerning these statistics that we have found investigates their asymptotic properties under the assumption of a normal population, whereas our bounds have only mild moment restrictions on $Y$. 
We mention here that 
%IP16 Corollary
Theorem~\ref{cor:iid} %RM12.27 Theorem~\ref{thm:iid} 
could be applied to several other popular statistics which are smooth functions of the sample covariance matrix $S$. 
For instance, our results can easily yield BE bounds for statistics proposed by Srivastava \cite{sri05} or Fisher et al.\ \cite{fisher10}; Chen et al.\ \cite{chen10} propose a statistic for the sphericity test which is a function of a $U$-statistic, for which the methods of this paper and \cite{chen07} could presumably be adapted. 
The reader is referred to 
\cite{muirhead09} %RM22 ,john71,nagao73} 
for other statistics used in testing for the equality of population covariances or independence between certain projections applied to $Y$.

\subsection{Principal component analysis (PCA)}
\label{subsec:pca}
It is well known that any simple eigenvalue of a (say, symmetric real matrix) and the orthoprojector onto the corresponding eigenspace are smooth functions of the matrix. Therefore, the delta method is almost universally applicable to PCA, and hence so are our results such as 
%IP16 Corollary
Theorem~\ref{cor:iid}. %RM12.27 Theorem~\ref{thm:iid}. 
The actual verification of the smoothness condition \eqref{eq:smooth} in PCA may involve operator perturbation theory and related tools, based on a representation of analytic functions of a linear operator as certain integrals of the resolvent. This representation largely reduces the problem of the smoothness of a general analytic function of an operator to the obvious smoothness of the map $A\mapsto A^{-1}$ on the set of all bounded invertible linear operators $A$ (cf.\ \eqref{eq:P=} and \eqref{eq:la=}). 
Whereas this idea is rather transparent, its execution may in some cases be rather nontrivial, and it may result in complicated expressions for $\ep$ and $M_\ep$ in \eqref{eq:smooth}. 

As an illustration of these general theses, let us consider here a statistic rather recently introduced by Cupidon et al.\ \cite{cupidon07,cupidon08}. 
Let $Y,Y_1,\dotsc,Y_n$ be iid r.v.'s taking values in a separable real Hilbert space $H$ with inner product $\ip\cdot\cdot$ and the corresponding norm $\|\cdot\|$. 
Assume at this point that $\E\|Y\|^2<\infty$, $\E Y=0$,    
and the covariance operator
\begin{equation*}
R:=\cov Y=\E(Y\otimes Y)
\end{equation*}
of $Y$ 
is (strictly) positive definite. 
Here, as usual, $\otimes$ denotes the tensor product on $H$, so that $Rx=\E\ip xY Y$ for all $x\in H$. 
Given the condition $\E\|Y\|^2<\infty$, the covariance operator $R$ is known to be compact, which allows its spectral decomposition -- see e.g.\ \cite[Theorem~2.10, page~260]{kato}; a short proof of the compactness of 
$R$ %RM12.15 !! $\R$ !! 
is presented in Appendix~\ref{sec:compact} %RMejs Supplement~\ref{compact} 
% \cite{supp}, %RM22\cite{compact}, %RM5 
for the readers' convenience. 

Next suppose that $H=H_1\oplus H_2$, where $H_1,H_2$ are closed orthogonal subspaces of $H$; for $j,k\in\{1,2\}$, let $\Pi_j$ denote the orthoprojector onto $H_j$, $R_{jk}:=\Pi_j R\Pi_k$, and also let $I_j$ denote the identity operator on $H_j$. 
Then, for any fixed $\al>0$, the \emph{regularized squared principal canonical correlation}, RSPCC or $\rho^2$, is defined by the formula 
\begin{equation}\label{eq:rspcc}
 \rho^2:=\rho^2(\al):=\mathop{\max_{x\in H_1\setminus\{0\}}}_{y\in H_2\setminus\{0\}}\frac{\langle x,R_{12}y\rangle^2}{\langle x,(\al I_1+R_{11})x\rangle\langle y,(\al I_2+R_{22})y\rangle};
\end{equation}
that this is a well-defined quantity is proved in \cite{cupidon07}. 
Define the sample RSPCC, $\hat\rho^2$, by replacing $R_{jk}$ in \eqref{eq:rspcc} with $S_{jk}$, where
\begin{equation*}
 S_{jk}=\Pi_jS\,\Pi_k,\quad S:=\YYtensbar-\bar Y\otimes\bar Y,\quad\bar Y:=\tfrac1n\tsum\nolim_iY_i,\quad\text{and}\quad\YYtensbar:=\tfrac1n\tsum\nolim_iY_i\otimes Y_i; 
\end{equation*}
thus, $S$ is the sample covariance operator of the random vector $Y$. 
See e.g.\ \cite{%RM22he03,
he04,eubank08} for discussion and results on the use of canonical correlations in functional data. 

Next define the (bounded self-adjoint nonnegative-definite linear) operators
\begin{equation}\label{eq:R1}
\begin{aligned}
 R_1:=&(\al I_1+R_{11})^{-1/2}R_{12}(\al I_2+R_{22})^{-1}R_{21}(\al I_1+R_{11})^{-1/2}, \\ 
 R_2:=&(\al I_2+R_{22})^{-1/2}R_{21}(\al I_1+R_{11})^{-1}R_{12}(\al I_2+R_{22})^{-1/2},
\end{aligned}
\end{equation}
and similarly let $\hat R_j$ denote the sample analogues of $R_j$ (obtained by replacing $R_{jk}$ with $S_{jk}$); under the  assumption that  $\E\|Y\|^2<\infty$ (which implies that $R$ is compact), we see that $R_1$ and $R_2$ are also compact. Moreover, by \cite[Theorem 2.4]{cupidon08}, $\|R_1\|=\|R_2\|=\rho^2$ and $\|\hat R_1\|=\|\hat R_2\|=\hat\rho^2$, where 
$\rho^2$ is as in \eqref{eq:rspcc} and 
$\|\cdot\|$ denotes the operator norm, so that $\|R_j\|$ is the largest eigenvalue of $R_j$. 

Fix any $j\in\{1,2\}$ and 
assume that $\rho^2$ is a simple nonzero eigenvalue of $R_j$, and then let $P$ denote the orthoprojector onto the corresponding (one-dimensional) eigenspace of $R_j$. 
Let $B(H)$ and $B(H_j)$ denote the Hilbert spaces of all bounded linear operators on $H$ and $H_j$, respectively, equipped with the corresponding operator norms.   

Let $g(x):=\|x+R_j\|-\|R_j\|$ for any $x\in B(H_j)$, so that $g(\hat R_j-R_j)=\hat\rho^2-\rho^2$. By formulas (3.6)--(3.8) on page 89, (2.32) on page 79, and (3.4) on page 88 in \cite{kato} (with $n=1$, $\varkappa=1$, $\hat\la(\vka)=\|x+\hat R_j\|$, $\la=\|R_j\|=\rho^2$, $\hat\la^{(1)}=\tr(xP)$,  
$T^{(1)}=x$, $T^{(2)}=T^{(3)}=\dots=0$, $a=\|x\|$, $c=0$, and $0^0:=1$), the smoothness condition \eqref{eq:g 2smooth} will be satisfied with $\ep_g=\be/\m$, $\be\in(0,1)$, $L(x)=\hat\la^{(1)}=\tr(xP)$, and $M_g=2\frac{\varrho\m^2}{1-\be}$, where $\varrho:=\max_{z\in\Ga}\abs{z-\la}$, $\m:=\max_{z\in\Ga}\|\res_j(z)\|$, $\res_j(z):=(R_j-zI)^{-1}$ is the resolvent of $R_j$, and $\Ga$ is the boundary of any open disc $D$ in $\C$ such that $\la\in D$ but the closure of $D$ does not contain $0$ or 
any eigenvalue of $R_j$ other than $\la$.  

\rule{0pt}{0pt}\big(The results from \cite{kato} referred to in the above paragraph were stated there for the case when the Hilbert space $H$ is finite-dimensional. All those results carry verbatim to the ``infinite-dimensional'' case. Such information can be extracted from other chapters in \cite{kato}. However, for readers' convenience, in Appendix~\ref{sec:inf-dim} %RMejs Supplement~\ref{inf-dim} 
% \cite{supp} %RM22 \cite{inf-dim} %RM5 
we provide the few necessary stepping stones to make the transition to the infinite dimension.\big)

By \cite[Theorem~2.1]{gilliam09}, condition \eqref{eq:h 2smooth} holds for 
the function $y\mapsto h_j(y):=(\al I_j+R_{jj}+y)^{-1/2}\break
-(\al I_j+R_{jj})^{-1/2}$ in place of $h$ 
for some real $\ep_{h_j}>0$ and all $y\in H_j$ with $\|y\|\le\ep_{h_j}$. 
So, in view of definitions \eqref{eq:R1} of $R_j$, their counterparts for $\hat R_j$, and Remark~\ref{re:compos}, one can set up a function $h\colon H\times B(H)\to H_j$ in a straightforward manner so that condition \eqref{eq:h 2smooth} holds and $h(\bar V)=\hat R_j-R_j$, with the zero-mean vector $V=(Y,Y\otimes Y-R)$. 
Using Remark~\ref{re:compos} once again, one sees that the function $f=g\circ h$ satisfies the smoothness condition \eqref{eq:smooth}, and at that $f(\bar V)=\hat\rho^2-\rho^2$.  
Thus, 
%IP16 Corollary
Theorem~\ref{cor:iid} %RM12.27 Theorem~\ref{thm:iid} 
yields

\begin{thm}\label{thm:rspcc}
Assume that $\tsi>0$ and $v_3<\infty$, for $\tsi$ and $v_p$ defined in 
\eqref{eq:tsi,ga3,v_al}. %RM12.27 \eqref{eq:tsi,v_p,vsi_p}. 
Then for all $n\in\N$ and $z\in\R$, 
\begin{equation}\label{eq:rspcc.uni}
 \Bigabs{\Bigprob{\frac{\hat\rho^2-\rho^2}{\tsi/\sqrt{n}}\le z}-\Phi(z)}\le\frac{\CC}{\sqrt n},
\end{equation}
where $\CC$ is a finite expression depending only on the distribution of $Y$; also, for any $\om>0$ and all real $z>0$ and $n\in\N$ satisfying \eqref{eq:z,iid},  
\begin{equation}\label{eq:rspcc.nonuni}
 \Bigabs{\Bigprob{\frac{\hat\rho^2-\rho^2}{\tsi/\sqrt{n}}\le z}-\Phi(z)}\le\frac{\CC}{z^3\,\sqrt n}, 
\end{equation}
where $\CC$ is a finite expression depending only on $\om$ and the distribution of $Y$. 
\end{thm}

Expressions for $\tsi$ can be obtained from \cite[(4.20), (5.1)]{cupidon07}. 
We see the recurring theme that $\norm{Y}_4<\infty$ is used to establish asymptotic normality of $\hat\rho^2$ (cf.\ \cite[(2.1), Theorem~4.2]{cupidon07}), while the moment restriction $\norm{Y}_6<\infty$ (equivalent to $v_3<\infty$ in Theorem~\ref{thm:rspcc}) is needed here to bound the rate of convergence on the order $\bigO(1/\sqrt{n})$. 
Again, it appears that the bounds in Theorem~\ref{thm:rspcc} are entirely new to the literature.

In Subsection~\ref{subsec:cov}, we considered various smooth functions of the determinant and trace of the sample covariance matrix for finite-dimensional r.v.'s $Y$, and in the present subsection we have a function of the largest eigenvalue of some smooth function of a sample covariance operator. 
Other statistics which are functions of eigenvalues from a sample covariance operator (be it constructed from a finite-dimensional or infinite-dimensional population) may of course lie in the class of statistics to which 
%IP16 Corollary
Theorem~\ref{cor:iid} %RM12.27 Theorem~\ref{thm:iid} 
could be applied; the primary problem to the practitioner is the demonstration of the smoothness condtion \eqref{eq:smooth}. 
The use of perturbation theory, as was done above, appears to be valuable for many such potential applications; we mention here statistics proposed in \cite{ji08,gaines-etal}, concerning the testing of equality of two covariance operators, as %IP2 yet 
further examples. 
Yet another potential application of our results would be to the empirical {W}asserstein
distance, for which central limit theorems were recently given in \cite{rippl-etal15}; cf.\ \cite{olkin-pukel82,dowson-landau,givens-shortt} (as noted by Dudley in his review MR0752258 on MathSciNet, the normality assumption is not actually needed there). 

%IP1
\newpage

\subsection{Maximum likelihood estimators (MLEs)}
\label{subsec:mle}
%IP16 copied from EJS draft 
Bounds on the closeness of the distribution of the MLE to normality in the so-called bounded Wasserstein distance, $d_{\bW}$, were recently obtained in \cite{anast-reinert} under certain regularity conditions. 
In \cite{anast-ley}, these bounds were improved in the rather common case when the MLE $\hat\th$ satisfies the condition 
\begin{equation}\label{eq:q(th)}
 q(\hat\th)=\frac1n\,\sum_{i=1}^n g(Y_i), %RM12.19 X_i), 
\end{equation}
where $q\colon\Th\to\R$ is a twice continuously differentiable one-to-one mapping, $g\colon\R\to\R$ is a Borel-measurable function, and the $Y_i$'s %RM12.19 $X_i$'s 
are i.i.d.\ %RM12.19 iid 
real-valued r.v.'s. 

%RM12.19 %%replace 'X' with 'Y' throughout
It is noted in \cite[Proposition~2.1]{anast-reinert} that for any r.v.\ $Y$ and $Z\sim N(0,1)$ one has $d_\Ko(Y,Z)\le2\sqrt{d_{\bW}(Y,Z)}$, where $d_\Ko$ denotes the Kolmogorov distance. 
This bound on $d_\Ko$ in terms of $d_{\bW}$ is the best possible one, up a constant factor. Indeed, for each real $\vp>0$, define a r.v.\ $Y_\vp$ as follows: $Y_\vp=\vp$ if $0<Z<\vp$ and $Y_\vp=Z$ otherwise. Then for any Lip$(1)$ function $h\colon\R\to\R$ one has 
$|\E h(Y_\vp)-\E h(Z)|\le\E|h(Y_\vp)-h(Z)|\le\E|Y_\vp-Z|=\int_0^\vp(\vp-z)\vpi(z)\,dz\le\vpi(0)\vp^2/2$. %RM12.20, 
%\E\ind{0<Z<\vp}\break
%\times|h(\vp)-h(Z)|\le\E\ind{0<Z<\vp}|\vp-Z|=\int_0^\vp(\vp-z)\vpi(z)\,dz\le\vpi(0)\vp^2/2$, 
%RM12.20 where $\vpi$ is the standard normal pdf. 
So, $d_{\bW}(Y_\vp,Z)\le d_\W(Y_\vp,Z)\le\vpi(0)\vp^2/2$, where $d_\W$ is the Wasserstein distance: $d_\W(X,Y):=\sup\{|\E h(X)-\E h(Y)|\colon h\in\text{Lip}(1), h\text{ bounded}\}$ for any r.v.'s $X$ and $Y$. 
On the other hand, $d_\Ko(Y_\vp,Z)\ge\P(Z<\vp)-\P(Y_\vp<\vp)=\Phi(\vp)-1/2\sim\vpi(0)\vp$, so that $d_\Ko(Y_\vp,Z)\ge\sqrt{2\vpi(0)-o(1)}\,\sqrt{d_{\bW}(Y_\vp,Z)}$ as $\vp\downarrow0$. 

Therefore, even though the bounds on $d_{\bW}$ obtained in \cite{anast-reinert,anast-ley} are of the optimal order $O(1/\sqrt n)$, the resulting bounds on the Kolmogorov distance are only of the order $O(1/n^{1/4})$. 

In this subsection, as an application of our general results, we shall obtain bounds of the optimal order $O(1/\sqrt n)$ on the closeness of the distribution of the MLE to normality in the Kolmogorov distance assuming a somewhat relaxed version of the condition \eqref{eq:q(th)}. 
In addition, we shall present a corresponding nonuniform bound. 
At that, our regularity conditions appear simpler than those in \cite{anast-reinert,anast-ley}. 

%\newpage

Indeed, let here 
$%IP16 
Y,Y_1,Y_2,\dots$ %RM12.20 $X_1,X_2,\dots$ 
be r.v.'s mapping a measurable space $(\Om,\A)$ to another measurable space $(\XXX,\B)$ and let $(\P_\th)_{\th\in\Th}$ be a parametric family 
of probability measures on $(\Om,\A)$ such that the r.v.'s 
$Y_1,Y_2,\dots$ %RM12.20 $X_1,X_2,\dots$ 
are 
i.i.d.\ %RM12.20 iid 
with respect to each of the probability measures $\P_\th$ with $\th\in\Th$; here the parameter space $\Th$ is assumed to be a subset of $\R$. 
As usual, let $\E_\th$ denote the expectation with respect to the probability measure $\P_\th$. 
Suppose that for each $\th\in\Th$ the distribution %IP16 $\P_\th$ 
$\P_\th Y^{-1}$ of the r.v.\ $Y$ with respect to the probability measure $\P_\th$ has a density $p_\th$ with respect to a measure $\mu$ on %IP16 $\F$
$\B$. 
For each point $\xx=(x_1,\dots,x_n)\in\XXX^n$ such that the likelihood function $\Th\ni\th\mapsto L_\xx(\th):=\prod_{i=1}^n p_\th(x_i)$ has a unique maximizer, denote this maximizer by $\hat\th_n(\xx)$; otherwise, assign to $\hat\th_n(\xx)$ any value in $\Th$. Let us then refer to $\hat\th_n(\X)$ as the MLE of $\th$, where 
$\X:=(Y_1,Y_2,\dots)$. %RM12.20 $\X:=(X_1,X_2,\dots)$. 
Clearly, this is a more general definition of the MLE than usual, and we can even allow the function $\hat\th_n$ to be non-measurable. So, the MLE $\hat\th_n(\X)$ does not have to be a r.v. 
Let $\th_0\in\Th$ be the ``true'' value of the unknown parameter $\th$, such that  
$\Th_0:=(\th_0-\vp,\th_0+\vp)\subseteq\Th$ for some real $\vp>0$. 

We assume the following relaxed version of the condition \eqref{eq:q(th)}:  
for some real constant $C>0$ and each natural $n$ there exists a set $E_n\in\B^{\otimes n}$ such that 
\begin{equation}\label{eq:notin E_n}
 \P_{\th_0}(\X\notin E_n)\le C/\sqrt n
\end{equation}
and 
for each point $\xx=(x_1,\dots,x_n)\in E_n$ the value $\hat\th_n(\xx)$ of the MLE belongs to the neighborhood $\Th_0$ of the point $\th_0$ and satisfies the condition 
\begin{equation}\label{eq:q(th),s}
 q\big(\hat\th_n(\xx)\big)=\frac1n\,\sum_{i=1}^n g(x_i), 
\end{equation} 
for some measurable function $g\colon\XX\to\R$ and some twice continuously differentiable %one-to-one 
mapping $q\colon\Th_0\to\R$ with $q'(\th)\ne0$ for all $\th\in\Th_0$, so that the mapping $q$ is one-to-one. 
Suppose also that the MLE $\hat\th_n(\X)$ is consistent at the point $\th_0$, that is, $\hat\th_n(\X)\underset{n\to\infty}\longrightarrow\th_0$ in probability with respect to the probability measure $\P_{\th_0}$; since the MLE $\hat\th_n(\X)$ does not have to be a r.v., the precise meaning of this consistency is that $(\P_{\th_0})^*(|\hat\th_n(\X)-\th_0|>\de)\underset{n\to\infty}\longrightarrow0$ for each real $\de>0$, where $(\P_{\th_0})^*$ denotes the outer measure induced by the probability measure $\P_{\th_0}$. Then, under the condition $\E_{\th_0}|g(Y_1%RM12.20 X_1
)|<\infty$, it follows from \eqref{eq:q(th),s} by the law of large numbers that $q(\th_0)=\mu_g:=\E_{\th_0}g(Y%RM12.20 X
_1)$ or, equivalently, $\th_0=q^{-1}(\mu_g)$, where $q^{-1}$ stands for the inverse of the function $q$. 

Assuming further that $\si_g:=\sqrt{\var_{\th_0}g(Y%RM12.20X
_1)}\in(0,\infty)$, let us introduce 
\begin{equation*}
 V%RM12.20 Y
 _i:=\frac{g(%IP16 X_i
 Y_i)-\mu_g}{\si_g}
\end{equation*}
for $i=1,\dots,n$ and 
\begin{equation*}
 f(v):=q^{-1}(\mu_g+\si_g v)-q^{-1}(\mu_g)=q^{-1}(\mu_g+\si_g v)-\th_0
 %RM12.20f(y):=q^{-1}(\mu_g+\si_g y)-q^{-1}(\mu_g)=q^{-1}(\mu_g+\si_g y)-\th_0
\end{equation*}
for real $v$ such that $\mu_g+\si_g v\in q(\Th_0)$ and $f(v)=0$ (say) for the other real values of $v$. Then, in view of \eqref{eq:q(th),s}, on the event $\{\X\notin E_n\}$ one has $f(\bar V%RM12.20Y
_n)=\hat\th_n(\X)-\th_0$, and at that $f(0)=0$, $f'(0)=\si_g\,(q^{-1})'(\mu_g)=\si_g/q'\big(q^{-1}(\mu_g)\big)=\si_g/q'(\th_0)$, and $f$ is twice continuously differentiable in a neighborhood of $0$. 
So, Theorem~\ref{thm:iid,uni-dim} immediately yields  

\begin{thm}\label{thm:MLE}
In addition to the conditions specified above, 
assume that $\E_{\th_0}|g(Y%RM12.20 X
_1)|^3<\infty$.  
Then for all $n\in\N$ and $z\in\R$
\begin{equation}\label{eq:MLE.bound.uni}
\Big|\P_{\th_0}\Big(\frac{\hat\th_n(\X)-\th_0}{\si_g/\sqrt n}\le\frac z{|q'(\th_0)|}\Big)-\Phi(z)\Big|\le\frac{C+\CC}{\sqrt n},
\end{equation} 
where $C$ is as in \eqref{eq:notin E_n} and $\CC$ is a finite expression depending only on the $\P_\th$-distributions of $Y%RM12.20 X
_1$ for $\th$ in a neighborhood of $\th_0$. % and the function $q$; 
Also, if in \eqref{eq:notin E_n} one can replace $\sqrt n$ by $n^2$, then 
for any $\om>0$ and for all real $z>0$ and $n\in\N$ satisfying \eqref{eq:z,iid}, 
\begin{equation}\label{eq:MLE.bound}
 \Big|\P_{\th_0}\Big(\frac{\hat\th_n(\X)-\th_0}{\si_g/\sqrt n}\le\frac z{|q'(\th_0)|}\Big)-\Phi(z)\Big|\le\frac{C+\CC}{z^3\,\sqrt n}, 
\end{equation}
where $\CC$ is a finite expression depending only on $\om$ and the $\P_\th$-distributions of $Y%RM12.20 X
_1$ for $\th$ in a neighborhood of $\th_0$. 
\end{thm}

As was noted, the MLE $\hat\th_n(\X)$ does not have to be a r.v., and so, the $\P_{\th_0}$-probability in \eqref{eq:MLE.bound.uni} and \eqref{eq:MLE.bound} does not have to be defined. Thus, strictly speaking, one should understand this probability as the corresponding outer or inner probability, $(\P_{\th_0})^*$ or $(\P_{\th_0})_*$ -- each one of the two versions will do in each of the two inequalities, \eqref{eq:MLE.bound.uni} and \eqref{eq:MLE.bound}.   

Let us show that, under certain mild and natural conditions, \eqref{eq:q(th),s} is fulfilled if the densities $p_\th$ form an exponential family with a natural parameter (cf.\ \cite{anast-ley}), so that  
\begin{equation}\label{eq:exp-fam}
 p_\th(x)=e^{\th g(x)-c(\th)}
\end{equation}
for some function $c\colon\Th\to\R$ and 
all $\th\in\Th$ and $x\in\XXX$. Here, as before, $g\colon\XXX\to\R$ is a  measurable function. The natural choice of the parameter space here is $\Th:=\{\th\in\R\colon\EE(\th):=\int_\XXX e^{\th g(x)}\mu(dx)<\infty\}$, and then of course $c(\th)=\ln\EE(\th)$ for all $\th\in\Th$. As before, assume that $\Th_0:=(\th_0-\vp,\th_0+\vp)\subseteq\Th$ for some real $\vp>0$. In fact, by decreasing $\vp$ if necessary, we may and shall assume that $[\th_0-\vp,\th_0+\vp]\subseteq\Th$. 
If $\mu(\{x\in\XXX\colon g(x)\ne a\})=0$ for some real $a$, then for all $\th\in\Th$ one has $\EE(\th)=e^{\th a}\mu(\XXX)<\infty$, whence $p_\th(x)=1/\mu(\XXX)$ for $x\in\XXX$, so that the densities $p_\th$ are the same for all $\th\in\Th$, and therefore parameter $\th$ is not identifiable. Let us exclude this trivial case. 
Note that the function $c$ is infinitely many times differentiable (and even real-analytic) on $\Th_0=(\th_0-\vp,\th_0+\vp)$. 
Moreover, its derivative $c'$ is (strictly) increasing and hence $c$ is strictly convex on $\Th_0$, because $c''(\th)=(\ln\EE)''(\th)=\var_\th  g(Y%RM12.20X
_1)>0$ for $\th\in\Th_0$, since the trivial case of the non-identifiability of $\th$ has just been excluded. 
In particular, it follows that the condition $\si_g:=\sqrt{\var_{\th_0}g(Y%RM12.20X
_1)}\in(0,\infty)$ holds. 
At that, $\mu_g:=\E_{\th_0}g(%IP16 X
Y_1)=c'(\th_0)$. 

Let now 
\begin{equation}\label{eq:E_n:=}
 E_n:%=\Big\{\xx\in\XXX^n\colon\sum_{i=1}^n g(x_i)\notin c'(\Th_0)\Big\}
 =\Big\{\xx\in\XXX^n\colon c'(\th_0-\vp)<\frac1n\,\sum_{i=1}^n g(x_i)<c'(\th_0+\vp)\big)\Big\}. 
\end{equation}
By Markov's inequality, 
\begin{align*}
 \P_{\th_0}\Big(\frac1n\,\sum_{i=1}^n g(Y%RM12.20X
 _i)\le c'(\th_0-\vp)\Big)
 &=\P_{\th_0}\Big(\exp\Big\{-\vp\,\sum_{i=1}^n g(Y%RM12.20X
 _i)\Big\}\ge\exp\big\{-n\vp\,c'(\th_0-\vp)\big\}\Big) \\ 
 &\le\exp\big\{n\vp\,c'(\th_0-\vp)\big\}\E_{\th_0}\exp\Big\{-\vp\,\sum_{i=1}^n g(Y%RM12.20X
 _i)\Big\} \\ 
 &=\exp\big\{n\vp\,c'(\th_0-\vp)+nc(\th_0-\vp)-nc(\th_0)\big\}
 =e^{-n\de(\vp)},  
\end{align*}
where $\de(\vp):=c(\th_0)-c(\th_0-\vp)-c'(\th_0-\vp)\vp>0$; the latter inequality holds because (i) the function $c$ is strictly convex and (ii) one has $h(u+v)>h(u)+h'(u)v$ for any strictly convex differentiable function $h$, any $u$, and any nonzero $v$. 
Quite similarly, $\P_{\th_0}\big(\frac1n\,\sum_{i=1}^n g(Y%RM12.20X
_i)\ge c'(\th_0+\vp)\big)\le e^{-n\de(-\vp)}$, with $\de(-\vp)>0$. 
So, 
\begin{equation}\label{eq:<e+e}
 \P_{\th_0}(\X\notin E_n)\le e^{-n\de(\vp)}+e^{-n\de(-\vp)}, 
\end{equation}
so that condition \eqref{eq:notin E_n} holds, even with $n^2$ in place of $\sqrt n$. 
On the other hand, in view of \eqref{eq:E_n:=} and because $c'$ is continuous and increasing on $\Th_0$, we see that \eqref{eq:q(th),s} holds for all $\xx\in E_n$, with $q(\th)=c'(\th)$ for all $\th\in\Th_0$. 
Now the consistency of the MLE at point $\th_0$ follows because (i) by \eqref{eq:<e+e}, $\P_{\th_0}(\X\notin E_n)\underset{n\to\infty}\longrightarrow0$ and (ii) by the law of large numbers, $\frac1n\,\sum_{i=1}^n g(Y%RM12.20X
_i)\underset{n\to\infty}\longrightarrow\E_{\th_0}g(Y%RM12.20X
_1)=\mu_g=c'(\th_0)$ in $\P_{\th_0}$-probability. 

Note finally that the condition $\E_{\th_0}|g(Y%RM12.20X
_1)|^3<\infty$ in Theorem~\ref{thm:MLE} holds as well, since 
$\E_{\th_0}\exp\{\vp\,|g(Y%RM12.20X
_1)|\}<\E_{\th_0}\exp\{\vp\,g(Y%RM12.20X
_1)\}+\E_{\th_0}\exp\{-\vp\,g(Y%RM12.20X
_1)\}
=c(\th_0+\vp)+c(\th_0-\vp)<\infty$. 

We have verified all the conditions needed in order to apply Theorem~\ref{thm:MLE}. In addition to this, note that in the present context of exponential families, $q'(\th)=c''(\th)=-\frac{\partial^2}{\partial\th^2}\ln p_\th(x)$ does not depend on $x$, whence for each $\th\in\Th_0$ one has $q'(\th)=-\E_\th\frac{\partial^2}{\partial\th^2}\ln p_\th(Y%RM12.20X
_1)=I(\th)$, the Fisher information contained in $Y%RM12.20X
_1$. Also, recall that $\si_g=\sqrt{\var_{\th_0}g(Y%RM12.20X
_1)}=\sqrt{c''(\th_0)}=\sqrt{I(\th_0)}$. 
Thus, we have 

\begin{cor}\label{cor:MLE,exp}
Suppose that the conditions introduced above starting with the exponential family condition \eqref{eq:exp-fam} hold. 
Then for all $n\in\N$ and $z\in\R$
\begin{equation}\label{eq:MLE.bound.uni.exp}
\Big|\P_{\th_0}\Big(\hat\th_n(\X)-\th_0\le\frac z{\sqrt{nI(\th_0)}}\Big)-\Phi(z)\Big|\le\frac{\CC}{\sqrt n},
\end{equation} 
where $\CC$ is a finite expression depending only on the $\P_\th$-distributions of $Y%RM12.20X
_1$ for $\th$ in a neighborhood of $\th_0$. % and the function $q$; 
Also,  
for any $\om>0$ and for all real $z>0$ and $n\in\N$ satisfying \eqref{eq:z,iid}, 
\begin{equation}\label{eq:MLE.bound.exp}
 \Big|\P_{\th_0}\Big(\hat\th_n(\X)-\th_0\le\frac z{\sqrt{nI(\th_0)}}\Big)-\Phi(z)\Big|\le\frac{\CC}{z^3\,\sqrt n}, 
\end{equation}
where $\CC$ is a finite expression depending only on $\om$ and the $\P_\th$-distributions of $Y%RM12.20X
_1$ for $\th$ in a neighborhood of $\th_0$. 
\end{cor}

\begin{ex}\label{ex:exp}
Let here $\XXX=\R$ and let $\B$ be the Borel $\si$-algebra over $\R$. 
Let the measure $\mu$ on $\B$ be defined by the formula $\mu(dx)=(x+1)^{-3}\ind{x\ge0}\,dx$, and let $g(x)=x$ for all real $x$. 
%the Lebesgue measure on $\B$, 
Let then $p_\th$ be as in \eqref{eq:exp-fam}, with $\Th=(-\infty,0]$. 
It follows that $c'$ increases on $\Th$, with $c'(0-)=\int_0^\infty x(x+1)^{-3}\,dx\big/\break
\int_0^\infty (x+1)^{-3}\,dx=1<\infty$. On the other hand, for each natural $n$, with nonzero $\P_\th$-probability for each $\th\in\Th$, the r.v.\ $\frac1n\,\sum_{i=1}^n g(Y%RM12.20X
_i)=\frac1n\,\sum_{i=1}^n Y%RM12.20X
_i$ may take arbitrarily large values, in particular values exceeding $1=c'(0-)=\sup_{\th\in\Th}c'(\th)$. So, the equality \eqref{eq:q(th)} will be violated with nonzero $\P_\th$-probability for each $\th\in\Th$ and for each natural $n$.  
However, Theorem~\ref{thm:MLE} and Corollary~\ref{cor:MLE,exp} will hold in this situation. 
This shows the usefulness of the relaxed version \eqref{eq:notin E_n}--\eqref{eq:q(th),s} of the condition \eqref{eq:q(th)}. 
\end{ex}

%IP16
As shown in \cite{MLE-CLT}, with more effort one can utilize the ``multivariate'' Theorem~\ref{cor:iid} (rather than the ``univariate'' Theorem~\ref{thm:iid,uni-dim}, used in this subsection) to obtain bounds of optimal order $O(1/\sqrt n)$ on the Kolmogorov distance for MLEs in general, without assuming \eqref{eq:q(th)} or \eqref{eq:notin E_n}--\eqref{eq:q(th),s}. It is also shown in \cite{MLE-CLT} that, again without assuming \eqref{eq:q(th)} or \eqref{eq:notin E_n}--\eqref{eq:q(th),s}, one can obtain the corresponding nonuniform bounds of the optimal orders in $n$ and $z$. All these results can be extended to the more general case of $M$-estimators or, even more generally, to the estimators that are zeros of estimating functions; see e.g.\ \cite{heyde97}. Indeed, the condition that $p_\th$ is a pdf for $\th\ne\th_0$ is used in our proofs only in order to state that $\E_\th\ell'_X(\th)=0$ and $\E_\th \ell'_X(\th)^2=-\E_\th \ell''_X(\th)=I(\th)\in(0,\infty)$. In the case of $M$-estimators or zeros of estimating functions, the corresponding conditions will have to be just assumed, with some other expressions in place of the Fisher information $I(\th)$, as it is done e.g.\ in \cite{pfanzagl71,pfanzagl73}, where 
uniform (but not nonuniform) bounds of optimal order $O(1/\sqrt n)$ for $M$-estimators were obtained (via different, specialized methods): in \cite{pfanzagl71} for a one-dimensional parameter space $\Th$ and in \cite{pfanzagl73} in the multidimensional case.       
%!! multivariate case

\section{Proofs}\label{sec:proofs}

All necessary proofs of the theorems and corollaries stated in the previous sections are provided here -- except for Corollaries~\ref{cor:centralT}, \ref{cor:centralT.nub}, and \ref{cor:Runif}, whose proofs are given in Appendix~\ref{app:num.proofs}. 

\subsection{Proofs of results from Section~\ref{sec:chen.mod}}
\label{sec:chen.proofs}

\begin{proof}[Proof of Theorem~\ref{thm:ub}]
As noted in Remark~\ref{re:ub}, the assertion of Theorem~\ref{thm:ub} is very similar to that of \cite[Theorem~2.1]{chen07}. 
From the condition that $|\De|\ge|T-W|$ (cf.\ \cite[(5.1)]{chen07})
\begin{equation}\label{eq:conc}
 -\P(z-|\De|\le W\le z)\le\P(T\le z)-\P(W\le z)\le\P(z\le W\le z+|\De|)
\end{equation}
for all $z\in\R$. 
The inequality 
\begin{equation*}
 \bigprob{z\le W\le z+|\bar\De|}
  \le \frac1{2\cc}\Bigp{4\de+\E\bigabs{W\bar\De}+\tsum\limits_{i=1}^n\E\bigabs{\xi_i(\bar\De-\De_i)}}
\end{equation*}
is proved by modifying the proof of \cite[Theorem~2.1]{chen07} -- replacing their $\De$ with our $\bar\De$ and their condition (2.2) with our \eqref{eq:de}. 
Recalling the condition \eqref{eq:barDe} on $\bar\De$, one has
\begin{equation}\label{eq:p+Geta(w)}
 \bigprob{z\le W\le z+\abs{\De}}
  \le\bigprob{z\le W\le z+\abs{\bar\De}}+\bigprob{\max\nolim_i\eta_i>w}.
\end{equation}
Then $\P(z-\abs{\De}\le W\le z)$ can be bounded in a similar fashion, using $z-\abs{\De}$ in place of $z$, and \eqref{eq:ub} follows.

In order to prove Remark~\ref{re:chen.improvement}, note that \cite[(5.6)]{chen07} still remains valid when $H_{1,2}$ there is replaced by 
\begin{equation*}
 H_{1,2}=\E\bigind{z\le W\le z+\abs{\bar\De}}\Bigabs{\tsum\nolim_i\bigp{\check\xi_i-\E\check\xi_i}},\quad
 \text{with}\quad
 \check\xi_i:=\abs{\xi_i}\bigp{\de\wedge\abs{\xi_i}};
\end{equation*} 
here, in distinction with the definition of $H_{1,2}$ in \cite{chen07}, 
the notation $\check\xi_i$ is used in place of $\eta_i$. 
Then the Cauchy-Schwarz inequality yields
\begin{equation*}
 H_{1,2}\le\sqrt{\E\bigind{z\le W\le z+\abs{\bar\De}}}\sqrt{\tsum\nolim_i\E\check\xi_i^2}\le\de\sqrt{\pp},\quad\text{where}\quad\pp:=\bigprob{z\le W\le z+\abs{\bar\De}};
\end{equation*}
cf.\ \cite[(5.8)]{chen07}. 
Following through with the remainder of the proof of \cite[Theorem~2.1]{chen07}, we have
\begin{equation*}
 \cc\pp-\de\pp^{1/2}\le\bb:=\frac12\Bigp{2\de+\E\bigabs{W\bar\De}+\tsum\nolim_i\E\bigabs{\xi_i(\bar\De-\De_i)}}. 
\end{equation*}
So, 
\begin{equation*}
 \pp\le\biggp{\frac{\de+\sqrt{\de^2+4\cc\bb}}{2\cc}}^2=\frac{2\de^2+4\cc\bb+2\de\sqrt{\de^2+4\cc\bb}}{4\cc^2}
  =\frac1{2\cc}\bigg(2\bb+\frac{\de^2}{\cc}
  +2\de\sqrt{\frac1{2\cc}
  \Big(2\bb+\frac{\de^2}{2\cc}\Big)}\,\bigg)
  ;
\end{equation*} 
in view of \eqref{eq:conc} and \eqref{eq:p+Geta(w)}, this verifies the improvement provided in  Remark~\ref{re:chen.improvement}. 
\end{proof}

\begin{proof}[Proof of Theorem~\ref{thm:nub}]
The proof of Theorem~\ref{thm:nub} largely follows the lines of that of \cite[Theorem~2.2]{chen07}; for the ease of comparison between the two proofs, we shall use notation similar to that in \cite{chen07}. 
The extension to $p$ other than 2 is obtained using a Cram\'er-tilt absolutely continuous transformation of measure along with the mentioned Rosenthal-type and exponential bounds. 
Introduce the Winsorized r.v.'s
\begin{equation}\label{eq:barW}
 \bar\xi_i:=\xi_i\wedge w\quad\text{and their sum,}\quad\bar W:=\sum_{i=1}^n\bar\xi_i.
\end{equation}
Note that in the statement of \cite[Lemma~5.1]{chen07} the $\bar\xi_i$'s are defined as the truncated r.v.'s $\xi_i\ind{\xi_i\le w}$ (with $w=1$). 
A problem with this definition arises on page~596 in \cite{chen07} concerning the assertion there that $\sum_i\E|\xi_i|(\de\wedge|\bar\xi_i|)=\sum_i\E|\xi_i|(\de\wedge|\xi_i|)$ whenever $\de\le0.07$; indeed, by letting $\xi_i$ take values $\pm2$ each with probability $\frac1{8n}$ and the value $0$ with probability $1-\frac1{4n}$, the assertion is seen to be false when $\bar\xi_i=\xi_i\ind{\xi_i\le 1}$ (while true if $\de\le w$ and $\bar\xi_i=\xi_i\wedge w$). 
See \cite{pin09.winsor} for a general discussion on comparative merits of the Winsorization vs.\ truncation, especially in regard to the Cram\'er tilt transformation. 

Recalling the definition \eqref{eq:Phat} of the measure $\hat\P$, one has 
\begin{equation}\label{eq:chenlem51}
\begin{split}
 &\hat\P\bigp{z-|\De|\le W\le z}\\
 &\qquad
   =\bigprob{z-|\De|\le W\le z,|\De|\le\pi_1z}\\
 &\qquad
   \le\tsum_{i=1}^n\bigprob{W\ge(1-\pi_1)z,\eta_i>w} + \bigprob{z-|\De|\le W\le z,|\De|\le\pi_1z,\max\nolim_i\eta_i\le w}\\
 &\qquad
   \le\tsum_{i=1}^n\bigprob{\xi_i>\pi_2z} + \tsum_{i=1}^n\bigprob{W-\xi_i\ge(1-\pi_1-\pi_2)z}\bigprob{\eta_i>w}\\
 &\qquad\qquad 
   + \bigprob{z-|\bar\De|\le\bar W\le z,|\bar\De|\le\pi_1z}\\
 &\qquad
   =\ga_z + \bigprob{z-|\bar\De|\le\bar W\le z,|\bar\De|\le\pi_1z}
   ;
\end{split}
\end{equation}
here the second inequality follows from the independence of $W-\xi_i$ and $\eta_i$, the condition \eqref{eq:barDe} on $\bar\De$, and the definition \eqref{eq:barW} of $\bar W$ (recall also the condition that  $\xi_i\le\eta_i$), and the second equality follows from the definitions of $\ga_z$ and $\pi_3$ in \eqref{eq:ga_z} and \eqref{eq:pi's}; cf.\ \cite[Lemma~5.1]{chen07}.

We must next establish the inequality
\begin{equation}\label{eq:conc<=tau}
 \P(z-|\bar\De|\le\bar W\le z,|\bar\De|\le\pi_1z)\le\tau e^{-(1-\pi_1)z/\th}; 
\end{equation}
cf.\ \cite[Lemma~5.2]{chen07}. 
Consider two cases: 
\begin{equation*}
	\text{(i) $\de>\de_0$\quad and\quad (ii) $0<\de\le\de_0\le w$}
\end{equation*}
(recall the restriction on the number $\de_0$ in \eqref{eq:pi's}). 
In the first case, when $\de>\de_0$, 
\begin{multline*}
 \bigprob{z-|\bar\De|\le\bar W\le z,|\bar\De|\le\pi_1z}
 \le\bigprob{\bar W\ge(1-\pi_1)z}
  \le\E e^{\bar W/\th}e^{-(1-\pi_1)z/\th} \\
  \le\tfrac{\de}{\de_0}\PUexp\bigp{\tfrac{1}{\th},w,1,\vp_1}e^{-(1-\pi_1)z/\th}
  \le c_3\de e^{-(1-\pi_1)z/\th}
  \le\tau e^{-(1-\pi_1)z/\th};
\end{multline*}
here \eqref{eq:PUexp} and \eqref{eq:PUcond} are used for the third inequality above \big(as well as the definitions \eqref{eq:vp1} and \eqref{eq:si_p} of $\vp_1$ and $\si_p$\big), and the definitions \eqref{eq:c3} and \eqref{eq:tau} of $c_3$ and $\tau$ are used for the last two inequalities there. 
Thus, \eqref{eq:conc<=tau} is established when $\de>\de_0$.

Consider now the second case, when $0<\de\le\de_0\le w$. 
Let
\begin{equation*}
 f_{\bar\De}(u):=
 \begin{cases}
  0&\text{if }u<z-|\bar\De|-\de,\\
  e^{u/\th}(u-z+|\bar\De|+\de)&\text{if }z-|\bar\De|-\de\le u<z+\de,\\
  e^{u/\th}(|\bar\De|+2\de)&\text{if }u\ge z+\de
 \end{cases}
\end{equation*}
be defined similarly to \cite[(5.16)]{chen07}. 
Then, by the independence of $(\De_i,\bar W-\bar\xi_i)$ and $\xi_i$,  
\begin{equation}\label{eq:G1+G2}
 \E Wf_{\bar\De}(\bar W)=G_1+G_2, 
\end{equation}
where 
\begin{equation*}
 G_1:=\tsum_{i=1}^n\E\xi_i\bigp{f_{\bar\De}(\bar W)-f_{\bar\De}(\bar W-\bar\xi_i)}
 \quad\text{and}\quad 
 G_2:=\tsum_{i=1}^n\E\xi_i\bigp{f_{\bar\De}(\bar W-\bar\xi_i)-f_{\De_i}(\bar W-\bar\xi_i)}.
\end{equation*}
Also, using an obvious modification of the arguments associated with \cite[(5.17)--(5.19)]{chen07}, one has 
\begin{equation}\label{eq:G11-G12}
 G_1\ge G_{1,1}-G_{1,2},
\end{equation}
where
\begin{equation}\label{eq:G11}
 G_{1,1}:=\cc\exp\bigb{\tfrac1{\th}\bigp{(1-\pi_1)z-\de}}\P(z-|\bar\De|\le\bar W\le z,|\bar\De|\le\pi_1z), 
\end{equation}
\begin{equation*}
 G_{1,2}:=\E\int_{|t|\le\de}e^{(\bar W-\de)/\th}\bigabs{\bar M(t)-\E\bar M(t)}\,dt, 
\end{equation*}
\begin{equation*}
 \bar M(t):=\tsum_{i=1}^n\bar M_i(t), 
 \quad\text{and}\quad 
 \bar M_i(t):=\xi_i\bigp{\ind{-\bar\xi_i\le t\le 0}-\ind{0<t\le-\bar\xi_i}};
\end{equation*}
in particular, the factor $c_*$ in the expression \eqref{eq:G11} for $G_{1,1}$ arises when one uses the relations $\int_{|t|\le\de}\E\bar M(t)\,dt=\sum_i\E|\xi_i|(\de\wedge|\xi_i|)\ge\cc$, which in turn follow by the condition $\de\le\de_0\le w$ of case (ii) and \eqref{eq:de}; cf.\ \cite[(5.19)]{chen07}. 
Further, 
\begin{equation*}
 \int_{|t|\le\de}\E\bigp{\bar M(t)-\E\bar M(t)}^2\,dt
  \le\sum_{i=1}^n\E\int_{|t|\le\de}\bar M_i(t)^2\,dt
  =\sum_{i=1}^n\E\xi_i^2\bigp{\de\wedge\abs{\bar\xi_i}}
  \le\de,
\end{equation*}
so that two applications of the Cauchy-Schwarz inequality yield
\begin{align}
\notag
 G_{1,2}&
  \le\E\biggp{\int_{|t|\le\de}e^{2(\bar W-\de)/\th}\,dt}^{1/2}\biggp{\int_{|t|\le\de}\bigp{\bar M(t)-\E\bar M(t)}^2\,dt}^{1/2}
  \le\Bigp{2\de\E e^{2(\bar W-\de)/\th}}^{1/2}\sqrt\de\\
&\label{eq:G12}
  \le\Bigp{2\PUexp\bigp{\tfrac2\th,w,1,\vp_1}}^{1/2}e^{-\de/\th}\de
  =\sqrt{2}\PUexp\bigp{\tfrac{2}{\th},w,\tfrac1{\sqrt{2}},\vp_1}e^{-\de/\th}\de,
\end{align}
where the last inequality follows from \eqref{eq:PUexp} and \eqref{eq:PUcond} (recalling also the definitions \eqref{eq:barW} and \eqref{eq:vp1} of $\bar W$ and $\vp_1$); the equality in \eqref{eq:G12} follows from the easily verified identity
\begin{equation}\label{eq:PUexp,pow}
 \PUexp\bigp{\la,y,B,\vp}^{\al}=\PUexp\Bigp{\la,y,\al^{1/2}B,\vp}\quad\text{for any }\al>0.
\end{equation}
Next (cf.\ \cite[(5.21)]{chen07}),
\begin{align}
\notag
 |G_2|&
  \le\tsum_{i=1}^n\E\bigabs{\xi_ie^{(\bar W-\bar\xi_i)/\th}(\bar\De-\De_i)}
  \le\tsum_{i=1}^n\bignorm{\xi_ie^{(\bar W-\bar\xi_i)/\th}}_p\,\bignorm{\bar\De-\De_i}_q
  =\tsum_{i=1}^n\E^{1/p}e^{\frac p\th(\bar W-\bar\xi_i)}\,\|\xi_i\|_p\|\bar\De-\De_i\|_q
\\\label{eq:G2}
& \le\PUexp\bigp{\tfrac{p}{\th},w,\tfrac1{\sqrt{p}},\vp_1}\tsum_{i=1}^n\norm{\xi_i}_p\norm{\bar\De-\De_i}_q.
\end{align}
Also, 
\begin{equation}\label{eq:Wf(W)}
 \E Wf_{\bar\De}(\bar W)\le\E\bigp{|\bar\De|+2\de}|W|e^{\bar W/\th}
  \le \bigp{\bignorm{\bar\De}_q+2\de}\bignorm{We^{\bar W/\th}}_p.
\end{equation}
Chen and Shao \cite{chen07} bounded $\E W^2e^{\bar W}$ (corresponding to the case when $p=2$ and $\th=2$ in \eqref{eq:Wf(W)}) with an absolute constant; in our case, more work is required to bound the last factor in \eqref{eq:Wf(W)} for the general $p$. 
Specifically, we apply Cram\'er's tilt transform to the $\xi_i$'s, using at that results of \cite{pin09.winsor,pin11tilt,pin12.tilt.symm}.

Let $\bs\xi:=(\xi_1,\dotsc,\xi_n)$, and for any real $c>0$ let $\hat{\bs\xi}=:(\hat{\xi}_1,\dotsc,\hat{\xi}_n)$ be a random vector such that 
\begin{equation*}
 \P(\hat{\bs\xi}\in E)=\frac{\E e^{c\bar W}\ind{\bs\xi\in E}}{\E e^{c\bar W}}
\end{equation*}
for all Borel sets $E\subseteq\R^n$. 
Then the $\hat\xi_i$'s are necessarily independent r.v.'s; moreover, if $f\colon\R^n\to\R$ is any nonnegative Borel function, then
\begin{equation}\label{eq:hatxi}
 \E f(\hat{\bs\xi})
  =\frac{\E f(\bs\xi)e^{c\bar W}}{\E e^{c\bar W}}.
\end{equation}
By \cite[Proposition~2.6,(I)]{pin11tilt}, $\E\hat\xi_i$ is nondecreasing in $c$, so that $\E\hat\xi_i\ge\E\xi_i=0$, and so, by \cite[Corollary~2.7]{pin11tilt},
\begin{equation*} \bigabs{\tsumnl_i\E\hat\xi_i}=\tsumnl_i\E\hat\xi_i\le\frac{e^{cw}-1}{w}\tsumnl_i\E\xi_i^2=\frac{e^{cw}-1}{w}. 
\end{equation*}
If the $\xi_i$'s are assumed to have symmetric distributions, then \cite[Theorem~1]{pin12.tilt.symm} allows for the factor $(e^{cw}-1)/w$ above to  be replaced by $\sinh(cw)/w$; cf.\ Remark~\ref{re:symm}. 
Choose now 
\begin{equation*}
	c=\frac p\th. 
\end{equation*}
Then, by \cite[Theorem~2.1]{pin09.winsor}, 
\begin{equation*}
 \E e^{c\bar\xi_i}=\E e^{c(\xi_i\wedge w)}=\E e^{cw(1\wedge\xi_i/w)}\ge L_{W;\,cw,\norm{\xi_i}_2/w}\ge L_{W;\,cw,\max_i\norm{\xi_i}_2/w}=a_1^{-1},
\end{equation*}
where $a_1$ is as defined in \eqref{eq:a1}; the last inequality above follows because $L_{W;\,c,\si}$ in \cite[(2.9)]{pin09.winsor} is nonincreasing in $\si$; the condition $c=\frac p\th$ was used here in the above display only for the last equality. 
So, 
\begin{equation*}
 \E\abs{\hat\xi_i}^p
  =\frac{\E\abs{\xi_i}^pe^{c\bar\xi_i}}{\E e^{c\bar\xi_i}}
  \le a_1e^{cw}\E\abs{\xi_i}^p,
\end{equation*}
with $\sum_i\E\hat\xi_i^2\le a_1e^{cw}$ a consequence of this. 
Next, 
\begin{equation}\label{eq:mod4.1}
\begin{split}
 \bignorm{\tsum\nolim_i\hat\xi_i}_p
 &\le\bignorm{\tsum\nolim_i(\hat\xi_i-\E\hat\xi_i)}_p
   +\bigabs{\tsum\nolim_i\E\hat\xi_i}
\\
%RM 02-21-13
%  &\le\AA_\R(p)\bigp{\tsum\nolim_i\E\abs{\hat\xi_i-\E\hat\xi_i}^p}^{1/p}
%    +\BB_\R(p)\bigp{\tsum\nolim_i\E(\hat\xi_i-\E\hat\xi_i)^2}^{1/2}
%    +(e^{pw/\th}-1)/w
% \\
 &\le
 \AA_{\R,\nc}(p)%RM 02-21-13 \AA_\R(p)
 \bigp{
 %RM 02-21-13 1.32
 \tsum\nolim_i\E\abs{\hat\xi_i}^p}^{1/p}+
 \BB_{\R,\nc}(p)%RM 02-21-13 \BB_\R(p)
 \bigp{\tsum\nolim_i\E\hat\xi_i^2}^{1/2}+(e^{pw/\th}-1)/w
\\
 &\le
 \AA_{\R,\nc}(p)%RM 02-21-13 \AA_\R(p)
 \bigp{
 %RM 02-21-13 1.32
 a_1e^{pw/\th}\si_p^p}^{1/p}+
 \BB_{\R,\nc}(p)%RM 02-21-13 \BB_\R(p)
 \bigp{a_1e^{pw/\th}}^{1/2}+(e^{pw/\th}-1)/w,
\end{split}
\end{equation}
where 
\eqref{eq:rosen.noncentral} %RM 02-21-13 \eqref{eq:rosen1} 
is used for the second inequality above. %RM 02-21-13 , and \cite[Theorem~2.3(v)--(vi)]{re-center} is used for the third inequality. 
Letting $f(x_1,\dotsc,x_n)\allowbreak\equiv |\sum_ix_i|^p$ in \eqref{eq:hatxi} and using \eqref{eq:PUexp}, \eqref{eq:PUcond}, and \eqref{eq:PUexp,pow} once more, one has
\begin{equation}\label{eq:mod4.2}
 \bignorm{We^{\bar W/\th}}_p
  = \Bigp{\E\bigabs{\tsum\nolim_i\xi_i}^pe^{p\bar W/\th}}^{1/p}
  = \Bigp{\E e^{p\bar W/\th}\E\bigabs{\tsum\nolim_i\hat\xi_i}^p}^{1/p}
  \le \PUexp\bigp{\tfrac{p}{\th},w,\tfrac1{\sqrt{p}},\vp_1}\bignorm{\tsum\nolim_i\hat\xi_i}_p.
\end{equation}

Thus, recalling the case condition $\de\le\de_0$, we have
\begin{equation*}
\begin{split}
 \bigprob{z-|\bar\De|\le\bar W\le z,|\bar\De|\le\pi_1z}
 &=\tfrac1{\cc}\,e^{-(1-\pi_1)z/\th}e^{\de/\th}G_{1,1}\\
 &\le\tfrac1{\cc}\,e^{-(1-\pi_1)z/\th}e^{\de/\th}\bigp{G_{1,2}+|G_2|+\E Wf_{\bar\De}(\bar W)}\\
 &\le\bigp{c_1\tsum\nolim_i\norm{\xi_i}_p\norm{\bar\De-\De_i}_q+c_2\norm{\bar\De}_q+c_3\de}e^{-(1-\pi_1)z/\th},
\end{split}
\end{equation*}
where the equality comes from the definition \eqref{eq:G11} of $G_{1,1}$, the first inequality follows from \eqref{eq:G1+G2} and \eqref{eq:G11-G12}, and the second inequality follows from \eqref{eq:G12}, \eqref{eq:G2}, \eqref{eq:Wf(W)}, \eqref{eq:mod4.2}, and \eqref{eq:mod4.1}, along with the definitions \eqref{eq:c1}, \eqref{eq:c2}, and \eqref{eq:c3} of $c_1$, $c_2$, and $c_3$. 
Thus, in view of the definition \eqref{eq:tau} of $\tau$, the inequality \eqref{eq:conc<=tau} is proved for the other case, $\de\le\de_0$.

Replace now $\P$ with $\hat\P$ in \eqref{eq:conc}, so that \eqref{eq:chenlem51} and \eqref{eq:conc<=tau} imply
\begin{equation*}
 \hat\P(W\le z)-\hat\P(T\le z)\le\ga_z+\tau e^{-(1-\pi_1)z/\th}.
\end{equation*}
In a similar fashion, one bounds $\hat\P(T\le z)-\hat\P(W\le z)$ from above, establishing \eqref{eq:nub}.
\end{proof}

\subsection{Proofs of results from Section~\ref{sec:f(S)}}
\label{sec:f(S).proofs}

%RM12.26
\begin{proof}[Proof of Remark~\ref{re:compos}]
In view of \eqref{eq:h 2smooth}, there exists $m_h\in(0,\infty)$ such that
\begin{equation}\label{eq:h 1smooth}
 \norm{h(x)}_\YY\le m_h\norm{x}_\XX\quad\text{for all $x\in\XX$ with $\norm{x}_\XX\le\ep_h$;}
\end{equation}
indeed, we may let $m_h:=\norm{L_h}+M_h\ep_h/2$. 
Assume that $\ep_h$ is chosen small enough to ensure $m_h\ep_h\le\ep_g$. 

Take any $x\in\XX$ with $\norm{x}_\XX\le\ep_h$. 
Then, by \eqref{eq:h 2smooth}, there is some $y_x\in\YY$ such that 
$\norm{y_x}_\YY\le1$ and $h(x)=L_h(x)+\frac12\,M_h\norm{x}_\XX^2y_x$. 
By \eqref{eq:h 1smooth}, $\norm{h(x)}_\YY\le m_h\ep_h\le\ep_g$, and so, by \eqref{eq:g 2smooth}, there is some $z_x\in\ZZ$ such that $\norm{z_x}_\ZZ\le1$ and 
\begin{align*}
 g(h(x))&=L_g(h(x))+\tfrac12\,M_g\norm{h(x)}_\YY^2z_x\\
  &=L_g(L_h(x))+\tfrac12\,M_h\norm{x}_\XX^2L_g(y_x)+\tfrac12\,M_g\norm{h(x)}_\YY^2z_x.
\end{align*}
Thus, by \eqref{eq:h 1smooth} (recall also $\norm{y_x}_\YY\le1$ and $\norm{z_x}_\ZZ\le1$), 
\begin{equation*}
 \norm{(g\circ h)(x)-(L_g\circ L_h)(x)}_\ZZ\le\tfrac12\bigp{M_h\norm{L_g}+M_gm_h^2}\norm{x}_\XX^2\quad\text{for all $x\in\XX$ with $\norm{x}_\XX\le\ep_h$;}
\end{equation*}
that is, \eqref{eq:smooth} with $\ZZ$ in place of $\R$ holds for $f=g\circ h$ with 
$L=L_g\circ L_h$, $\Mf=M_h\norm{L_g}+M_gm_h^2$, and $\ep=\ep_h$. 
\end{proof}

The uniform and nonuniform BE type bounds in Theorems~\ref{thm:f(S).ub} and \ref{thm:f(S).nub} rely on the corresponding bounds of Section~\ref{sec:chen.mod}. 
Let $f$ be a function satisfying \eqref{eq:smooth}, and also let $X_1,\dotsc,X_n$ be independent zero-mean $\XX$-valued random vectors. 
Further let $\si=\norm{L(S)}_2$, as in \eqref{eq:si}, and for $i=1,\dotsc,n$ let
\begin{equation*}
 g_i(x):=\frac{L(x)}{\si}\quad\text{and}\quad\xi_i=g_i(X_i)=\frac{L(X_i)}{\si},
\end{equation*}
in accordance with \eqref{eq:xi,eta}. 
The choices for the functions $h_i$ (used to define the r.v.'s $\eta_i$) will depend on the value of $p$ and the type of bound (uniform or nonuniform) being derived (cf.\ \eqref{eq:eta.ub} and \eqref{eq:eta.nub}). 
Next, let
\begin{equation*}
 T:=\frac{f(S)}{\si},\quad W:=\tsum\nolim_i\xi_i=\frac{L(S)}{\si},
\end{equation*}
and also
\begin{equation}\label{eq:tT}
 \tilde T:=T\,\ind{\norm{S}\le\ep}+W\ind{\norm{S}>\ep}.
\end{equation}
Finally, let
\begin{equation}\label{eq:De}
 \De:=\frac{\Mf}{2\si}\,\norm{S}^2.
\end{equation}
Then, by \eqref{eq:smooth},
\begin{align*}
 |\tilde T-W|
 & =\si^{-1}\bigabs{f(S)-L(S)}\ind{\norm{S}\le\ep}
   \le\tfrac{\Mf}{2\si}\norm{S}^2=\De.
\end{align*}
Adopt some more notation:
\begin{equation}\label{eq:tX}
 \tX_i:=X_i\bigind{\eta_i\le w},\quad
 \tS:=\tsum\nolim_i\tX_i,
\end{equation}
\begin{equation}\label{eq:De-bar}
 \bar\De:=\tfrac{\Mf}{2\si}\Bigp{\norm{S}^2\ind{p=3}+\norm{\tS}^2\ind{p<3}},
\end{equation}
\begin{equation}\label{eq:De_i}
 \De_i:=\tfrac{\Mf}{2\si}\Bigp{\norm{S-X_i}^2\ind{p=3}+\norm{\tS-\tX_i}^2\ind{p<3}}.
\end{equation}
Then the assumptions of Theorems~\ref{thm:ub} and \ref{thm:nub} are satisfied for the nonlinear statistic $\tilde T$ (in place of $T$) and its linear approximation $W$; particularly, $\E\xi_i=0$, $\var W=1$, $|\De|\ge|\tilde T-W|$, $\bar\De$ satisfies \eqref{eq:barDe}, and $\De_i$ satisfies the condition that $X_i$ and $(\De_i,(X_j\colon j\ne i))$ are independent (which further implies that $X_i$ and $(\De_i,W-\xi_i)$ are independent).

\begin{lem}\label{lem:barDe}
Under the conditions of Theorem~\ref{thm:f(S).ub}, $\norm{\bar\De}_q\le\uu,$ where $\uu$ is as defined in \eqref{eq:u}.
\end{lem}

\begin{lem}\label{lem:xi(De-De_i)}
Under the conditions of Theorem~\ref{thm:f(S).ub}, $\tsum_{i=1}^n\norm{\xi_i}_p\norm{\bar\De-\De_i}_q\le\si_p\vv$, where $\si_p$ and $\vv$ are as defined in \eqref{eq:si_p} and \eqref{eq:v}, respectively.
\end{lem}

The proofs of these lemmas (and subsequent ones) are deferred to the end of this subsection.

\begin{proof}[Proof of Theorem~\ref{thm:f(S).ub}]
Recall that the conditions of Theorem~\ref{thm:ub} hold, with $\tilde T$ in place of $T$, so that \eqref{eq:tT} and \eqref{eq:ub} imply
\begin{align}\notag
 \bigabs{\P(T\le z)-\P(W\le z)}
 &\le\P(\norm{S}>\ep)+\bigabs{\P(\tilde T\le z)-\P(W\le z)}
\\\label{eq:f(S),ub,1}
 &\le\P(\norm{S}>\ep)+\frac1{2\cc}\Bigp{4\de+\bignorm{W}_p\bignorm{\bar\De}_q+\tsum_{i=1}^n\E\norm{\xi_i}_p\norm{\bar\De-\De_i}_q}+G_\eta(w)
\end{align}
for all $z\in\R$. 
%RM 12.21.12
% The use of \cite[Remark~2.1]{chen07} allows us to choose
% \begin{equation}\label{eq:de.to.si}
%  \de=c_4\si_p^\tq,
% \end{equation}
% in accordance with \eqref{eq:de}, where $\tq$ and $c_4$ are as defined in \eqref{eq:c4}. 
Along with \eqref{eq:f(S),ub,1},
%RM 02-21-13 choose $\de$ as in \eqref{eq:c4}, %RM 12.21.12
use Lemmas~\ref{lem:barDe} and \ref{lem:xi(De-De_i)}, and apply the Rosenthal-type inequality \eqref{eq:rosen} to obtain $\norm{W}_p\le\AA_\R(p)\si_p+\BB_\R(p)$. Then \eqref{eq:f(S).ub} follows, and the proof of Theorem~\ref{thm:f(S).ub} is complete. 
\end{proof}

\begin{proof}[Proof of Theorem~\ref{thm:f(S).nub}]
Recall that the conditions of Theorem~\ref{thm:nub} hold with $\tilde T$ in place of $T$. 
Also, by \eqref{eq:De}, \eqref{eq:z}, and \eqref{eq:om}, 
\begin{equation*}
 \bigb{\abs{\De}\le\pi_1z}=\bigb{\norm{S}\le(2\pi_1\si z/\Mf)^{1/2}}\subseteq\bigb{\norm{S}\le(2\pi_1\om/\Mf)^{1/2}}\subseteq\bigb{\norm{S}\le\ep}.
\end{equation*}
Thus, by Remark~\ref{re:Phat}, \eqref{eq:Phat}, \eqref{eq:tT}, and \eqref{eq:nub}, 
\begin{equation*}
\begin{split}
 \bigabs{\P(T\le z)-\P(W\le z)}&\le\bigabs{\hat\P(T\le z)-\hat\P(W\le z)}+\bigprob{\abs{\De}>\pi_1z}\\
  &=\bigabs{\hat\P(\tilde T\le z)-\hat\P(W\le z)}+\bigprob{\abs{\De}>\pi_1z}\\
  &\le \tilde\ga_z+\tau e^{-(1-\pi_1)z/\th}
\end{split}
\end{equation*}
for all $z$ as in \eqref{eq:z}, where $\tilde\ga_z$ is as in \eqref{eq:tga_z}. 
Recall the definitions \eqref{eq:tau} and \eqref{eq:tau2} of $\tau$ and $\tilde\tau$, respectively, to see that $\tau\le\tilde\tau$ follows from 
%RM 02-21-13 \eqref{eq:c4} %RM12.21.12 \eqref{eq:de.to.si} 
% and 
Lemmas~\ref{lem:barDe} and \ref{lem:xi(De-De_i)}. 
Then \eqref{eq:f(S).nub} is proved.
\end{proof}

The following lemma provides two bounds on $\tilde\ga_z$ in \eqref{eq:tga_z} which will be used in the proofs of %IP16 Corollary
Theorem~\ref{cor:iid} and Theorem~\ref{thm:iid,p=3,nonunif}. 

\begin{lem}\label{lem:tga_z}
Assume that the conditions of Theorem~\ref{thm:f(S).nub} hold. 
Take any real numbers $\ka_2>0$ and $\ka_3>0$, and let 
\begin{equation}\label{eq:x2,y2,vp2}
 x_2:=\Bigp{\frac{2\pi_1}{\Mf}\,\si z}^{1/2},\quad 
 y_2:=\frac{x_2}{\ka_2},\quad
 \vp_2:=\frac{s_p^p}{s_2^2y_2^{p-2}}\wedge1,\quad 
 S_{y_2}:=\tsum_{i=1}^nX_i\ind{\norm{X_i}\le y_2},
\end{equation}
\begin{equation}\label{eq:x3,y3,vp3}
 x_3:=\pi_3z,\quad y_3:=\frac{x_3}{\ka_3},\quad\vp_3:=\frac{\si_p^p}{y_3^{p-2}}\wedge1,
\end{equation}
\begin{equation}\label{eq:PU2,PU3}
 \PUtwo:=\PU\bigp{x_2,y_2,s_2^2,\E\norm{S_{y_2}},\vp_2},
 \quad\text{and}\quad\PUthr:=\PU\bigp{x_3,y_3,1,0,\vp_3},
\end{equation}
where $\PU$ is as in \eqref{eq:PU}. 
Then 
\begin{align}
\label{eq:ga<=tga}
 &\tilde\ga_z\le G_X(y_2)+\PUtwo+G_\xi(\pi_2z)+\bigp{G_\xi(y_3)+\PUthr}G_\eta(w)
\end{align}
for all $z>0$, where $\tilde\ga_z$ is as in \eqref{eq:tga_z}. 

One consequence of \eqref{eq:ga<=tga} is that
\begin{equation}\label{eq:ga_z2}
\begin{split}
 \tilde\ga_z
 &\le
   G_X\biggp{\Bigp{\frac{\pi_1}{2p^2\om\Mf}}^{1/2}\,\si z}+\Bigp{\frac{2ep\Mf D^2}{\pi_1}\,\frac{s_2^2}{\si z}}^p
  + G_{\xi}\bigp{\pi_2z} + \biggp{
     G_{\xi}\Bigp{\frac{2\pi_3}{p}\,z}+\frac{(ep/2)^{p/2}}{(\pi_3z)^p}}
    G_\eta(w)
\end{split}
\end{equation}
for all $z$ as in \eqref{eq:z}.
\end{lem}

In the proof of %IP16 Corollary
Theorem~\ref{cor:iid}, let us write $a\O b$ if $|a|\le\CC b$ for some $\CC$ as in Corollary~\ref{cor:iid}. 
Let us then write $a\asymp b$ if $a\O b$ and $b\O a$. 

\begin{proof}[Proof of %IP16 Corollary
Theorem~\ref{cor:iid}]
Set $\cc=\frac12$, $w=\de_0=1$, $\pi_1=(\Mf\ep^2/(2\om))\wedge\frac13$, $\pi_2=\pi_3=\frac12(1-\pi_1)$, and $\th=\tth(1-\pi_1)$ in the statements of Theorems~\ref{thm:f(S).ub} and \ref{thm:f(S).nub}, so that \eqref{eq:pi's} and \eqref{eq:om} be satisfied. 
Further let $X_i=\frac1n V_i$. 
Then $S=\sum_{i=1}^nX_i=\bar V$ and, by the definitions \eqref{eq:si}, \eqref{eq:s_al}, \eqref{eq:si_p}, and \eqref{eq:la_p},
\begin{equation}\label{eq:iid,orders}
 \si=\frac{\tsi}{n^{1/2}},\quad 
 s_\al=\frac{\norm{V}_\al}{n^{1-1/\al}},\quad 
 \si_\al=\frac{\norm{L(V)}_\al}{\tsi n^{1/2-1/\al}},\quad\text{and}\quad
 \la_\al=\frac{\norm{L}\norm{V}_\al}{\tsi n^{1/2-1/\al}}
\end{equation}
for any $\al\ge1$. 
Letting $\de$ be as in \eqref{eq:c4}, and recalling %RM 02-21-13 Recalling 
also the definitions \eqref{eq:u}, \eqref{eq:v}, 
%RM 02-21-13 \eqref{eq:c4}, 
and \eqref{eq:tau2}, as well as Remark~\ref{re:tau}, one has
\begin{equation}\label{eq:iid,u,v,ttau}
 \uu\asymp n^{-1/2},\quad \si_p\vv\asymp n^{-1/2},\quad 
 \de%IP 03.12.13\asymp
 \asymp%RM 03.20.13 \O 
 n^{-1/2},\quad %RM 02-21-13 \si_p^{\tq}\asymp n^{-1/2},\quad
 \text{and hence}\quad \tilde\tau
 \asymp%RM 03.20.13 \O 
 n^{-1/2}
\end{equation}
for all $p\in(2,3]$; moreover, it is clear that the above expressions depend on the distribution of $V$ only through $\tsi$, $\norm{L(V)}_p$, $\norm{V}_q$, $\norm{V}_2$, and $\norm{V}_p$. 
Also, for any $t>0$, \eqref{eq:eta.ub} and \eqref{eq:eta.nub} imply 
\begin{align}
\label{eq:G_X.iid}
 G_X(t)   &=n\bigprob{\norm{V}>nt}\le\frac{\norm{V}_p^p}{n^{p-1}t^p}\\
 \text{and }G_\xi(t) &\le G_\eta(t)
  \le n\bigprob{\norm{L}\norm{V}>\sqrt{n}\tsi t}
 \O\frac1{n^{p/2-1}t^p}. 
\label{eq:G_eta.iid}
\end{align}

By \eqref{eq:|S|>ep}, $\P(\norm{S}>\ep)\le\norm{V}_2^2/(\ep^2n)$. 
Next, there exists a positive absolute constant $A$ such that
\begin{equation*}
 \sup_{z\in\R}\bigabs{P(\sqrt{n}L(\bar V)/\tsi\le z)-\Phi(z)}\le A\,\frac{\norm{L(V)}_p^p}{n^{p/2-1}},
\end{equation*}
which follows from, say, Theorem~6 of \cite[Chapter V]{pet75}. 
Then \eqref{eq:f(S).ub}, \eqref{eq:iid,u,v,ttau}, and \eqref{eq:G_eta.iid} yield \eqref{eq:f(S).iid}. 

Using \eqref{eq:iid,u,v,ttau} and recalling that $\th=(1-\pi_1)\tth$, one has $\tilde\tau e^{-(1-\pi_1)z/\th}\O 1/(\sqrt ne^{z/\tth})$. 
In view of \eqref{eq:f(S).nub}, \eqref{eq:G_X.iid}, \eqref{eq:G_eta.iid}, and \eqref{eq:ga_z2}, one obtains \eqref{eq:f(S).iid.nu} with $\Phi(z)$ there replaced by $\P(\sqrt{n}L(\bar V)/\tsi\le z)$. 
To obtain \eqref{eq:f(S).iid.nu} as stated, note that 
%IP12.13.12 \cite[Corollary~1.3]{pin11.BE} implies
\begin{equation}\label{eq:os}
 \Bigabs{\Bigprob{\frac{L(\bar V)}{\tsi/\sqrt{n}}\le z}-\Phi(z)}\O\frac{\si_p^p}{e^{z/\tth}}+G_\xi\Bigp{\frac{z}{1+p/2}}+\frac{G_\xi(1)}{z^p}
\end{equation} 
for all $z>0$%IP12.13.12 
; this follows by \cite[Corollary~1.3]{pin11.BE} with $v=w=1$, $c=0$, and $\la=1/\tth$ (in notation therein), using at that the inequalities $\be_v\le\mu_p/v^p$ (displayed right after \cite[(1.2)]{pin11.BE}) and $P_1\wedge\dots\wedge P_5\le P_4$. 
Combining \eqref{eq:os} with \eqref{eq:iid,orders} and \eqref{eq:G_eta.iid}, one completes the proof.
\end{proof}

\begin{proof}[Proof of Lemma~\ref{lem:barDe}]
Suppose first that $p=3$, so that, in accordance with \eqref{eq:De-bar}, $\bar\De=\frac{\Mf}{2\si}\norm{S}^2$. 
Then, by the Rosenthal-type inequality \eqref{eq:rosen} and the definitions \eqref{eq:la_p} and \eqref{eq:u} of $\la_\al$ and $\uu$, respectively, 
\begin{equation*}
 \norm{\bar\De}_q
  =\tfrac{\Mf}{2\si}\,\norm{S}_{2q}^2
  \le\tfrac{\Mf}{2\si}\bigp{
  \AA_\XX(2q)^2%RM 03.20.13 \AA_\XX^2(2q)
  s_{2q}^2+
  \BB_\XX(2q)^2%RM 03.20.13 \BB_\XX^2(2q)
  s_2^2}
  =\tfrac{\Mf\si}{2\norm{L}^2}\bigp{
  \AA_\XX(2q)^2%RM 03.20.13 \AA_\XX^2(2q)
  \la_{2q}^2+
  \BB_\XX(2q)^2%RM 03.20.13 \BB_\XX^2(2q)
  \la_2^2}
  =\uu,
\end{equation*}
which proves the lemma when $p=3$.

Now suppose that $p\in(2,3)$.
By \eqref{eq:tX}, \eqref{eq:la_p} and \eqref{eq:eta.ub},
\begin{equation}\label{eq:||EtS||}
\begin{split}
 \norm{\E\tS}
 &=\bignorm{\tsum\nolim_i\E X_i\ind{\eta_i>w}}
  \le\tsum\nolim_i\E\norm{X_i}\ind{\eta_i>w}
  \le w^{-(p-1)}\tsum\nolim_i\E\norm{X_i}\eta_i^{p-1}
  =\tfrac{\si w}{\norm{L}}\,\la_p^p.
\end{split}
\end{equation}
Let
\begin{equation*}
 \hat X_i:=\tX_i-\E\tX_i
 \quad\text{and}\quad
 \hat S:=\tsum_i\hat X_i=\tS-\E\tS,
\end{equation*}
so that
\begin{equation}\label{eq:hatX->tX}
 \norm{\hat X_i}_\al
  \le\norm{\tilde X_i}_\al+\norm{\E\tilde X_i}
  \le 2\norm{\tilde X_i}_\al
  \le 2\norm{X_i}_\al
\end{equation}
for all $\al\ge1$, and also
\begin{equation}\label{eq:hatX->X}
 \norm{\hat X_i}_\al\le 2\norm{\tX_i}_\al\le2\bigp{\tfrac{\si w}{\norm{L}}}^{1-p/\al}\,\norm{X_i}_p^{p/\al}
\end{equation}
for all $\al\ge p$. 
Then 
\begin{equation*}
\begin{split}
 \tfrac{2\si}{\Mf}\,\norm{\bar\De}_q
  =\norm{\tS}_{2q}^2
 &\le\bigp{\norm{\hat S}_{2q}+\norm{\E\tS}}^2\le\tfrac54\norm{\hat S}_{2q}^2+5\norm{\E\tS}^2
\\
 &\le\tfrac54\Bigp{
 \AA_\XX(2q)^2%RM 03.20.13 \AA_\XX^2(2q)
 \bigp{\tsum\nolim_i\norm{\hat X_i}_{2q}^{2q}}^{1/q}+
 \BB_\XX(2q)^2%RM 03.20.13 \BB_\XX^2(2q)
 \tsum\nolim_i\norm{\hat X_i}_2^2}+5\bigp{\tfrac{\si w}{\norm L}\la_p^p}^2
\\
 &\le5\Bigp{
 \AA_\XX(2q)^2%RM 03.20.13 \AA_\XX^2(2q)
 \bigp{\tfrac{\si w}{\norm{L}}}^{2-p/q}s_p^{p/q}+
 \BB_\XX(2q)^2%RM 03.20.13 \BB_\XX^2(2q)
 s_2^2+\bigp{\tfrac{\si w}{\norm L}\la_p^p}^2}
\\
 &=5\bigp{\tfrac{\si w}{\norm{L}}}^2\bigp{
 \AA_\XX(2q)^2%RM 03.20.13 \AA_\XX^2(2q)
 \la_p^{p-1}+
 \BB_\XX(2q)^2%RM 03.20.13 \BB_\XX^2(2q)
 \la_2^2+\la_p^{2p}}
  =\tfrac{2\si}{\Mf}\,\uu,
\end{split}
\end{equation*}
where the easily verified inequality $(x+y)^2\le\frac54x^2+5y^2$ is used in the first line above, the Rosenthal-type inequality \eqref{eq:rosen} and \eqref{eq:||EtS||} are used in the second line, \eqref{eq:hatX->tX} and \eqref{eq:hatX->X} are used in the third line, and the definitions \eqref{eq:la_p} and \eqref{eq:u} of $\la_\al$ and $\uu$, respectively, are used in the last line. 
This completes the proof of the lemma.
\end{proof}

\begin{proof}[Proof of Lemma~\ref{lem:xi(De-De_i)}]
Suppose first that $p=3$. 
Then, by \eqref{eq:De-bar} and \eqref{eq:De_i}, for each $i=1,\dotsc,n$, 
\begin{equation*}
\begin{split}
 \tfrac{2\si}{\Mf}\,\bigabs{\bar\De-\De_i}
 &=\bigabs{\norm{S}^2-\norm{S-X_i}^2}
  =\bigabs{\norm{S}-\norm{S-X_i}}\bigp{\norm{S}+\norm{S-X_i}}
\\
 &\le\norm{X_i}\bigp{\norm{X_i}+2\norm{S-X_i}}
  =\norm{X_i}^2+2\norm{X_i}\norm{S-X_i}.
\end{split}
\end{equation*}
Also, by \eqref{eq:D}, $\norm{S-X_i}_q\le\norm{S-X_i}_2\le Ds_2$. It follows that 
\begin{equation}\label{eq:barDe-De_i}
 \bignorm{\bar\De-\De_i}_q
  \le\tfrac{\Mf}{2\si}\bigp{\norm{X_i}_{2q}^2+2\norm{X_i}_q\norm{S-X_i}_q}
  \le\tfrac{\Mf}{2\si}\bigp{\norm{X_i}_{2q}^2+2Ds_2\norm{X_i}_q},
\end{equation}
So, 
\begin{equation*}
\begin{split}
 \tsum_{i=1}^n\bignorm{\xi_i}_p\bignorm{\bar\De-\De_i}_q
 &\le\tfrac{\Mf}{2\si}\tsum_{i=1}^n\bignorm{\xi_i}_p\bigp{\norm{X_i}_{2q}^2+2Ds_2\norm{X_i}_q}
  \le\tfrac{\Mf}{2\si}\,\si_p\,\bigp{s_{2q}^2+2Ds_2s_q}
  =\si_p\vv, 
\end{split}
\end{equation*}
where H\"older's inequality is used for the last inequality, and the definitions \eqref{eq:la_p} and \eqref{eq:v} are used for the equality. 
This proves the lemma when $p=3$.

Suppose now that $p\in(2,3)$. Similarly to \eqref{eq:barDe-De_i} and using the truncation in the definition \eqref{eq:tX},
\begin{equation*}
 \norm{\bar\De-\De_i}_q
  \le\tfrac{\Mf}{2\si}\bigp{\norm{\tX_i}_{2q}^2+2\norm{\tX_i}_q\norm{\tS-\tX_i}_q}
  \le\tfrac{\Mf}{2\si}\bigp{\bigp{\tfrac{\si w}{\norm L}}^{2-p/q}\norm{X_i}_{p}^{p/q}+2\norm{X_i}_q\norm{\tS-\tX_i}_2};
\end{equation*}
also using \eqref{eq:D} and \eqref{eq:||EtS||}, and reasoning as in \eqref{eq:hatX->tX}, one has
\begin{align*}
 \bignorm{\tS-\tX_i}_2
  \le\bignorm{\hat S-\hat X_i}_2+\bignorm{\E\tS-\E\tX_i}
  \le 2Ds_2+\tfrac{\si w}{\norm{L}}\,\la_p^p
  =\tfrac{\si w}{\norm{L}}\bigp{2D\la_2+\la_p^p}.
\end{align*}
So, 
\begin{equation*}
\begin{split}
 \tsum_{i=1}^n\E\norm{\xi_i}_p\norm{\bar\De-\De_i}_q
 &\le\tfrac{\Mf}{2\si}\tsum_{i=1}^n\norm{\xi_i}_p\Bigp{
  \bigp{\tfrac{\si w}{\norm L}}^{2-p/q}\norm{X_i}_p^{p/q}
   +2\tfrac{\si w}{\norm L}\norm{X_i}_q\bigp{2D\la_2+\la_p^p}}
\\
 &\le\tfrac{\Mf\si w^2}{2\norm{L}^2}\,\si_p\Bigp{
   \la_p^{p/q}+2\la_q\bigp{2D\la_2+\la_p^p}}
  =\si_p\vv.
\end{split}
\end{equation*}
Thus, the lemma is proved for $p\in(2,3)$ as well.
\end{proof}

\begin{proof}[Proof of Lemma~\ref{lem:tga_z}]
By \eqref{eq:x3,y3,vp3}, for each $i=1,\dots, n$ 
\begin{equation*}
 \bigprob{W-\xi_i\ge\pi_3z}\le\Bigprob{\max_{j\ne i}\xi_j>y_3}+\Bigprob{\tsum_{j\ne i}\xi_j\bigind{\xi_j\le y_3}\ge x_3}
  \le G_\xi(y_3)+\PUthr,
\end{equation*}
with the last inequality following from \eqref{eq:PU}, \eqref{eq:PUcond}, and the definition of $\PUthr$ in \eqref{eq:PU2,PU3}. 
A similar use of truncation, together with \eqref{eq:PU,vector}, \eqref{eq:x2,y2,vp2}, and \eqref{eq:PU2,PU3},  
yields
\begin{equation*}
 \biggprob{\norm{S}>\Bigp{\frac{2\pi_1\si z}{\Mf}}^{1/2}}
  =\bigprob{\norm{S}>x_2}
  \le G_X(y_2)+\bigprob{\norm{S_{y_2}}>x_2}
  \le G_X(y_2)+\PUtwo.
\end{equation*}
Then \eqref{eq:ga<=tga} follows from the definitions \eqref{eq:ga_z} and \eqref{eq:tga_z} of $\ga_z$ and $\tilde\ga_z$. 

By \eqref{eq:incr} and the definition of $\BH$ right after \eqref{eq:PUcond},
\begin{equation}\label{eq:PU<=BH}
\begin{split}
 \PU(x,y,B,m,\vp) &\le\BH(x,y,B,m)=\exp\Bigb{\tfrac{(x-m)_+}{y}\Bigp{1-\bigp{1+\tfrac{B^2}{(x-m)_+y}}\ln\bigp{1+\tfrac{(x-m)_+y}{B^2}}}}\\
  &\le\Bigp{\frac{eB^2}{(x-m)_+y}}^{(x-m)_+/y}\wedge1,
\end{split}
\end{equation}
where the equality is implied by \cite[(2.9)]{hoeff63}. 
Now let $\ka_2=2p$ and $\ka_3=p/2$. 
Since $G_X(y_2)\le G_X(y_2(\si z/\om)^{1/2})$ whenever \eqref{eq:z} is satisfied, \eqref{eq:ga_z2} follows from \eqref{eq:ga<=tga} and \eqref{eq:PU<=BH} once it is demonstrated that
\begin{equation}\label{eq:PU2<=}
 \PUtwo\le\Bigp{\La_1\,\frac{s_2^2}{\si z}}^p,\quad\text{where }\La_1:=\frac{2pe\Mf D^2}{\pi_1}.
\end{equation}
Assume now that $\La_1s_2^2\le\si z$, since otherwise \eqref{eq:PU2<=} trivially holds. 
Then
\begin{align*}
 \E\norm{S_{y_2}}
 &\le\E\norm{S}+\E\norm{S-S_{y_2}}
  \le\norm{S}_2+\E\bignorm{\tsum\nolim_iX_i\ind{\norm{X_i}>y_2}}
  \le Ds_2+\frac{s_2^2}{y_2}\\
 &=\frac{x_2}{4}\biggp{\Bigp{\frac{16D^2s_2^2}{2\pi_1\si z/\Mf}}^{1/2}+\frac{8ps_2^2}{2\pi_1\si z/\Mf}}
  <\frac{x_2}4\biggp{\Bigp{\La_1\,\frac{s_2^2}{\si z}}^{1/2}+\La_1\,\frac{s_2^2}{\si z}}\le\frac{x_2}2,
\end{align*}
where \eqref{eq:D} is used in the first line above, the definitions in \eqref{eq:x2,y2,vp2} are used for the equality, and the inequalities $8<2pe$ and $4<2eD^2$ (which follow since $p\ge2$ and $D\ge1$) are used for the penultimate inequality. 
Thus, $\PUtwo\le\PU(x_2,y_2,s_2^2,x_2/2,\vp_2)$ follows -- cf.\ \eqref{eq:PUcond}; \eqref{eq:PU2<=} is then seen to hold after an application of \eqref{eq:PU<=BH}. 
\end{proof}

\subsection{Proofs of results from Section~\ref{sec:apps}}
\label{sec:app.proofs}

\begin{proof}[Proof of Theorem~\ref{thm:iid,p=3,unif}]
Note that the conditions of Theorem~\ref{thm:f(S).ub} hold when we set 
\begin{equation*}
 X_i=V_i/n 
\end{equation*}
and take any real $w>0$. 
Then, recalling also that $p=3$
and \eqref{eq:eta.ub} and \eqref{eq:G}, one has $G_\eta(w)=0$. %RM 02-21-13
% \begin{equation*}
%  q=\tfrac32,\quad 
%  G_\eta(w)=0 \text{ (cf.\ \eqref{eq:eta.ub} and \eqref{eq:G})},\quad
%  \tq=3,\quad\text{and}\quad 
%  c_4=\tfrac1{4(1-\cc)}\text{ (cf.\ \eqref{eq:c4})}.
% \end{equation*}
By \eqref{eq:iid,orders}, and in accordance with the notation \eqref{eq:tsi,ga3,v_al}, 
\begin{equation}\label{eq:si,iid}
 \si=\frac{\tsi}{\sqrt{n}},\quad
 \si_p=\frac{\vsi_3}{n^{1/6}},\quad
 s_\al=\frac{v_\al}{n^{1-1/\al}},\quad\text{and}\quad
 \la_\al=\frac{\norm{L}v_\al}{\tsi n^{1/2-1/\al}}
\end{equation}
for any $\al\ge1$. 
Further, use 
%RM 02-21-13 \eqref{eq:A1,B1} and 
the assumption that $\XX$ is a Hilbert space to let $D=1$
as well as use the constants in \eqref{eq:A1,B1}. %RM 02-21-13 , $\AA_\R(p)=\AA_\XX(3)=1$, and  $\BB_\R(p)=\BB_\XX(3)=
% 2^{1/3}$. %RM 02-21-13 3^{1/3}$.
Then, in view of \eqref{eq:u}, \eqref{eq:v}, and the inequality $v_{3/2}\le v_2$, 
\begin{equation}\label{eq:u,v,p=3}
 \uu=\frac1{\sqrt{n}}\,\frac{\Mf}{2\tsi}\Bigp{\frac{v_3^2}{n^{1/3}}+
 2^{2/3} %RM 02-21-13 3^{2/3}
 v_2^2}\quad\text{and}\quad
 \si_p\vv\le\frac{1}{\sqrt{n}}\,\frac{\Mf}{2\tsi}\,\vsi_3\Bigp{\frac{v_3^2}{n^{1/2}}+2v_2^2}.
\end{equation}
One also has $\P(\norm{S}>\ep)\le\KKu{\ep}$ by Remark~\ref{re:|S|} and \eqref{eq:A1,B1}. 
%RM 02-21-13 %% added following sentence %%
Concerning the choice of $\de$, since $\si_1=\vsi_1\sqrt{n}\le\sqrt{n}$ and $\cc\in[\frac12,1)$, by \eqref{eq:de-imp} we may choose
\begin{equation}\label{eq:de.iid}
 \de=\frac{\vsi_3^3-(2\cc-1)^2}{4(1-\cc)\sqrt{n}}.
\end{equation}
Then \eqref{eq:f(S).ub}, combined with \eqref{eq:tyurin,michel} and the above substitutions and inequalities, yields \eqref{eq:iid,p=3}. 
Using now Young's inequality 
\begin{equation}\label{eq:young's}
 \vsi_3^iv_\al^2\le\frac{1}{3}\,\frac{\vsi_3^{3i}}{\ka_{\al,i}^3}+\frac{2}{3}\,\ka_{\al,i}^{3/2}v_\al^3\text{ for }(\al,i)\in\{2,3\}\times\{0,1\},
\end{equation}
one deduces \eqref{eq:iid,p=3,young's} from \eqref{eq:iid,p=3}. 
\end{proof}

\begin{proof}[Proof of Corollary~\ref{cor:asymp}] 
Let $n\to\infty$. 
%RM 02-21-13 In accordance with \eqref{eq:de.to.si} and \eqref{eq:c4}, let $\de=\frac{\si_3^3}{4(1-\cc)}=\frac{\vsi_3^3}{4(1-\cc)\sqrt{n}}$, so that $\de\to0$. 
Following the lines of the proof of \eqref{eq:iid,p=3}, one can see that the bound there 
equals 
\begin{equation}\label{eq:iid,p=3,bound}
%IP 03.12.13
\frac{0.13925+0.33554\vsi_3^3}{\sqrt{n}}
%	0.13925+0.33554%RM12.21.12 0.4748
%	\frac{\vsi_3^3}{\sqrt{n}}
	+\frac{4\de}{2\cc}+\frac{\mathfrak C}{\sqrt n}
	+\frac{\KKu{\ep}}{\sqrt{n}}, %RM 02-21-13 
\end{equation}
where $\mathfrak C:=(\KKu{20}+\KKu{21}\vsi_3)v_2^2+(\KKu{30}+\KKu{31}\vsi_3)v_3^2
%RM 02-21-13 +\KKu{\ep}
$ is an upper bound on 
$\frac1{2\cc}\big(\E\bigabs{W\bar\De} + \tsum_i\E\bigabs{\xi_i(\bar\De-\De_i)}\big)$ -- cf.\ \eqref{eq:ub}. 
Restricting $\cc$ to be in $[\frac12,1)$ and then letting $\de$ be as in \eqref{eq:de.iid}, so that $\de\to0$, %RM 02-12-13 Hence, 
by Remark~\ref{re:chen.improvement} the term $4\de$ in the bound \eqref{eq:iid,p=3,bound} may be replaced by  
\begin{equation*}
	2\de+\frac{\de^2}{\cc}+2\de\sqrt{\frac\de{\cc}+\frac{\de^2}{4\cc^2}+\frac{\mathfrak C}{\sqrt n}}\sim2\de. 
\end{equation*}  
So, the %RM 03.22.13 
terms $\KKu{1}=0.33554 %RM12.21.12 0.4748
+\frac1{2\cc(1-\cc)}$ %RM 03.22.13 
and $\KKu{0}=0.13925-\frac{(2\cc-1)^2}{2\cc(1-\cc)}$ in \eqref{eq:iid,p=3} can be replaced by ones asymptotic to $0.33554 %RM12.21.12 0.4748
+\frac1{4\cc(1-\cc)}$
and %RM 03.22.13 similarly the term $\KKu{0}$ may be replaced by 
$0.13925-\frac{(2\cc-1)^2}{4\cc(1-\cc)}$%RM 03.22.13 . %RM 02-21-13
, respectively. 

Let now $\ep=\ep_n=n^{-1/8}$; the assumed continuity of $f''$ implies $\Mf\downarrow\norm{f''(0)}$, and from \eqref{eq:Kuep} we see that $\KKu{\ep}\downarrow0$. 
Moreover, then $\bigp{\KKu{20},\KKu{21},\KKu{30},\KKu{31}}\to
\frac{\norm{f''(0)}}{2\cc\tsi}((2/\pi)^{1/6},1,0,0)$. %RM 02-21-13 \frac{\Mf}{4\cc\tsi}(3,2,0,0)$. 

Thus, 
%IP 03.12.13  
\begin{equation*}
 \limsup_{n\to\infty}\,\sup_{z\in\R}
 %RM 03.20.13 |z|^3
 \sqrt{n}\bigabs{\bigprob{\tfrac{f(\bar V)}{\tsi/\sqrt{n}}\le z}-\Phi(z)} \le0.13925+0.33554\vsi_3^3+\frac{\vsi_3^3-(2\cc-1)^2}{4\cc(1-\cc)}+\frac{\norm{f''(0)}}{2\cc\tsi}\bigp{(\tfrac2\pi)^{1/6}+\vsi_3}v_2^2.
%RM 02-21-13
%  \limsup_{n\to\infty}\,\sup_{z\in\R}\sqrt{n}\bigabs{\bigprob{\tfrac{f(\bar V)}{\tsi/\sqrt{n}}\le z}-\Phi(z)}
%   \le0.13925+\bigp{0.33554+%RM12.21.12\bigp{0.4748+
%   \tfrac{1}{4\cc(1-\cc)}}\vsi_3^3+\tfrac{\norm{f''(0)}}{4\cc\tsi}\bigp{3+2\vsi_3}v_2^2.
\end{equation*}
%\begin{equation*}
% \limsup_{n\to\infty}\,\sup_{z\in\R}\sqrt{n}\bigabs{\bigprob{\tfrac{f(\bar V)}{\tsi/\sqrt{n}}\le z}-\Phi(z)}
%  \le0.13925+0.33554\vsi_3^3+\frac{\vsi_3^3-(2\cc-1)^2}{4\cc(1-\cc)}+\frac{\norm{f''(0)}}{2\cc\tsi}\bigp{1+\vsi_3}v_2^2.
%%RM 02-21-13
%%  \limsup_{n\to\infty}\,\sup_{z\in\R}\sqrt{n}\bigabs{\bigprob{\tfrac{f(\bar V)}{\tsi/\sqrt{n}}\le z}-\Phi(z)}
%%   \le0.13925+\bigp{0.33554+%RM12.21.12\bigp{0.4748+
%%   \tfrac{1}{4\cc(1-\cc)}}\vsi_3^3+\tfrac{\norm{f''(0)}}{4\cc\tsi}\bigp{3+2\vsi_3}v_2^2.
%\end{equation*}
Since 
\begin{equation*}
 \min_{\cc\in[1/2,1)}\Bigp{\frac{\vsi_3^3-(2\cc-1)^2}{4\cc(1-\cc)}+\frac{y_*}{2\cc}}=\frac{1+\vsi_3^3+y_*+\sqrt{(\vsi_3^3-1)(\vsi_3^3-1+2y_*)}}2,
\end{equation*}
the inequality \eqref{eq:asymp.unif} follows. 
%RM 02-21-13 $\min_{\cc\in(0,1)}\big(\frac{a}{\cc(1-\cc)}+\frac{b}{\cc}\big)=2a+b+2\sqrt{a(a+b)}$ whenever $a>0$ and $b>0$, the first inequality in \eqref{eq:asymp.unif} follows. %RM12.21.12; the second inequality there follows by letting $\cc$ take the suboptimal value of $\frac12$.

To prove \eqref{eq:asymp.nonunif}, fix any real $\tth>0$ and let $z_0=g(n)$, $\om=\tsi/g(n)$, $K_1=\sqrt{n}/\vsi_3^3$, $K_2=\tsi^3z_0^3\sqrt{n}/v_2^4$, and $K_3=\tsi^3z_0^3\sqrt{n}/v_3^3$, so that \eqref{eq:z.n.conds} holds for all $z\in[g(n),\sqrt{n}/g(n)]$. 
Then, for $z\ge z_0$ and large enough $n$ we have $z^3e^{-z/\tth}\le z_0^3e^{-z_0/\tth}\to0$. 
Concerning the pre-constants in Theorem~\ref{thm:iid,p=3,nonunif} in Appendix~\ref{app:nonunif}, one can clearly choose values for the corresponding parameters so that 
(i) $\Cnecon,\dotsc,\Cnevthr$ be absolutely bounded; 
(ii) $\Cnvtwoa$, $\Cnvtwob$, $\Cnvthra$, and $\Cnvthrb$ all vanish in the limit (since $\om\downarrow0$); and 
(iii) $\Cnvsi\to30.2211+\pi_2^{-3}$. 
Moreover, one can replace the factor $\vsi_3^3$ in the second inequality (and, if so desired, in the other two inequalities) in \eqref{eq:G<} by the asymptotically much smaller expression $\E\big(\frac{L(V)}\tsi\big)^3\I\{L(V)>\pi_2\tsi z\sqrt n\}=o(\vsi_3^3)$. 
Then the limit of the corresponding improved expression for $\Cnvsi$ becomes just $30.2211$, instead of $30.2211+\pi_2^{-3}$. 
Now \eqref{eq:asymp.nonunif} follows by Theorem~\ref{thm:iid,p=3,nonunif}.
\end{proof}

\begin{proof}[Proof of Corollary~\ref{cor:centralT2}] % RM 03-05-13 %% supporting calculations in notebook `03-08-13.unifbound.4thmoments.nb' %%
Take any natural number $N_0\ge1$ and any real numbers $\ep\in(0,1)$, 
$\cc\in[\frac12,1)$, %RM 02-21-13 $\cc\in(0,1)$, 
$\ka_1>0$, and $\ka_2>0$, and let 
%RM 03-05-13 $X_i:=Y_i$, 
$\xi_i:=Y_i/\sqrt n$, $W:=\sum_i\xi_i$, and $\bar{Y^2}:=
%RM 03-05-13 %IP 03.23.13 
\frac1n\sum_iY_i^2=
\sum_i\xi_i^2$. 
Further let 
\begin{equation*} %RM 03-05-13 $ 
\tilde T:=T_1\ind{\abs{\bar{Y^2}-1}\le\ep}+W\ind{\abs{\bar{Y^2}-1}>\ep}, 
\end{equation*} %RM 03-05-13 $, 
where $T_1=W/\sqrt{\bar{Y^2}}$ is the self-normalized sum as defined in \eqref{eq:T1}. 
Then
\begin{equation*}
 \begin{split}
  \abs{\tilde T-W}
   &=\Bigabs{W\Bigp{\tfrac1{\sqrt{\bar{Y^2}}}-1}}\bigind{\abs{\bar{Y^2}-1}\le\ep}
    =\bigabs{W(\bar{Y^2}-1)}\,\frac{\ind{\abs{\bar{Y^2}-1}\le\ep}}
    {\bar{Y^2}+\sqrt{\bar{Y^2}}} %RM 03-05-13 {\sqrt{\bar{Y^2}}\bigp{1+\sqrt{\bar{Y^2}}}}
    \le\check\Mf\bigabs{W(\bar{Y^2}-1)},
 \end{split}
\end{equation*}
where
\begin{equation*}
 \check\Mf:=\frac1{1-\ep+\sqrt{1-\ep}}. %RM 03-05-13 \frac1{\sqrt{1-\ep}\,(1+\sqrt{1-\ep})}.
\end{equation*}
Accordingly, let
\begin{equation*}
 \bar\De:=\De:=
 \check\Mf W\bigp{\bar{Y^2}-1}\quad\text{and}\quad\De_i:=
 \check\Mf W_{(i)}\bigp{\bar{Y^2}_{(i)}-1},
\end{equation*}
where $W_{(i)}=W-\xi_i$ and $\bar{Y^2}_{(i)}=\bar{Y^2}-\xi_i^2$. 
Then the conditions of Theorem~\ref{thm:ub} hold with $\tilde T$ in place of $T$ if we let $\eta_i=0$ for $i=1,\dotsc,n$ (and then allow $w$ to take any positive value). 

Recall that $\norm{Y}_2=1$ is being assumed, whence $\norm{Y^2-1}_2^2=\norm{Y}_4^4-1$, 
$\norm{\bar{Y^2}-1}_2=\sqrt{\norm{Y}_4^4-1}/\sqrt{n}$ and $\norm{W}_4^4=n\E\xi_1^4+3n(n-1)(\E\xi_1^2)^2=\frac1n(3(n-1)+\norm{Y}_4^4)$. %RM 03-05-13; also, Young's inequality implies
% \begin{equation*}
%  \norm{Y}_4^2\norm{Y^2-1}_2\le\frac12\Bigp{\frac{\norm{Y}_4^4}{\ka_1^2}+\ka_1^2\norm{Y^2-1}_2^2}=\frac{\ka_1^2+1/\ka_1^2}{2}\,\norm{Y}_4^4-\frac{\ka_1^2}{2}
% \end{equation*}
% and
% \begin{equation*}
%  \norm{Y^2-1}_2\le\frac12\Bigp{\frac1{\ka_2^2}+\ka_2^2\norm{Y^2-1}_2^2}=\frac{\ka_2^2}{2}\,\norm{Y}_4^4-\frac{\ka_2^2-1/\ka_2^2}{2}.
% \end{equation*}
% From $\norm{\bar{Y^2}-1}_2=\norm{Y^2-1}_2/\sqrt{n}$ and $\norm{W}_4^4=n\E\xi_1^4+3n(n-1)(\E\xi_1^2)^2\le3+\norm{Y}_4^4/n$, it follows that
Then we have %RM 03-05-13
\begin{equation*}
 \E\bigabs{W
 %RM 03-05-13 \bar
 \De}
  =\check\Mf\E W^2\bigabs{\bar{Y^2}-1}
  \le\check\Mf\norm{W}_4^2\bignorm{\bar{Y^2}-1}_2
  \le\frac{\check\Mf}{\sqrt n}\,\sqrt{\norm{Y}_4^4-1}\sqrt{3+\norm{Y}_4^4/n}. %RM 03-05-13  \le\tfrac{\check\Mf}{\sqrt n}\bignorm{Y^2-1}_2\bigp{\sqrt{3}+\tfrac{1}{\sqrt{n}}\norm{Y}_4^2}.
\end{equation*}
Also, $%RM 03-05-13 \bar
\De-\De_1=\check\Mf(\xi_1^2W_{(1)}+\xi_1(\bar{Y^2}-1))$, whence  
\begin{equation*}
 \begin{split}
  \E\bigabs{\xi_1(\De-\De_1)}
  &\le 
  \check\Mf\bigp{\E\abs{\xi_1}^3\E\abs{W_{(1)}}+\E\xi_1^2\abs{\bar{Y^2}-1}}
   \le 
   \check\Mf\bigp{\norm{\xi_1}_3^3\norm{W_{(1)}}_2+\norm{\xi_1}_4^2\bignorm{\bar{Y^2}-1}_2}\\
  &\le\frac{\check\Mf}{n^{3/2}}\Bigp{\norm{Y}_3^3+\norm{Y}_4^2\sqrt{\norm{Y}_4^4-1}\,}. %RM 03-05-13 &\le\tfrac{\check\Mf}{n\sqrt{n}}\bigp{\norm{Y}_3^3+\norm{Y}_4^2\norm{Y^2-1}_2}.
 \end{split}
\end{equation*}

In the case where $n\ge N_0$, combine \eqref{eq:tyurin,michel} and \eqref{eq:ub} (use also 
\eqref{eq:de.iid}) %RM 02-21-13 \eqref{eq:de.to.si} with $p=3$) 
to obtain
\begin{align} %RM 03-05-13
\notag
&\bigabs{\P(T_1\le z)-\Phi(z)}
 \le\bigprob{\bigabs{\bar{Y^2}-1}>\ep}+\bigabs{\P(\tilde T\le z)-\P(W\le z)}+\bigabs{\P(W\le z)-\Phi(z)}\\
\notag
&\quad
 \le\frac{\E\abs{\bar{Y}^2-1}^2}{\ep^2}+\frac{4\de+\E\abs{W\De}+n\E\abs{\xi_1(\De-\De_1)}}{2\cc}+0.33554n\bigp{\norm{\xi_1}_3^3+0.415\norm{\xi_1}_2^3}\\
\label{eq:4thmoms}
\begin{split}
 &\quad\le\frac1{\sqrt{n}}\biggp{0.33554\bigp{\norm{Y}_3^3+0.415}+\frac{\norm{Y}_4^4-1}{\ep^2\sqrt{n}}+\frac{\norm{Y}_3^3-(2\cc-1)^2}{2\cc(1-\cc)}\\
 &\quad\qquad+\frac{\check\Mf}{2\cc}\biggp{\norm{Y}_3^3+\sqrt{\norm{Y}_4^4-1}\Bigp{\norm{Y}_4^2+\sqrt{3+\norm{Y}_4^4/n}\,}}}
\end{split}
\\
\notag
&\quad\le\frac{A_3\norm{Y}_3^3+A_4\norm{Y}_4^4-A_0}{\sqrt{n}},
\end{align}
% \begin{equation*}
%  \begin{split}
%   &\bigabs{\P(T_1\le z)-\Phi(z)}
%    \le\bigprob{\bigabs{\bar{Y^2}-1}>\ep}+\bigabs{\P(\tilde T\le z)-\P(W\le z)}+\bigabs{\P(W\le z)-\Phi(z)}\\
%   &\quad
%    \le\tfrac{\norm{Y}_4^4-1}{\ep^2n}
%    +\tfrac1{2\cc}\Bigp{\tfrac{\norm{Y}_3^3
%    -(2\cc-1)^2%RM 02-21-13 
%    }{(1-\cc)\sqrt n}+\E\abs{W\bar\De}+n\E\abs{\xi_1(\bar\De-\De_1)}}
%    +\tfrac{0.13925+0.33554%RM12.21.12 0.4748
%    \norm{Y}_3^3}{\sqrt{n}}\\
%   &\quad\tfrac1{\sqrt{n}}\Bigp{0.33554\bigp{\norm{Y}_3^3+0.415}+\tfrac{\norm{Y}_4^4-1}{\ep^2\sqrt{n}}+\tfrac{\norm{Y}_3^3-(2\cc-1)^2}{2\cc(1-\cc)}+\tfrac{\check\Mf}{2\cc}\bigp{\sqrt{\norm{Y}_4^4-1}\bigp{\sqrt{3+\norm{Y}_4^4/n}+\norm{Y}_4^2}+\norm{Y}_3^3}}\\
%   &\quad
%    \le\tfrac1{\sqrt{n}}\bigp{A_3\norm{Y}_3^3+A_4\norm{Y}_4^4-A_0},
% \end{split}
% \end{equation*}
where
\begin{equation*}
 \begin{split}
 A_3&:=0.33554%RM12.21.12 0.4748
 +\frac1{2\cc(1-\cc)}+\frac{\check\Mf}{2\cc},\\ 
 A_4&:=\frac{\check\Mf}{4\cc}
  \Bigp{\ka_1^2+\frac1{\ka_1^{2}}+\ka_2^2+\frac1{\ka_2^2N_0}} %RM 03-05-13 
%   \Big(\Big(1+\frac1{\sqrt{N_0}}\Big)
%   \Big(\ka_1^2+\frac1{\ka_1^2}\Big)
%   +\sqrt3\ka_2^2\Big)
  +\frac1{\ep^2\sqrt{N_0}},\\
A_0&:=\frac{\check\Mf}{4\cc}
 \Bigp{\ka_1^2+\ka_2^2-\frac{3}{\ka_2^2}} %RM 03-05-13 
% \Bigp{\Big(1+\frac1{\sqrt{N_0}}\Big)\ka_1^2+\sqrt3\Big(\ka_2^2-\frac1{\ka_2^2}\Big)}+\frac1{\ep^2\sqrt{N_0}}
 +\frac{(2\cc-1)^2}{2\cc(1-\cc)} %RM 02-21-13 
 -0.13925; %RM12.21.12 
 \end{split}
\end{equation*}
Young's inequality, specifically $\sqrt{a\,b}\le\frac12(\ka^2a+b/\ka^2)$ for any positive $\ka$ and any nonnegative $a$ and $b$, is used on the last two terms in \eqref{eq:4thmoms}. %RM 03-05-13 
Then the inequality \eqref{eq:centralT.4th} holds for any of the triples in \eqref{eq:Ajlist.4th}, in the case where $n\ge N_0$, when the parameter values in the table below (to be interpreted as rational numbers) are substituted into the expressions for $A_3$, $A_4$, and $A_0$ above: 
\begin{center}
 \begin{tabular}{ccc|ccccc}
 $A_3$&$A_4$&$A_0$&$N_0$&$\ep$   &$\cc$  &$\ka_1$&$\ka_2$\\\hline
 3.00 &4.66 &4.33 &12   &0.335595&0.613  &2.1149 &1.656\\
 3.17 &2.04 &1.07 &18   &0.4944  &0.64847&1.12925&1.114\\
 3.48 &1.27 &-1.43&39   &0.5878  &0.7211 &0.6066 &1 %RM 03-05-13 
%   3.00 &5.02 & 4.69  &12 &0.332&0.618&1.696&1.639\\
%   3.18 &2.09 & 1.01  &19 &0.490&0.657&1.094&0.806\\
%   3.60 &1.30 &$-1.44$&41 &0.559&0.749&1.019&0.424 %RM 02-21-13
%   3.01 &5.16 &  4.75 &12   &0.357&0.606&1.821  &1.778\\
%   3.20 &2.20 &  1.14 &19   &0.494&0.664&1.147  &0.930\\
%   3.65 &1.31 &$-1.45$&41   &0.565&0.755&1.020  &0.479
%RM 12.21.12  
%   3.14 &5.16 &  4.89 &12   &0.357&0.606&1.821  &1.778\\
%   3.34 &2.20 &  1.28 &19   &0.494&0.664&1.147  &0.930\\
%   3.79 &1.31 &$-1.32$&41   &0.565&0.755&1.020  &0.479
 \end{tabular}
\end{center}
In the case where $n<N_0$ (or hence $n\le N_0-1$), it suffices to use the trivial bound $|\P(T_1\le z)-\Phi(z)|\le\sqrt{N_0-1}/\sqrt{n}$ and then note that $\sqrt{N_0-1}\le A_3+A_4-A_0\le A_3\norm{Y}_3^3+A_4\norm{Y}_4^4-A_0$ for any of the three triples $(A_3,A_4,A_0)$ in the table above. 
\end{proof}

\appendix

\clearpage

\section{An explicit nonuniform bound}
\label{app:nonunif}

In this appendix, we state and prove Theorem~\ref{thm:iid,p=3,nonunif}, which presents an explicit nonuniform BE-type bound for the normalized statistic $\sqrt{n}f(\bar V)/\tsi$ when the summands $V_i$ are i.i.d. 
The following lemma quotes expressions found in \cite{pinutev86,pinBH} for the exponential bound $\PU$ on the tail probability defined in \eqref{eq:PU}. 
These expressions will be needed in applications of Theorem~\ref{thm:iid,p=3,nonunif}, wherein $\PU$ enters the expressions for several pre-constants. 

\begin{lem}\label{lem:PU}
For any real $x\in\R$, $y>0$, $B>0$, $m$, and $\vp\in(0,1]$, let
\begin{equation*}
 u:=\frac{(x-m)_+y}{B^2}\quad\text{and}\quad\ka:=\frac{(x-m)_+}{y}.
\end{equation*}
Then 
\begin{align}\label{eq:PUchar1}
 \PU(x,y,B,m,\vp)=\PU(u,\ka,\vp):=
 \begin{cases}
  1&\text{if }u=0,\\
  \PUalt(u,\ka,\vp)&\text{if $u>0$ and $\vp<1$},\\ 
  \BHalt(u,\ka)&\text{if $u>0$ and $\vp=1$},
 \end{cases}
\end{align}
where
\begin{equation*}
 \BHalt(u,\ka):=\exp\Bigb{\ka\Bigp{1-\Big(1+\frac1{u}\Big)\ln(1+u)}},
\end{equation*}
\begin{equation}\label{eq:PU1} 
 \PUalt(u,\ka,\vp):=\exp\biggb{\frac{\ka}{2(1-\vp)u}\biggp{(1-\vp)^2\Big[1
  +\Lam\Big(\frac{\vp}{1-\vp}\,\exp\frac{\vp+u}{1-\vp}\Big)\Big]^2-(\vp+u)^2-(1-\vp^2)}},
\end{equation}
and $\Lam$ is Lambert's product-log function with domain restricted to the positive real numbers (so that for positive $w$ and $z$ one has $\Lam(z)=w$ if and only if $z=we^w$); 
in \eqref{eq:PUchar1}, we allowed ourselves the slight abuse of notation, by using the same symbol, $\PU$, to denote two different functions, represented by two expressions, which take the same values but expressed using two different sequences of arguments: $(x,y,B,m,\vp)$ and $(u,\ka,\vp)$. 

One also has the alternative identity
\begin{equation}\label{eq:PUchar2}
 \PU(u,\ka,\vp)=\inf_{0<\al<1}\exp\bigb{L_1\vee L_2},
\end{equation}
where
\begin{equation}\label{eq:L1} 
 L_1:=L_1(\al,u,\ka,\vp):=\ka\biggp{1-\al-\al\,\frac{\vp}{1-\vp}-\frac{\al(2-\al)}{2(1-\vp)}\,u},\quad\text{with}\quad L_1(\al,u,\ka,1):=-\infty,
\end{equation}
and
\begin{equation}\label{eq:L2} 
 L_2:=L_2(\al,u,\ka,\vp):=\ka\biggp{1-\al-\Bigp{1-\frac{\al}{2}+\frac{\vp}{u}}\ln\Bigp{1+(1-\al)\frac{u}{\vp}}},\quad\text{with}\quad L_2(\al,0,\ka,\vp):=0.
\end{equation}
\end{lem}

Indeed, \eqref{eq:PUchar1} is essentially \cite[Proposition~3.1]{pinBH}, with the ``boundary'' case $\vp=1$ resulting in the Bennett--Hoeffding bound $\BHalt(u,\ka)$. 
Next, \eqref{eq:PUchar2} (for $\vp<1$) is established in \cite[Corollary~1]{pinutev86} and, again, immediately follows for $\vp=1$ using $\BHalt(u,\ka)$. 

\begin{thm}\label{thm:iid,p=3,nonunif}
Assume that the conditions of Theorem~\ref{thm:iid,p=3,unif} hold, and let
\begin{equation}\label{eq:nonunif.params}
 \cc,\ \th,\ w,\ \de_0,\ \pi_1,\ \pi_2,\ \pi_3,\ z_0,\ \om,\ \ka_{2,0},\ \ka_{3,0},\ \ka_{2,1},\ \ka_{3,1},\ \ka_2,\ \ka_3,\ \al,\ \vp_\ast,\ K_1,\ K_2,\ \text{and}\ K_3
\end{equation}
all be positive real numbers satisfying the constraints 
\begin{equation}\label{eq:z,om,Ki}
 \cc<1,\ \de_0\le w,\ \pi_1+\pi_2+\pi_3=1,\ \om\le\frac{\Mf\ep^2}{2\pi_1},\ \ka_3\ge\tfrac32,\ \al<1,\ \vp_\ast<1,\ \hat\ka_2\ge2,\ \text{and}\ \hat\ga<1,
\end{equation}
where
\begin{equation}\label{eq:gahat<1}
 \hat\ga:=\Bigp{\frac{\Mf^2\om}{4\pi_1^2K_2}}^{1/4}+\frac{\ka_2^2}{K_3}\Bigp{\frac{\Mf\om}{2\pi_1}}^{3/2}
\end{equation}
and
\begin{equation}\label{eq:ka1hat}
 \hat\ka_2:=(1-\hat\ga)\ka_2.
\end{equation}
Also introduce 
\begin{equation}\label{eq:t2,t3,u0}
 t_2:=\frac{\pi_1\al(2-\al)(1-\hat\ga)^2}{\Mf(1-\vp_{\ast})}\Bigp{\frac{K_2}{\om}}^{1/2},\quad
 t_3:=\frac{\ka_2^{2}}{(1-\hat\ga)K_3}\Bigp{\frac{\Mf\om}{2\pi_1}}^{3/2},\quad
 u_0:=\frac{2\pi_1(1-\hat\ga)}{\Mf\ka_2}\,\Bigp{\frac{K_2}{\om}}^{1/2},
\end{equation}
\begin{gather}\label{eq:ta1}
 \tvp_1:=\tfrac1{K_1w},\quad\text{and}\quad\ta_1:=1/L_{W;\,3w/\th,\tvp_1},
\end{gather}
where $L_{W;\,c,B}$ is as in \eqref{eq:L_W}; 
further let $\tc_1$, $\tc_2$, and $\tc_3$ be obtained from $c_1$, $c_2$, and $c_3$ in \eqref{eq:c1}--\eqref{eq:c3} by replacing there $a_1$, $\vp_1$, and $\si_p$ by $\ta_1$, $1\wedge\tvp_1$, and $K_1^{-1/3}$, respectively.  
%RM 03.20.13 and also letting $(\AA_{\R,\nc}(p),\BB_{\R,\nc}(p))=(1.316^{1/3},2^{1/3})$ (as in \eqref{eq:A1,B1}). %RM 02-21-13
Recall also the definition of $\PU$ in \eqref{eq:PUchar1}. 
Then for all $z\in\R$ and $n\in\N$ such that
\begin{equation}\label{eq:z,iid,p=3}
 z_0\le z\le\frac{\om}{\tsi}\,\sqrt{n},
\end{equation}
\begin{equation}\label{eq:K2,K3,ga}
 \frac{K_1\vsi_3^3}{\sqrt{n}}\le1,\quad
 \frac{K_2v_2^4}{\tsi^3z^3\sqrt{n}}\le1,\quad\text{and}\quad
 \frac{K_3v_3^3}{\tsi^3z^3\sqrt{n}}\le1
\end{equation}
one has
\begin{multline}\label{eq:iid,p=3,nu}
 \Bigabs{\Bigprob{\frac{f(\bar V)}{\tsi /\sqrt{n}}\le z}-\Phi(z)}\\
 \le 
  \frac{\Cnvsi\vsi_3^3+\bigp{(\Cnvtwoa\vee \Cnvtwob)v_2^4}\vee\bigp{\Cnvthra v_3^3}+\Cnvthrb v_3^3}{z^3\sqrt n}
  +\frac{\Cnecon+\Cnevsi\vsi_3^3+\Cnevtwo v_2^3+\Cnevthr v_3^3}{e^{(1-\pi_1)z/\th}\sqrt n},
\end{multline}
where
\begin{equation}\label{eq:C7}
 \Cnvsi:=30.2211+\frac1{\pi_2^3}+\frac{\ka_3^{3/2}}{(w\pi_3)^3}\biggp{\frac{\ka_3^{3/2}}{K_1}+\sup_{u\ge\pi_3^2z_0^2/\ka_3}u^{3/2}\PU\Bigp{u,\ka_3,\frac{\ka_3}{K_1\pi_3z_0}\wedge1}},
\end{equation}
\begin{equation}\label{eq:C81}
 \Cnvtwoa:=\frac{\om\exp\{\hat\ka_2(1-\al-\frac{\al\vp_{\ast}}{1-\vp_\ast})\}}{\tsi^3}\Bigp{\frac{\Mf(1-\vp_{\ast})}{\pi_1\al(2-\al)(1-\hat\ga)^2}}^{2}\sup_{t\ge t_2}t^{2}e^{-t},
\end{equation}
\begin{equation}\label{eq:C82}
 \Cnvtwob:=\frac{\om}{\tsi^3}\Bigp{\frac{\Mf\ka_2}{2\pi_1(1-\hat\ga)}}^{2}\sup_{u\ge u_0}u^{2}\PU(u,\hat\ka_2,\vp_{\ast}),
\end{equation}
\begin{equation}\label{eq:C83}
 \Cnvthra :=\frac{\ka_2^2e^{\hat\ka_2(1-\al)}}{\tsi^3(1-\hat\ga)}\Bigp{\frac{\Mf\om}{2\pi_1}}^{3/2}\sup_{t\in(0,t_3]}\frac1t \exp\biggb{-\hat\ka_2\Bigp{1-\frac\al2+t}\ln\Bigp{1+\frac{1-\al}t}},
\end{equation}
\begin{equation}\label{eq:C9}
 \Cnvthrb:=\Bigp{\frac{\ka_2}{\tsi}}^3\Bigp{\frac{\Mf\om}{2\pi_1}}^{3/2},
\end{equation}
\begin{equation}\label{eq:C10}
 \Cnecon:=\frac{\Mf\tc_2}{6\tsi}\Bigp{\frac1{\ka_{3,0}^3K_1^{2/3}}+\frac{
 2%RM 02-21-13 3
 ^{2/3}}{\ka_{2,0}^3}}
 -\frac{(2\cc-1)_+^2}{4(1-\cc)}\,\tc_3, %RM 02-21-13
\end{equation}
\begin{equation}\label{eq:C11}
 \Cnevsi:=
  d(\cc)\tc_3 %RM 02-21-13 \frac{\tc_3}{4(1-\cc)}
  +\frac{\Mf\tc_1}{6\tsi}\Bigp{\frac1{\ka_{3,1}^3K_1}+\frac{2}{\ka_{2,1}^3}},
\end{equation}
\begin{equation}\label{eq:C12}
 \Cnevtwo:=\frac{\Mf}{3\tsi}\bigp{2\tc_1\ka_{2,1}^{3/2}+
 2%RM 02-21-13 3
 ^{2/3}\tc_2\ka_{2,0}^{3/2}},
\end{equation}
\begin{equation}\label{eq:C13}
 \Cnevthr:=\frac{\Mf}{3\tsi}\Bigp{\frac{\tc_1\ka_{3,1}^{3/2}}{K_1}+\frac{\tc_2\ka_{3,0}^{3/2}}{K_1^{2/3}}},
 %RM 02-21-13 ;
\end{equation}
\begin{equation}\label{eq:d(cc)} %RM 02-21-13 %% added following definition of d(\cc) %%
 d\colon(0,1)\to\R\text{ is defined by }d(\cc)=\begin{cases}\cc&\text{if $\cc\in(0,\frac12]$,}\\\frac{1}{4(1-\cc)}&\text{if $\cc\in(\frac12,1)$;}\end{cases}
\end{equation}
moreover, each of the expressions in \eqref{eq:C7}--\eqref{eq:C13} is finite. 
\end{thm}

\begin{remark}\label{re:symm.nonunif}
Suppose here that $L(V)$ is symmetric. 
Then the statement of Theorem~\ref{thm:iid,p=3,nonunif} holds when the replacement mentioned in Remark~\ref{re:symm} is made in the expression \eqref{eq:c2} for the pre-constant $c_2$ and, accordingly, in the expression for $\tc_2$ defined right after \eqref{eq:ta1}. 
Also, one can take $\Cnvsi$ in \eqref{eq:iid,p=3,nu} to be defined as
\begin{equation} 
\Cnvsi:=30.2211+\frac1{2\pi_2^3}+\frac{\ka_3^{3/2}}{2(w\pi_3)^3}\biggp{\frac{\ka_3^{3/2}}{2K_1}+\sup_{u\ge\pi_3^2z_0^2/\ka_3}u^{3/2}\PU\Bigp{u,\ka_3,\frac{\ka_3}{K_1\pi_3z_0}\wedge1}},
\end{equation}
because one can then use   
$G_\eta(t)=G_\xi(t)\le\vsi_3^3/(2t^3\sqrt{n})$ in place of  $G_\eta(t)=G_\xi(t)\le\vsi_3^3/(t^3\sqrt{n})$ to improve the bounds in \eqref{eq:G<} (in the proof of Theorem~\ref{thm:iid,p=3,nonunif}). 
\end{remark}

\begin{remark}\label{re:ka3,ka1hat}
That all the pre-constants in Theorem~\ref{thm:iid,p=3,nonunif} are finite is easily verifiable by inspection, except perhaps for the pre-constants $\Cnvsi$, $\Cnvtwob$, and $\Cnvthra$, whose expressions in \eqref{eq:C7}, \eqref{eq:C82}, and \eqref{eq:C83} involve comparatively complicated suprema. 
However, Lemma~\ref{lem:incr-decr} in Appendix~\ref{app:num.proofs} provides the sufficient conditions $\ka_3\ge\frac32$ and $\hat\ka_2\ge2$ in \eqref{eq:z,om,Ki} for these three suprema, and hence for the pre-constants $\Cnvsi$, $\Cnvtwob$, and $\Cnvthra$, to be finite.  
\end{remark}

One can substantially improve the bound on $\tilde\ga_z$ in \eqref{eq:ga_z2}. 
The following lemma is key to that, and its proof will be given after the proof of Theorem~\ref{thm:iid,p=3,nonunif}. 

\begin{lem}\label{lem:exp.iid.p=3}
Assume that the conditions of Theorem~\ref{thm:iid,p=3,nonunif} hold. 
Then, for all $z\in\R$ and $n\in\N$ satisfying the inequalities in \eqref{eq:z,iid,p=3} and \eqref{eq:K2,K3,ga},
\begin{equation}\label{eq:exp3}
 \PUtwo\le\frac{((\Cnvtwoa\vee \Cnvtwob)v_2^4)\vee(\Cnvthra v_3^3)}{z^3\sqrt{n}},
\end{equation}
where $\PUtwo$ is as defined in \eqref{eq:PU2,PU3} and $\Cnvtwoa$, $\Cnvtwob$, and $\Cnvthra $ are as defined in \eqref{eq:C81}, \eqref{eq:C82}, and \eqref{eq:C83}, respectively.
\end{lem}

\begin{proof}[Proof of Theorem~\ref{thm:iid,p=3,nonunif}]
Take any $z\in\R$ and $n\in\N$ such that \eqref{eq:z,iid,p=3} and \eqref{eq:K2,K3,ga} hold. 
The conditions of Theorem~\ref{thm:f(S).nub} are met when we let $p=3$ and $X_i=V_i/n$, so that \eqref{eq:tyurin,michel} and \eqref{eq:f(S).nub} imply
\begin{equation}\label{eq:Michel}
 \Bigabs{\Bigprob{\frac{f(\bar V)}{\tsi/\sqrt{n}}\le z}-\Phi(z)}\le\frac{30.2211\vsi_3^3}{z^3\sqrt{n}}+\tilde\ga_z+\tilde\tau e^{-(1-\pi_1)z/\th}.
\end{equation}

Recall \eqref{eq:eta.nub} to see that $\eta_i=\xi_i=L(V_i)/(\tsi\sqrt{n})$; then, for any $t>0$, $G_\eta(t)=G_\xi(t)\le\vsi_3^3/(t^3\sqrt{n})$ and $G_X(t)\le v_3^3/(t^3n^2)$ (cf.\ the inequalities \eqref{eq:G_X.iid} and \eqref{eq:G_eta.iid}). 
Using these inequalities and also the first inequality of \eqref{eq:K2,K3,ga}, \eqref{eq:z,iid,p=3}, \eqref{eq:incr}, and Lemma~\ref{lem:PU}, one has 
\begin{gather}
 G_\eta(w)\le\frac{1}{w^3}\,\frac{\vsi_3^3}{\sqrt{n}},\quad 
 G_\xi(\pi_2z)\le\frac{1}{\pi_2^3}\,\frac{\vsi_3^3}{z^3\sqrt{n}},\quad
 G_\xi(y_3)\le\frac{\ka_3^3}{\pi_3^3}\,\frac{\vsi_3^3}{z^3\sqrt{n}}\le\frac{\ka_3^3}{K_1\pi_3^3}\,\frac1{z^3},
 \label{eq:G<}\\
 G_X(y_2)\le\frac{v_3^3}{y_2^3n^2}\,\frac{\om^{3/2}n^{3/4}}{(\tsi z)^{3/2}}=\frac{\ka_2^3\om^{3/2}}{\tsi^3(2\pi_1/\Mf)^{3/2}}\,\frac{v_3^3}{z^3\sqrt{n}}, \notag\\
\PUthr=\PU\bigp{x_3y_3,\ka_3,\vp_3}\le\frac{\ka_3^{3/2}}{\pi_3^3z^3}\sup_{u\ge\pi_3^2z_0^2/\ka_3}u^{3/2}\PU\Bigp{u,\ka_3,\frac{\ka_3}{K_1\pi_3z_0}\wedge1}.
\notag
\end{gather}
Then \eqref{eq:ga<=tga} and Lemma~\ref{lem:exp.iid.p=3} yield
\begin{equation}\label{eq:tga1}
 \frac{30.2211\vsi_3^3}{z^3\sqrt{n}}+\tilde\ga_z
 \le\frac1{z^3\sqrt{n}}\Bigp{\Cnvsi\vsi_3^3+\bigp{(\Cnvtwoa\vee \Cnvtwob)v_2^4}\vee\bigp{\Cnvthra v_3^3}+\Cnvthrb v_3^3};
\end{equation}
where $\Cnvsi,\dotsc,\Cnvthrb$ are as in \eqref{eq:C7}--\eqref{eq:C9}. 

Next, in the definitions \eqref{eq:c1}--\eqref{eq:c3} and \eqref{eq:a1}, set $p=3$, 
$\AA_{\R,\nc}(p)=1.316^{1/3}$, and $\BB_{\R,\nc}(p)=2^{1/3}$ %RM 03.20.13 $\AA_\R(p)=1$, and $\BB_\R(p)=3^{1/3}$ 
-- recall here \eqref{eq:A1,B1}. 
Also, by the first inequality of \eqref{eq:K2,K3,ga} and \eqref{eq:vp1}, $\si_p=\vsi_3n^{-1/6}\le K_1^{-1/3}$, $\vp_1\le1\wedge(K_1w)^{-1}$, and $\norm{\xi_i}_2/w=1/(w\sqrt{n})\le\vsi_3^3/(w\sqrt{n})\le1/(K_1w)$. 
Then, referring to \eqref{eq:ta1}, we see that $a_1\le\ta_1$ (as $L_{W;\,c,\si}$ is nonincreasing with respect to $\si$) and $c_j\le\tc_j$ for $j=1,2,3$. 
%RM 02-21-13 %% added following sentence %%
By Remark~\ref{re:de} and \eqref{eq:de.iid}, we see that \eqref{eq:de} is satisfied when
\begin{equation*}
 \de=\frac{d(\cc)\vsi_3^3-(2\cc-1)_+^2/(4(1-\cc))}{\sqrt{n}},
\end{equation*}
where $d$ is as in \eqref{eq:d(cc)}. 
Using the definition \eqref{eq:tau2} of $\tilde\tau$, as well as \eqref{eq:u,v,p=3} and \eqref{eq:young's}, one obtains the inequalities 
\begin{equation}\label{eq:ttau<=}
\begin{split}
 \tilde\tau&\le\frac{\Mf}{2\tsi\sqrt{n}}\biggp{\tc_1\vsi_3\Bigp{\frac{v_3^2}{\sqrt{n}}+2v_2^2}+\tc_2\Bigp{\frac{v_3^2}{n^{1/3}}+
 2%RM 02-21-13 3
 ^{2/3}v_2^2}}+
 \frac{\tc_3}{\sqrt{n}}\Bigp{d(\cc)\vsi_3^3-\frac{(2\cc-1)_+^2}{4(1-\cc)}}\\ %RM 02-21-13\frac{\tc_3}{4(1-\cc)}\,\frac{\vsi_3^3}{\sqrt{n}}\\
 &\le\frac1{\sqrt{n}}\Bigp{\Cnecon+\Cnevsi\vsi_3^3+\Cnevtwo v_2^3+\Cnevthr v_3^3},
\end{split}
\end{equation}
where $\Cnecon,\dotsc,\Cnevthr$ are as in \eqref{eq:C10}--\eqref{eq:C13}; here, the first inequality of \eqref{eq:K2,K3,ga} is again used to see that $n\ge K_1^2\vsi_3^6\ge K_1^2$. 

Combine now the inequalities \eqref{eq:Michel}, \eqref{eq:tga1}, and \eqref{eq:ttau<=}; then \eqref{eq:iid,p=3,nu} follows.
\end{proof}

\begin{proof}[Proof of Lemma~\ref{lem:exp.iid.p=3}]
As we have let $X_i=V_i/n$ and $p=3$ in Theorem~\ref{thm:iid,p=3,nonunif}, \eqref{eq:si,iid} holds. 
Let now 
\begin{equation*}
 c_x:=\Bigp{\frac{2\pi_1}{\Mf}}^{1/2},\text{ so that }x_2=\frac{c_x(\tsi z)^{1/2}}{n^{1/4}}\text{ and }y_2=\frac{c_x(\tsi z)^{1/2}}{\ka_2n^{1/4}},
\end{equation*}
by \eqref{eq:x2,y2,vp2}. 
Then
\begin{align*} \E\norm{S_{y_2}}&\le\E\norm{S}+\E\norm{S_{y_2}-S}\le\norm{S}_2+\tsum\nolim_i\E\bignorm{X_i\ind{\norm{X_i}>y_2}}\le s_2+\frac{s_3^3}{y_2^2}=\frac{v_2}{\sqrt{n}}+\frac{v_3^3}{y_2^2n^2}\\
 &=x_2\Bigp{\frac{v_2}{c_x(\tsi z)^{1/2}n^{1/4}}+\frac{\ka_2^2v_3^3}{c_x^3(\tsi z)^{3/2}n^{5/4}}}
  =x_2\biggp{\frac1{c_x}\Bigp{\frac{v_2^4}{\tsi^3z^3\sqrt{n}}}^{1/4}\Bigp{\frac{\tsi z}{\sqrt{n}}}^{1/4}+\frac{\ka_2^2}{c_x^3}\,\frac{v_3^3}{\tsi^3z^3\sqrt{n}}\Bigp{\frac{\tsi z}{\sqrt{n}}}^{3/2}}\\
 &\le x_2\Bigp{\frac1{c_x}\Bigp{\frac{\om}{K_2}}^{1/4}+\frac{\ka_2^2\om^{3/2}}{c_x^3K_3}}=\hat\ga x_2,
\end{align*}
where \eqref{eq:z,iid,p=3} and \eqref{eq:K2,K3,ga} are used to obtain the last inequality above, and the definition \eqref{eq:gahat<1} of $\hat\ga$ is used for the last equality. 
Then, since $\hat\ga<1$ is assumed in \eqref{eq:z,om,Ki}, Lemma~\ref{lem:PU} yields  
\begin{equation}\label{eq:PUargs.iid}
 \PUtwo\le\PU(\hat u,\hat\ka_2,\vp_2),\text{ where }
 \hat u:=\frac{(1-\hat\ga)x_2y_2}{s_2^2}=\frac{c_x^2(1-\hat\ga)}{\ka_2}\,\frac{\tsi z}{v_2^2}\,\sqrt{n}
 \text{ and }\hat\ka_2:=(1-\hat\ga)\ka_2
\end{equation}
(recall \eqref{eq:PUcond}). 
Also, in accordance with \eqref{eq:x2,y2,vp2}, $\vp_2=\dfrac{s_3^3}{s_2^2y_2}\wedge1=\dfrac{\ka_2v_3^3}{c_xv_2^2(\tsi z)^{1/2}n^{3/4}}\wedge1$.

The inequality in \eqref{eq:exp3} is proved by taking any $\vp_\ast\in(0,1)$, as in Theorem~\ref{thm:iid,p=3,nonunif}, and considering two cases: (i) $\vp_2\in(\vp_\ast,1]$ and (ii) $\vp_2\in(0,\vp_\ast]$. 
Assume first that $\vp_2\in(\vp_\ast,1]$. 
By \eqref{eq:PUargs.iid} and \eqref{eq:PUchar2}, 
\begin{equation}\label{eq:L1veeL2}
 \PUtwo\le\PU\bigp{\hat u,\hat\ka_2,\vp_2}\le\exp\{L_1(\al,\hat u,\hat\ka_2,\vp_2)\}\vee\exp\{L_2(\al,\hat u,\hat\ka_2,\vp_2)\}
\end{equation}
for any $\al\in(0,1)$. 
Now introduce 
\begin{align}\label{eq:r2<=}
 r_2^2:=\frac1{\hat u}=\frac{\ka_2}{c_x^2(1-\hat\ga)}\,\frac{v_2^2}{\tsi z\sqrt{n}}
  =  \frac{\ka_2}{c_x^2(1-\hat\ga)}\Bigp{\frac{v_2^4}{\tsi^3z^3\sqrt{n}}}^{1/2}\Bigp{\frac{\tsi z}{\sqrt{n}}}^{1/2}
  &\le\frac{\ka_2\om^{1/2}}{c_x^2(1-\hat\ga)}\Bigp{\frac{v_2^4}{\tsi^3z^3\sqrt{n}}}^{1/2}\\
\label{eq:r2<=,2}
 &\le\frac{\ka_2\om^{1/2}}{K_2^{1/2}c_x^2(1-\hat\ga)}=\frac1{u_0}
\end{align}
and
\begin{align}\label{eq:rp<=}
 r_3^3:=\frac{\vp_2}{\hat u}\le\frac{\ka_2v_3^3}{c_xv_2^2(\tsi z)^{1/2}n^{3/4}}\,\frac{\ka_2v_2^2}{(1-\hat\ga)c_x^2\tsi z\sqrt{n}}
  =   \frac{\ka_2^2}{(1-\hat\ga)c_x^3}\,\frac{v_3^3}{\tsi^3z^3\sqrt{n}}\Bigp{\frac{\tsi z}{\sqrt{n}}}^{3/2}
  &\le \frac{\ka_2^2\om^{3/2}}{(1-\hat\ga)c_x^3}\,\frac{v_3^3}{\tsi^3z^3\sqrt{n}}\\
\label{eq:rp<=,2}
  & \le\frac{\ka_2^2\om^{3/2}}{(1-\hat\ga)c_x^3K_3}=t_3,
\end{align}
where \eqref{eq:z,iid,p=3} is used to establish the inequalities in \eqref{eq:r2<=} and \eqref{eq:rp<=}, and \eqref{eq:K2,K3,ga} and \eqref{eq:t2,t3,u0} are used for \eqref{eq:r2<=,2} and \eqref{eq:rp<=,2}. 

Next, in view of \eqref{eq:r2<=,2}, \eqref{eq:ka1hat}, and \eqref{eq:t2,t3,u0}, one has 
\begin{equation*}
 \frac{\hat\ka_2\al(2-\al)}{2(1-\vp_2)}\,\uhat
  \ge\frac{\hat\ka_2\al(2-\al)}{2(1-\vp_{\ast})}\,\frac{c_x^2(1-\hat\ga)}{\ka_2}\Bigp{\frac{K_2}{\om}}^{1/2}
  =\frac{\pi_1\al(2-\al)(1-\hat\ga)^2}{\Mf(1-\vp_\ast)}\Bigp{\frac{K_2}{\om}}^{1/2}
  =t_2. 
\end{equation*}
So, the case condition $\vp_2\in(\vp_\ast,1]$ together with the definitions of \eqref{eq:L1} and \eqref{eq:r2<=} of $L_1$ and $r_2^2$ imply
\begin{equation}\label{eq:L1<=}
 e^{L_1}\le e^{\hat\ka_2(1-\al-\al\vp_{\ast}/(1-\vp_\ast))}\Bigp{\frac{2(1-\vp_{\ast})}{\hat\ka_2\al(2-\al)}}^{2}\biggp{\sup_{t\ge t_2}t^{2}e^{-t}}r_2^{4}\le \Cnvtwoa\,\frac{v_2^4}{z^3\sqrt{n}},
\end{equation}
where the last inequality follows by the definition \eqref{eq:C81} of $\Cnvtwoa$ and \eqref{eq:r2<=} (on recalling also that $\hat\ka_2=(1-\hat\ga)\ka_2$). 
Note that if $\vp_2=1$ then, by the definition, $L_1=-\infty$, which makes \eqref{eq:L1<=} trivial (using the convention $\exp\{-\infty\}:=0$).

Again by the case condition $\vp_2\in(\vp_\ast,1]$, now together with \eqref{eq:L2} and \eqref{eq:rp<=,2}, 
\begin{equation}\label{eq:L2<=}
e^{L_2}\le e^{\hat\ka_2(1-\al)}\biggp{\sup_{t\in(0,t_3]}\frac1t\exp\biggb{-\hat\ka_2\Bigp{1-\frac\al2+t}\ln\Bigp{1+\frac{1-\al}t}}}r_3^3\le \Cnvthra \,\frac{v_3^3}{z^3\sqrt{n}},
\end{equation}
where the last inequality follows by the definition \eqref{eq:C83} of $\Cnvthra $ and \eqref{eq:rp<=}. 
Now, upon combining \eqref{eq:L1veeL2}, \eqref{eq:L1<=}, and \eqref{eq:L2<=}, we obtain the result \eqref{eq:exp3} in the case $\vp_2\in(\vp_\ast,1]$.

Consider the remaining case, when $\vp_2\in(0,\vp_{\ast}]$. 
Then, by \eqref{eq:incr}, \eqref{eq:PUargs.iid}, \eqref{eq:r2<=}, \eqref{eq:r2<=,2}, and the definition \eqref{eq:C82} of $\Cnvtwob$, 
\begin{equation}\label{eq:P(||S_y||>x)}
\begin{split}
 \PUtwo\le\PU\bigp{\hat u,\hat\ka_2,\vp_2}\le\PU\bigp{\hat u,\hat\ka_2,\vp_\ast}\le r_2^{4}\Bigp{\sup_{u\ge u_0}u^{2}\PU\bigp{u,\hat\ka_2,\vp_\ast}}\le \Cnvtwob\,\frac{v_2^4}{z^3\sqrt{n}}.
\end{split}
\end{equation}
Thus, \eqref{eq:P(||S_y||>x)} yields \eqref{eq:exp3} in the case $\vp_2\in(0,\vp_{\ast}]$ as well, and the lemma is proved.
\end{proof}

\clearpage
\section{Optimality of the restriction \texorpdfstring{$z=O(\sqrt{n})$}{z=O(sqrt\{n\})} for the nonuniform bound}
\label{app:z<=sqrt.n}

The following proposition shows that the upper bound on $z$ in \eqref{eq:z,iid}, and hence in \eqref{eq:z}, is in general optimal, up to the choice of the constant factor $\om$.

\begin{prop}\label{prop:z}
Let $\XX=\R$ and $f(x)\equiv x+x^2$, so that \eqref{eq:smooth} is satisfied when $L(x)\equiv x$, $\Mf=2$, and $\ep=1$. 
For any $p\in(2,3]$, let $V,V_1,\dots,V_n$'s be real-valued symmetric i.i.d.\ r.v.'s with density $|v|^{-p-1}\ln^{-2}|v|$ for all $|v|\ge v_0$, where the real number $v_0>1$ and the density values on $(-v_0,v_0)$ are chosen so that $\norm{V}_2=1$; note that then $\|V\|_p<\infty$. 
For any triple $b:=(b_1,b_2,b_3)$ of positive real numbers, let $\NZ(b)$ denote the set of all pairs $(n,z)\in\N\times(0,\infty)$ for which the inequality \eqref{eq:f(S).iid.nu} with $b_1,b_2,b_3$ in place of the three instances of $\CC$ holds. Then there exists a constant $\om(b)\in(0,\infty)$ depending only on $b$ such that \eqref{eq:z,iid} holds for all pairs $(n,z)\in\NZ$. 
\end{prop}

\begin{remark}\label{re:zconds.altbounds}
Let $r\in(0,p)$. Then an application of Chebyshev's inequality to the first two terms in the bound of \eqref{eq:f(S).iid.nu} yields
\begin{equation} \label{eq:f(S).iid.nu.cheby1}
\begin{split}
 &\Bigabs{\Bigprob{\frac{f(\bar V)}{\tsi /\sqrt n}\le z}-\Phi(z)}\\
 &\qquad
   \le \CC\Bigp{\frac{\E\norm{V}^r\ind{\norm{V}>\CC
   z\sqrt{n}%RM 03.20.13 z
   }}{z^rn^{r/2-1}%RM 03.20.13 z^r
   }+\frac{\E\norm{V}^r\ind{\norm{V}>\CC\sqrt{n}}}{z^rn^{r/2-1}%RM 03.20.13 z^r
   }+\frac{1}{(z\sqrt{n}%RM 03.20.13 z
   )^p}+\frac{1}{e^{z/\tth}n^{p/2-1}%RM 03.20.13 e^{z/\tth}
   }}
\end{split}
\end{equation}
for any $z$ satisfying \eqref{eq:z,iid}.  
The arguments of the proof of Proposition~\ref{prop:z} can be used to demonstrate that the bound of \eqref{eq:f(S).iid.nu.cheby1} (larger than that in \eqref{eq:f(S).iid.nu}) generally fails to hold if $z/\sqrt{n}\to\infty$. 
Using Chebyshev's inequality when $r=p$ yields 
\begin{equation}\label{eq:f(S).iid.nu.cheby2}
 \Bigabs{\Bigprob{\frac{f(\bar V)}{\tsi /\sqrt n}\le z}-\Phi(z)}\le\frac{\CC}{z^pn^{p/2-1}%RM 03.20.13 z^p
 }. 
\end{equation}
One might hope that a bound of the form in \eqref{eq:f(S).iid.nu.cheby2} could hold for all $f$ satisfying the smoothness condition \eqref{eq:smooth} and for all $z>0$. 
However, another modification of the proof of Proposition~\ref{prop:z} (which will be given in Section~\ref{sec:proofs}) demonstrates that \eqref{eq:f(S).iid.nu.cheby2} fails to be true whenever
\begin{equation}\label{eq:z.iid.cheby}
 \frac{z}{\sqrt{n}\ln^\al n}\to\infty,\quad\text{where $\al$ is any fixed number such that $\al p>1$;}
\end{equation}
the extra log factor above is needed because the bound in \eqref{eq:f(S).iid.nu.cheby2} is worse than that in \eqref{eq:f(S).iid.nu.cheby1}. 
\end{remark}

\begin{proof}[Proof of Proposition~\ref{prop:z}]
Let $S=\bar V$, so that $\si=\norm{L(S)}_2=1/\sqrt{n}$, $T=f(S)/\si=\sqrt{n}(S+S^2)$, and $W=L(S)/\si=\sqrt{n}S$. 
To obtain a contradiction, assume that Proposition~\ref{prop:z} is false. Then for some triple $b\in(0,\infty)^3$ and each value of $\om\in\N$ there is a pair $(n,z)=(n_\om,z_\om)\in\NZ(b)$ such that $z>\frac{\om}{\tsi}\,\sqrt n$. 
Now, for the rest of the proof of Proposition~\ref{prop:z}, let $\om\to\infty$, so that   
\begin{equation*}
 \zeta:=z/\sqrt{n}\to\infty;
\end{equation*} 
further let
\begin{equation*}
 \vt:=\zeta^{1/2}n=z^{1/2}n^{3/4},
\end{equation*}
so that $\vt/n=\zeta^{1/2}\to\infty$.
Note that for $v>v_0$ 
\begin{equation*}
 \P(V>v)=\int_v^\infty\frac{du}{u^{p+1}\ln^2u}
  \asymp\frac{1}{v^p\ln^2v}
\end{equation*}
as $v\to\infty$, which follows by l'Hospital's rule. 

So, 
\begin{gather*}
 n\bigprob{\norm{V}>\CC z\sqrt{n}}
  \asymp\frac{n}{z^pn^{p/2}\ln^2(z\sqrt{n})}
  =\frac{\ln^2(\zeta^{1/2}n)}{\zeta^{p/2}\ln^2(\zeta n)}\,\frac{n}{\vt^p\ln^2\vt}
  =o\bigp{n\P(V>\vt)},
\\
 \frac{n\P(\norm{V}>\CC\sqrt{n})}{z^p}
  \asymp\frac{n}{z^pn^{p/2}\ln^2\sqrt{n}}
  =\frac{\ln^2(\zeta^{1/2}n)}{\zeta^{p/2}\ln^2\sqrt{n}}\,\frac{n}{\vt^p\ln^2\vt}
  =o\bigp{n\P(V>\vt)},
\\
 \frac{1}{(z\sqrt{n}%RM 03.20.13 z
 )^p}
  =\frac{\ln^2(\zeta^{1/2}n)}{\zeta^{p/2}n}\,\frac{n}{\vt^p\ln^2\vt}
  =o\bigp{n\P(V>\vt)},
\\
 \frac{1}{e^{z/\tth}n^{p/2-1}%RM 03.20.13 e^{z/\tth}
 }
  =\frac{\zeta^{p/2}n^{p/2}\ln^2(\zeta^{1/2}n)}{e^{\zeta\sqrt{n}/\tth}}\,\frac{n}{\vt^p\ln^2\vt}
  =o\bigp{n\P(V>\vt)}, 
\end{gather*}
and 
\begin{equation}\label{eq:1-Phi} 
 1-\Phi(z)
  \asymp\frac{1}{ze^{z^2/2}}
  =\frac{\zeta^{p/2-1}n^{p-3/2}\ln^2(\zeta^{1/2}n)}{e^{\zeta^2n/2}}\,\frac{n}{\vt^p\ln^2\vt}
  =o\bigp{n\P(V>\vt)}. 
\end{equation}
Then \eqref{eq:f(S).iid.nu} and \eqref{eq:os} imply that $|\P(T\le z)-\Phi(z)|$ and $|\P(W\le z)-\Phi(z)|$ are both $o(n\P(V>\vt))$. 
Now let $\De=T-W=\sqrt{n}S^2$, so that 
\begin{equation}\label{eq:De=o(nP(V>v))}
 \P(\De>2z)
  \le\P(T>z)+\P(-W>z)
  =\P(T>z)+\P(W>z)
  =o\bigp{n\P(V>\vt)}, 
\end{equation}
by \eqref{eq:1-Phi}. 

On the other hand, by \cite[Lemma~2.3]{acosta79},
\begin{equation*}
 \P(\De>2z)
  =\P(\sqrt nS^2>2z)
  =\bigprob{\bigabs{\tsum\nolim_iV_i}>\sqrt{2}\vt}
 \ge\tfrac12\,(1-e^{-\psi})
\end{equation*}
for large enough $n$, where
\begin{equation*}
 \psi:=n\bigprob{|V|>\sqrt{2}\vt}=2n\bigprob{V>\sqrt{2}\vt}.
\end{equation*}
Since $\vt/n=\zeta^{1/2}\to\infty$, one has $\psi=o(n^{-p+1})\to0$, whence 
\begin{equation*}
 \P(\De>2z)\ge\tfrac\psi3>\tfrac2{3\cdot2^p}\,n\P(V>\vt)
\end{equation*}
for large enough $n$, which contradicts \eqref{eq:De=o(nP(V>v))}.

The statements of Remark~\ref{re:zconds.altbounds} are proved with only a few modifications to the above arguments, using the relation 
\begin{equation*}
 \E\norm{V}^r\ind{\norm{V}>v}\asymp\frac{1}{v^{p-r}\ln^2v}
\end{equation*}
as $v\to\infty$, for any $r\in(0,p)$. 
In order to show that \eqref{eq:f(S).iid.nu.cheby2} fails to hold simultaneously with \eqref{eq:z.iid.cheby}, let $V$ have density $1/(|v|^{p+1}\ln^{\al p}|v|)$ for $|v|\ge v_0>1$ (and still assume that $V$ is symmetric, with $v_0$ and density on $(-v_0,v_0)$ chosen to ensure that $\norm{V}_2=1$), $\zeta:=z/(\sqrt{n}\ln^\al n)$, and $\vt:=\zeta^{1/2}n=z^{1/2}n^{3/4}/\ln^{\al/2}n$. 
After these redefinitions, it is easy to verify that
\begin{equation*}
 \frac{1}{z^pn^{p/2-1}%RM 03.20.13 z^p
 }
  =\frac{\ln^{\al p}(\zeta^{1/2}n)}{\zeta^{p/2}\ln^{\al p}n}\,\frac{n}{\vt^p\ln^{\al p}\vt}
  \asymp\frac{\ln^{\al p}(\zeta^{1/2}n)}{\zeta^{p/2}\ln^{\al p}n}\,n\P(V>\vt)
  =o\bigp{n\P(V>\vt)},
\end{equation*}
from which \eqref{eq:De=o(nP(V>v))} follows and the contradiction is derived as done previously.
\end{proof}

\clearpage

\section{Proofs of bounds with explicit numerical constants, using a computer algebra system (CAS)}
\label{app:num.proofs}

This appendix contains proofs of Corollaries~\ref{cor:centralT}, \ref{cor:centralT.nub}, and \ref{cor:Runif}.  
The numerical computations that arise in these proofs are easily performed with a CAS; of course the calculations could, in principle, be done without the aid of a computer, but the amount of time required for such a task makes the use of a CAS practically indispensable. 

\begin{proof}[Proof of Corollary~\ref{cor:centralT}]
Consider the i.i.d.\ r.v.'s $V:=(Y,Y^2-1)$ and $V_i:=(Y_i,Y_i^2-1)$, taking values in $\XX=\R^2$ with the standard Euclidean norm, and let $f(\x):=x_1/\sqrt{1+x_2}$ for $\x=(x_1,x_2)\in\XX$ with $x_2>-1$ (also let $f(\x)$ take an arbitrary value for all other $\x\in\XX$). 
Further let $L=f'(0)$, so that $\norm{L}=1$, $L(V)=Y$, and $\tsi=\norm{L(V)}_2=1$. 
Then $\sqrt{n}f(\bar V)/\tsi=T_1$ a.s., by \eqref{eq:T1}. 
On recalling \eqref{eq:L,M,iid}, it is clear that $f$ satisfies the smoothness condition \eqref{eq:smooth} whenever $\ep<1$, whence the conditions of Theorem~\ref{thm:iid,p=3,unif} hold. 

For any $\x\in\XX$ such that $\norm{\x}\le\ep<1$, the spectral norm of the Hessian matrix $f''(\x)$ is 
\begin{equation*}
 \bignorm{f''(\x)}=\Bigabs{\tfrac{3x_1+\sqrt{9x_1^2+16(1+x_2)^2}}{8(1+x_2)^{5/2}}}\vee\Bigabs{\tfrac{3x_1-\sqrt{9x_1^2+16(1+x_2)^2}}{8(1+x_2)^{5/2}}}.
\end{equation*}
It is easy to see that $\norm{f''(\x)}$ is symmetric with respect to $x_1$; moreover, $\norm{f''(\x)}$ is increasing in $x_1\ge0$ and decreasing in $x_2$. 
Hence, 
\begin{equation}\label{eq:C1.selfnorm}
 \Mf=\sup_{\norm{\x}\le\ep}\bignorm{f''(\x)}=\sup_{\norm{\x}=\ep}\bignorm{f''(\x)}=\sup_{-\ep\le x_2\le0}\tfrac{3\sqrt{\ep^2-x_2^2}+\sqrt{9(\ep^2-x_2^2)+16(1+x_2)^2}}{8(1+x_2)^{5/2}};
\end{equation}
given some specific rational $\ep$, a CAS can be used to obtain an algebraic expression for $\Mf$. 

Next, introduce
\begin{equation}\label{eq:y3,y4,y6}
 y_3:=\norm{Y}_3,\quad y_4:=\norm{Y}_4,\quad\text{and}\quad y_6:=\bignorm{Y^2-1}_3^{1/2};
\end{equation} 
then \eqref{eq:tsi,ga3,v_al} yields
\begin{equation}\label{eq:ga3,v2,v3}
 \vsi_3=y_3,\quad v_2=y_4^2,\quad\text{and}\quad v_3=\bignorm{Y^2+(Y^2-1)^2}_{3/2}^{1/2}.
\end{equation}
For any nonnegative numbers $\tw_0$, $\tw_3$, and $\tw_4$, let
\begin{equation}\label{eq:nu3}
 \nu_3:=\nu_3(\tw_0,\tw_3,\tw_4):=\sup_{y\in\R}\frac{(y^2+(y^2-1)^2)^{3/2}}{\tw_0(1-y^2)+\tw_3|y|^3+\tw_4y^4+|y^2-1|^3},
\end{equation}
so that \eqref{eq:y3,y4,y6} and \eqref{eq:ga3,v2,v3} imply
\begin{equation}\label{eq:v3split}
 v_3^3\le\nu_3\cdot(\tw_3y_3^3+\tw_4y_4^4+y_6^6);
\end{equation} 
note that, whenever the numbers $\tw_0$, $\tw_3$, and $\tw_4$ happen to be such that the denominator in \eqref{eq:nu3} is negative for some $y\in\R$, then necessarily  $\nu_3(\tw_0,\tw_3,\tw_4)=\infty$ and the inequality in \eqref{eq:v3split} is trivially satisfied. 

Introduce arbitrary positive parameters $N_0\in\N$, $w_4$, and $w_6$. Consider two cases: (i) $n\le N_0-1$ and (ii) $n\ge N_0$. 
In the first case, when $n\le N_0-1$, use the inequalities $y_4\ge y_3\ge1$ to see that 
\begin{equation}\label{eq:smallN,unif}
 \bigabs{\P(T_1\le z)-\Phi(z)}\le1
  \le\frac{\sqrt{N_0-1}}{\sqrt{n}}
  \le\frac1{\sqrt{n}}\bigp{A_{3,1}y_3^3+A_{4,1}y_4^6+A_{6,1}y_6^6},
\end{equation}
where
\begin{equation*}
 A_{3,1}:=\frac{\sqrt{N_0-1}}{1+w_4},\quad A_{4,1}:=\frac{w_4\sqrt{N_0-1}}{1+w_4},\quad\text{and}\quad A_{6,1}:=0.
\end{equation*}

Consider then the case when $n\ge N_0$, and let $\cc$, $\ka_{2,0}$, $\ka_{3,0}$, $\ka_{2,1}$, and $\ka_{3,1}$ be as in \eqref{eq:ka's}. 
Further let $w_{6,2}:=1$, take any nonnegative numbers $w_{0,2}$, $w_{3,2}$, and $w_{4,2}$ (to be specified later), and let
\begin{equation}\label{eq:nujk,unif}
 \nu_{j,2}:=\nu_3\bigp{w_{0,2},w_{3,2},w_{4,2}}w_{j,2}\quad\text{for }j\in\{3,4,6\},
 \quad\text{so that}\quad
 v_3^3\le\nu_{3,2}y_3^3+\nu_{4,2}y_4^6+\nu_{6,2}y_6^6,
\end{equation}
by \eqref{eq:v3split}. 
Then \eqref{eq:iid,p=3,young's} and \eqref{eq:nujk,unif} imply 
\begin{equation}\label{eq:largeN,unif}
 \bigabs{\P(T_1\le z)-\Phi(z)}
 \le\frac1{\sqrt{n}}\bigp{\Cucon+\Cuvsi y_3^3+\Cuvtwo y_4^6+\Cuvthr v_3^3}\\
 \le\frac1{\sqrt{n}}\bigp{A_{3,2}y_3^3+A_{4,2}y_4^6+A_{6,2}y_6^6}
\end{equation}
where $\Cucon,\dotsc,\Cuvthr$ are as in \eqref{eq:tKu}, but with $N_0$ replacing each instance of $n$ in those expressions, 
\begin{equation*}\label{eq:Aj2,unif}
 A_{3,2}:=\pi
 \bigp{\Cucon}_+ %RM 02-21-13 \Cucon
 +\Cuvsi+\nu_{3,2}\Cuvthr,\quad A_{4,2}:=(1-\pi)
 \bigp{\Cucon}_+ %RM 02-21-13 \Cucon
 +\Cuvtwo +\nu_{4,2}\Cuvthr,\quad A_{6,2}:=\nu_{6,2}\Cuvthr,
\end{equation*}
and $\pi$ is any number in the interval $[0,1]$. 
Now choose $\pi$ to minimize $A_{3,2}\vee(A_{4,2}/w_4)$ subject to the constraint that $\pi\in[0,1]$; that is, let
\begin{equation*}
 \pi:=1\wedge\Bigp{\frac{\Cucon+\Cuvtwo +\nu_{4,2}\Cuvthr-w_4(\Cuvsi+\nu_{3,2}\Cuvthr)}{\Cucon(1+w_4)}}_+\,.
\end{equation*}
Of course, if $(\Cucon)_+=0$ then we may let $\pi$ be arbitrary. %RM 02-21-13

Referring now to \eqref{eq:smallN,unif} and \eqref{eq:largeN,unif}, we see that \eqref{eq:centralT} holds when
\begin{equation*}
 A_j:=A_{j,1}\vee A_{j,2}\quad\text{for}\quad j\in\{3,4,6\}.
\end{equation*} 
As mentioned before, the two triples $(A_3,A_4,A_6)$ in \eqref{eq:Ajlist.Tunif} are the result of trying to approximately minimize $A_3\vee(A_4/w_4)\vee(A_6/w_6)$, with weights $(w_4,w_6)\in\{1,0.25\}\times\{0.05\}$. 
Using a CAS to find exact expressions for $\Mf$ in \eqref{eq:C1.selfnorm} and $\nu_3$ in \eqref{eq:nu3}, and substituting the parameter values given in the table below (which should be interpreted as exact, rational numbers), one can verify that \eqref{eq:centralT} indeed holds with the specific values of the triples $(A_3,A_4,A_6)$ listed in \eqref{eq:Ajlist.Tunif}.  
\begin{center}
 \begin{tabular}{cc|*{9}{c@{\hspace{7pt}}}c||ccc}  $w_4$&$w_6$&$\ep$&$N_0$&$\cc$&$\ka_{2,0}$&$\ka_{3,0}$&$\ka_{2,1}$&$\ka_{3,1}$&$w_{0,2}$&$w_{3,2}$&$w_{4,2}$  &$A_3$&$A_4$&$A_6$\\\hline\hline
  1   &0.05&0.360&36&0.69 &1   &0.77&1.16&0.85&0.39&0&1  &2.99&2.99&0.15\\
  0.25&0.05&0.378&32&0.842&0.99&0.97&1.04&0.86&0.32&0&0.7&4.46&1.12&0.22
%RM 02-21-13 
%   1    &0.05 &0.346&44   &0.684&1          &0.870      &1          &0.874      &0.392    &0        &1          &3.33 &3.33 &0.17 \\
%   0.25 &0.05 &0.352&53   &0.744&0.762      &0.919      &0.768      &0.939      &0        &0.629    &0          &5.79 &1.45 &0.26
 \end{tabular}
\end{center}
%RM 02-21-13 
% To prove that the set in \eqref{eq:Ajlist.symm} can replace that in \eqref{eq:Ajlist.Tunif} when the symmetry of $Y$ is assumed, one need only redefine $\KKu{20},\dotsc,\KKu{31}$ according to Remark~\ref{re:symm.unif}, and then use the set of parameter values tabulated below.
% \begin{center}
%  \begin{tabular}{cc|*{9}{c@{\hspace{7pt}}}c||ccc}
%   $w_4$&$w_6$&$\ep$&$N_0$&$\cc$&$\ka_{2,0}$&$\ka_{3,0}$&$\ka_{2,1}$&$\ka_{3,1}$&$w_{0,2}$&$w_{3,2}$&$w_{4,2}$  &$A_3$&$A_4$&$A_6$\\\hline\hline
%   1    &0.05 &0.370&39   &0.669&1          &0.905      &1          &0.904      &0.392    &0        &1          &3.09 &3.09 &0.16 \\
%   0.25 &0.05 &0.377&45   &0.726&0.771      &0.886      &0.781      &0.931      &0        &0.498    &0          &5.32 &1.33 &0.27
%  \end{tabular}
% \end{center}
\end{proof}

\begin{proof}[Proof of Corollary~\ref{cor:centralT.nub}]
Adopt the notation used in the proof of Corollary~\ref{cor:centralT}; particularly recall \eqref{eq:y3,y4,y6} and \eqref{eq:ga3,v2,v3}). 
Recall also the positive parameters in \eqref{eq:nonunif.params} satisfying the constraints in \eqref{eq:z,om,Ki}; we shall specify their values later in the proof. 
In addition, take any
\begin{equation*}
 \ep\in(0,1),\ c\in(0,1),\ \pi_4\in[0,1],\ w_{j,k}\ge 0\text{ for }(j,k)\in\{0,3,4\}\times\{2,3\},\text{ and }w_{g,k}\ge0\text{ for }k\in\{1,2,3\}.
\end{equation*}
Then let (cf.\ \eqref{eq:nu3} and \eqref{eq:v3split})
\begin{equation}\label{eq:nu3k}
 \nu_{j,k}:=w_{j,k}\,\nu_3(w_{0,k},w_{3,k},w_{4,k})\text{ for }(j,k)\in\{3,4,6\}\times\{2,3\},\text{ so that }v_3^3\le\nu_{3,k}y_3^3+\nu_{4,k}y_4^8+\nu_{6,k}y_6^6,
\end{equation}
where $w_{6,k}:=1$ for $k\in\{2,3\}$, and also let (cf.\ \eqref{eq:g})
\begin{equation}\label{eq:g_j}
 g_k(z):=\frac1{z^3}+\frac{w_{g,k}}{e^{z/\tth}}\quad\text{for }k\in\{1,2,3\},\quad\text{where }\tth:=\frac{\th}{1-\pi_1}.
\end{equation}

Similarly to the proof of \cite[Theorem~1.1]{pin11}, consider three cases. 

\emph{Case 1 (``small $z$''): $0<z<z_0$.} 

Let $(A_3,A_4,A_6)$ be any triple of constants such that \eqref{eq:pin,T} holds; we shall provide specific values for the triple $(A_3,A_4,A_6)$ at the end of the proof, using general expressions obtained in  \cite[Theorem~1.2]{pin11}. 
Since $g_1$ in \eqref{eq:g_j} is decreasing on $(0,\infty)$, \eqref{eq:pin,T} and the case condition $0<z<z_0$ then imply
\begin{equation}\label{eq:case1}
 \bigabs{\P(T_1\le z)-\Phi(z)}\le\frac1{\sqrt{n}}\,\bigp{A_3y_3^3+A_4y_4^8+A_6y_6^6}\le\frac{g_1(z)}{\sqrt{n}}\,\bigp{A_{3,1}y_3^3+A_{4,1}y_4^8+A_{6,1}y_6^6},
\end{equation}
where
\begin{equation}\label{eq:Ak1}
 A_{j,1}:=\frac{A_j}{g_1(z_0)}\quad\text{for }j\in\{3,4,6\}.
\end{equation}

\emph{Case 2 (``large $z$, small $n$''): $z\ge z_0$ and \eqref{eq:K2,K3,ga} fails to hold.} 

Recall the definition \eqref{eq:T1} of $T_1$ and also that $c\in(0,1)$, and then note that
\begin{equation*}
 \P(T_1>z)
  \le\bigprob{\sqrt{n}\,\bar{Y}>\sqrt{c}z}+\bigprob{\bar{Y^2}< c}.
\end{equation*}
By \eqref{eq:tyurin,michel}, 
\begin{equation*}
 \bigprob{\sqrt{n}\,\bar{Y}>\sqrt cz}\le1-\Phi(\sqrt cz)+\frac{30.2211}{c^{3/2}}\,\frac{y_3^3}{z^3\sqrt n}.
\end{equation*}
Next, by \cite[Theorem~7]{pinutev86} with $\xi_i:=-Y_i^2$, 
\begin{equation*}
 \bigprob{\bar{Y^2}\le c}
  \le\exp\Bigb{-\frac n{y_4^4}\,(1-c+c\ln c)}
  \le\Bigp{\frac{2}{e(1-c+c\ln c)}}^2\frac{y_4^8}{n^2}
  \le\om^3\Bigp{\frac{2}{e(1-c+c\ln c)}}^2\frac{y_4^8}{z^3\sqrt n},
\end{equation*}
where $\sup_{x>0}x^2e^{-x}=(2/e)^2$ is used for the penultimate inequality above, and the restriction on $z$ \eqref{eq:z,T} is used for the last inequality. 
Thus, since $1-\Phi(z)<1-\Phi(z\sqrt{c})$ and $1/z^3\le g_2(z)$,
\begin{equation}\label{eq:case2,1}
 \abs{\P(T_1\le z)-\Phi(z)}\le h(z)+\frac{\tA_{3,2}y_3^3}{z^3\sqrt{n}}+\frac{\tA_{4,2}y_4^8}{z^3\sqrt{n}}\le h(z)+\frac{g_2(z)}{\sqrt{n}}\bigp{\tA_{3,2}y_3^3+\tA_{4,2}y_4^8},
\end{equation}
where
\begin{equation}\label{eq:h,tA}
 h(z):=1-\Phi(\sqrt{c}z),\quad\tA_{3,2}:=\frac{30.2211}{c^{3/2}},\quad\text{and}\quad\tA_{4,2}:=\om^3\Bigp{\frac{2}{e(1-c+c\ln c)}}^2.
\end{equation}
By the assumed conditions of Case~2, at least one of the inequalities in \eqref{eq:K2,K3,ga} fails to hold. Therefore and in view of \eqref{eq:nu3k},
\begin{equation}\label{eq:hle}
\begin{split}
 h(z)&\le h(z)\Bigp{\frac{K_1y_3^3}{\sqrt{n}}\vee\frac{K_2y_4^8}{z^3\sqrt{n}}\vee\frac{K_3v_3^3}{z^3\sqrt{n}}}\\
 &\le\frac{g_2(z)}{\sqrt{n}}\,
 \max\big(K_1S_{3,2}y_3^3,\,K_2S_{4,2}y_4^8,\,K_3S_{4,2}(\nu_{3,2}y_3^3+\nu_{4,2}y_4^8+\nu_{6,2}y_6^6)\big),  
\end{split}
\end{equation}
where 
\begin{equation}\label{eq:S32,S42}
 S_{3,2}:=\sup_{z\ge z_0}\frac{h(z)}{g_2(z)}\quad\text{and}\quad S_{4,2}:=\sup_{z\ge z_0}\frac{h(z)}{z^3g_2(z)}. 
\end{equation}
Thus, by \eqref{eq:case2,1} and \eqref{eq:hle}, 
\begin{equation}\label{eq:case2}
 \bigabs{\P(T_1\le z)-\Phi(z)}\le\frac{g_2(z)}{\sqrt{n}}\bigp{A_{3,2}y_3^3+A_{4,2}y_4^8+A_{6,2}y_6^6},
\end{equation}
where
\begin{equation*}
\begin{split}
 A_{3,2}&:=\tA_{3,2}+\max(K_1S_{3,2},\, K_3S_{4,2}\nu_{3,2}),\\
 A_{4,2}&:=\tA_{4,2}+S_{4,2}\max(K_2,\, K_3\nu_{4,2}),\\
 \text{and }A_{6,2}&:=K_3S_{4,2}\nu_{6,2}.
\end{split}
\end{equation*}

\emph{Case 3 (``large $z$, large $n$''): $z\ge z_0$ and \eqref{eq:K2,K3,ga} is true.}

In this final case, the assumptions of Theorem~\ref{thm:iid,p=3,nonunif} all hold when $\Mf$ is as in \eqref{eq:C1.selfnorm}. 
Recall now the definition \eqref{eq:g_j} of $\tth$, the inequality in \eqref{eq:nu3k}, and also note that $\Cnecon\le \Cnecon(\pi_4y_3^3+(1-\pi_4)y_4^8)$ (which follows because $1\le y_3\le y_4$). 
Then \eqref{eq:iid,p=3,nu} yields
\begin{equation}\label{eq:BEp+e}
\bigabs{\P(T_1\le z)-\Phi(z)}\le
\frac{z^{-3}\bigp{\be_3y_3^3+\be_4y_4^8+\be_6y_6^6}+e^{-z/\tth}\bigp{\be_{3,e}y_3^3+\be_{4,e}y_4^8+\be_{6,e}y_6^6}}{\sqrt n},
\end{equation}
where
\begin{equation*}
\begin{aligned}
 \be_3&:=\Cnvsi+\nu_{3,3}(\Cnvthra +\Cnvthrb),&\be_4&:=\Cnvtwoa\vee \Cnvtwob\vee(\nu_{4,3}\Cnvthra)+\nu_{4,3}\Cnvthrb,&\be_6&:=\nu_{6,3}(\Cnvthra+\Cnvthrb),\\ \be_{3,e}&:=\pi_4
 (\Cnecon)_+ %RM 02-21-13 \Cnecon
 +\Cnevsi+\nu_{3,3}\Cnevthr,&\be_{4,e}&:=(1-\pi_4)
 (\Cnecon)_+ %RM 02-21-13 \Cnecon
 +\Cnevtwo+\nu_{4,3}\Cnevthr,&\be_{6,e}&:=\nu_{6,3}\Cnevthr. 
\end{aligned}
\end{equation*}
Next, let
\begin{equation}\label{eq:ez}
 e_z:=\sup_{z\ge z_0}z^3e^{-z/\tth}.
\end{equation}
Then, by the definition \eqref{eq:g_j} of $g_3(z)$, for any $j\in\{3,4,6\}$
\begin{equation}\label{eq:b+be}
  \frac{\be_j}{z^3}+\frac{\be_{j,e}}{e^{z/\tth}}
   \le g_3(z)\sup_{z\ge z_0}\frac{\be_jz^{-3}+\be_{j,e}e^{-z/\tth}}{z^{-3}+w_{g,3}e^{-z/\tth}} 
    =g_3(z)\sup_{r\in(0,e_z]}\frac{\be_j+\be_{j,e}r}{1+w_{g,3}r}
   =A_{j,3}\,g_3(z), 
\end{equation}
where
\begin{equation*}
 A_{j,3}:=\be_j\vee\frac{\be_j+\be_{j,e}e_z}{1+w_{g,3}e_z}\quad\text{for }j\in\{3,4,6\}.
\end{equation*}
Now, by \eqref{eq:BEp+e} and \eqref{eq:b+be}, 
\begin{equation}\label{eq:case3}
 \bigabs{\P(T_1\le z)-\Phi(z)}\le\frac{g_3(z)}{\sqrt{n}}\,\bigp{A_{3,3}y_3^3+A_{4,3}y_4^8+A_{6,3}y_6^6}. 
\end{equation}

Now combine the inequalities in \eqref{eq:case1}, \eqref{eq:case2}, and \eqref{eq:case3}, and recall also the definitions \eqref{eq:g_j} of the functions $g_k$, to see that
\begin{equation}\label{eq:finalcase}
 \bigabs{\P(T_1\le z)-\Phi(z)}\le\frac{z^{-3}+w_ge^{-z/\tth}}{\sqrt{n}}\,\bigp{\hat A_3y_3^3+\hat A_4y_4^8+\hat A_6y_6^6},
\end{equation}
where
\begin{equation}\label{eq:Ahat}
 \hat{A}_j:=\max_{k\in\{1,2,3\}}A_{j,k}\text{ for }j\in\{3,4,6\},\quad\text{and}\quad w_g:=\max_{(j,k)\in\{3,4,6\}\times\{1,2,3\}}\frac{A_{j,k}}{\hat A_j}\,w_{g,k}.
\end{equation}
In view of \eqref{eq:finalcase} and \eqref{eq:centralT2}, the proof will be complete upon demonstrating the existence of a set of parameters such that the constants listed in Table~\ref{tab:Tbound} are in accordance with the definitions in \eqref{eq:Ahat}. 

Similarly to the proof of \cite[Theorem~1.1]{pin11}, those constants are obtained by trying to minimize the value of $\hat A_3\vee(\hat A_4/w_4)\vee(\hat A_6/w_6)$ for each of the points $\bigp{\om,w_g,(w_4,w_6)}\in\{0.1,0.5\}\times\{0,1\}\times\{(1,1),(0.5,0.2)\}$. Note that treating $w_g$ in \eqref{eq:Ahat} as an arbitrarily fixed constant introduces the restriction that $w_{g,k}\le w_g\min_j\hat A_j/A_{j,k}$ for each $k\in\{1,2,3\}$, and so $w_{g,k}=0$ when $w_g=0$; further, 
on recalling the definitions \eqref{eq:g} and \eqref{eq:g_j} of $g$ and $\tth$ along with the bound in \eqref{eq:finalcase}, one has the additional restriction that $\tth=2$, whence $\th=2(1-\pi_1)$. 

The parameters used to obtain the constants $\hat A_j$ are tabulated in Tables~\ref{tab:Tparms.wg=0} and \ref{tab:Tparms.wg=1} below. 
There are a few remarks to be made concerning the verification that the values listed in those tables indeed prove the statement of Corollary~\ref{cor:centralT.nub}. 
First, it is a practical necessity to use a sufficiently powerful CAS; we performed the calculations with the Mathematica software. 
In order to skirt any issue of rounding error in intermediate calculations, the values in Tables~\ref{tab:Tparms.wg=0} and \ref{tab:Tparms.wg=1} should be interpreted as being exact rational numbers; in this way, $\hat A_j$ (and the expressions upon which the $\hat A_j$'s depend) can be calculated to within any prescribed precision. 

Some care must be taken in order to implement the expressions for the $\hat A_j$'s. 
Note that $\nu_3$ in \eqref{eq:nu3} (used in the definition of $\nu_{j,k}$ in \eqref{eq:nu3k}) and $\Mf$ in \eqref{eq:C1.selfnorm} are algebraic expressions and therefore can be calculated exactly in a CAS. 
Concerning the numbers $A_{j,1}$ in \eqref{eq:Ak1}, the triples $(A_3,A_4,A_6)$ are obtained by similar calculations (with exact rational numbers) as directed by the proof of \cite[Theorem~1.2]{pin11}; one should also replace the absolute constant 0.4785 in the proof there (due to Tyurin \cite{tyurin09}) with the smaller constant 0.4748 (due to Shevtsova \cite{shev11}). 
For each of the two pairs $(w_4,w_6)\in\{(1,1),(0.5,0.2)\}$ considered here, the parameters used to obtain the triples $(A_3,A_4,A_6)$ are listed below (using the notation of \cite{pin11}):
\begin{equation*}
 \begin{array}{cc|ccc|ccccccc}
 w_4&w_6&A_3   &A_4   &A_6   &\al  &\vp_2&\vp_3&\vp_4&\th_3&\th_4&\ka  \\\hline
 1  &1  &1.5175&1.4852&1.4814&0.080&0.206&3.187&0.135&0.415&2.898&0.173\\
 0.5&0.2&1.9946&0.9996&0.1897&0.216&0.369&0.761&0.278&0.408&3.532&0.275
 \end{array}
\end{equation*}

Also note that \eqref{eq:C7}--\eqref{eq:C83}, \eqref{eq:S32,S42}, and \eqref{eq:ez} contain expressions of the general form $\sup_{x\ge x_0}k(x)$ or $\sup_{0<x\le x_0}k(x)$ for some function $k$ and positive number $x_0$. For the specific values of the parameters listed in Tables~\ref{tab:Tparms.wg=0} and \ref{tab:Tparms.wg=1}, one can use Lemma~\ref{lem:incr-decr} below to see that these suprema are all attained at the boundary point $x_0$. 
Finally, bounding $\ta_1$ in \eqref{eq:ta1} involves estimating the root of the equation in \cite[(2.3)]{pin09.winsor}; as noted at the end of the paragraph containing formula \eqref{eq:L_W}, 
$L_{W;\,c,\si}$ is nonincreasing in $\si>0$, and hence any upper bound on the mentioned root results in an upper bound on $\ta_1$. 
Implementation of the expressions $\hat A_j$ in accordance with the above remarks and the parameter values listed in Tables~\ref{tab:Tparms.wg=0} and \ref{tab:Tparms.wg=1} will then demonstrate that \eqref{eq:centralT2} holds. 

To prove that the statement of Corollary~\ref{cor:centralT.nub} holds when $Y$ is assumed to be symmetric and Table~\ref{tab:Tbound.symm} is used in place of Table~\ref{tab:Tbound}, one need only amend the definitions of $\tc_2$ and $\Cnvsi$ as prescribed by Remark~\ref{re:symm.nonunif}, and then use the parameter values given below in Tables~\ref{tab:Tparms.wg=0.symm} and \ref{tab:Tparms.wg=1.symm}.
\end{proof}

\begin{table}[ht]\small 
\begin{center}
\caption{Parameters associated with Corollary~\ref{cor:centralT.nub}, for $w_g=0$}\label{tab:Tparms.wg=0}
\begin{tabular}{c||cc|cc}
\multicolumn{5}{p{0.52\textwidth}}{For all columns below, set $w_{g,1}=w_{g,2}=w_{g,3}=0$, $\ka_3=1.5$, $\pi_2=1-\pi_1-\pi_3$, $\vp_\ast=0.001$, $\ka_{2,0}=\ka_{3,0}=\ka_{3,1}=1$, and $\pi_4=0$.}
\\\hline
$\om$&\multicolumn{2}{c|}{0.5}&\multicolumn{2}{c}{0.1}\\
$(w_4,w_6)$&(1,1)&(0.5,0.2)&(1,1)&(0.5,0.2)\\\hline
$\ep$ & 0.232 & 0.301 & 0.054 & 0.073  \\
$z_0$ & 4.782 & 4.855 & 4.629 & 4.390  \\
$c$   & 0.757 & 0.759 & 0.900 & 0.891  \\
$K_1$ & $6.9\times 10^4$ & $1.3\times 10^5$ & $2.0\times 10^5$ & $9.2\times 10^4$  \\
$K_2$ & $6.3\times 10^6$ & $4.0\times 10^6$ & $2.2\times 10^7$ & $4.1\times 10^6$  \\
$K_3$ & $6.9\times 10^6$ & $3.4\times 10^6$ & $2.3\times 10^7$ & $1.6\times 10^6$  \\
$w_{0,2}$& 0.156 & 0.380 & 0.206 & 0.147  \\
$w_{3,2}$& 0.400 & 0.036 & 0.600 & 1.000  \\
$w_{4,2}$& 0.380 & 1.000 & 0.742 & 0.778  \\
$\cc$ & 0.536 & 0.621 & 0.500 & 0.514  \\
$\th$ & 0.861 & 0.880 & 0.875 & 0.978  \\
$w$   & 0.360 & 0.316 & 0.376 & 0.398  \\
$\de_0$& 0.007 & 0.009 & 0.007 & 0.010  \\
$\pi_1$& 0.042 & 0.083 & 0.008 & 0.015  \\
$\pi_3$& 0.645 & 0.635 & 0.660 & 0.660  \\
$\ka_2$& 2.108 & 2.093 & 2.102 & 2.116  \\
$\ka_{2,1}$& 1.570 & 0.800 & 6.050 & 1.612  \\
$\al$  & 0.070 & 0.050 & 0.075 & 0.080  \\
$w_{0,3}$& 0.278 & 0.275 & 0.216 & 0.392  \\
$w_{3,3}$& 0     & 0.365 & 0     & 0          \\
$w_{4,3}$& 0.595 & 0.980 & 0.45 & 1       \\\hline
$\hat A_3$& 166   & 229   & 151  & 169\\%RM 03.20.13 170     \\
$\hat A_4$& 166   & 115   & 148  & 85      \\
$\hat A_6$& 165   & 45    & 147  & 29     
\end{tabular}
\end{center}
\end{table}

\begin{table}[ht] \small 
\begin{center}
\caption{Parameters associated with Corollary~\ref{cor:centralT.nub}, for $w_g=1$}\label{tab:Tparms.wg=1}
\begin{tabular}{c||cc|cc}
\multicolumn{5}{p{0.52\textwidth}}{For all columns below, set $w_{g,1}=w_{g,3}=1$, $w_{g,2}=0$, $\ka_3=1.5$, $\pi_2=1-\pi_1-\pi_3$, $\th=2(1-\pi_1)$, $\al=0.05$, $\vp_\ast=0.001$, and $\ka_{3,0}=\ka_{3,1}=1$.}\\\hline
$\om$      &\multicolumn{2}{c|}{0.5}&\multicolumn{2}{c}{0.1}\\
$(w_4,w_6)$&(1,1)&(0.5,0.2)&(1,1)&(0.5,0.2)\\\hline
$\ep$      & 0.363 & 0.438 & 0.066 & 0.112  \\
$z_0$      & 6.800 & 7.175 & 6.550 & 6.074  \\
$c$        & 0.738 & 0.708 & 0.885 & 0.874  \\
$K_1$ & $4.0\times 10^5$ & $5.0\times 10^7$ & $7.9\times 10^6$    & $1.2\times 10^6$  \\
$K_2$ & $9.0\times 10^7$ & $3.5\times 10^9$ & $6.6\times 10^{10}$ & $1.6\times 10^9$  \\
$K_3$ & $1.0\times 10^8$ & $5.5\times 10^9$ & $1.3\times 10^{10}$ & $3.5\times 10^8$  \\
$w_{0,2}$  & 0.040 & 0.133 & 0.263 & 0.142  \\
$w_{3,2}$  & 0     & 1     & 0.100 & 0.590  \\
$w_{4,2}$  & 0.080 & 0.600 & 0.588 & 0.396  \\
$\cc$      & 0.490 & 0.741 & 0.500 & 0.552  \\
$w$        & 1.160 & 0.940 & 1.655 & 1.530  \\
$\de_0$    & 0.039 & 0.027 & 0.016 & 0.018  \\
$\pi_1$    & 0.108 & 0.257 & 0.012 & 0.038  \\
$\pi_3$    & 0.422 & 0.409 & 0.415 & 0.423  \\
$\ka_2$    & 2.095 & 2.012 & 2.011 & 2.017  \\
$\ka_{2,0}$& 1     & 0.799 & 1     & 1.046  \\
$\ka_{2,1}$& 0.983 & 1.496 & 4.750 & 1.104  \\
$\pi_4$    & 0     & 0.467 & 0     & 0.220  \\
$w_{0,3}$  & 0.392 & 0.318 & 0.392 & 0.392  \\
$w_{3,3}$  & 0     & 0.224 & 0     & 0      \\
$w_{4,3}$  & 1     & 1     & 1     & 1      \\\hline
$\hat A_3$ & 48    & 66    & 38    & 39     \\
$\hat A_4$ & 48    & 33    & 36    & 20     \\
$\hat A_6$ & 42    & 13    & 36    & 7      \\
\end{tabular}
\end{center}
\end{table}

\begin{table}[ht]\small 
\begin{center}
\caption{Parameters associated with Remark~\ref{re:symm.T1}, for $w_g=0$}\label{tab:Tparms.wg=0.symm}
\begin{tabular}{c||cc|cc}
\multicolumn{5}{p{0.52\textwidth}}{For all columns below, set $w_{g,1}=w_{g,2}=w_{g,3}=0$, $\ka_3=1.5$, $\pi_2=1-\pi_1-\pi_3$, $\vp_\ast=10^{-4}$, and $\ka_{3,0}=\ka_{3,1}=1$.}
\\\hline
$\om$&\multicolumn{2}{c|}{0.5}&\multicolumn{2}{c}{0.1}\\
$(w_4,w_6)$&(1,1)&(0.5,0.2)&(1,1)&(0.5,0.2)\\\hline
$\ep$& 0.264 & 0.310 & 0.072 & 0.082 \\
$z_0$& 4.527 & 4.679 & 4.328 & 4.170 \\
$c$  & 0.750 & 0.762 & 0.900 & 0.918 \\
$K_1$& $2.3\times 10^4$ & $3.0\times 10^4$ & $5.4\times 10^4$ & $3.4\times 10^4$ \\
$K_2$& $1.8\times 10^6$ & $1.4\times 10^6$ & $5.1\times 10^6$ & $8.0\times 10^5$ \\
$K_3$& $2.0\times 10^6$ & $1.5\times 10^6$ & $2.3\times 10^6$ & $6.0\times 10^5$ \\
$w_{0,2}$& 0.274 & 0.173 & 0.144 & 0.153 \\
$w_{3,2}$& 0.214 & 0.852 & 0.100 & 1 \\
$w_{4,2}$& 0.688 & 0.916 & 0.300 & 1 \\
$\cc$& 0.565 & 0.643 & 0.510 & 0.581 \\
$\th$& 0.849 & 0.894 & 0.890 & 1.060 \\
$w$  & 0.320 & 0.381 & 0.430 & 0.344 \\
$\de_0$& 0.010 & 0.048 & 0.009 & 0.038 \\
$\pi_1$& 0.054 & 0.090 & 0.009 & 0.019 \\
$\pi_3$& 0.655 & 0.601 & 0.664 & 0.655 \\
$\ka_2$& 2.137 & 2.119 & 2.143 & 2.159 \\
$\ka_{2,0}$& 1 & 1.127 & 0.848 & 1 \\
$\ka_{2,1}$& 1.310 & 0.868 & 3.819 & 1.142 \\
$\al$& 0.200 & 0.150 & 0.120 & 0.150 \\
$w_{0,3}$& 0.276 & 0.220 & 0.280 & 0.392 \\
$w_{3,3}$& 0 & 0.595 & 0 & 0 \\
$w_{4,3}$& 0.590 & 1 & 0.600 & 1 \\\hline
$\hat A_3$ & 141    & 205%RM 03.20.13 206
& 124    & 145     \\
$\hat A_4$ & 138    & 103    & 123%RM 03.20.13 124
& 73     \\
$\hat A_6$ & 138    & 42    & 121    & 22      \\
\end{tabular}
\end{center}
\end{table}

\begin{table}[ht]\small 
\begin{center}
\caption{Parameters associated with Remark~\ref{re:symm.T1}, for $w_g=1$}\label{tab:Tparms.wg=1.symm}
\begin{tabular}{c||cc|cc}
\multicolumn{5}{p{0.52\textwidth}}{For all columns below, set $w_{g,1}=w_{g,3}=1$, $w_{g,2}=0$, $\ka_3=1.5$, $\pi_2=1-\pi_1-\pi_3$, $\th=2(1-\pi_1)$,  $\vp_\ast=10^{-4}$, and $\ka_{3,0}=\ka_{3,1}=1$.}
\\\hline
$\om$&\multicolumn{2}{c|}{0.5}&\multicolumn{2}{c}{0.1}\\
$(w_4,w_6)$&(1,1)&(0.5,0.2)&(1,1)&(0.5,0.2)\\\hline
$\ep$& 0.365 & 0.456 & 0.153 & 0.131 \\
$z_0$& 6.800 & 6.885 & 6.200 & 6.015 \\
$c$& 0.738 & 0.677 & 0.910 & 0.894 \\
$K_1$& $1.0\times 10^5$ & $8.2\times 10^4$ & $4.0\times 10^5$ & $4.0\times 10^5$\\
$K_2$& $3.0\times 10^6$ & $1.0\times 10^9$ & $5.0\times 10^7$ & $2.0\times 10^8$ \\
$K_3$& $1.0\times 10^7$ & $5.0\times 10^8$ & $9.0\times 10^7$ & $1.0\times 10^8$ \\
$w_{0,2}$& 0 & 0.392 & 0.224 & 0.018 \\
$w_{3,2}$& 0.030 & 0 & 0.481 & 0.514 \\
$w_{4,2}$& 0 & 1 & 0.704 & 0.041 \\
$\cc$& 0.760 & 0.703 & 0.470 & 0.625 \\
$w$& 0.692 & 0.913 & 1.612 & 1.163 \\
$\de_0$& 0.124 & 0.078 & 0.055 & 0.282 \\
$\pi_1$& 0.144 & 0.291 & 0.023 & 0.052 \\
$\pi_3$& 0.453 & 0.393 & 0.432 & 0.461 \\
$\ka_2$& 2.082 & 2.015 & 2.053 & 2.024 \\
$\ka_{2,0}$& 1.588 & 1.101 & 1.476 & 1.313 \\
$\ka_{2,1}$& 0.838 & 0.796 & 2.474 & 3.073 \\
$\pi_4$& 0.487 & 0.950 & 0 & 0.368 \\
$\al$& 0.067 & 0.150 & 0.103 & 0.137 \\
$w_{0,3}$& 0.363 & 0.251 & 0.239 & 0.383 \\
$w_{3,3}$& 0 & 0.461 & 0 & 0.026 \\
$w_{4,3}$& 0.856 & 1 & 0.500 & 1\\\hline
$\hat A_3$& 48 & 57 & 35 & 37 \\
$\hat A_4$& 48 & 29 & 32 & 19 \\
$\hat A_6$& 41 & 12 &31 & 5
\end{tabular}
\end{center}
\end{table}

\begin{lem}\label{lem:incr-decr}
Say that a function $k$ is $\nearrow\searrow$ on $(0,\infty)$ whenever there exists a point $x_\ast\in(0,\infty)$ such that $k$ is increasing on $(0,x_\ast)$ and decreasing on $(x_\ast,\infty)$. 
Also say that ``the supremum of a function $k$ is attained at the finite (or positive) boundary point'' if $\sup_{x\ge x_0}k(x)=k(x_0)$ (or $\sup_{0<x\le x_0}k(x)=k(x_0)$). 
Then the following statements are all true:
\begin{enumerate}[(i)]
 \item For any $c\in(0,1)$, the function $h$ as defined in \eqref{eq:h,tA} is decreasing on $(0,\infty)$.
 \item For any $p>0$ and $\ka>0$, the function $x\mapsto x^pe^{-\ka x}$ is $\nearrow\searrow$ on $(0,\infty)$.
 \item For any $0<p\le\ka$ and $\vp\in(0,1)$, the function $x\mapsto x^p\PU(x,\ka,\vp)$ is $\nearrow\searrow$ on $(0,\infty)$.
 \item For any $\ka\ge2$ and $\al\in(0,1)$, the function $x\mapsto \frac1x\,\exp\bigb{-\ka(1-\al/2+x)\ln(1+(1-\al)/x)}$ is $\nearrow\searrow$ on $(0,\infty)$.
 \item For any $c\in(0,1)$, the function $x\mapsto x^3h(x)$ is $\nearrow\searrow$ on $(0,\infty)$, where $h$ is as in \eqref{eq:h,tA}.
\end{enumerate}
The suprema in the expressions \eqref{eq:C7}--\eqref{eq:C83}, \eqref{eq:S32,S42}, and \eqref{eq:ez} are all attained at the respective finite (or positive) boundary points whenever the values in Tables~\ref{tab:Tparms.wg=0} and \ref{tab:Tparms.wg=1} are substituted in those expressions.
\end{lem}

\begin{proof}[Proof of Lemma~\ref{lem:incr-decr}]
Statements (i) and (ii) are trivial to verify by differentiation. 

By \eqref{eq:PU1}, to prove statement (iii), it suffices to show that
\begin{equation*}
 x\mapsto p\ln x+\frac{\ka}{2(1-\vp)x}\biggp{(1-\vp)^2\Bigp{1+\Lam\bigp{\tfrac{\vp}{1-\vp}\exp\bigb{\tfrac{\vp+x}{1-\vp}}}}^2-(\vp+x)^2-(1-\vp^2)}
\end{equation*}
is $\nearrow\searrow$ on $(0,\infty)$. 
Now let $w:=\frac{1-\vp}{\vp}\Lam\bigp{\frac{\vp}{1-\vp}\exp\{\frac{\vp+x}{1-\vp}\}}$, whence $x=(1-\vp)(\frac{\vp}{1-\vp}\,w+\ln w)-\vp$, and note that $w$ continuously increases from 1 to $\infty$ as $x$ increases from 0 to $\infty$. 
Thus, it suffices to show that
\begin{equation*}
 k(w):=p\ln\Bigp{(1-\vp)\bigp{\tfrac{w\vp}{1-\vp}+\ln w}-\vp}+\frac{\ka\bigp{(1-\vp)^2\bigp{1+\frac{w\vp}{1-\vp}}^2-(1-\vp)^2\bigp{\frac{w\vp}{1-\vp}+\ln w}^2-(1-\vp^2)}}{2(1-\vp)\bigp{(1-\vp)(\frac{w\vp}{1-\vp}+\ln w)-\vp}}
\end{equation*}
is $\nearrow\searrow$ on $(1,\infty)$. Next, introduce 
\begin{equation*}
 k_1(w):=\frac{2w((w-1)\vp+(1-\vp)\ln w)^2}{1-\vp+\vp w}\,k'(w)=2\bigp{p+\vp(\ka-p)}\ln w-2(w-1)(\ka-p)\vp-(1-\vp)\ka\ln^2w,
\end{equation*}
and note that $k_1$ and $k'$ have the same sign on $(1,\infty)$. Also introduce 
\begin{equation*}
 k_2(w):=\frac w2\,k_1'(w)=p-\vp(w-1)(\ka-p)-(1-\vp)\ka\ln w.
\end{equation*}
Then $k_2$ and $k_1'$ share the same sign on $(1,\infty)$ and $k_2$ is decreasing on $(1,\infty)$. 
Further, since $k_2(1)=p>0$ and $k_2(\infty)=-\infty$, we see that $k_2$ and hence $k_1'$ change sign once from $+$ to $-$ on $(1,\infty)$; that is, $k_1$ is $\nearrow\searrow$ on $(1,\infty)$. 
As $k_1(1)=0$ and $k_1(\infty)=-\infty$, it follows that $k_1$ and hence $k'$ change sign once from $+$ to $-$ on $(1,\infty)$. 
That is, $k$ is $\nearrow\searrow$ on $(1,\infty$), and thus statement (iii) is proved.

To prove (iv), let
\begin{equation*}
\begin{split}
 k(x)&:=-\ka\Bigp{1-\tfrac\al2+x}\ln\Bigp{1+\tfrac{1-\al}x}-\ln x,\\
 k_1(x)&:=k'(x)=\frac{\ka(1-\al)(2+2x-\al)}{2x(1-\al+x)}-\frac1x-\ka\ln\Bigp{1+\tfrac{1-\al}x},\\
 k_2(x)&:=2x^2(1-\al+x)^2k_1'(x)=2x^2-2x(1-\al)(\ka-2)-(1-\al)^2(\ka(2-\al)-2).
\end{split}
\end{equation*}
Then $k_2$ is decreasing on $(0,x_\ast)$ and increasing on $(x_\ast,\infty)$, where $x_\ast:=\frac12(\ka-2)(1-\al)$. 
Since $k_2(0)=-(1-\al)^2(\ka(2-\al)-2)<0$ and $k_2(\infty)=\infty$, it follows that $k_2$ and hence $k_1'$ change sign once from $-$ to $+$ on $(0,\infty)$. 
So, $k_1$ is $\searrow\nearrow$ on $(0,\infty)$; as $k_1(0+)=\infty$ and $k_1(\infty)=0$, we see that $k_1$ changes sign once from $+$ to $-$, and hence $k$ is $\nearrow\searrow$ on $(0,\infty)+$. 
Thus, $x\mapsto\exp\{k(x)\}$ is $\nearrow\searrow$ on $(0,\infty)$, proving statement (iv). 

The proof of part (v) is easily done by using the l'Hospital-type rule for monotonicity, as in the proof of Lemma~3 in \cite{pin07}. 

To finish the proof, make the various substitutions from Tables~\ref{tab:Tparms.wg=0} and \ref{tab:Tparms.wg=1} into the respective expressions of \eqref{eq:C7}--\eqref{eq:C83}, \eqref{eq:S32,S42}, and \eqref{eq:ez}; note that, since $w_{g,2}=0$ in all of the parameter sets, $g_2(z)=z^{-3}$ and hence $S_{4,2}=h(z_0)$ follows from statement (i) and \eqref{eq:S32,S42}. 
Next estimate the unique positive critical point $x_\ast$ of each of the functions in statements (ii)--(v) by finding rational numbers $x_1<x_2<x_3$ such that $k(x_1)<k(x_2)$ and $k(x_2)>k(x_3)$; then we shall know that $x_\ast\in(x_1,x_3)$. 
So, it will follow that $\sup_{x\ge x_0}k(x)$ is attained at the boundary point $x_0$ by checking that $x_0\ge x_3$, and that $\sup_{0<x\le x_0}k(x)$ is attained at $x_0$ by checking that $x_0<x_1$. 
Thus, one completes the proof. 
\end{proof}

\begin{proof}[Proof of Corollary~\ref{cor:Runif}]
For $\al\ge1$, let
\begin{gather*}
 y_\al:=\norm{Y}_\al\quad\text{and}\quad z_\al:=\norm{Z}_\al.
\end{gather*}
Also adopt the notation of Theorem~\ref{thm:R}, with $\rho=0$, so that $V=(Y,Z,Y^2-1,Z^2-1,YZ)$, $L(V)=YZ$, and $\tsi=\norm{YZ}_2$. 
Take any natural number $N_0$ and any real number $b_3>0$, and consider the two cases: \break (i) $n\le N_0-1$ and (ii) $n\ge N_0$. 

In the first case, when $n\le N_0-1$, note that $1\le(y_6^6+z_6^6)/2$ (since $1=y_2\le y_6$ and $1=z_2\le z_6$) and $\tsi^3\le(y_4z_4)^3\le y_6^3z_6^3\le(y_6^6+z_6^6)/2$ (which follows by H\"older's and Young's inequalities). 
Then
\begin{equation}\label{eq:smallN,Runif}
 \bigabs{\P(\sqrt{n}R/\tsi\le z)-\Phi(z)}\le1\le\frac{\sqrt{N_0-1}}{\sqrt{n}}\le\frac{y_6^6+z_6^6}{\sqrt{n}}\Bigp{B_{0,1}+\frac{B_{3,1}}{\tsi^3}},
\end{equation}
where
\begin{equation}\label{eq:Bj1,Runif}
 \bigp{B_{0,1},B_{3,1}}:=\frac{\sqrt{N_0-1}}{2(1+b_3)}\,\bigp{1,b_3}.
\end{equation}

Suppose then that $n\ge N_0$. 
Take any $\ep\in(0,\sqrt{3}/2)$ and 
%IP 03.11.13 $\cc\in(0,1)$, 
$\cc\in[\frac12,1)$ so that the conditions of Theorem~\ref{thm:iid,p=3,unif} are satisfied (cf.\ the discussion following \eqref{eq:f.R}); also introduce the 
parameter $\ka>0$. %RM 02-21-13 parameters $\ka>0$ and $\pi\in[0,1]$, so that $\KKu{0}\le\pi\KKu{0}+(1-\pi)\KKu{0}\vsi_3^3$. %RM 12.21.12 $\ka>0$.  
Recall the notation in \eqref{eq:tsi,ga3,v_al}, so that
\begin{gather*}
\vsi_3=\norm{YZ}_3/\tsi\le y_6z_6/\tsi,\quad\vsi_3^3\le\tfrac1{2}\bigp{y_6^6+z_6^6}/\tsi^3,\\
1\le v_2^3\le v_3^3\le\sup_{(y,z)\in\R^2}\frac{(y^2+z^2+(y^2-1)^2+(z^2-1)^2+y^2z^2)^{3/2}}{1-y^2+1-z^2+y^6+z^6}\,\bigp{y_6^6+z_6^6}=\tfrac{3^{3/2}}{2}\bigp{y_6^6+z_6^6},\\
v_2^2\le v_3^2\le1+\tfrac2{3^{3/2}}v_3^3\le\tfrac12\bigp{y_6^6+z_6^6}+\tfrac{2}{3^{3/2}}v_3^3\le\tfrac32\bigp{y_6^6+z_6^6},\\ 
\vsi_3v_2^2\le\vsi_3v_3^2\le y_6z_6v_3^2/\tsi\le\bigp{y_6^3z_6^3+\tfrac2{3^{3/2}}v_3^3}/\tsi\le\tfrac32\bigp{y_6^6+z_6^6}/\tsi;
\end{gather*}
in the last two lines we use the following instance of Young's inequality: $ab\le a^3+2(b/3)^{3/2}$ for $a\ge0$ and $b\ge0$. 
Then \eqref{eq:iid,p=3} implies
\begin{equation}\label{eq:largeN,Runif} 
 \Bigabs{\Bigprob{\frac{R}{\tsi/\sqrt{n}}\le z}-\Phi(z)}
  \le\frac{y_6^6+z_6^6}{\sqrt{n}}\Bigp{A_0+\frac{A_1}{\tsi}+\frac{A_2}{\tsi^2}+\frac{A_3}{\tsi^3}}
  \le\frac{y_6^6+z_6^6}{\sqrt{n}}\Bigp{B_{0,2}+\frac{B_{3,2}}{\tsi^3}},
\end{equation}
where
\begin{equation}\label{eq:Aj,Runif}
A_0:=
\tfrac12\,\bigp{\KKu{0}}_+ + %RM 02-21-13
 \tfrac{3}{2\ep^2\sqrt{N_0}}\wedge\tfrac{3^{3/2}(
 2%RM 03.20.13 3
 +1/\sqrt{N_0})}{2\ep^3N_0},\quad
A_1:=\tfrac32(\KKu{20}+\KKu{30})\tsi,\quad
A_2:=\tfrac32(\KKu{21}+\KKu{31})\tsi,\quad
A_3:=\tfrac12\,\KKu{1},
\end{equation}
with $N_0$ replacing $n$ in the expressions $\KKu{1},\dotsc,\KKu{3,1}$, 
\begin{equation}\label{eq:Bj2,Runif}
 B_{0,2}:=A_0+\tfrac23\,\ka^{-3/2}A_1+\tfrac13\,\ka^{-3}A_2,\quad\text{and}\quad
 B_{3,2}:=A_3+\tfrac13\,\ka^3A_1+\tfrac23\ka^{3/2}A_2.
\end{equation}

Then \eqref{eq:smallN,Runif} and \eqref{eq:largeN,Runif} yield the desired inequality \eqref{eq:Runif} if we let
\begin{equation}\label{eq:B0,B3}
 B_0:=B_{0,1}\vee B_{0,2}\quad\text{and}\quad B_3:=B_{3,1}\vee B_{3,2}.
\end{equation}
We shall show that, for $f$ as in \eqref{eq:f.R}, 
\begin{equation}\label{eq:M,Runif}
 \text{\eqref{eq:smooth} holds for any pair }(\ep,\Mf)\in\bigb{(0.06,1.094),(0.17,1.365),(0.25,1.688),(0.30,1.962)}.
\end{equation} 
Then, substituting the values of the parameters $b_3$, $N_0$, $\ep$, 
%RM 02-21-13 $\Mf$, 
$\cc$, and $\ka$ given in the table below into the expressions for $B_0$ and $B_3$ in \eqref{eq:B0,B3} (which depend on the expressions in \eqref{eq:Bj1,Runif}, \eqref{eq:Bj2,Runif}, \eqref{eq:Aj,Runif}, and \eqref{eq:Ku}), one will see that \eqref{eq:Runif} holds for any of the pairs $(B_0,B_3)$ listed in \eqref{eq:Bjlist,Runif}.  

\begin{center}
\begin{tabular}{c|cccc||cc}
$b_3$&$N_0$&$\ep$&$\cc$ &$\ka$  &$B_0$&$B_3$\\\hline\hline
1    &209  &0.25 &0.77  &0.983  &3.61 &3.61\\
8    &405  &0.3  &0.877 &1.745  &1.12 &8.94\\
1/8  &900  &0.17 &0.6115&0.4416 &13.33&1.69\\
27   &965  &0.3  &0.909 &2.339  &0.56 &14.97\\
1/27 &5674 &0.06 &0.5635&0.28273&36.32&1.37 %RM 02-21-13
%RM 12.21.12 
% \begin{tabular}{c|ccccc||cc}
% $b_3$&$N_0$&$\ep$&$\Mf$&$\cc$&$\ka$&$B_0$&$B_3$\\\hline\hline
% 27   &1580 &0.3  &1.962  &0.922&2.324&0.71 &19.16\\
% 8    &614  &0.3  &1.962  &0.893&1.745&1.38 &11.02\\
% 1    &267  &0.25 &1.688 &0.794&0.978&4.08 &4.08 \\
% 1/8  &1098 &0.17 &1.365&0.644&0.452&14.73&1.85 \\
% 1/27 &6618 &0.06 &1.094 &0.570&0.290&39.22&1.47
\end{tabular}
\end{center}

%RM 02-21-13
% To prove Remark~\ref{re:R.symm}, one only needs to redefine the pre-constants $\KKu{20},\dotsc,\KKu{31}$ as directed by Remark~\ref{re:symm.unif} and use the parameter values tabulated below: 
% 
% \begin{center}
% \begin{tabular}{c|cccccc||cc}
% $b_3$&$N_0$&$\ep$&$\Mf$&$\cc$&$\ka$&$\pi$&$B_0$&$B_3$\\\hline\hline
% 27   &1436 &0.3  &1.962&0.922&2.374&0&0.68 &18.27\\
% 8    &555  &0.3  &1.962&0.892&1.807&0&1.31 &10.47\\
% 1    &224  &0.25 &1.688&0.792&1.0148&0&3.74&3.74 \\
% 1/8  &881  &0.17 &1.365&0.627&0.452&1&13.19&1.65 \\
% 1/27 &5699 &0.06 &1.094&0.559&0.285&1&36.40&1.35
%RM 12.21.12 
% \begin{tabular}{c|ccccc||cc}
% $b_3$&$N_0$&$\ep$&$\Mf$&$\cc$&$\ka$&$B_0$&$B_3$\\\hline\hline
% 27   &1436 &0.3  &1.962&0.922&2.374&0.68 &18.27\\
% 8    &555  &0.3  &1.962&0.892&1.807&1.31 &10.47\\
% 1    &224  &0.25 &1.688&0.792&1.0148&3.74&3.74 \\
% 1/8  &940  &0.17 &1.365&0.625&0.443&13.63&1.71 \\
% 1/27 &6243 &0.06 &1.094&0.558&0.278&38.10&1.42
% \end{tabular}
% \end{center}

To complete the proof of Corollary~\ref{cor:Runif},
%RM 02-21-13 (and Remark~\ref{re:R.symm}), 
it now remains to verify \eqref{eq:M,Runif}. 
Toward that end, take any $\ep\in(0,\sqrt{3}/2)$, and recall the definition \eqref{eq:f.R} of $f$ (with $\rho=0$) to see that
\begin{equation*}
 f(x_1,x_2,x_3,x_4,x_5)\equiv f(-x_1,-x_2,x_3,x_4,x_5)\equiv -f(-x_1,x_2,x_3,x_4,-x_5)\equiv -f(x_1,-x_2,x_3,x_4,-x_5)
\end{equation*}
and
\begin{equation*}
 f(x_1,x_2,x_3,x_4,x_5)\equiv f(x_2,x_1,x_4,x_3,x_5)
\end{equation*}
for any $\x\in\R^5$ such that $\norm{\x}\le\ep$. 
The above identities then imply
\begin{equation*}
 \Mf^*:=\sup_{\norm{\x}\le\ep}\norm{f''(\x)}=\sup\bigb{\norm{f''(\x)}\colon\x\in\tB_\ep\cap\tilde\R^5};
\end{equation*}
here $\tB_\ep$ denotes the open $\ep$-ball about the origin and 
\begin{equation*}
 \tilde\R^5:=\bigb{\x\in\R^5\colon\sgn(x_1)=\sgn(x_2)=\sgn(x_5)\text{ and }x_3\le x_4},\quad\text{where}\quad\sgn(x):=\ind{x\ge0}-\ind{x<0}.
\end{equation*}

Next take any positive $m\in\N$, and let $\de_\ep:=\ep/m$. 
For any $\bu=(u_1,\dotsc,u_5)\in\Z^5$, let
\begin{equation*}  C_\bu:=\prod_{j=1}^5[u_j\de_\ep,(u_j+1)\de_\ep],\quad\text{and}\quad\bc_\bu:=\bigp{\regp{u_1+\tfrac12}\de_\ep,\dotsc,\regp{u_5+\tfrac12}\de_\ep};
\end{equation*}
that is, $C_\bu$ is the cube of side length $\de_\ep$ with its ``southwest'' corner at the point $\de_\ep\bu$ and center at $\bc_\bu$. 
Introduce also the set 
\begin{equation*}
 U:=\Bigb{\bu\in\Z^5\cap\tilde\R^5\colon\tB_\ep\cap C_\bu\ne\emptyset}=\Bigb{\bu\in\Z^5\cap\tilde\R^5\colon\tsum_{i=1}^5\bigp{u_j+\tfrac12-\tfrac12\,\sgn(u_j)}^2< m^2},
\end{equation*}
so that $\tB_\ep\cap\tilde\R^5\subseteq\bigcup_{\bu\in U}C_\bu$. 
Then 
\begin{align}
\notag
 \Mf^*&\le\max_{\bu\in U}\sup_{\x\in C_\bu}\norm{f''(\x)}\le\max_{\bu\in U}\Bigp{\bignorm{f''(\bc_\bu)}+\sup_{\x\in C_\bu}\bignorm{f''(\x)-f''(\bc_\bu)}_F}\\
 \label{eq:S_ep<=}
 &\le\max_{\bu\in U}\Bigp{\bignorm{f''(\bc_\bu)}+\sqrt{5}\,\frac{\de_\ep}2\sup_{\x\in C_\bu}\bignorm{f'''(\x)}_F},
\end{align}
where
\begin{equation*}
 \bignorm{f'''(\x)}_F:=\biggp{\tsum_{i,j,k=1}^5\bigp{f_{ijk}(\x)}^2}^{1/2} 
\end{equation*}
and $f_{ijk}=\partial^3f/(\partial x_i\partial x_j\partial x_k)$; here we assume that $m$ is chosen large enough (whence $\de_\ep$ is small enough) so as to ensure $f_{ijk}$ exists and is continuous on each cube $C_\bu$ (i.e.\ $\min_{\bu\in U}\inf_{\x\in C_\bu}(1+x_3-x_1^2)(1+x_4-x_2^2)>0$). 

Take now any $\bu\in U$, and then take any $\x\in \intr C_\bu$, so that $x_j\ne0$ for any $j\in\{1,\dotsc,5\}$. 
It is easy to see with a CAS that 
\begin{equation}\label{eq:f'''.to.poly}
 \norm{f'''(\x)}_F^2=\frac{3\tx_3\tx_4}{64}\,p(\tilde\x),\quad\text{where}\quad\tilde\x:=\bigp{\tx_1,\dotsc,\tx_5}:=\Bigp{x_1,x_2,\frac{1}{1+x_3-x_1^2},\frac{1}{1+x_4-x_2^2},x_5},
\end{equation}
and $p$ is a polynomial, namely, the sum of 172 monomials with integer coefficients; note that $\tx_3$ and $\tx_4$ are both positive. 
Further, $p(\tilde\x)$ can be bounded from above by bounding each of the 172 monomials. 
To do that, for $j\in\{1,2,5\}$ introduce
\begin{equation*}
 \tx_{j,1}:=\bigp{u_j+\tfrac12+\tfrac12\,\sgn(u_j)}\de_\ep\quad\text{and}\quad\tx_{j,-1}:=\bigp{u_j+\tfrac12-\tfrac12\,\sgn(u_j)}\de_\ep,
\end{equation*}
so that $|\tx_{j,-1}|\le|\tx_j|\le|\tx_{j,1}|$; 
also, for $j\in\{3,4\}$ let
\begin{equation*}
 \tx_{j,1}:=\frac{1}{1+u_{j}\de_\ep-\tx_{j-2,1}^2}\quad\text{and}\quad\tx_{j,-1}:=\frac{1}{1+(u_j+1)\de_\ep-\tx_{j-2,-1}^2},
\end{equation*}
so that $0<\tx_{j,-1}\le\tx_j\le\tx_{j,1}$. 
Then, for any nonnegative integers $d_1,\dotsc, d_5$, any integer $a$, and $s:=\sgn(a)\sgn(u_1)^{d_1}\sgn(u_2)^{d_2}\sgn(u_5)^{d_5}$,
\begin{equation}\label{eq:sup.mon}
 a\tx_1^{d_1}\tx_2^{d_2}\tx_3^{d_3}\tx_4^{d_4}\tx_5^{d_5}
 =s|a||\tx_1|^{d_1}\dotsb|\tx_5|^{d_5}
 \le s|a||\tx_{1,s}|^{d_1}\dotsb|\tx_{5,s}|^{d_5}
 =a\tx_{1,s}^{d_1}\tx_{2,s}^{d_2}\tx_{3,s}^{d_3}\tx_{4,s}^{d_4}\tx_{5,s}^{d_5},
\end{equation}
which follows since $\tx_j\ge0$ whenever $u_j\ge0$ (and $\tx_j\le0$ whenever $u_j<0$) for $j\in\{1,2,5\}$. 
Replacing each of the monomial summands in $p(\tilde\x)$ with their upper bound in \eqref{eq:sup.mon}, we see from \eqref{eq:f'''.to.poly} that
\begin{equation}\label{eq:p1,p-1}
 \norm{f'''(\x)}_F\le\frac{\sqrt{3\tx_{3,1}\tx_{4,1}}}{8}\,\sqrt{p_{\sgn(u_1)}(\tx_{1,1},\dotsc,\tx_{5,1},\tx_{1,-1},\dotsc,\tx_{5,-1})},
\end{equation}
where $p_1$ and $p_{-1}$ are each polynomials in the 10 variables (in fact, $p_{-1}$ is a polynomial in only the five variables $\tx_{1,1},\dotsc,\tx_{5,1}$, as it turns out that $s=1$ for each of the monomials of $p(\tx)$ for $\bu\in U$ with $u_1<0$). 

Thus, combining \eqref{eq:S_ep<=} and \eqref{eq:p1,p-1}, one has
\begin{equation*}
 \Mf^*\le\max_{\bu\in U}\Bigp{\norm{f''(\bc_\bu)}+\frac{\ep\sqrt{15\tx_{3,1}\tx_{4,1}}}{16m}\,\sqrt{p_{\sgn(u_1)}(\tx_{1,1},\dotsc,\tx_{5,1},\tx_{1,-1},\dotsc,\tx_{5,-1}}\,}.
\end{equation*}
One can then write a program in a CAS which will give an algebraic number for the latter upper bound (and then to bound that algebraic number with a rational). 
In particular, upon letting $m=19$ and implementing the bound above for $\ep\in\{\frac6{100},\frac{17}{100},\frac{25}{100},\frac{30}{100}\}$, \eqref{eq:M,Runif} follows. 
\end{proof}

\clearpage

\section{On Fisher's \texorpdfstring{$z$}{z} transform}\label{app:z} 
A statistic closely related to Pearson's $R$ is commonly known as the Fisher $z$ transform, defined by the formula $R_z:=\tanh^{-1}(R)=\frac12\ln\bigp{\frac{1+R}{1-R}}$. 
An advantage to using $R_z$ (as opposed to $R$) in making statistical inferences about $\rho$ follows from its variance-stabilizing property in normal populations; that is, $n\var(R_z)\to1$ for all $\rho\in(-1,1)$ as $n\to\infty$, as opposed to $n\var(R)\to(1-\rho^2)^2$, whenever $(Y,Z)$ has a bivariate normal distribution. 
Moreover, the distribution of $R_z$ converges to normality more rapidly than does the distribution of $R$ (especially for non-zero values of $\rho$) when the pair $(Y,Z)$ comes from a normal population; see e.g.\ Fisher \cite{fisher21}, David \cite{david38}, and Hotelling \cite{hot53}. 
In his discussion of Hotelling's paper, Kendall provides heuristics suggesting that such variance stabilization of the distribution of a statistic may often result in it being closer to normality. 
Namely, if an approximate constancy of the variance of a statistic were the same as an approximate constancy of its distribution itself, and if the distribution is close to normality at least for one value of the parameter (say, $\rho$, as in the present case), then it would be close to normality for all values of $\rho$. 

However, it is well known that the closeness of the distribution of a statistic to normality is usually mainly determined, not by the variance, but by the third moments of the underlying distribution. 
It is therefore natural to wonder whether or to what extent the nice properties of the $z$ transform hold for non-normal populations. 
For moderate sample sizes $n$, Gayen \cite{gayen51} observed that the convergence to normality for both $R$ and $R_z$ is lessened for non-normal populations with $\rho\ne0$, and Monte Carlo sampling performed by Berry and Mielke \cite{berry-mielke} suggests that the presence of skewness or heavy tails in the population of $(Y,Z)$ significantly reduces the accuracy of a normal approximation to $R_z$ when $\rho\ne0$. 
In \cite{corr.asymp}, explicit expressions for $\De_R=\lim_{n\to\infty}\sqrt{n}\abs{F_R-\Phi}_K$ and $\De_{R_z}=\lim_{n\to\infty}\sqrt{n}\abs{F_{R_z}-\Phi}_K$ are derived, where $F_R$ and $F_{R_z}$ are the d.f.'s of $R$ and $R_z$ and $\abs{\cdot}_K$ denotes the Kolmogorov distance. 
These ``asymptotic distances'' generally depend on up to the sixth moments of $Y$ and $Z$ when $\rho\ne0$, and it is demonstrated in \cite{corr.asymp} that, if the distribution of $(Y,Z)$ is not bivariate normal, $\De_{R_z}$ can be just as easily greater than $\De_R$ as less.   

In light of the above considerations, we now briefly investigate how any of the BE-type bounds of Section~\ref{sec:f(S)}, when applied to the statistic $R_z$, would fare in a comparison with corresponding bounds associated with $R$. 
Aside from the choice of parameter values, the only differences between the applications of our bounds to $R$ and $R_z$ are those arising from the choice of $f$; namely, upon letting $g(\x):=\tanh^{-1}(f(\x)+\rho)-\tanh^{-1}\rho$ for all $\x$ with $f$ as in \eqref{eq:f.R}, one has $g(\bar V)=R_z-\tanh^{-1}\rho$. 
In the case when $\rho=0$, we see that $f'(0)=g'(0)$ and $g''(0)=f''(0)$; moreover, in view of results in \cite{bhatt78}, one can see that 
an asymptotic expansion up to $\bigO(1/\sqrt{n})$ of the d.f.\ of $R$ is identical to that of $R_z$, whether or not the population of $(Y,Z)$ is Gaussian. 

Despite these similarities between $R$ and $R_z$, it appears that $M_g:=\sup_{\|\x\|\le\ep}\|g''(\x)\|>M_f:=\break\sup_{\|\x\|\le\ep}\|f''(\x)\|$ for $\ep>0$, at least when $\rho=0$. 
In particular, we showed (in the proof of Corollary~\ref{cor:Runif}) that $M_f\le1.962$ when $\ep=\frac3{10}$; on the other hand, one can see that 
\begin{equation*}
\norm{g''(\x)}>2.104\text{ and }\norm{\x}<\tfrac3{10}\text{ when }\x=-\bigp{\tfrac{28269}{200000},\tfrac{28269}{200000},\tfrac{45081}{500000},\tfrac{45081}{500000},\tfrac{183801}{1000000}},  
\end{equation*}
so that $M_g>2.104>1.962\ge M_f$, which will result, at least using the method presented in this paper, in a worse BE-type bound for $R_z$ as compared with that for $R$. 
In view of these points, one can conclude that, at least for $\rho=0$, the use of Fisher's $z$ transform $R_z$ in place of Pearson's $R$ will hardly yield better BE-bounds. 

%RM12.26
\clearpage%RMejs
\section{Compactness of the covariance operator} 
\label{sec:compact}%RMejs 
\newcommand{\HH}{\mathbb{H}}

Here we give a short proof that the covariance operator of a r.v.\ with finite second moment is compact. 
Let $X$ be a r.v.\ taking values in a separable Hilbert space $\HH$ such that $\E\norm{X}^2<\infty$ and $\E X=\mu$. 
Then the covariance operator $R\colon\HH\to\HH$ is defined by
\begin{equation*}
 Rx:=\E\langle x,X-\mu\rangle (X-\mu)=\E\ol{\langle X-\mu,x\rangle}(X-\mu);
\end{equation*}
let us assume w.l.o.g.\ that $\mu=0$. 
Note that $R$ is both self-adjoint and nonnegative-definite: for all $x,y\in\HH$
\begin{equation*}
 \langle Rx,y\rangle=\E\langle x,X\rangle\langle X,y\rangle=\E\ol{\langle y,X\rangle}\,\ol{\langle X,x\rangle}=\ol{\langle Ry,x\rangle}=\langle x,Ry\rangle
\end{equation*}
and
\begin{equation*}
 \langle Rx,x\rangle=\E\langle x,X\rangle\langle X,x\rangle=\E\bignorm{\langle x,X\rangle}^2\ge0.
\end{equation*}

Now let $(e_j)_{j\in\N}$ be any orthonormal basis of $\HH$, so that $X=\sum_j\langle X,e_j\rangle e_j$. 
Further take any $x\in\HH$, so that $Rx=\E\langle x,X\rangle\sum_j\langle X,e_j\rangle e_j$. 
For $n\in\N$, define the operator $R_n$ by $R_nx=\E\langle x,X\rangle\sum_{j=1}^n\langle X,e_j\rangle e_j$, and note that the range of $R_n$ is finite-dimensional. 
Moreover, if $\norm{x}\le 1$, then
\begin{align*}
  \bignorm{(R-R_n)x}&=\Bignorm{\E\langle x,X\rangle\tsum_{j=n+1}^\infty\langle X,e_j\rangle e_j}\\
  &\le\E\Bignorm{\langle x,X\rangle\tsum_{j=n+1}^\infty\langle X,e_j\rangle e_j}\\
  &\le\E\norm{X}\sqrt{\tsum_{j=n+1}^\infty\langle X,e_j\rangle^2}\mathop{\to}_{n\to\infty}0;
\end{align*}
the limit holds by dominated convergence, since $\sqrt{\sum_{j=n+1}^\infty\langle X,e_j\rangle^2}\le\sqrt{\sum_{j=1}^\infty\langle X,e_j\rangle^2}=\norm{X}$. 
As $x$ was arbitrary and the above majorant of $\norm{(R-R_n)x}$ does not depend on $x$, it follows that $\norm{R-R_n}\to0$; that is, $R$ is the limit (in the operator norm) of a sequence of finite-dimensional linear operators on $\HH$, and so is compact.

\clearpage%RMejs
\section{On the spectral decomposition of a covariance operator of a random vector in an arbitrary separable Hilbert space} 
\label{sec:inf-dim}%RMejs 

\newcommand{\rs}{\operatorname{res}}
\newcommand{\RR}{\mathcal{R}}
\renewcommand{\d}{\mathrm{d}}
\renewcommand{\sp}{\operatorname{sp}}

Let $X$ be a random vector in a separable Hilbert space $(H,\ip\cdot\cdot)$ with $\E\|X\|^2<\infty$. 
Let $R$ be the covariance operator of $X$. So, $R$ is self-adjoint.  
Obviously, any self-adjoint operator is normal. 
Hence, by \cite[Theorem~2.10, page~260]{kato}, 
\begin{equation}\label{eq:sum}
 R=\sum_{\la\in\La}\la P_\la, 
\end{equation}
where $\La$ is the (necessarily at most countable) set of all (necessarily nonnegative) eigenvalues of $R$;
(in the case when the set $\La$ is infinite) the sum converges in the operator norm; and, for each $\la\in\La$, $P_\la$ is the orthoprojector onto the eigenspace (say $E_\la$) of $\la$, which is necessarily of a finite dimension $n_\la:=\dim E_\la=\tr P_\la$ if $\la\ne0$. At that, 
\begin{equation}\label{eq:total}
\sum_{\la\in\La}P_\la=I,  
\end{equation}
the identity operator, and the eigenspaces $E_\la$ are pairwise mutually orthogonal:  
\begin{equation}\label{eq:ortho}
P_\la P_\mu=\ind{\la=\mu}P_\la   
\end{equation}
for all $\la$ and $\mu$ in $\La$. 

Moreover, for each $\la\in\La$, let $B_\la$ be any orthonormal basis of $E_\la$, so that $B:=\bigcup_{\la\in\La}B_\la$ is an orthonormal basis of $H$. 
Then
$\tr R=\sum_{\la\in\La}\la n_\la=\sum_{\la\in\La}\sum_{e\in B_\la}\ip{Re}e
=\sum_{e\in B}\E|\ip eX|^2=\E\sum_{e\in B}|\ip eX|^2=\E\|X\|^2<\infty$, so that $\sum_{\la\in\La}\la n_\la<\infty$. 
So, the set $\La$ of all eigenvalues of $R$ may have at most one limit point, and any limit point of $\La$ must be $0$. 

The spectrum $\sp R$ of $R$ is defined as the set of all $z\in\C$ such that the linear operator $R-zI$ does not have a bounded inverse. 
It follows that $\sp R$ coincides with $\La$ if $\dim H<\infty$ and with $\La\cup\{0\}$ if $\dim H=\infty$. The complementary set $\rs R:=\C\setminus\sp R$ is called the resolvent set. 
Let $B(H)$ denote the Banach space of all bounded linear operators $A\colon H\to H$. 

One can now define the resolvent $\RR\colon\rs R\to B(H)$ by the formula 
\begin{equation}\label{eq:res}
 \RR(z):=(R-zI)^{-1}=\sum_{\la\in\La}\frac1{\la-z}\,P_\la;  
\end{equation}
the latter equality can be easily verified in view of \eqref{eq:sum}, \eqref{eq:ortho}, and \eqref{eq:total}, because $R-zI=\sum_{\la\in\La}(\la-z)\,P_\la$. 

Take now any nonzero $\la\in\La$, which is necessarily an isolated point of the set $\La$. 
So, there is an open disc $D_\la$ in $\C$ such that $\la\in D_\la$ but no other point of the set $\La\cup\{0\}$ is in the closure of $D_\la$. 
Let  
$\Ga_\la$ be the boundary of $D_\la$. 
Then, by \eqref{eq:res} and the Cauchy integral theorem, 
\begin{equation}\label{eq:P=}
 P_\la=-\frac1{2\pi i}\,\int_{\Ga_\la}\RR(z)\,\d z,  
\end{equation}
whence 
\begin{equation}\label{eq:la=}
 \la=\frac1{n_\la}\tr P_\la=\frac1{n_\la}\tr RP_\la
 =-\frac1{2\pi i n_\la}\,\int_{\Ga_\la}\tr R\RR(z)\,\d z. 
\end{equation}
Formulas \eqref{eq:P=} and \eqref{eq:la=} are important, because it is comparatively easy to analyze the resolvent.

\clearpage
%IP
\textbf{Acknowlegment.} We are pleased to thank the Referees for their stimulating comments, which resulted in improved and more explicit bounds, as well as in better exposition. 

%RM12.26
\bibliographystyle{imsart-number}
\bibliography{citations.nodoi}
% \bibliographystyle{abbrv}
% \bibliography{citations.02.21.13} %RM 02-21-13 12.13.12}
%RM 12.21.12 \bibliography{citations}

%IP 03.12.13 edited {chen-fang}, {pin13.rosen}, and {pin12-2smooth} in citations.02.21.13 for capitalization

\end{document}